\documentclass{amsart}
\usepackage{amssymb,latexsym,graphics}
\usepackage[mathscr]{eucal}

\newtheorem{thm}{Theorem}[section]

\newtheorem{lem}[thm]{Lemma}
\newtheorem{prop}[thm]{Proposition}

\theoremstyle{remark}
\newtheorem{rem}[thm]{Remark}

\newtheorem{defn}[thm]{Definition}
\theoremstyle{definition}

\newcommand{\lp}[2]{\Vert \, #1 \, \Vert_{#2}}
\newcommand{\llp}[1]{ \Vert \, #1 \, \Vert }
\newcommand{\uF}{\underline{F}}
\newcommand{\uD}{\underline{D}}
\newcommand{\RR}{\mathbb{R}^n}

\newcommand{\td}{\widetilde}

\newcommand{\ou}{{}^\omega \! u }
\newcommand{\oL}{{}^\omega \! L }
\newcommand{\uL}{ {}^\omega \! \underline{L}}
\newcommand{\oA}{{}^\omega \! A }
\newcommand{\og}{{}^\omega \! g }
\newcommand{\oB}{{}^\omega \! B }
\newcommand{\uA}{ \underline{A}_{\, \bullet \ll 1} }
\newcommand{\uAl}{ \underline{A}_{\, \bullet \ll \lambda} }
\newcommand{\uAlg}{ \underline{A}_{\, \lambda\lesssim \bullet} }
\newcommand{\Co}{ \mathbb{C}\otimes \mathfrak{o}(m) }
\newcommand{\oPi}{ {}^\omega \! \Pi }
\newcommand{\ooPi}{ {}^\omega \! \overline{\Pi} }
\newcommand{\sd}{ {\, {\slash\!\!\! d  }} }
\newcommand{\suA}{ {\, {\slash\!\!\!\! \underline{A}_{\, \bullet \ll 1}  }} }
\newcommand{\oC}{{}^\omega \! C }
\newcommand{\uoC}{{}^\omega \! \underline{C} }
\newcommand{\uoA}{{}^\omega \! \underline{A} }
\newcommand{\oh}{{}^\omega \! h }

\newcommand{\ret}{\vspace{.3cm}}

\begin{document}

\title[Regularity for 6D YM]
{Global Regularity for the Yang--Mills Equations
on High Dimensional Minkowski Space}
\author{Joachim Krieger}
\address{}
\email{}
\author{Jacob Sterbenz}
\address{}
\email{}
\thanks{}
\date{}
\dedicatory{}
\commby{}


\begin{abstract}
We study here the global Cauchy problem for the Yang-Mills equations
on $(6+1)$ and higher dimensional Minkowski space, when the initial
data sets  are small in the critical gauge covariant Sobolev space
$\dot{H}_A^{(n-4)/{2}}$. Regularity is obtained through a certain 
``microlocal geometric renormalization'' of the equations which is
constructed via a family of approximate Coulomb gauge--null
Cr\"onstrom gauge transformations. The proof is then reduced to
controlling Hodge systems and degenerate elliptic equations on high
index and non-isotropic $L^p$ spaces, as well as the proof of some
bilinear estimates in auxiliary square-function spaces.
\end{abstract}

\maketitle


\section{Introduction}

In this work we investigate the global in time regularity properties
of the Yang-Mills equations on high dimensional Minkowski space with compact
semi-simple gauge group $G$. Specifically, we show that if a certain gauge
covariant Sobolev norm is small, the so called \emph{critical} regularity
$\dot{H}_A^\frac{n-4}{2}$, and the dimension satisfies $6\leqslant n$,
then a global solution exists and remains regular for all times given
that the initial data is regular. This is in the same spirit as the recent
result \cite{RT_MKG} for the Maxwell-Klein-Gordon system, as well as
earlier results for high dimensional wave-maps
(see \cite{T_wm1}, \cite{KR_wm},
\cite{SS_wm}, and \cite{Uetal_wm}). Our approach shares many
similarities with those works, whose underlying philosophy in basically
the same. That is, to introduce Coulomb type gauges in order to treat
a specific potential term as a quadratic error. In our setup, we use a
non-abelian variant of the remarkable
parametrix construction contained in \cite{RT_MKG}, in conjunction
with a version of the Uhlenbeck lemma \cite{U_cg} on the existence of
global Coulomb gauges. This latter result has been used for
high dimensional wave-maps to globally ``renormalize'' the equation
so that the existence theory can be treated directly through
Strichartz estimates applied to multi-linear expressions.
In the present situation, as was the case with the
Maxwell-Klein-Gordon system, the corresponding
renormalization procedure is necessarily more involved because
it needs to be done separately for each distinct direction
in phase space. That is, we provide a renormalization of the
Yang-Mills equations through the construction of a Fourier
integral operator with $G$-valued phase. The construction and
estimation of such an object relies heavily on elliptic-Coulomb
theory, primary due to the difficulty one faces in that
the $G$-valued phase function cannot be localized within a
neighborhood of
any given point on the group due to the critical nature of the problem
(if you like, there is a logarithmic ``twisting'' of the group element as
one moves around in physical space; fortunately the group is compact
so this doesn't ruin things).\\

To get things started, we now give a simple gauge covariant
description of the equations we are considering. The (hyperbolic)
Yang-Mills equations arise as the evolution equations for a
connection on the bundle $V=\mathcal{M}^n\times \mathfrak{g}$, where
$\mathcal{M}^n$ is some $n$ (spatial) dimensional Minkowski space,
with metric $g := (-1,1,\ldots,1)$ in inertial coordinates $(x^0 ,
x^i)$, and $\mathfrak{g}$ is the Lie algebra of some compact
semi-simple Lie group $G$. Here we are considering $V$ with the
$Ad(G)$ gauge structure. If $\phi$ is any section to $V$ over
$\mathcal{M}$, then a connection assigns to every vector-field $X$
on the base $\mathcal{M}^n$,  a derivative which we denote as $D_X$,
such that the following Leibniz rule is satisfied for every scalar
field $f$:
\begin{equation}
        D_X( f\phi) \ = \ X(f)\phi + f D_X\phi \ . \notag
\end{equation}
In this setup, we assume that $V$ is equipped with an $Ad(G)$ invariant metric
$\langle \cdot , \cdot \rangle$ which respects the action of $D$.
That is, one has the formula:
\begin{equation}
        d \langle \phi , \psi \rangle \ = \
        \langle D \phi , \psi \rangle + \langle \phi , D \psi \rangle
        \ . \label{connection_compat}
\end{equation}
In the present situation we will take $\langle \cdot , \cdot \rangle$
to be the Killing form on $\mathfrak{g}$.
The curvature associated to $D$ is the $\mathfrak{g}$ valued two-form $F$
which arises from the commutation of
covariant derivatives and is defined via the formula:
\begin{equation}
        D_X D_Y \phi - D_Y D_X \phi - D_{[X,Y]}\phi \ = \
    [ F(X,Y) , \phi] \ . \notag 
\end{equation}
We say that the connection $D$ satisfies the \emph{Yang-Mills}
equations if its curvature is a (formal) local minima of the
following Maxwell type functional:
\begin{equation}
        \mathcal{L}[F] \ = \
    -\frac{1}{4} \ \int_{\mathcal{M}^{n}}
    \langle F_{\alpha\beta} , F^{\alpha\beta}
        \rangle\ DV_{\mathcal{M}^{n}} \ . \label{YM_funct}
\end{equation}
The Euler-Lagrange equations of \eqref{YM_funct}
read:
\begin{equation}
        D^\beta F_{\alpha\beta} \ = \ 0 \ . \label{YM_eq1}
\end{equation}
Also, from the fact that $F$ arises as the curvature of some
connection, we have that the following identity known as ``Bianchi'' is
satisfied:\begin{equation}
        D_{[ \alpha} F_{\beta\gamma]} \ = \ 0 . \label{YM_eq2}
\end{equation}
From now on we will refer to the system \eqref{YM_eq1}--\eqref{YM_eq2}
as the first order Yang--Mills equations (FYM).\\

As we have already mentioned, our aim is to study the regularity properties of
the Cauchy problem for the (FYM) system. To describe this in
a geometrically invariant way, we make use of the
following splitting of the connection-curvature pair $(F,D)$:\ Foliating
$\mathcal{M}$ into the standard Cauchy hypersurfaces $t=const.$, we
decompose:
\begin{equation}
        (F,D) \ = \ (\uF,\uD) \oplus (E , D_0) \ , \notag 
\end{equation}
where $(\uF,\uD)$ denotes the portion of $(F,D)$ which is tangent
to the surfaces $t=const.$ (i.e. the induced connection), and
$(E , D_0)$ denotes respectively the interior product of $F$ with
the foliation generator $T=\partial_t$, and the normal portion
of $D$. In inertial coordinates we have:
\begin{equation}
        E_i \ = \ F_{0i} \ . \notag
\end{equation}
On the initial Cauchy hypersurface $t=0$ we call a set
$(\uF(0),\uD(0) , E(0))$ \emph{admissible Cauchy data}\footnote{Of
course, this set is overdetermined as the curvature $\uF$ depends
completely on the connection $\uD$. Also, it is perhaps not
completely obvious at first that the set $(\uF(0),\uD(0) , E(0))$
determined \emph{uniquely} a solution $(F,D)$ to
\eqref{YM_eq1}--\eqref{YM_eq2}. For example, the initial normal
derivative $D_0(0)$ does not need to be specified. We will show this
is the case in the sequel (in particular see Proposition
\ref{local_coulomb_prop}).} if it satisfies the following compatibility
condition:
\begin{equation}
        \uD^i E_i(0) \ = \ 0 \ . \label{YM_compat_cond}
\end{equation}
We define the \emph{Cauchy problem} for the Yang-Mills equation to
be the task of construction a connection $(F,D)$ which solves \eqref{YM_eq1},
and has Cauchy data equal to $(\uF(0),\uD(0) , E(0))$. \\

In order to understand
what the appropriate condition on the initial data
should be (and what we would like it to be!), it is necessary to
consider the following two basic mathematical features
of the system \eqref{YM_eq1}--\eqref{YM_eq2}. The first is
\emph{conservation}. From the Lagrangian nature of the field
equations  \eqref{YM_eq1}--\eqref{YM_eq2}, we have the tensorial
conservation law:
\begin{align}
        Q_{\alpha\beta}[F] \ &= \ \langle F_{\alpha\gamma} ,
        F_{\beta}^{\ \gamma} \rangle - \frac{1}{4} g_{\alpha\beta}
    \langle F_{\gamma\delta} ,
        F^{\gamma\delta} \rangle \ , \notag\\
    \nabla^\alpha Q_{\alpha\beta}[F] \ &= \ 0 \ , \notag
\end{align}
where $\nabla$ is the covariant derivative on $\mathcal{M}^n$. In
particular, contracting $Q$ with the vector-field $T=\partial_t$,
we arrive at the following constant of motion for the system
\eqref{YM_eq1}--\eqref{YM_eq2}:
\begin{equation}
        \int_{\RR} Q_{00} \, dx \ = \
    \frac{1}{2} \
    \int_{\RR} \big(|E|^2 + |\uF |^2\big)
    \, dx \ . \label{consv_energy}
\end{equation}
The second main aspect of the system \eqref{YM_eq1}--\eqref{YM_eq2} is
that of \emph{scaling}. If we perform the transformation:
\begin{equation}
        (x^0,x^i) \ \rightsquigarrow \ (\lambda x^0,\lambda x^i) \ ,
        \label{spatial_scale_trans}
\end{equation}
on $\mathcal{M}^n$, then an easy calculation shows that:
\begin{align}
        D \ &\rightsquigarrow \ \lambda D \ ,
    	&F \ &\rightsquigarrow \ \lambda^2 F 
	\ . \label{connection_scale_trans}
\end{align}
If we now define the gauge covariant (integer) Sobolev spaces:
\begin{equation}
        \lp{F}{\dot{H}_A^s}^2 \ := \
    \sum_{|I| = s} \ \lp{\underline{D}^I
    F}{L^2(\mathbb{R}^n)}^2        \ , \label{gauge_cov_sob}
\end{equation}
where for each multiindex $I=(i_1,\ldots, i_s)$ we have that
$D^I = D_{\partial_{i_1}}\ldots D_{\partial_{i_s}}$ is the repeated
covariant differentiation with respect to the translation
invariant spatial vector-fields $\{\partial_1 , \ldots , \partial_n\}$,
then for even\footnote{For odd spatial dimensions, the above
discussion needs to be modified somewhat because we will not make an
attempt here to define fractional powers of the spaces $\dot{H}_A^s$.
Instead, what one should do is to simply
put things in a Coulomb gauge and then use the usual fractional
Sobolev spaces. This later approach is what we will take in the
sequel, although for sake of concreteness we will only discuss the
case of even dimensions. We have
opted for the covariant approach in the introduction because it makes
stating our main result a bit easier, and has an appealing
simplicity. Also, since we shall need many specifics
on how Coulomb gauges are constructed in order to create and control our
parametrix, we will explain how the Coulomb gauge relates to the
Cauchy problem in detail in the following two sections.}
spatial dimensions, the norm
$\dot{H}_A^\frac{n-4}{2}$ is invariant with respect to
the scaling transformation \eqref{connection_scale_trans}. In
particular, the conserved quantity \eqref{consv_energy} is invariant
when $n=4$ and this is called the \emph{critical} dimension.\\

Now, based on numerical evidence as well as analytical arguments, it
is suspected that in general the Cauchy problem for
\eqref{YM_eq1}--\eqref{YM_eq2} with smooth initial data will not be
well behaved without size control of the critical regularities $s_c
= \frac{n-4}{2}$ in high dimensions. What we will take this
statement to mean here is simply that if $4 \leqslant  n$ and the
$\dot{H}_A^{s_c}$ norm of the initial data is not sufficiently
small, then one can expect the existence of regular (i.e.
$C_A^\infty$) sets $(\uF(0),\uD(0) , E(0))$ such that the
corresponding solution to \eqref{YM_eq1}--\eqref{YM_eq2} will
develop a singularity in finite time. By singularity development, we
mean  that some higher norm of the type \eqref{gauge_cov_sob} will
fail to be bounded at a later time, given that it was initially; or
even more specifically, that the $L^\infty$ norm of the curvature
$F$ will blow up in finite time for some regular initial data sets.
Since these norms are gauge covariant, this type of singularity
development would correspond to an intrinsic geometric breakdown of
the equations, and could not be an artifact of poorly chosen local
coordinates (gauge) on $V$. This has been rigorously demonstrated in
the equivariant category for the \emph{supercritical} dimensions
$5\leqslant n$ (see \cite{CST_blowup}). In the critical dimension
things are much less clear, although there is numerical evidence
that on still has blowup for large initial data (see
\cite{B_blowup}). This is thought to be connected with the existence
of large static solutions (instantons). One possible conjecture is
that there is global regularity when the  norm \eqref{consv_energy}
is below the ground state energy.\\

Going in the other direction, it is expected that if the critical
norm $\dot{H}_A^\frac{n-4}{2}$ is sufficiently small, then regular
initial data will remain regular for all times. This can be seen as
an easier preliminary step toward understanding in detail the issue
of large data for dimension $n=4$, and is furthermore an interesting
problem in its own right. A central difficulty in the  demonstration
of this conjecture is to construct a stable set coordinates on the
bundle $V$ such that the Christoffel symbols of $D$ are well behaved
in the sense that they obey the natural range of estimates one
expects for this type of problem. This is precisely what we shall do
in dimensions $6\leqslant n$ through the well known procedure of
using (spatial) Coulomb gauges. Unfortunately, this preliminary
gauge construction is far from sufficient to close the regularity
argument, and it will in fact be necessary for us to go much further
and control infinitely many Coulomb gauges, each of which correspond
to a distinct polarized plane wave solution to the usual (flat) wave
equation $\Box = \nabla^\alpha\nabla_\alpha$. However, this does not
effect the statement of our main result which is in fact quite
simple:\\

\begin{thm}[Critical regularity for high dimensional
Yang-Mills]\label{main_Th}
Let the number of spatial dimensions be even and such that
$6\leqslant n$. Then there
exists fixed constants $0 < \varepsilon_0 , C$ such that if
$(\uF(0),\uD(0) , E(0))$ is an admissible data set which satisfies the
smallness condition:
\begin{equation}
        \lp{(\uF(0) , E(0))}{\dot{H}^\frac{n-4}{2}_A } \ \leqslant \
        \varepsilon_0 \ , \label{initial_smallness}
\end{equation}
and there exists constants $M_k < \infty $, $\frac{n-4}{2} <
k \in \mathbb{N}$ such that:
\begin{equation}
        \lp{(\uF(0) , E(0))}{\dot{H}^k_A } \
    = \ M_k \ , \label{main_th_higher_noms}
\end{equation}
then there exists a unique global solution to the field equations
\eqref{YM_eq1}--\eqref{YM_eq2}
with this initial data, and furthermore one has that the
following inductive norm bounds hold:
\begin{align}
        \lp{F}{\dot{H}^\frac{n-4}{2}_A } \ &\leqslant \
    C \varepsilon_0 \ , \notag \\
    \lp{F}{\dot{H}^k_A } \ &\leqslant \
    C(M_\frac{n-4}{2},\ldots,M_{k-1})\, M_k \ . \notag
\end{align}
In particular, in this case $F$ remains smooth (in the gauge
covariant sense) and bounded for all times.
\end{thm}\ret

\begin{rem}
As alluded to above, we will more specifically prove the existence of a global
(in space and time) spatial Coulomb gauge such that the
coefficient functions of the
curvature $F$, as well as the Christoffel symbols (gauge potentials)
of the connection $D$ are in the classical Sobolev spaces $\dot{H}^s$,
and such that they satisfy appropriate angularly and spatially
microlocalized Strichartz estimates. We have elected to
eliminate a discussion of this from the
statement from the main theorem in favor of the simpler geometric language
so that the reader can at a first glance gain an idea of the
content of our result without being confronted with too many technical
details.
\end{rem}\ret


\subsection*{Acknowledgements}
First and foremost, we would like to thank our advisors Sergiu
Klainerman and Matei Machedon for their continuing support
and encouragement. This subject matter as well as our own point
of view owes much to them.
We would also like to thank Igor Rodnianski,
Terry Tao, and Daniel Tataru for many interesting and helpful
conversations. This work began at the Institute for Advanced Study  
during the Fall 2003 semester when both authors were in attendance. 
The second author would like to
thank Harvard University for for its hospitality during the Spring 
of 2004 and Winter 2005. The first author was partially supported
under NSF grant DMS-0401177. 
The second author was supported in part by an
NSF postdoctoral fellowship.\\

\ret

\section{Some Basic Notation}\label{basic_not_sect}

We list here some of the basic conventions used in this work, as
well as some constants which will be needed in the sequel. We use
the usual notation $a \lesssim b$, to denote that $a\leqslant C\cdot
b$ for some (possibly large) constant $C$ which may change from line
to line. Likewise we write $a\ll b$ to mean $a \leqslant C^{-1}\cdot
b$ for some \emph{large} constant $C$. In general, $C$ will denote a
large constant, but at times we will also call $C$ a connection. The
difference should be clear from context. Overall, we will have use
for a family of small constants, which satisfy the hierarchy:
\begin{equation}
    0 \ < \ \varepsilon_0 \ \ll \ \epsilon_0 \ \ll \ \td{\epsilon_0} \
    \ll \ \mathcal{E} \ \ll \ \gamma \ \ll \ \delta \ \ll \ 1 \ .
    \label{constants_line}
\end{equation}

\ret

\section{Some gauge-theoretic preliminaries}

In this paper, we are working with a compact semi-simple group Lie
$G$. However, all of our calculations will be carried out in a
somewhat larger context. Firstly, we will assume that $G$ is
embedded as a subgroup of matrices of some (possibly) larger
orthogonal group $O(m)$. In particular, we can identify the Lie
algebra $\mathfrak{g}$ with an appropriate sub-algebra of
$\mathfrak{o}(m)$. This allows us to perform all of our calculations
on a specific collection of matrices. Since our main computation
involves complex valued integral operators, we will further need to
work in the complexified algebra $\mathbb{C}\otimes
\mathfrak{o}(m)$.  The Killing form $\langle \cdot , \cdot \rangle$
on $\mathfrak{g}$ extends easily to this context to yield the
bilinear form:
\begin{equation}
        \langle A , B \rangle \ = \  \hbox{trace} (AB^*)
        \ . \label{explicit_killing}
\end{equation}
Notice that this is a positive definite form when restricted to the real
vector space $\mathfrak{o}(m)$, and is a sesquilinear form on the
corresponding complexified algebra $\mathbb{C}\otimes \mathfrak{o}(m)$.
More importantly, $\langle \cdot , \cdot \rangle$ is $Ad(O(m))$
invariant, and in fact the more general identity holds:
\begin{equation}
        \langle g_1^{-1} A h_1 , g_2^{-1} B h_2 \rangle \ = \
    \langle g_2 g_1^{-1} A h_1 h_2^{-1} , B \rangle \ , \label{Killing_iden}
\end{equation}
for $A,B\in \mathbb{C}\otimes \mathfrak{o}(m)$ and $g_i , h_i \in
O(m)$. In fact, it is not difficult to see that the form 
\eqref{explicit_killing} extends to a  sesquilinear form on all
complex matrices in $M(m\times m)$, and that it can be identified
with the usual matrix inner product: 
\begin{equation}
    \langle A , B\rangle \ = \ \sum_{i,j}\ a_{ij}\overline{b}_{ij}
    \ , \label{matrix_inner}
\end{equation}
which come from considering these matrices as vectors in
$\mathbb{C}^{m^2}$. Furthermore, it is easy to see that the
general adjoint formula \eqref{Killing_iden} continues to hold in
this context. This will be of fundamental importance in the sequel.
In general, we will use the notation $\lp{A}{}^2 = \langle A , A
\rangle$ to denote the action of this norm on any matrix. Also,
notice that directly from \eqref{Killing_iden} one has the isometric 
identity:
\begin{align}
        \lp{gA}{} \ &= \ \lp{A}{} \ ,
    &g \ &\in \ O(m) \ . \label{matrix_isom}
\end{align}
These are all very simple algebraic identities, but our method is
incredibly sensitive to them and would collapse entirely
if they did not hold.\\

In the context of matrices, we may compute the action of the
connection $D$ on sections $F$ to $V$ as follows:
\begin{equation}
        D_X F \ = \ X^\alpha \big( \nabla_\alpha (F) + [ A_\alpha ,
        F]\big) \ , \notag
\end{equation}
Here, the gauge potentials $\{A_\alpha\}$ are $\mathfrak{g}$-valued,
and are defined via the equation:
\begin{equation}
        D_\alpha \textbf{1}_V \ = \ [A_\alpha ,\textbf{1}_V ] \ , \notag
\end{equation}
where $\textbf{1}_V$ denotes some chosen orthonormal frame
in $V$, and we are abusively writing $F = F \textbf{1}_V$. In
shorthand notation, we write:
\begin{equation}
        D \ = \ d + A \ , \notag 
\end{equation}
where $d$ is the usual exterior derivative on matrix valued
functions. Likewise, in
this notation we have the well known identity for the curvature of
$D$:
\begin{equation}
        F^A \ = \ dA + [A,A] \ . \notag 
\end{equation}
In this last formula, we use the superscript notation to
emphasize the fact that the
curvature is \emph{not} gauge invariant, but transforms according to
the $Ad(G)$ action:
\begin{equation}
        F \ \rightsquigarrow \ g F g^{-1} \ , \notag
\end{equation}
whenever one performs the change of frame
$\textbf{1}_V \rightsquigarrow g \textbf{1}_V g^{-1}$. As is well known,
the potentials $\{A_\alpha\}$ themselves do not transform
according to $Ad(G)$, but instead take on an affine group of
transformations:
\begin{equation}
        B \ = \ gAg^{-1} + g\, dg^{-1} \ , \label{gauge_trans}
\end{equation}
where $\{B_\alpha\}$ represents the connection $D$ in the frame
$g \textbf{1}_V g^{-1}$. In particular, the difference of two
connections obeys the $Ad(G)$ structure, a fact we will have use for
in a moment. For instance, any connection $\{C_\alpha\}$ with $F^C=0$
obeys $Ad(G)$. Furthermore, as is the basic fact of gauge theory, such
connections always lead to a globally\footnote{Of course this ODE is
non-linear, but in the present context it also satisfies the
conservation law $g g^\dagger = I$.} integrable ODE:
\begin{equation}
        dg \ = \ g \, C \ , \notag
\end{equation}
where the solution $g$ belongs to $G$. Thus, we may identify flat
connections $C$ with infinitesimal gauge transformations, and it is
easy to see that every gauge transformation  \eqref{gauge_trans}
leads to a flat connection which we may define as $C = g^{-1} dg$.
This completes our discussion of elementary gauge theory.\\

It will also be necessary for us to make use of the basic facts from
(non-gauge-covariant) Hodge theory. Even though the connections we
work  with in this
paper are on the full space-time $\mathcal{M}^n$, our use of Hodge
theory will always be restricted to  time slices
$\{t\}\times\mathbb{R}^n$.
In particular we use the general notation $d,
d^*$ for the exterior derivative and its adjoint acting on
$\mathfrak{g}$ (and more generally $M(m\times m)$) valued
differential forms on $\mathbb{R}^n$. To emphasize this restriction,
we will use Latin indices when computing these operators. For example:
\begin{align}
        (dA)_{ij} \ &= \ \nabla_{\{ i} A_{j\}} \ ,
        &(dF)_{ijk} \ &= \ \nabla_{[ i} F_{jk]} \ , \notag
\end{align}
where $\{\ldots\}$ and $[\ldots]$ denote anti-symmetric and symmetric
cyclic summing respectively.
Also, the adjoint here is taken with respect to the
Killing form \eqref{explicit_killing}. In particular, we have the Hodge
Laplacean:\begin{equation}
        \Delta \ = \ -( d d^* + d^* d) \ , \label{Hodge_lap}
\end{equation}
which in our context is simply the usual scalar Laplacean acting
component-wise on matrices. Finally, we have the Hodge decomposition
which we write as $A = A^{df} + A^{cf}$ where:
\begin{align}
        A^{df} \ &=\ - d^* d \Delta^{-1} A \ , \notag\\
    A^{cf} \ &=\ - d d^* \, \Delta^{-1} A \ . \notag
\end{align}
This decomposition is bounded on $L^p$ spaces for $1 < p < \infty$ as
the operators involved are SIO's. Also, since these operators are all
real,  this decomposition respects the Lie algebra structure
of $\mathfrak{g}$ inside of $\mathbb{C}\otimes \mathfrak{o}(m)$.\\

The last topic we cover here is the basic underpinning of much of
analysis in the context of compact gauge groups. This is the
remarkable Uhlenbeck lemma, which allows one to ``straighten out'' a
connection as long as its curvature satisfies appropriate bounds. The
important thing for us is that these bounds are precisely at the level
of the critical regularity $\dot{H}_A^\frac{n-4}{2}$. This result is:\\

\begin{lem}[Classical Uhlenbeck lemma]\label{Uhl_lem}
Let $D^A = d + A$ be a connection with compact (matrix) group on
$\mathbb{R}^n$. Then there is a pair of constants 
$\epsilon_0,C$ which only depend on the dimension
$n$ such that if  the curvature $F^A$ of $D^A$ satisfies the
bound:
\begin{equation}
        \lp{F^A}{L^\frac{n}{2}} \ \leqslant \ \epsilon_0 \ , \notag
\end{equation}
then $D^A$ is gauge equivalent to a connection $D^B = d + B$ where the
potentials $\{B_i\}$ satisfy the condition:
\begin{equation}
        d^*B \ = \ 0 , \notag
\end{equation}
and such that the following estimate holds:
\begin{equation}
	\lp{B}{L^n} \ \leqslant \ C\epsilon_0 \ . 
	\label{B_size_bounds}
\end{equation}
\end{lem}\ret

\noindent In the sequel, it will be useful for us to have a somewhat
more refined version of Lemma \ref{Uhl_lem} which does not make
reference to the size of the curvature, but rather to the size of the
connection $\{A_\alpha\}$ itself in a critical norm which does \emph{not}
involve derivatives. This will allow us to prove certain
connections exist more directly. Furthermore, since the basic formulas used
in the proof of this result will be important in
constructing our parametrix, it will set the pace for much of what follows.
Finally, we mention here that our proof is a bit
different from that of \cite{U_cg} in that it
does not rely on any implicit function theorem type
arguments, and is instead completely explicit being based on a
simple Picard iteration.\\

\begin{lem}[Uhlenbeck lemma for small $L^n$ perturbations of Coulomb
potentials with small $L^\frac{n}{2}$ curvature.]\label{A_Uhl_lem}
Let $D^A= d + A$ be a connection on $\mathbb{R}^n\times V$ with compact
(matrix) gauge group $G$. Then there exists constants $\epsilon_0,C$
such that if:
\begin{equation}
	\lp{F^A}{L^\frac{n}{2}} \ \leqslant \ \epsilon_0
	\ , \label{A_Uhl_lem_curv_smallness}
\end{equation}
and such that $d+A$ is gauge equivalent to 
$d+B$ with $d^* B = 0$ \ , where one has the bounds:
\begin{equation}
	\lp{A}{L^n} \ \leqslant \ C\epsilon_0 
	\ , \label{A_size_control}
\end{equation}
then for every connection $d + \td{A}$ such that:
\begin{equation}
        \lp{\td{A} - A}{L^n} \ \leqslant \ 
	\sqrt{C}\epsilon_0 \ ,
    	\label{A_tdA_Ln_smallness}
\end{equation}
there exists a gauge equivalent connection $d+\td{B}$ such that
$d^* \td{B} = 0$, and one has the same size control:
\begin{equation}
	\lp{\td{B}}{L^n} \ \leqslant \ C\epsilon_0 
	\ . \label{tdB_size_control}
\end{equation}
\end{lem}\ret

\begin{rem}
Before continuing with proof, let us remark here that Lemma
\ref{A_Uhl_lem} is in fact more general that the classical Uhlenbeck
Lemma. Specifically, \ref{A_Uhl_lem} easily implies \ref{Uhl_lem}
with smallness condition $\epsilon_0/2$ (where $\epsilon_0$ is
determined by Lemma \ref{A_Uhl_lem}) through a straightforward
induction procedure which we outline now.\\

First of all, from Lemma \ref{A_Uhl_lem} we see that the set of all
connections $d+A$ with curvature such that:
\begin{equation}
        \lp{F^A}{L^\frac{n}{2}} \ \leqslant \ \frac{\epsilon_0}{2} \ ,
        \label{small_curv_bnd}
\end{equation}
and such that $d+A$ is equivalent to $d+B$ with $d^*B=0$, and such that
one has the bounds \eqref{B_size_bounds}, is an open
set in the intersection of $L^n$ with the set determined by
\eqref{small_curv_bnd} (in the sense of distributions). Therefore,
if the conclusion of Lemma \eqref{Uhl_lem} were to be violated, it
must then be the case that there is a smallest number $r^*$ such
that the sphere of radius $r^*$ contains a connection $d+A$ with the
property that it \emph{cannot} be put in the Coulomb gauge
(with $L^n$ bounds), even
though the bound \eqref{small_curv_bnd} is valid for this
connection. Now, consider the set of connections $d+\lambda A$ where
$0 < (1-\lambda)\ll 1$. A quick calculation shows that these have
curvature:
\begin{equation}
        F^{\lambda A} \ = \ \lambda F^A + \lambda(\lambda-1)[A,A] \
        . \notag
\end{equation}
Choose $\lambda$ such that:
\begin{equation}
        (1-\lambda) \ \leqslant \ (1+r^*)^{-2}\cdot \frac{\epsilon_0}{2} \ .
        \notag
\end{equation}
By the triangle and H\"olders inequality, and the definition of
$r^*$, we have that:
\begin{equation}
        \lp{F^{\lambda A}}{L^\frac{n}{2}} \ \leqslant \ \epsilon_0 \
        . \notag
\end{equation}
Therefore, by the minimality of $r^*$ we have that $d+\lambda A$ can
be Coulomb gauged. Again, by the definition of $\lambda$, we have
that:
\begin{equation}
        d+A \ = \ d + \lambda A + \td{A} \ , \notag
\end{equation}
where we easily have the bound (we may assume $1 \leqslant C$):
\begin{equation}
        \lp{\td{A}}{L^n} \ \leqslant \ \sqrt{C}\epsilon_0 \ . \notag
\end{equation}
Therefore, by an application of Lemma \ref{A_Uhl_lem} we have that
$d+A$ can be put in the Coulomb gauge with the \eqref{B_size_bounds} 
holds. This contradicts the
minimality of $r^*$ as was to be shown.
\end{rem}\ret

\begin{proof}[Proof of Lemma \ref{A_Uhl_lem}]
It suffices to show that $d+\td{A}$ is gauge equivalent to $d+\td{B}$,
where $d^* \td{B}=0$ and with the bound \eqref{tdB_size_control}, 
provided that:
\begin{equation}
        \lp{\td{A}}{L^n} \ \leqslant \ 2\sqrt{M}\epsilon_0 
	\ , \label{tdA_Ln_smallness}
\end{equation}
when $\epsilon_0$ chosen sufficiently small, and where $M$ is some 
sufficiently large fixed
constant which will be determined in a moment, and which will be
chosen to be our $C$ in the estimates \eqref{A_size_control} and
\eqref{tdB_size_control} (the reason for the notation switch will
become clear in a moment). 
To see this, notice that the smallness condition
\eqref{A_tdA_Ln_smallness} is gauge invariant because the difference
of two connections transforms according to the $Ad(G)$ action which
fixes the Killing form used to compute $\lp{\cdot}{L^n}$.
Therefore, we may assume from the
start that the original connection $A$ is in the Coulomb gauge with
size control \eqref{A_size_control}. In particular, the connection $d+A$
satisfies the div-curl system:
\begin{align}
        dA \ &= \ F^A - [A,A] \ , &d^* A \ &= \ 0 \ , \notag
\end{align}
which we can integrate to form the equation:
\begin{equation}
	A \ = \ -d^*\Delta (F^A - [A,A]) . \label{FA_coulomb_int}
\end{equation}
Everything we do now will be based on the Riesz operator bounds:
\begin{align}
        \nabla^2 \Delta^{-1} \ : \ L^n \ &\hookrightarrow \
    	L^n \ , \label{SIO_bound}\\
    	\nabla \Delta^{-1} \ : \ L^\frac{n}{2}
    	\ &\hookrightarrow \ L^n \ . \label{fractional_int_bound}
\end{align}
We choose our constant $C$ such that $\frac{\sqrt{M}}{8}$ is the
constant in (the various vector analogs) of the embeddings
\eqref{SIO_bound}--\eqref{fractional_int_bound}. Using these 
bounds and the integral equation \eqref{FA_coulomb_int} in
conjunction with the assumed smallness conditions 
\eqref{A_Uhl_lem_curv_smallness} and \eqref{A_size_control},
and a round of H\"older's inequality, we have the following
improved bounds for $d+A$:
\begin{equation}
        \lp{A}{L^n} \ \leqslant \ \sqrt{M}\epsilon_0 \ ,
\end{equation}
as long as $\epsilon_0$ is chosen sufficiently small.
In particular, using \eqref{A_tdA_Ln_smallness} and some addition
and subtraction we have the bound \eqref{tdA_Ln_smallness}.\\

We now construct by hand the gauge transformation:
\begin{equation}
        dg \ = \ g \td{A} - \td{B} g \ , \label{tdA_gauge_trans}
\end{equation}
with $d^* \td{B}=0$. This will be done by constructing
the infinitesimal gauge transformation $C = g^{-1}dg$. A quick
calculation shows that this must satisfy the following div-curl system:
\begin{subequations}\label{C_dc}
\begin{align}
        dC \ &= \ - \, [C,C] \ , \label{C_dc1}\\
    d^* C \ &= \ d^* \td{A} + [\td{A},C] \ .  \label{C_dc2}
\end{align}
\end{subequations}
Unfortunately, the system \eqref{C_dc} cannot be solved
constructively, say through an iteration scheme. This is because
implicit in its structure is the compatibility condition $d\, [C,C] =
0$, which gets destroyed through (at least the usual) Picard
iteration. This could be side-stepped by using an implicit function
theorem type argument, but since we prefer to do things explicitly we
proceed as follows:\ We first write the system \eqref{C_dc} in terms
of integral equations:
\begin{subequations}\label{C_dc_int}
\begin{align}
        C^{df} \ &= \ \frac{d^*}{\Delta} [C,C] \ , \label{C_dc_int1}\\
    C^{cf} \ &= \ \frac{d}{\Delta} \left( - d^* \td{A} -
    [\td{A},C]\right) \ . \label{C_dc_int2}
\end{align}
\end{subequations}
Here $C = C^{df} + C^{cf}$ denotes the Hodge decomposition of the
matrix valued one-form $C$. A solution to system \eqref{C_dc_int} can
now be constructed from scratch via Picard iteration starting with
$C^{(0)} = 0$. The condition \eqref{tdA_Ln_smallness} and the
embeddings \eqref{SIO_bound}--\eqref{fractional_int_bound}
guarantee convergence to a solution. Furthermore, because it is true
for each iterate, one has the following bounds on the solution:
\begin{equation}
        \lp{C}{L^n} \ \leqslant \ 2\cdot \frac{\sqrt{M}}{8}
	\lp{\td{A}}{L^n} \ \leqslant \ 
        \frac{M}{2} \epsilon_0 \ . \label{C_Ln_smallness}
\end{equation}
Also, since each iterate belongs pointwise to $\mathfrak{g}$, the
solution does also due to the fact that $\mathfrak{g}$ is a linear
(and hence closed) subspace of the matrices $M(m\times m)$. We now
need to show that this $C$ is indeed a solution to the original
system \eqref{C_dc}. That is, we need to establish the identity:
\begin{equation}
        d \, d^* \Delta^{-1} [C,C] \ = \ -\, [C,C] \ . \label{C_brack_iden}
\end{equation}
Notice that this does not follow \emph{algebraically} from the form of the
integral system \eqref{C_dc_int}, because it is not a-priori clear that
in fact $d [C,C] = 0$. However, this is the case, which is a
consequence of the following a-priori \emph{estimate} for solutions to
\eqref{C_dc_int}:
\begin{equation}
        \lp{d\, d^* \Delta^{-1} [C,C] + [C,C]}{L^\frac{n}{2}}
    \ \lesssim \ \lp{C}{L^n}\cdot\lp{d\, d^* \Delta^{-1} [C,C] + [C,C]
    }{L^\frac{n}{2}} \ . \label{C_apriori_Jacobi}
\end{equation}
Notice that \eqref{C_Ln_smallness} and \eqref{C_apriori_Jacobi} taken
together immediately imply the identity \eqref{C_brack_iden}.\\

In order to show \eqref{C_apriori_Jacobi}, we first use the Hodge
Laplacean \eqref{Hodge_lap} to write:
\begin{equation}
        d\, d^*\Delta^{-1} [C,C] + [C,C] \ =\
    - d^*\Delta^{-1}\left(d [C,C] \right) \ . \notag
\end{equation}
Next, we compute that:
\begin{align}
        (d [C,C])_{ijk} \ &= \ \nabla_{[i} [ C_j , C_{k]} ]
        \ , \notag\\
    &= \ [ \nabla_{[i} C_j , C_{k]}   ]  \ - \
    [\nabla_{[i}  C_k , C_{j]}   ] \ , \notag\\
    &= \ -\, [ C_{[i} , (dC)_{jk]} ] \ . \notag
\end{align}
Therefore, using this last identity in conjunction with fractional
integration, and using the identity from line
\eqref{C_dc_int1} above, we have that:
\begin{align}
        \lp{d\, d^* \Delta^{-1} [C,C] + [C,C]}{L^\frac{n}{2}} \ &= \
    \lp{d^*\Delta^{-1} [C , dC] }{L^\frac{n}{2}} \ , \notag\\
    &\lesssim \ \lp{[C, dC] }{L^\frac{n}{3}} \ , \notag\\
    &\leqslant \ \lp{\big[ C, (d d^*\Delta^{-1} [C,C] + [C,C])
    \big] }{L^\frac{n}{3}} \ + \
    \lp{\big[C , [C,C] \big]}{L^\frac{n}{3}} \ , \notag\\
    &\leqslant \ 2 \, \lp{C}{L^n}\cdot \lp{d d^*\Delta^{-1} [C,C] + [C,C]
    }{L^\frac{n}{2}} \ . \notag
\end{align}
Notice that the last inequality here follows simply from the Jacobi identity
$\big[C, [C,C]\big] = 0$.\\

To wrap things up, we only need to establish the existence of $g$ on
line \eqref{tdA_gauge_trans} above with $d^* \td{B}=0$, and
such that we have the size control \eqref{tdB_size_control} (with 
constant $M$). Now, by design we have that $F^C = 0$, so we may 
integrate the equation:
\begin{equation}
        dg \ = \ g C \ , \notag
\end{equation}
with initial conditions $g(0)=I$ on all of $\mathbb{R}^n$.
Defining now:
\begin{equation}
        \td{B} \ = \ g \td{A} g^{-1} + g\, dg^{-1} \ , \notag
\end{equation}
we have that:
\begin{align}
        - \, d^* \td{B} \ &= \ D^{\td{B}}_i \td{B}^i \ , \notag\\
    &= \ g\,  D^{\td{A}}_i ( g^{-1} \td{B}^i g )\, g^{-1} \ ,
    \notag\\
    &= \ g\,  D^{\td{A}}_i (  \td{A}^i - C^i )\, g^{-1} \ ,
    \notag\\
    &= \ g\left( -\, d^* \td{A} + d^*\, C - [\td{A},
    C]\right)g^{-1} \ , \notag\\
    &= 0 , \notag
\end{align}
as was to be shown. Finally, by using the bounds 
\eqref{tdA_Ln_smallness} and
\eqref{C_Ln_smallness} and the definition of the 
potentials $\{\td{B}\}$  and $\{C\}$ we have the bound:
\begin{equation}
	\lp{\td{B}}{L^n} \ \leqslant \ 
	\lp{\td{A}}{L^n} + \lp{C}{L^n} \ \leqslant \
	M\epsilon_0 \ . \notag
\end{equation}
This completes the proof of Lemma \ref{A_Uhl_lem}.
\end{proof}

\ret

\section{Some analytic preliminaries}\label{basic_anal_sect}

We record here some useful formulas, mostly from elementary harmonic
analysis, which will be used many times in the sequel. Firstly, we
define the Fourier transform on $\mathbb{C}\otimes \mathfrak{o}(m)$,
which is merely the usual scalar Fourier transform acting
component-wise on matrices:
\begin{equation}
        \widehat{A}(\xi) \ = \ \int_{\mathbb{R}^n}\ e^{-2\pi i
        x\cdot\xi}\ A(x)\ dx \ . \label{spatial_FT}
\end{equation}
The Plancherel theorem with respect to the Killing form
\eqref{explicit_killing} reads:
\begin{equation}
        \int_{\mathbb{R}^n_x}\ \langle A , B \rangle \ dx \ =
        \int_{\mathbb{R}^n_\xi}\ \langle \widehat{A} , \widehat{B}
        \rangle\ d\xi \ . \notag
\end{equation}
This follows simply from definition of the inner product
\eqref{matrix_inner}. While the constructions we make in the sequel
are almost explicitly based on the spatial transform
\eqref{spatial_FT}, it will in certain places be convenient for us
to work with the space-time Fourier transform:
\begin{equation}
    \widehat{A}(\tau,\xi) \ = \ \int_{\mathbb{R}^{n+1}}\
    e^{-2\pi i(t\tau + x\cdot\xi)}\ A(t,x)\ dt\, dx \ . \notag
\end{equation}\ret

In the sequel, we will have much use for dyadic frequency
decompositions with respect to the spatial variable. For the most
part, we will use a fairly loose and heuristic notation for this
operation. This will help us to avoid having to come up with
different symbols for multipliers which are basically the same.
First of all, we let $\chi(\xi)$ denote some smooth bump function
adapted to the unit frequency annulus $\{2^{-a} \leqslant |\xi|
\leqslant 2^a \}$, where $1 \leqslant a$ is some constant used to define
$\chi$ which may change from line to line. For a dyadic number
$\mu\in \{2^i \ \big| \ i\in\mathbb{Z}\}$, we define the rescaled
cutoffs:
\begin{equation}
        \chi_\mu(\xi) \ = \ \chi(\mu^{-1}\xi) \ , \notag
\end{equation}
and the associated Fourier multipliers $\widehat{P_\mu A} = \chi_\mu
\widehat{A}$. The two main facts we will need about these
multipliers is the Bernstein inequality:
\begin{equation}
        \lp{P_\mu A}{L^p} \ \lesssim \ \mu^{n(\frac{1}{q} -
        \frac{1}{p})} \ \lp{A}{L^q} \ , \label{Bernstein1}
\end{equation}
which holds for all $1\leqslant q \leqslant p \leqslant \infty$, and
the Littlewood-Paley equivalence:
\begin{equation}
        \lp{( \sum_\mu |P_\mu A|^2 )^\frac{1}{2}}{L^p} \ \sim \
        \lp{A}{L^p} \ , \label{LP_thm}
\end{equation}
which holds under the restriction $1 < p < \infty$. All of the norms
above can be taken with respect to \eqref{explicit_killing}.\\

There are two simple analysis lemmas involving derivatives and
multipliers which will come in useful in the sequel. The first of
these is the low frequency (operator) commutator estimate:
\begin{equation}
    \lp{[A,P_{1}]\cdot F}{L^p} \ \lesssim \
    \lp{\nabla_x A}{L^q}\cdot\lp{F}{L^r} \ , \label{Taos_est}
\end{equation}
where $\frac{1}{p} = \frac{1}{q} + \frac{1}{r}$ (see \cite{RT_MKG}).
The second is the homogeneous paraproduct estimate:
\begin{equation}
    \lp{\nabla_x^k(A\cdot F)}{L^p} \
    \lesssim \ \lp{\nabla^k_x A}{L^{q_1}}\cdot\lp{F}{L^{r_1}} +
    \lp{ A}{L^{q_2}}\cdot\lp{\nabla^k_x F}{L^{r_2}} \ ,
    \label{para_est}
\end{equation}
for $1 <  p , q_i , r_i  < \infty$,  $\frac{1}{p}
= \frac{1}{q_1} + \frac{1}{r_1}$, and
$\frac{1}{p} = \frac{1}{q_2} + \frac{1}{r_2}$ whenever
$0 < k$. This estimate is true even for non-integer $0\leqslant
k$ by a simple Littlewood-Paley argument. We note here that we only use it
the integer case, and there it is only employed as a convenience.
For a proof of this, see e.g. Chapter 2 of \cite{T_tools}.\\

We would now like to set up a system to formalize many of the dyadic
estimates which will appear in this paper. This is most easily done
using the language of Besov spaces. Since we have a specific purpose
for these in mind, we introduce the following notation:
\begin{equation}
    \lp{A}{\dot{B}_2^{p,(q,s)}}^2 \ = \ \sum_\mu \
    \mu^{2s - 2n(\frac{1}{q} - \frac{1}{p})}
    \lp{P_\mu A}{L^p}^2 \ , \label{Besov_norm}
\end{equation}
This notation may seem a bit mysterious at first, but the thing to
keep in mind here is that the first index $p$ in some sense controls
the decay, while the second double index $(q,s)$ controls the
scaling, which is the same as $\dot{W}^{s,q}$ (homogeneous $L^q$
Sobolev space). In general, the second index will be fixed, so we
will strive to have $p$ as low as possible (see Remark
\ref{freq_loc_rem} below). This notation has the following simple
significance:\ $\dot{B}_2^{p,(q,s)}$ is the $\ell^2$ Besov space of
Lebesgue index $p$ which \emph{contains} the standard Besov space
$\dot{B}_2^{q,s}$ defined by:
\begin{equation}
    \lp{A}{\dot{B}_2^{q,s}}^2 \ = \ \sum_\mu \
    \mu^{2s}\ \lp{P_\mu A}{L^q}^2 \ . \notag
\end{equation}
This identification is a direct consequence of the Bernstein
embedding \eqref{Bernstein1}. In general, one has the
inclusions:
\begin{align}
    \dot{B}_2^{p_1,(q,s)} \ &\subseteq \
    \dot{B}_2^{p_2,(q,s)} \ ,
    &q\leqslant p_1 \leqslant p_2 \leqslant \infty \ . \label{Besov_nesting}
\end{align}
Furthermore, a quick application of the Littlewood-Paley identity
\eqref{LP_thm} gives the Lebesgue space inclusion:
\begin{align}
    \dot{B}_2^{p,\big(q,n(\frac{1}{q} - \frac{1}{p})\big)}
     \ &\subseteq \ L^p \ ,
    &2\leqslant p  < \infty \ . \label{Lebesgue_besov_incl}
\end{align}
The reason we prefer to use this more involved notation, instead of
the usual Besov norm convention is that ours allows one to tell at
first glance which norms are critical, which is particularly useful
in a scale invariant problem like the one of this paper.
Specifically, the norms $\dot{B}_2^{p,(2,\frac{n-2}{2})}$ will play
a prominent role in what follows.\\

It will also be necessary for us to employ
the $\ell^1$ summing version of the norm
\eqref{Besov_norm}, which we label by $\dot{B}_1^{p,(q,s)}$. This
will essentially be used for one purpose only, and that is that the
$L^\infty$ end-point of \eqref{Lebesgue_besov_incl} is true for this
space:
\begin{align}
    \dot{B}_1^{\infty,(q,\frac{n}{q})} \ &\subseteq \  L^\infty
    \ , &1\leqslant q  \leqslant  \infty \ . \label{Linfty_besov_incl}
\end{align}\\

Besov spaces are particularly well behaved with respect to the
action of Riesz operators, which is exactly why we use them. In
general, we define the operator $|D_x|^{-\sigma}$ to be the Fourier
multiplier with symbol $|\xi|^{-\sigma}$. The basic embedding we
will
use in the sequel is the following:\\

\begin{lem}\label{Besov_lem}
One has the following bilinear estimate for Besov spaces for
$0\leqslant \sigma$:
\begin{equation}
        |D_x|^{-\sigma} \ : \ \dot{B}_2^{p,(2,s_1)}\cdot\dot{B}_2^{q,(2,s_2)}
        \ \hookrightarrow \ \dot{B}_1^{r,(2,s_3)} \ ,
        \label{general_besov_embed}
\end{equation}
where the indices
$1\leqslant p,q,r \leqslant \infty$ and
$\sigma,s_i$ satisfy the following conditions:
\begin{align}
        s_3 \ &= \ s_1 + s_2 + \sigma -\frac{n}{2} \ ,
    &\hbox{(scaling)}\label{sc_cond}\ , \\
        \sigma + \frac{n}{2} - s_3 \ &< \ n(\frac{1}{p} + \frac{1}{q}) \
        , &\hbox{(High $\times$High)}\label{gap_cond}\ , \\
        s_{1} \ &< \ \frac{n}{2} + \min\{n( \frac{1}{q} - \frac{1}{r} )\ ,\ 0\} \ ,
        &\hbox{(Low $\times$High)}\ , \label{pos_cond1}\\
        s_{2} \ &< \ \frac{n}{2} + \min\{n( \frac{1}{p} - \frac{1}{r} )\ ,\ 0\} \ ,
        &\hbox{(High $\times$Low)}\ , \label{pos_cond2}\\
        \frac{1}{r} \ &\leqslant \ \frac{1}{p} + \frac{1}{q} \ ,
        &\hbox{(Lebesgue)} \ . \label{Lb_cond}
\end{align}
\end{lem}\ret

\begin{rem}\label{freq_loc_rem}
As will become apparent in the proof, it is possible to show
frequency localized versions of the embedding
\eqref{general_besov_embed} such that not all of the conditions
\eqref{gap_cond}--\eqref{pos_cond2} need to be satisfied. Indeed, we
will show the following two frequency localized ``improvements'' are
possible:
\begin{align}
        |D_x|^{-\sigma} \ : \ P_{\bullet\ll\lambda}(
        \dot{B}_2^{p,(2,s_1)})\cdot P_\lambda(\dot{B}_2^{q,(2,s_2)})
        \ &\hookrightarrow \ P_\lambda(\dot{B}_1^{r,(2,s_3)}) \ ,
        \label{freq_loc_general_besov_embed1}\\
        |D_x|^{-\sigma} \ : \ P_{\lambda}(
        \dot{B}_2^{p,(2,s_1)})\cdot P_\lambda(\dot{B}_2^{q,(2,s_2)})
        \ &\hookrightarrow \ \left(\frac{\mu}{\lambda}
        \right)^\delta P_\mu(\dot{B}_1^{r,(2,s_3)}) \ ,
        \label{freq_loc_general_besov_embed2}
\end{align}
where $\delta=n(\frac{1}{p} + \frac{1}{q}) + s_3 - \sigma -
\frac{n}{2}$ in estimate \eqref{freq_loc_general_besov_embed2}.
Estimate \eqref{freq_loc_general_besov_embed1} holds whenever
\eqref{sc_cond}, \eqref{pos_cond1}, and \eqref{Lb_cond} are
satisfied. The second estimate \eqref{freq_loc_general_besov_embed2}
is valid whenever we have \eqref{sc_cond}, \eqref{gap_cond}, and
\eqref{Lb_cond}. In particular, notice that for larger $\sigma$ this
estimate requires lower values of $p,q$. This fact will have an
immense bearing on the estimates we prove in the sequel, and seems
to be one of the most difficult factors in lowering the dimension of
the overall argument from $n=6$ (apart from even more difficult
things such as null-form estimates).
\end{rem}\ret

\begin{proof}[Proof of estimate \eqref{general_besov_embed}]
The proof is a simple matter of the standard technique of
trichotomy. That is, we start with two test matrices $A$ and $C$,
and we run a frequency decomposition on the product:
\begin{equation}
        A\cdot C \ = \ \sum_{\lambda,\mu_i}\
        P_\lambda \left(P_{\mu_1}A\cdot P_{\mu_2}C\right) \ . \notag
\end{equation}
Setting now:
\begin{equation}
        \gamma \ = \ \min\Big\{\ \frac{n}{2}-s_1 \ , \ \frac{n}{2}-s_2 \ , \
        n(\frac{1}{p} + \frac{1}{q}) + s_3 -
        \sigma - \frac{n}{2}  \ , \
        \frac{n}{2} + n( \frac{1}{q} - \frac{1}{r}) -
        s_1 \ , \
        \frac{n}{2} + n( \frac{1}{p} - \frac{1}{r} ) -
        s_2 \ \Big\} \ , \notag
\end{equation}
we have from the conditions \eqref{gap_cond}--\eqref{pos_cond2}
that $0 < \gamma$. To prove \eqref{general_besov_embed}
it suffices to show that:
\begin{align}
   \begin{split}
        &\sum_{\substack{\mu_1 \ : \\ \mu_1\ll\mu_2\\
        \lambda \sim \mu_2}}\
        \lambda^{s_3 - n(\frac{1}{2} - \frac{1}{r})-\sigma}\
        \lp{P_\lambda \left(P_{\mu_1}A\cdot P_{\mu_2}C\right)}{L^r}
        \ \lesssim \notag\\
    &\hspace{1.75in}
    \sum_{\substack{\mu_1 \ : \\ \mu_1\ll\mu_2\\
        \lambda \sim \mu_2}}\
    \left(\frac{\mu_1}{\mu_2}\right)^\gamma\
    \lp{P_{\mu_1} A}{\dot{B}^{p,(2,s_1)}}\cdot
        \lp{P_{\mu_2}C}{\dot{B}^{q,(2,s_2)}} \ ,
   \end{split} \notag\\
   \begin{split}
        &\sum_{\substack{\mu_2 \ : \\ \mu_2\ll\mu_1\\
        \lambda \sim \mu_1}}\
        \lambda^{s_3 - n(\frac{1}{2} - \frac{1}{r})-\sigma}\
        \lp{P_\lambda \left(P_{\mu_1}A\cdot P_{\mu_2}C\right)}{L^r}
        \ \lesssim \notag\\
    &\hspace{1.75in}
    \sum_{\substack{\mu_2 \ : \\ \mu_2\ll\mu_1\\
        \lambda \sim \mu_1}}\ \left(\frac{\mu_2}{\mu_1}\right)^\gamma\
        \lp{P_{\mu_1}A}{\dot{B}^{p,(2,s_1)}}\cdot
        \lp{P_{\mu_2} C}{\dot{C}^{q,(2,s_2)}} \ ,
   \end{split} \notag\\
        &\sum_{\substack{\lambda \ : \\ \mu_2\sim\mu_1\\
        \lambda \lesssim  \mu_i}}\
        \lambda^{s_3 - n(\frac{1}{2} - \frac{1}{r})-\sigma}\
        \lp{P_\lambda \left(P_{\mu_1}A\cdot P_{\mu_2}C\right)}{L^r}
        \ \lesssim \
        \lp{P_{\mu_1}A}{\dot{B}^{p,(2,s_1)}}\cdot
        \lp{P_{\mu_2}C}{\dot{B}^{q,(2,s_2)}} \ ,
        \notag
\end{align}
That \eqref{general_besov_embed} follows from these three estimates
is a simple consequence of Young's inequality and Cauchy-Schwartz
respectively. These estimates, in turn, are all a consequence of the
single fixed frequency bound:
\begin{multline}
        \lambda^{s_3 - n(\frac{1}{2} - \frac{1}{r})-\sigma}\
        \lp{P_\lambda \left(P_{\mu_1}A\cdot P_{\mu_2}C\right)}{L^r}
        \ \lesssim \\
        \left(\frac{\lambda}{\max\{\mu_i\}}\right)^\gamma\cdot
        \min\{ \left(\frac{\mu_1}{\mu_2}\right)^\gamma
        , \left(\frac{\mu_2}{\mu_1}\right)^\gamma\}\cdot
        \lp{P_{\mu_1}A}{\dot{B}^{p,(2,s_1)}}\cdot
        \lp{P_{\mu_2}C}{\dot{B}^{q,(2,s_2)}} \ .
        \label{gen_besov_fixed_freq}
\end{multline}
The proof of \eqref{gen_besov_fixed_freq} is a simple matter of
H\"olders and Bernstein's inequalities, and counting weights. There
are three cases corresponding to the three summing estimates above.
In the first case, we assume that $\lambda \ \lesssim \
\mu_1\sim\mu_2$. Since \eqref{gen_besov_fixed_freq} is scale
invariant, we may assume in this case that both $\mu_i\sim 1$. Using
now H\"olders inequality which is permissible by \eqref{Lb_cond},
followed by the Bernstein inequality, we have that:
\begin{equation}
        \lp{P_\lambda \left(P_{\mu_1}A\cdot P_{\mu_2}C\right)}{L^r}
        \ \lesssim \ \lambda^{n(\frac{1}{p} + \frac{1}{q} - \frac{1}{r})}
        \ \lp{P_{\mu_1}A}{L^p}\cdot\lp{P_{\mu_2}C}{L^q} \ . \notag
\end{equation}
Multiplying this last estimate by the weight $\lambda^{s_3 -
n(\frac{1}{2} - \frac{1}{r})-\sigma}$  we arrive at the bound:
\begin{equation}
        \hbox{(L.H.S.)}\eqref{gen_besov_fixed_freq} \ \lesssim \
        \lambda^{n(\frac{1}{p} + \frac{1}{q}) +s_3 -\sigma
        -\frac{n}{2}}\ \lp{P_{\mu_1}A}{L^p}\cdot\lp{P_{\mu_2}C}{L^q} \ . \notag
\end{equation}
Then \eqref{gen_besov_fixed_freq} follows in this case from the
definition of $\gamma$ and the fact that $\mu_i \sim 1$. The other
two cases, which correspond to $\mu_1 \ll \mu_2$ or vice versa are
similar, so it suffices to consider the first. In this case we
rescale to $\mu_2\sim \lambda \sim 1$. In the case where $r < q$
we set $\frac{1}{\td{p}} = \frac{1}{r} - \frac{1}{q}$,
and we again use H\"older and Bernstein to estimate:
\begin{equation}
        \lp{P_\lambda \left(P_{\mu_1}A\cdot P_{\mu_2}C\right)}{L^r}
        \ \lesssim \ \mu_1^{n(\frac{1}{p} - \frac{1}{\td{p}})}
        \ \lp{P_{\mu_1}A}{L^p}\cdot\lp{P_{\mu_2}C}{L^q} \ . \notag
\end{equation}
If it is the case that $q\leqslant r$, then we simply estimate:
\begin{equation}
    \lp{P_\lambda \left(P_{\mu_1}A\cdot P_{\mu_2}C\right)}{L^r}
        \ \lesssim \ \mu_1^{\frac{n}{p}}
        \ \lp{P_{\mu_1}A}{L^p}\cdot\lp{P_{\mu_2}C}{L^q} \ . \notag
\end{equation}
In either case, the claim \eqref{gen_besov_fixed_freq} follows
from the definition of $\gamma$. This completes the proof of
\eqref{general_besov_embed}.
\end{proof}\ret

Before continuing on, let us note here a slight refinement of the
Besov norms \eqref{Besov_norm} and the embedding
\eqref{general_besov_embed}. This involves taking into account functions
which live at frequency $\lesssim 1$. If we let $\langle D_x \rangle$
denote the multiplier with symbol $( 1 + |\xi|^2)^\frac{1}{2}$,
then we form the low frequency spaces:
\begin{equation}
    \lp{A}{\dot{B}_{2,10n}^{p,(q,s)}} \ = \
    \lp{\langle D_x\rangle^{10n} A}{\dot{B}_{2}^{p,(q,s)}} \ , \label{exp_Besov}
\end{equation}
with a similar definition for the $\ell^1$ version
$\dot{B}_{1,10n}^{p,(q,s)}$. By a straightforward adaptation of the previous
argument, it is easy to see that the embedding \eqref{general_besov_embed}
is equally valid for these  low frequency spaces. We leave the details
to the reader.\\

It will also be necessary for us to perform various dyadic
decompositions with respect to the angular frequency variable. For
each fixed direction $\omega$ in the frequency plane
$\mathbb{R}^n_\xi$, we decompose the unit sphere
$\mathbb{S}^{n-1}_\xi$ into dyadic conical regions:
\begin{equation}
        \mathcal{R}(\omega,\theta) \ = \ \{\eta\in\mathbb{S}^{n-1}_\xi
        \ \big| \ \angle(\omega,\eta) \sim \theta\} \ ,
        \label{angular_region}
\end{equation}
where $\theta\in \{\frac{\pi}{2}\cdot 2^i \ \big| \ i \in \mathbb{Z}
, i \leqslant 0\}$. Here we will not bother to fix the constant in
the $\sim$ notation used to define the regions
\eqref{angular_region}, but we will let it change from line to line
as we have done for the spatial multipliers above. We also define a
smooth partition of unity adapted to these regions, which we label
by $b^\omega_\theta$. These can always be chosen (e.g. by defining
them on a larger sphere and then rescaling) so that they satisfy the
differential bounds:
\begin{align}
        |(\omega\cdot\nabla_\xi)^k_\omega\, p_1
    b^\omega_\theta| \ &\lesssim \ 1 \ ,
    &|(\omega^\perp\cdot\nabla_\xi)^k \,
    p_1 b^\omega_\theta| \ &\lesssim \
        \theta^{-k} \ , \notag
\end{align}
where the implicit constants depend on $k$ but are uniform in
$\theta$. In particular, if we define the multipliers
$\widehat{\oPi_\theta A} = b^\omega_\theta \widehat{A}$, then the
operators $\oPi_\theta P_\mu$ are bounded on all $L^p$ spaces
uniformly in $\mu$ and $\theta$. In fact, the following refinement
of the inequality \eqref{Bernstein1} holds, which we also call
Bernstein:
\begin{equation}
        \lp{\oPi_\theta P_\mu A}{L^p} \ \lesssim \ \mu^{n(\frac{1}{q} -
        \frac{1}{p})}\theta^{(n-1)(\frac{1}{q} -
        \frac{1}{p})} \ \lp{A}{L^q} \ . \label{Bernstein2}
\end{equation}
In all of the above inequalities, we have kept $\omega$ as a fixed
directional value. However, it will also be necessary for us to have
an account of how our multipliers depend on this parameter. In
particular, we will need to have bounds for the operators
$\nabla_\omega \oPi_\theta$. This is easily achieved by
differentiating the associated multiplier. In fact, one has the
bounds for fixed $\xi$:
\begin{equation}
        |\nabla_\omega^k b^\omega_\theta| \ \lesssim \ \theta^{-k} \
        . \label{ang_symbol_bnds}
\end{equation}
The way we shall express this bound in calculations is through the
following heuristic operator identity:
\begin{equation}
        \nabla_\omega^k \oPi_\theta \ \approx \ \theta^{-k}
        \oPi_\theta \ , \label{heuristic_op_bnd}
\end{equation}
which we shall take to mean that \emph{the left hand side satisfies
all $L^p$ space bounds as the right hand side}. Notice that this
relation has a preferred direction (left$\Rightarrow$right).
In practice, this means that we
have the bound \eqref{Bernstein2} for the operator on the left hand
side of \eqref{heuristic_op_bnd} with the added factor of
$\theta^{-k}$.\\

Finally, let us end this section by making the following
conventions. Firstly, it will be convenient for us at times to write
$P_\mu A = A_\mu$ for a localized object. This should not be
confused with the $\mu^{th}$ component of $A$ in the case that it is
a one-form. This should usually be clear from context. Secondly, it
will be necessary for us to ensure that certain of our multipliers
have real symbol so that they respect the subalgebra
$\mathfrak{g}(m)\subseteq M(m\times m)$. This will be done by taking
their real part which simply symmetrizes their (real) symbols. In
particular, we will denote this by:
\begin{equation}
    \Re (\oPi_\theta) \ = \ {\ooPi}_\theta \ . \notag
\end{equation}
Secondly, we
use the following bulleted notation for the sum of various cutoffs
over a given range:
\begin{align}
        P_{\bullet < c} \ &= \ \sum_{\mu < c} P_\mu \ ,
        &\oPi_{\bullet < c} \ &= \ \sum_{\theta < c} \oPi_\theta \ ,
        \notag
\end{align}
etc. We will also use the notation $A_{\bullet < c}$ etc. for these
operators applied to tensors. Finally, we will set aside a special
notation here for cutting off on angles sectors whose width depends
on the frequency:
\begin{equation}
        \ooPi^{ (\sigma) } \ = \ \sum_\mu \ \ooPi_{\mu^\sigma < \bullet}
        \, P_\mu \ . \label{variable_angular_cutoff}
\end{equation}
Notice that this multiplier does not satisfy good bounds of the form
\eqref{ang_symbol_bnds}. However, it can be dealt with using the
Littlewood-Paley equivalence \eqref{LP_thm} if there is a little
extra room left to sum over fixed angular dyadics. This ends our
description of the basic analysis we will use in this paper.

\ret

\section{Gauge construction for the initial data; Reduction to a second
order system and the main a-priori estimate}

We now begin our proof of the main theorem \ref{main_Th}. As we have
already mentioned, one of the central components of the proof is to
construct a stable set of ``elliptic  coordinates'' on the bundle
$V$. The way we will do this is to construct the desired frame
on the $t=0$ slice $\RR\times \mathfrak{g}$. We will
then show that this frame propagates as the system evolves by
solving an auxiliary set of equations for the gauge potentials which
respects the chosen frame automatically. The regularity of this system
of equations will be provided in the usual translation 
invariant Sobolev spaces.
We then show that our auxiliary solution is in fact a true solution to the
system of equations \eqref{YM_eq1}--\eqref{YM_eq2} by employing a
bootstrapping procedure which is similar to that used in the proof
of Lemma \ref{A_Uhl_lem}. The desired gauge covariant regularity, which is
contained in the statement of Theorem \ref{main_Th}, will be provided
by a comparison principle. These constructions are all local in
time and are more or less standard. We have included them here for the
convenience of the reader, the sake of completeness, and the fact
that some of the formulas we develop along the way will be
central to what we do in later sections.\\

With the local theory established, the global conclusion of Theorem
\ref{main_Th} will then be a consequence of a certain a-priori
estimate on the (usual Sobolev) energy of solutions to
\eqref{YM_eq1}--\eqref{YM_eq2} in the gauge we construct. Our task
will then be to show that this a-priori estimate is true for all
solutions to yet another system of auxiliary equations, this time for
the curvature. This can be considered to be the main 
estimate of the paper. The proof turns out to be
quite involved, and will occupy the rest of the paper. In the next
section, we will prove the main a-priori estimate itself with the
help of a certain family of microlocalized space-time (Strichartz)
estimates for solutions to second order covariant wave equations on
bundles with connections satisfying estimates consistent with our
bootstrapping assumptions. The breakdown here is based on the
Smith-Tataru (see \cite{ST_quasi}) $\mathcal{E}$-parametrix
idea, which allows one to reduce the needed Strichartz estimates
to proving them for a suitable family of approximate frequency localized
fundamental solutions. Our rendition of this is essentially equivalent
to that contained
in the paper \cite{RT_MKG}.\\

Finally, in the remaining sections of
the paper we develop the linear theory. This is by far the most
involved portion of the present work, and requires the construction
of some fairly sophisticated oscillatory integrals and microlocal
function  spaces. This material can be read without reference to the
non-linear problem, as long as one is familiar with the algebraic
and analytic assumptions we make on the geometry (frequency localized
connection). While these come  from the non-linear problem, they are
of course a bit more general.\\


\subsection{Construction of the initial frame, and the comparison
principle}

The first thing we do here is to put the initial connection
$\underline{D}$ into the Coulomb gauge. Via the Uhlenbeck lemma
\eqref{Uhl_lem}, we simply need to show that:
\begin{equation}
        \lp{\underline{F}}{L^\frac{n}{2}} \ \lesssim 
	\  \varepsilon_0 \
        , \notag 
\end{equation}
for $\varepsilon_0$ the sufficiently small parameter
from line \eqref{initial_smallness} (which should not be confused
with the small constant from Lemma \ref{Uhl_lem} above). This $L^p$
bound follows immediately from the gauge covariant Sobolev embedding
(for $n$ even):
\begin{equation}
        \dot{H}^\frac{n-4}{2}_A \ \subseteq \ L^\frac{n}{2} \
        , \notag 
\end{equation}
which in turn follows from repeated application of the usual 
single derivative Sobolev
embeddings and the Kato estimate (which follows immediately
from \eqref{connection_compat} and Cauchy-Schwatrz):
\begin{equation}
        \big| d |F| \big| \ \leqslant \ |\underline{D}F| \ , \label{Kato_est}
\end{equation}
where $F$ is any section to $\mathcal{M}\times \mathfrak{g}$ and the
absolute norm $|\cdot |$ is taken with respect to the Killing inner
product \eqref{explicit_killing}.\\

We may now assume that we are dealing with an initial data set:
\begin{equation}
        (\uF(0),\uD(0) , E(0)) \ , \label{the_data}
\end{equation}
for the system which is such that connection $\uD(0) = d +
\underline{A}(0)$ satisfied the elliptic div-curl system:
\begin{align}
        d\underline{A}(0) + [\underline{A}(0),
        \underline{A}(0)] \ &= \ \uF(0) \ ,
        &d^*\underline{A}(0) \ = \ 0 \ , \label{initial_div_curl}
\end{align}
and such that the compatibility condition \eqref{YM_compat_cond}
is satisfied. Furthermore, from \eqref{B_size_bounds}
we have the bounds:
\begin{equation}
        \lp{\underline{A}(0)}{L^n} \ \lesssim  \ \varepsilon_0 \ .
        \notag
\end{equation}
We will now use this last bound to show that the initial data set
\eqref{the_data} is in fact in the classical Sobolev spaces
$\dot{H}^k$. This is a consequence of the following:\\

\begin{lem}[Comparison principle for Sobolev norms on $\RR$]\label{comp_lem}
Let $\underline{D}
=d+\underline{A}$ be a connection on $\mathbb{R}^n$, with $n$
even, such that one has the potential and curvature bounds:
\begin{align}
        \lp{\underline{A}}{L^n} \ , \ 
	\lp{\uF}{\dot{H}^\frac{n-4}{2}_A} \
        &\leqslant \ \epsilon_0 \ , \label{comp_lem_crit}\\
        \lp{\uF}{\dot{H}^k_A} \ &\leqslant \ M_k \ ,
        \label{comp_lem_large}
\end{align}
for $\frac{n-4}{2} < k$. Suppose also that $\underline{D}$ 
is in the gauge
$d^*\underline{A}=0$. Then we have the critical classical Sobolev
bounds:
\begin{align}
        \lp{\uF}{\dot{H}^\frac{n-4}{2}} \ &\leqslant \
        C\epsilon_0 \ , \label{F_cSob_bnds}\\
        \lp{\underline{A}}{\dot{H}^\frac{n-2}{2}} \ &\leqslant \
        C\epsilon_0 \ .
        \label{A_cSob_bnds}
\end{align}
Furthermore, if $G$ is any $\mathfrak{g}$ valued function, then we
have the following inductive comparison of norms:
\begin{align}
     \begin{split}
    C^{-1}(M_{\frac{n-4}{2}},\ldots,M_{k-1})\
        \lp{G}{H^{[k^*,k]}} \ \ \  &\leqslant \ \ \
        \lp{G}{H_A^{[k^*,k]}} \ ,  \\
    &\leqslant \ \ \
        C(M_{\frac{n-4}{2}},\ldots,M_{k-1})\
        \lp{G}{H^{[k^*,k]}} \ ,
     \end{split}\label{general_comp}
\end{align}
where the index $k^*$ is such that $\frac{n-4}{2} \leqslant k^* <
n$, and where we have set:
\begin{equation}
    \lp{G}{H_A^{[k^*,k]}}^2 \ = \ \sum_{k^* \leqslant m \leqslant k}
    \ \lp{\underline{D}^m G}{L^2}^2 \ , \notag
\end{equation}
to be the interval gauge-covariant Sobolev space. We use an analogous
definition for the space $H^{[k^*,k]}$. We also have the non-inductive
equivalence between $\nabla_x\underline{A}$ and $\uF$:
\begin{equation}
        N^{-1}_k\, \lp{\underline{A}}{\dot{H}^k} \ \leqslant \
        \lp{\uF}{\dot{H}^{k-1}} \ \leqslant \ N_k\,
        \lp{\underline{A}}{\dot{H}^k} \ , \label{AF_equiv}
\end{equation}
where $N_k$, $\frac{n-2}{2} \leqslant k$, is a set of constants
which depends only on the dimension and \emph{not} on the constant
$\epsilon_0$ once it is sufficiently small. In particular, combining
all of this, we have the following classical Sobolev bounds on the
pair $(\underline{A},\uF)$:
\begin{align}
        \lp{\uF}{\dot{H}^k} \ &\leqslant \
        C(M_{\frac{n-4}{2}},\ldots,M_{k-1})\, M_k
        \ , \label{higher_F_cSob_bnds}\\
        \lp{\underline{A}}{\dot{H}^{k+1}} \ &\leqslant \
        C(M_{\frac{n-4}{2}},\ldots,M_{k-1})\, M_k \ .
        \label{higher_A_cSob_bnds}
\end{align}
for $\frac{n-4}{2} < k$.
\end{lem}\ret

\begin{proof}[Proof of Lemma \ref{comp_lem}]
The proof will be accomplished via a series of inductions.
In what follows, we will assume the estimate \eqref{AF_equiv},
whose proof follows from simple analysis of the elliptic
system \eqref{initial_div_curl} in Besov spaces of the
kind $\dot{B}_2^{p,(2,s)}$. We will perform many reductions like
this in the sequel so we leave this one to the reader.\\

The first step is to prove the critical classical Sobolev
\eqref{F_cSob_bnds}. Note that the potential bounds
\eqref{A_cSob_bnds} follow from this and \eqref{AF_equiv}.
The inductive hypothesis that we make here is that:
\begin{equation}
    \lp{\nabla_x^l \underline{D}^m\uF}{L^\frac{n}{k}} \ \lesssim\
    \epsilon_0 \ , \label{first_indct_est}
\end{equation}
for $k=l+m +2\leqslant \frac{n}{2}$
whenever $0\leqslant l\leqslant l_0$. Notice that this
hypothesis is verified for $l_0=0$ on account of the assumption
\eqref{comp_lem_crit} and by applying the Kato estimate
\eqref{Kato_est} in conjunction with integer Sobolev embeddings.
Notice also that by applying Riesz operator estimates to the elliptic
system \eqref{initial_div_curl}, and using the product
estimate \eqref{para_est} along with Sobolev embeddings
we have the bounds:
\begin{align}
        \lp{\nabla_x^{l+1} \underline{A}}{L^\frac{n}{k}} \ &\lesssim \
        \lp{\nabla^l_x \uF}{L^\frac{n}{k}} +
        \lp{\nabla^{l}_x ([\underline{A},\underline{A}])}{L^\frac{n}{k}}
        \ , \notag\\
        &\lesssim \ \lp{\nabla^l_x \uF}{L^\frac{n}{k}} +
        \lp{\nabla^{l}_x \underline{A}}{L^\frac{n}{k-1}}
        \cdot \lp{\underline{A}}{L^n} \ , \notag\\
        &\lesssim \ \lp{\nabla^l_x \uF}{L^\frac{n}{k}} + \epsilon_0\cdot
        \lp{\nabla_x^{l+1} \underline{A}}{L^\frac{n}{k}} \ . \notag
\end{align}
Therefore, the inductive hypothesis \eqref{first_indct_est}
may be assumed to also contain the estimate:
\begin{equation}
    \lp{\nabla^{l+1}_x\underline{A}}{L^\frac{n}{k}}
    \ \lesssim \ \epsilon_0 \ , \label{first_indct_A}
\end{equation}
for $k=l+2\leqslant \frac{n}{2}$ and $l\leqslant l_0$. To show that \eqref{first_indct_est}
holds for all $l\leqslant l_0+1$, we start with $l\leqslant l_0$ and
we compute using \eqref{para_est} and Sobolev embeddings that:
\begin{align}
    &\lp{\nabla_x^{l+1} \underline{D}^{m-1} \uF}{L^\frac{n}{k}}
    \ , \notag\\
    \lesssim \
    &\lp{\nabla_x^l \underline{D}^m \uF}{L^\frac{n}{k}} +
    \lp{\nabla_x^l ([\underline{A},\underline{D}^{m-1} \uF])
    }{L^\frac{n}{k}} \ , \notag\\
    \lesssim \ &\epsilon_0 + \lp{\nabla_x^l \underline{A}}
    {L^\frac{n}{l+1}}\cdot\lp{\underline{D}^{m-1}
    \uF}{L^\frac{n}{k-l-1}} + \lp{\underline{A}}
    {L^n}\cdot
    \lp{\nabla_x^l \underline{D}^{m-1} \uF}{L^\frac{n}{k-1}}
    \ , \notag\\
    \lesssim \ &\epsilon_0 + \epsilon_0\cdot
    \lp{\nabla_x^{l+1} \underline{D}^{m-1} \uF}{L^\frac{n}{k}} \ . \notag
\end{align}
This inductively establishes \eqref{first_indct_est} and hence proves
\eqref{F_cSob_bnds}.\\

We now show \eqref{general_comp}. We first deal with the leftmost
inequality. Our inductive hypothesis here is that:
\begin{equation}
    \lp{\nabla_x^l \underline{D}^m G}{L^2}
    \ \lesssim \ C(M_{\frac{n-4}{2}},\ldots,M_{k-1})\
        \lp{G}{H_A^{[k^*,k]}} \ , \label{general_indct_hyp}
\end{equation}
where $l+m = k_0$ for $k_0=k$ or $k_0=k^*$, and for all
$l\leqslant l_0$. To compute $\nabla_x^{l+1} \underline{D}^{m-1} G$
in terms of this, we need to split into cases depending on whether
or not $l + 1 < \frac{n}{2}$. In the former case we compute that:
\begin{align}
    &\lp{\nabla_x^{l+1} \underline{D}^{m-1} G}{L^2} \ , \notag\\
    \lesssim \ &\lp{\nabla_x^{l} \underline{D}^{m} G}{L^2}
    + \lp{\nabla_x^{l}([\underline{A},
    \underline{D}^{m-1} G])}{L^2} \ , \notag\\
    \begin{split}
    \lesssim \ &C(M_{\frac{n-4}{2}},\ldots,M_{k-1})\
        \lp{G}{H_A^{[k^*,k]}} +
    \lp{\nabla_x^l\underline{A}}{L^\frac{n}{l+1}}
    \cdot\lp{\underline{D}^{m-1} G}{L^\frac{2n}{n-2l-2}}\\
    &\ \ \ \ \ \ \ \ \ \ \ \ + \lp{\underline{A}}{L^n}
    \cdot\lp{\nabla_x^l \underline{D}^{m-1} G}
    {L^\frac{2n}{n-2}} \ ,
     \end{split}\notag\\
    \lesssim \ &C(M_{\frac{n-4}{2}},\ldots,M_{k-1})\
        \lp{G}{H_A^{[k^*,k]}}
    + \epsilon_0 \cdot\lp{\nabla_x^{l+1} \underline{D}^{m-1} G}{L^2}
    \ . \notag
\end{align}
In the case where $\frac{n}{2}-1 \leqslant l$ we have the inequality:
\begin{align}
    &\lp{\nabla_x^{l+1} \underline{D}^{m-1} G}{L^2} \ , \notag\\
    \begin{split}
    \lesssim \ &C(M_{\frac{n-4}{2}},\ldots,M_{k-1})\
        \lp{G}{H_A^{[k^*,k]}} +
    \lp{\nabla_x^l\underline{A}}{L^\frac{2n}{n-2}}
    \cdot\lp{\underline{D}^{m-1} G}{L^n}\\
    &\ \ \ \ \ \ \ \ \ \ \ \ + \lp{\underline{A}}{L^n}
    \cdot\lp{\nabla_x^l \underline{D}^{m-1} G}
    {L^\frac{2n}{n-2}} \ ,
     \end{split}\notag\\
     \begin{split}
    \lesssim \ &C(M_{\frac{n-4}{2}},\ldots,M_{k-1})\
        \lp{G}{H_A^{[k^*,k]}}
    + \lp{\nabla_x^{l+1}\underline{A}}{L^2}
    \cdot\lp{\underline{D}^{\frac{n-2}{2} + m-1} G}{L^2}\\
    &\ \ \ \ \ \ \ \ \ \ \ \
    + \epsilon_0 \cdot\lp{\nabla_x^{l+1}
    \underline{D}^{m-1} G}{L^2} \ .
     \end{split}\notag
\end{align}
Notice that this last line above used the $L^2\hookrightarrow L^n$
gauge covariant Sobolev embedding.
To bound the second term on this  line, notice that since
$\frac{n}{2}-1 \leqslant l$ and we must assume that $1\leqslant m$ for
the induction to make sense, we have the bound
$k^*\leqslant \frac{n-2}{2} + m-1\leqslant k$. This allows us to bound:
\begin{equation}
    \lp{\underline{D}^{\frac{n-2}{2} + m-1} G}{L^2}
     \ \leqslant \ \lp{G}{H_A^{[k^*,k]}} \ . \notag
\end{equation}
Furthermore, by placing all of these
calculations within an induction on the value of $k$ itself,
and using the bound \eqref{AF_equiv} while
noting that $l\leqslant k-1$ we may assume the bound:
\begin{equation}
    \lp{\nabla_x^{l+1}\underline{A}}{L^2} \ \lesssim \
    \lp{\nabla^l_x \uF}{L^2} \ \lesssim \
    C(M_{\frac{n-4}{2}},\ldots,M_{k-1}) \ . \notag
\end{equation}
This completes our inductive proof of \eqref{general_indct_hyp} above.\\

The proof of the second inequality on line \eqref{general_comp}
follows from reasoning similar as that used to prove
\eqref{general_indct_hyp} inductively. We leave it to the reader to
set up the inductive hypothesis for this case and work out the details.
This completes our proof of Lemma \ref{comp_lem}.
\end{proof}\ret

Using Lemma \ref{comp_lem} and the assumed bounds
\eqref{initial_smallness}--\eqref{main_th_higher_noms}, we may
assume that our initial data \eqref{the_data} is such that:
\begin{align}
        \lp{(\uF(0),E(0))}{\dot{H}^\frac{n-4}{2}} \ &\leqslant \
        \td{\epsilon}_0
        \ , \label{c_initial_cSob_bnds}\\
        \lp{\underline{A}(0)}{\dot{H}^\frac{n-2}{2}} \ &\leqslant \
        \td{\epsilon}_0 \ ,
        \label{c_initial_A_cSob_bnds}\\
        \lp{(\uF(0),E(0))}{\dot{H}^k} \ &\leqslant \
        \td{M}_k
        \ , \label{h_initial_cSob_bnds}\\
        \lp{\underline{A}(0)}{\dot{H}^{k+1}} \ &\leqslant \
        \td{M}_k \ ,
        \label{h_initial_A_cSob_bnds}
\end{align}
where $\frac{n-4}{2} < k$, and the $\td{M}_k$ depend on the $M_k$ in
some inductive way, and we also have that $\td{\epsilon_0}\leqslant
C\varepsilon_0$ for some constant $C$ which depends only on the
dimension. Here $M_k$ and $\varepsilon_0$ refer to the constants
introduced in the statement of Theorem \ref{main_Th}.\\

We now decompose the initial field strength $\{E_i(0)\}$ in a way
that will be consistent with the evolution of the system
\eqref{YM_eq1}--\eqref{YM_eq2}. This will be convenient for
discussing the Cauchy problem. Our first step is to define the
following elliptic quantity:
\begin{equation}
        \Delta a_0 \ = \ -[a_i,\nabla^i a_0] + [a^i,E_i] \ .
        \label{initial_a0}
\end{equation}
where for convenience we have labeled $\{a_i\} =
\{\underline{A}_i(0)\}$. We then define the auxiliary set of
quantities:
\begin{equation}
        \dot{a}_i \ = \ E_i + \nabla_i a_0 -[a_0,a_i] \ .
        \label{initial_adoti}
\end{equation}
Notice that as an immediate consequence of the constraint equation
\eqref{YM_compat_cond}, the form of \eqref{initial_a0}, and the
Coulomb condition $d^* a=0$, we have the secondary Coulomb
condition:
\begin{equation}
        \nabla^i \dot{a}_i \ = \ 0 \ . \notag
\end{equation}
This will turn out to be important in a moment. Now, from the
definition of the quantities \eqref{initial_a0} and
\eqref{initial_adoti}, the already established bounds
\eqref{c_initial_cSob_bnds}--\eqref{h_initial_A_cSob_bnds}, and
several rounds of  Sobolev embeddings, we have the following
differential bounds on the quantities $\{\dot{a}_i\}$:
\begin{align}
        \lp{\dot{a}}{\dot{H}^\frac{n-4}{2}} \ \leqslant \ \td{\epsilon}_0 \ ,
        \label{c_initial_adoti_sob}\\
        \lp{\dot{a}}{\dot{H}^k} \ \leqslant \ \td{M}_k \ ,
        \label{h_initial_adoti_sob}
\end{align}
for $\frac{n-4}{2} < k$ (after a possible slight redefinition of the
constants $\td{\epsilon}_0 , \td{M}_k$ via multiplication by
some fixed dimensional constant). We now define a
\emph{Coulomb admissible} initial data set to be a collection
$(\uF,\{a_i\},\{\dot{a}_i\})$ such that:
\begin{align}
        d a + [a,a] \ &= \ \uF \ ,
        &d^* a \ &= \ 0 \ ,
        &d^* \dot{a} \ &= \ 0 \ . \  \label{coulomb_data}
\end{align}
Notice that $\uF$ is uniquely determined by the $\{a_i\}$, therefore
we do not need to include it in the definition of initial data. We
define the \emph{Coulomb-Cauchy problem} to be the task of finding a
space-time connection $D= d+A$ such that it satisfies the set of
equations:
\begin{subequations}\label{coulomb_YM}
\begin{align}
        D^\beta F_{\alpha\beta} \ = \ 0 \ , \label{coulomb_YM1}\\
        dA + [A,A] \ = \ F \ , \label{coulomb_YM2}\\
        d^*\underline{A} \ = \ 0 \ , \label{coulomb_YM3}
\end{align}
\end{subequations}
and such that at time $t=0$ we have that:
\begin{align}
        \underline{A}(0) \ &= \ a \ ,
        &\partial_t \underline{A}\, (0) \ = \ \dot{a} \ .
        \label{coulomb_initial}
\end{align}\ret

We remark briefly here that solving the problem
\eqref{coulomb_data}--\eqref{coulomb_initial} provides a solution to
the original Yang Mills system \eqref{YM_eq1}--\eqref{YM_eq2} with
Cauchy data \eqref{the_data} as long as we define the collection
$\{\dot{a}\}$ according to the equations
\eqref{initial_a0}--\eqref{initial_adoti}. All we need to do to
prove this assertion is to show that:
\begin{equation}
        F_{0i}(0) \ = \ E_i \ . \notag
\end{equation}
Our proof of this follows the same bootstrapping philosophy used to
show the equivalence \eqref{C_brack_iden} in the proof of Lemma
\ref{A_Uhl_lem}. The claim will follow at once from equation
\eqref{initial_adoti} if we can first establish that:
\begin{equation}
        A_0(0) \ = \ a_0 \ , \notag
\end{equation}
where $a_0$ is defined by \eqref{initial_a0}. Now, from the system
of equations \eqref{coulomb_YM} we have that the quantity $A_0$ is
elliptically determined by the equation:
\begin{equation}
        \Delta_{\underline{A}} A_0 \ = \ [A_i,\nabla_t A^i] \ ,
        \label{first_initial_ell}
\end{equation}
where $\Delta_{\underline{A}} = \underline{D}^i
\underline{D}_i$ is the gauge covariant
Laplacean. Furthermore, by using equation \eqref{initial_adoti} as
the definition of $E_i$, and substituting this into equation
\eqref{initial_a0}, we have that the quantity $a_0$ is elliptically
determined by the equation:
\begin{equation}
        \Delta_a a_0 \ = \ [a_i,\dot{a}^i] \ .
        \label{second_initial_ell}
\end{equation}
By subtracting \eqref{second_initial_ell} from
\eqref{first_initial_ell} at time $t=0$ we have that:
\begin{equation}
        \Delta_a (A_0(0) - a_0 ) \ = \ 0 \ . \notag
\end{equation}
Uniqueness now comes from the Sobolev type estimate:
\begin{equation}
        \lp{B}{L^{n}} \ \lesssim \ \lp{\Delta_a B}{L^\frac{n}{3}} \ , \notag
\end{equation}
which follows from the smallness condition
\eqref{c_initial_A_cSob_bnds} and the usual Sobolev estimates. The
details of the proof are left to the reader.\\

Keeping the equivalence we have just established in mind, and  the
first inequality contained in the comparison estimates
\eqref{general_comp} and \eqref{AF_equiv}, we have reduced the
demonstration of Theorem \ref{main_Th} to showing the following
non-gauge covariant global regularity theorem:\\

\begin{thm}[Global regularity in the Coulomb
gauge]\label{main_coulomb_Th} Let the number of spatial dimensions
be $6\leqslant n$. Then there exists a set of constants
$\td{\epsilon}_0$ and $C,C_k$, $\frac{n-2}{2}\leqslant k$ such that
if $(\uF, \{a_i\},\{\dot{a}_i\})$ is a Coulomb admissible initial
data set such that is satisfies the bounds:
\begin{subequations}\label{initial_coulomb_bounds}
\begin{align}
        \lp{\uF}{\dot{H}^\frac{n-4}{2}} \ &\leqslant \ \td{\epsilon}_0 \ ,
        &\lp{\uF}{\dot{H}^k} \ &\leqslant \ \td{M}_k \ ,
        \label{initial_coulomb_bounds1}\\
        \lp{a}{\dot{H}^\frac{n-2}{2}} \ &\leqslant \ \td{\epsilon}_0 \ ,
        &\lp{\dot{a}}{\dot{H}^\frac{n-4}{2}} \ &\leqslant \ \td{\epsilon}_0 \ ,
        \label{initial_coulomb_bounds2}\\
        \lp{a}{\dot{H}^k} \ &\leqslant \ \td{M}_{k-1} \ ,
        &\lp{\dot{a}}{\dot{H}^{k-1}} \ &\leqslant \ \td{M}_{k-1} \ ,
        \label{initial_coulomb_bounds3}
\end{align}
\end{subequations}
then if $\td{\epsilon}_0$ is sufficiently small there exists a
unique global solution $\{A_\alpha\}$ to the system
\eqref{coulomb_YM} with this initial data. Furthermore, this
solution obeys the following differential estimates:
\begin{subequations}\label{later_coulomb_bounds}
\begin{align}
        \lp{A}{\dot{H}^\frac{n-2}{2}} \ &\leqslant \ C\td{\epsilon}_0 \ ,
        &\lp{\partial_t A}{\dot{H}^\frac{n-4}{2}} \
        &\leqslant \ C\td{\epsilon}_0 , \ \label{later_coulomb_bounds1}\\
        \lp{A}{\dot{H}^k} \ &\leqslant \ C_{k-1}\td{M}_{k-1} \ ,
         &\lp{\partial_t A}{\dot{H}^{k-1}} \ &\leqslant \ C_{k-1}\td{M}_{k-1} \ ,
        \label{later_coulomb_bounds2}
\end{align}
\end{subequations}
\end{thm}\ret


\subsection{Local existence in the Coulomb gauge}
Our goal here is to reduce the proof of Theorem
\eqref{main_coulomb_Th} to a certain a-priori estimate involving the
energies of the field strength $F$. This amounts to proving a local
existence theorem for the system
\eqref{coulomb_data}--\eqref{coulomb_initial}. The proof of this
will allow us to set up a system of equations for the coulomb
potentials $\{A_\alpha\}$ which will be of central importance in the
sequel. We will show that:\\

\begin{prop}[Local existence in the Coulomb
gauge]\label{local_coulomb_prop} Let the number of spatial
dimensions be $6\leqslant n$. Then for every set of constants
$C,C_{k}$, $\frac{n-2}{2}\leqslant k$, there exists an
$\td{\epsilon}_0$
which only depends on $C$ with the following property: If
$(\{a_i\},\{\dot{a}_i\})$ is any set of Coulomb admissible initial
data such that:
\begin{align}
        \lp{a}{\dot{H}^\frac{n-2}{2}} \ &\leqslant \ C \td{\epsilon}_0 \ ,
        &\lp{\dot{a}}{\dot{H}^\frac{n-4}{2}} \ &\leqslant \ C \td{\epsilon}_0 \ ,
        \label{initial_small_coulomb}\\
        \lp{a}{\dot{H}^k} \ &\leqslant \ C_{k-1} \td{M}_{k-1} \ ,
        &\lp{\dot{a}}{\dot{H}^{k-1}} \ &\leqslant \ C_{k-1} \td{M}_{k-1} \ ,
        \label{initial_large_coulomb}
\end{align}
then for $\td{\epsilon}_0$ sufficiently small there exists a time $0
< T^*$, which only depends on the quantities $C \td{\epsilon}_0   ,
C_{\frac{n}{2}} \td{M}_{\frac{n}{2}} ,
C_{\frac{n+2}{2}}\td{M}_{\frac{n+2}{2}}$ such that there exists a
unique local solution $\{A_\alpha\}$ to the system
\eqref{coulomb_data}--\eqref{coulomb_initial} with this set of
initial data. Furthermore, on the time interval $[0,T^*]$ one has
the following norm bounds on the collection $\{A_\alpha\}$:
\begin{align}
        \sup_{0 \leqslant t \leqslant T^*}\
        \lp{A (t)}{\dot{H}^\frac{n-2}{2}} \ &\leqslant \ 2 C\td{\epsilon}_0 \ ,
        \label{local_coulomb_bounds1}\\
        \sup_{0 \leqslant t \leqslant T^*}\
        \lp{\partial_t A\, (t)}{\dot{H}^\frac{n-4}{2}} \ 
	&\leqslant \ 2 C\td{\epsilon}_0 \ ,
        \label{local_coulomb_bounds2}\\
        \sup_{0 \leqslant t \leqslant T^*}\
        \lp{A(t)}{\dot{H}^k} \ &\leqslant \ 2 
	C_{k-1}\td{M}_{k-1} \ ,
        \label{local_coulomb_bounds3}\\
        \sup_{0 \leqslant t \leqslant T^*}\
        \lp{\partial_t A\, (t)}{\dot{H}^{k-1}} \ &\leqslant 
	\ 2 C_{k-1}\td{M}_{k-1} \ .
        \label{local_coulomb_bounds4}
\end{align}
\end{prop}\ret

\begin{proof}[Proof of Proposition \ref{local_coulomb_prop}]
The proof will be reduced to the standard procedure of energy
estimates and Sobolev embeddings. Since we are assuming that the
initial data has enough smoothness to cover $L^\infty$, this is more
or less trivial. We start by plugging \eqref{coulomb_YM2} directly
into \eqref{coulomb_YM1}. After an application of the gauge
condition $d^*\underline{A}=0$ this yields a general second order
system of equations which we write as:
\begin{equation}
        \Box A_\beta\ = \ -\partial_\beta \partial_t A_0 +
        [\partial_t A_0 , A_\beta] - [A_\alpha,\partial^\alpha A_\beta]
        - [A^\alpha , F_{\alpha\beta}] \ . \label{generic_YM_coulomb}
\end{equation}
To split this into a hyperbolic-elliptic system, we decompose the
set of equations \eqref{generic_YM_coulomb} into its spatial and
temporal parts, and apply the Leray projection:
\begin{equation}
        \mathcal{P} \ = \ -\, \frac{d^* d}{\Delta}
        \ = \ \big(I - \nabla_x \frac{(\hbox{div})}{\Delta}\big) \ , \notag
\end{equation}
to the resulting spatial equation. After some rearrangement of the
elliptic equation this yields the coupled system:
\begin{subequations}\label{elliptic_2nd_YM}
\begin{align}
        \Box A_i \ &= \ \mathcal{P}\big( [\partial_t A_0 , A_i]
        - [A_\alpha,\partial^\alpha A_i]
        - [A^\alpha , F_{\alpha i}]\big) \ , \label{elliptic_2nd_YM1}\\
        \Delta A_0 \ &= \  -\,
        [A_i,\partial^i A_0]
        + [A^i , F_{0 i}] \ . \label{elliptic_2nd_YM2}
\end{align}
\end{subequations}
The above system of equations can be solved locally in time with the
bounds \eqref{local_coulomb_bounds1}--\eqref{local_coulomb_bounds1}
through a Picard iteration scheme. We leave this as an exercise for
the reader. Notice that the projection $\mathcal{P}$ can be removed
in energy estimates because it is an order zero operator. Notice
also that even though the smallness of the time interval $[0,T^*]$
will not make up for estimates involving the elliptic equation
\eqref{elliptic_2nd_YM2}, the critical smallness assumption
\eqref{initial_small_coulomb} allows one to obtain the bootstrapping
estimates
\eqref{local_coulomb_bounds1}--\eqref{local_coulomb_bounds1} if one
uses Littlewood-Paley decompositions and paraproducts to make sure
at least one factor in the non-linearity on the right hand side of
\eqref{elliptic_2nd_YM2} goes in a critical space. This same comment
goes for bounding terms on the right hand side of
\eqref{elliptic_2nd_YM1} in energy estimates when one is
bootstrapping the higher norm constants $C_k \td{M}_k$ for
$\frac{n+2}{2} < k$. Again, the smallness in time makes up for the
size of the first few constants $C \td{\epsilon}_0   ,
C_{\frac{n}{2}} \td{M}_{\frac{n}{2}} ,
C_{\frac{n+2}{2}}\td{M}_{\frac{n+2}{2}}$.\\

Having now produced a local solution to the system
\eqref{elliptic_2nd_YM} with the desired properties, we have shown
the conclusion of Proposition \eqref{local_coulomb_prop} once we
show that the spatial potentials which solve
\eqref{elliptic_2nd_YM1} are in fact solutions to the spatial
portion of the original second order equation
\eqref{generic_YM_coulomb}. This will be shown through  our general
strategy of coming up with a quantity which yields a critical
elliptic bootstrapping estimate which will force it to be zero. This
time, the desired quantity turns out to be related to the
conservation of electric charge for the Yang-Mills equations. We
first write the spatial portion of the non-linearity on the right
hand side of \eqref{generic_YM_coulomb} as a vector:
\begin{equation}
        \mathcal{N}_i \ = \ -\, \partial_i \partial_t A_0
        + [\partial_t A_0 , A_i]
        - [A_\alpha,\partial^\alpha A_i]
        - [A^\alpha , F_{\alpha i}] \ . \label{spat_N_def}
\end{equation}
We would like to show that the equations \eqref{elliptic_2nd_YM}
force $(I-\mathcal{P}) \mathcal{N} = 0$. We compute that:
\begin{equation}
        (I-\mathcal{P}) \mathcal{N} \ = \ \nabla_x\Delta^{-1}\left(
        -\partial_t \Delta A_0 - \partial^i\partial^\alpha [A_\alpha,A_i]
        - \partial^i [A^\alpha , F_{\alpha i}]
        \right) \ . \notag
\end{equation}
Now, using the equation \eqref{elliptic_2nd_YM} to compute
$\partial_t\Delta A_0$, this last line becomes:
\begin{align}
        (I-\mathcal{P}) \mathcal{N} \ &= \ \nabla_x\Delta^{-1}\left(
         - \partial^\beta \partial^\alpha [A_\alpha,A_\beta]
        - \partial^\beta [A^\alpha , F_{\alpha \beta}]
        \right) \ , \notag\\
        &= \ -\, \nabla_x\Delta^{-1}
        \partial^\beta [A^\alpha , F_{\alpha \beta}]
         \ . \notag
\end{align}
where the equality of the second line follows on account of skew
symmetry. We now isolate the interesting portion of the term on the
right hand side of the last line above and use the Jacobi identity
to compute that:
\begin{align}
        \partial^\beta [A^\alpha , F_{\alpha \beta}] \ &= \
        \frac{1}{2} [(dA)^{\alpha\beta}, F_{\alpha\beta} ]
        + [A^\alpha , \partial^\beta F_{\alpha\beta}] \ , \notag\\
        &= \ \frac{1}{2} \big[[A^\alpha,A^\beta], F_{\alpha\beta} \big]
        - \big[A^\alpha , [A^\beta , F_{\alpha\beta} ]\big]
        + [A^\alpha , D^\beta F_{\alpha\beta}] \ , \notag \\
        &= \ [A^\alpha , D^\beta F_{\alpha\beta}] \ . \notag
\end{align}
Now, again using equation \eqref{elliptic_2nd_YM2} we have that
$D^\beta F_{0\beta}=0$. Furthermore, from equation
\eqref{elliptic_2nd_YM1} we also have the identity:
\begin{equation}
        D^\beta F_{i\beta} \ = \ -\,
       (I-\mathcal{P})_i\mathcal{N} \ . \notag
\end{equation}
Combining all of this, we have the following equality:
\begin{equation}
        (I-\mathcal{P})\mathcal{N} \ = \ \nabla_x\Delta^{-1} [A^i ,
        (I-\mathcal{P})_i \mathcal{N}] \ . \label{N_trap_eq}
\end{equation}
Finally, from the form of \eqref{spat_N_def} and the already
established estimates
\eqref{local_coulomb_bounds1}--\eqref{local_coulomb_bounds4} as well
as the boundedness properties of the operator $(1-\mathcal{P})$ we
have that:
\begin{equation}
        \lp{(I-\mathcal{P})\mathcal{N}(t)}{L^\frac{n}{3}} \ < \ \infty
        \ , \notag
\end{equation}
for all times $t\in [0,T^*]$. However, from the smallness bound
\eqref{local_coulomb_bounds1}, the identity \eqref{N_trap_eq}, and a
Sobolev embedding we also have the fixed time bound:
\begin{align}
        \lp{ (I-\mathcal{P})\mathcal{N} }{L^\frac{n}{3}} \ &\lesssim \
        \lp{[A^i , (1-\mathcal{P})_i \mathcal{N}]}{L^\frac{n}{4}} \
        , \notag\\
        &\leqslant \ \lp{A}{L^n}\cdot\lp{(I-\mathcal{P})
        \mathcal{N}}{L^\frac{n}{3}} \ , \notag\\
        &\lesssim \ \td{\epsilon}_0 \cdot\lp{(I-\mathcal{P})
        \mathcal{N}}{L^\frac{n}{3}} \ . \notag
\end{align}
Therefore, for $\td{\epsilon}_0$ sufficiently small we see that we
must have $(I-\mathcal{P})\mathcal{N}=0$ as was to be shown. This
completes the proof that the solution to \eqref{elliptic_2nd_YM} is
a solution to the general system \eqref{generic_YM_coulomb}, and
therefore ends our proof of Proposition \ref{local_coulomb_prop}.
\end{proof}


\subsection{The second order curvature equation and the main
a-priori estimate}

Through a repeated application of the local existence theorem
\ref{local_coulomb_prop}, we may reduce the proof of the global
existence theorem \ref{main_coulomb_Th} to showing a-priori that any
solution to the Coulomb system
\eqref{coulomb_data}--\eqref{coulomb_initial} which exists on a time
interval $[0,T^*]$ (possibly large!), and such that it obeys the
\emph{both} the initial data bounds
\eqref{initial_coulomb_bounds1}--\eqref{initial_coulomb_bounds3}, as
well as the evolution bounds
\eqref{local_coulomb_bounds1}--\eqref{local_coulomb_bounds4}, in fact
obeys the improved evolution bounds
\eqref{later_coulomb_bounds1}--\eqref{later_coulomb_bounds2}.\\

Now, it turns out that the system of equations
\eqref{elliptic_2nd_YM} is by itself not so well
adapted\footnote{Strictly speaking, this is not entirely true. This
can be seen from the fact that if one looks at the localized
commutator $[\Box_A , \mathcal{P}]P_\lambda$, where the connection
$\{A_\alpha\}$ is assumed to be of much lower frequency than
$\lambda$, then this is essentially a ``derivative falls on low''
interaction which can be handled with the available Strichartz
estimates in $5\leqslant n$ dimensions. We have elected instead to
follow a formulation of the YM system which is based on the
curvature because of its conceptual appeal. However, in lower
dimensions, it may be best to work  directly with the connection
$\{A_\alpha\}$, in part to help mitigate bad $High\times High
\Rightarrow Low$ frequency interactions which come from the
quadratic term on the right hand side of
\eqref{2nd_order_curvature}.} to the proof of such an a-priori
estimate. This stems from the fact that these equations are not
covariant. This manifests itself in the projection operator
$\mathcal{P}$. If one were to try to write the hyperbolic system of
equations \eqref{elliptic_2nd_YM1} in terms of covariant wave
operator $\Box_A$ and a source term, the projection operator which
is non-local would end up causing problems in various commutator
terms. The way around this is to not only consider the system
\eqref{elliptic_2nd_YM}, but to also work directly with the
curvature in the equations \eqref{coulomb_YM1}--\eqref{coulomb_YM2}.
This is possible because we are not attempting to set up an
iteration scheme, but are instead merely trying to prove an a-priori
estimate, so we may safely assume that the quantities we work with
satisfy any equation which results from the system
\eqref{coulomb_YM}. We will in fact use several such elliptic and
hyperbolic equations. As a very rough description of this kind of
philosophy, the reader may find it useful to keep in mind the
following schematic:
\begin{align}
        \hbox{Weak control of the connection} \ &\Longrightarrow \
        \hbox{Improved control of the curvature}   \ , \notag\\
        &\Longrightarrow \ \hbox{Improved control of the connection}
        \ , \notag\\
        &\Longrightarrow \ \hbox{Weak control of the connection
        for longer times} \ . \notag
\end{align}
To provide the improved control on the curvature, we will employ a
second order equation for it. To derive this, we write the Bianchi
identities \eqref{coulomb_YM2} in the form \eqref{YM_eq2} and then
contract this expression with the covariant derivative $D$. This
yields the equations:
\begin{align}
        0 \ &= \ D^\gamma( D_\alpha F_{\beta\gamma} +
        D_\gamma F_{\alpha\beta} + D_\beta F_{\gamma\alpha} ) \ ,
        \notag \\
        &= \ \Box_A F_{\alpha\beta} + [F^\gamma_{\ \ \alpha} , F_{\beta\gamma}]
        + [F^\gamma_{\ \ \beta} , F_{\gamma\alpha}] \ , \notag\\
        &= \ \Box_A F_{\alpha\beta} -2[F_{\alpha\gamma} , F_\beta^{\
        \gamma}] \ . \label{2nd_order_curvature}
\end{align}
In addition to \eqref{2nd_order_curvature} and the system
\eqref{elliptic_2nd_YM}, it will also be useful for us to employ a
secondary elliptic equation. This will be for the quantity
$\partial_t A_0$:
\begin{equation}
        \partial_t A_0 \ =\ \Delta^{-1}\partial^i (
        -\, [A_i,\partial_t A_0] + [A_0 , \partial_t A_i]
        + [A^\alpha , F_{i\alpha}] ) \ . \label{dtA0_eq}
\end{equation}
This equation follows immediately from differentiating the equation
\eqref{elliptic_2nd_YM2} with respect to time, and then applying the
conservation law $\nabla^\alpha [A^\beta,F_{\alpha\beta}] =0$ to the
resulting expression. We are now ready to state our main a-priori
estimate:\\

\begin{thm}[Main a-priori estimate for the curvature of the Coulomb
system
\eqref{coulomb_data}--\eqref{coulomb_initial}]\label{apriori_thm}
Let the space-time connection $D=d+A$ on $\mathbb{R}^{(n+1)}$, where
$6\leqslant n$, be given such that it satisfies the following system
of equations on some finite time interval $[0,T^*]$:
\begin{subequations}\label{2nd_order_YM}
\begin{align}
        \Box_A F_{\alpha\beta} \ &= \ 2[F_{\alpha\gamma} , F_\beta^{\
        \gamma}] \ , \label{2nd_order_YM1}\\
        d{A} + [{A},{A}] \ &= \
        F \ , \label{2nd_order_YM2} \\
        d^* \underline{A} \ &= \ 0 \ , \label{2nd_order_YM3}\\
        \Box A_i \ &= \ \mathcal{P}\big( \, [\partial_t A_0 , A_i]
        - [A_\alpha,\partial^\alpha A_i]
        - [A^\alpha , F_{\alpha i}]\big) \ , \label{2nd_order_YM4}\\
        \Delta A_0 \ &= \
        \partial^i\, [A_0, A_i]
        + [A^i , F_{0 i}] \ , \label{2nd_order_YM5}\\
        \Delta (\partial_t A_0) \ &=\ \partial^i \big(
        -\, [A_i,(\partial_t A_0)] + [A_0 , \partial_t A_i]
        + [A^\alpha , F_{i\alpha}] \big) \ . \label{2nd_order_YM6}
\end{align}
\end{subequations}
Here we have split $\{A_\alpha\} = ( A_0 , \{\underline{A}_i\})$.
Let there also be given a set of fixed constants $L, N , L_k , N_k $
for the indices $\frac{n-2}{2}\leqslant k$, such that at time $t=0$
we have the initial bounds:
\begin{align}
        \lp{F(0)}{\dot{H}^\frac{n-4}{2}} \ &\leqslant \
        \td{\epsilon}_0 \ ,
        &\lp{\partial_t F\, (0)}{\dot{H}^\frac{n-6}{2}} \ &\leqslant \
        L\td{\epsilon}_0 \ , \label{main_ap_F_small}\\
        \lp{F(0)}{\dot{H}^k} \ &\leqslant \ \td{M}_k \ ,
        &\lp{\partial_t F\, (0)}{\dot{H}^{k-1}} \ &\leqslant \
        L_{k}\td{M}_k \ . \label{main_ap_F_large}
\end{align}
Then if $\td{\epsilon_0}$ is chosen as to be sufficiently small on
line \eqref{main_ap_F_small} above,  there exists a collection
constants $C , C_k$, which only depend on the dimension and the
collection $L,N,L_k ,N_k$ but \emph{not} on $\td{\epsilon}_0$ (once
it is small enough) or the collection $\td{M}_k$, such that if at
later times we have the bounds:
\begin{align}
        \sup_{0\leqslant t \leqslant T^*}
        \lp{\underline{A}(t)}{\dot{H}^\frac{n-2}{2}} \ &\leqslant \
        2N C \td{\epsilon}_0 \ ,
        &\sup_{0\leqslant t \leqslant T^*}
        \lp{\partial_t \underline{A}\, (t)}{\dot{H}^\frac{n-4}{2}} \ &\leqslant \
        2N C \td{\epsilon}_0 \ , \label{later_main_ap_A_small}\\
        \sup_{0\leqslant t \leqslant T^*}
        \lp{F(t)}{\dot{H}^\frac{n-4}{2}} \ &\leqslant \
        2N C \td{\epsilon}_0 \ ,
        &\sup_{0\leqslant t \leqslant T^*}
        \lp{\partial_t F\, (t)}{\dot{H}^\frac{n-6}{2}} \ &\leqslant \
        2N C \td{\epsilon}_0 \ , \label{later_main_ap_F_small}\\
        \sup_{0\leqslant t \leqslant T^*}
        \lp{\underline{A}(t)}{\dot{H}^k} \ &< \ \infty \ ,
        &\sup_{0\leqslant t \leqslant T^*}
        \lp{\partial_t \underline{A}\, (t)}{\dot{H}^{k-1}} \ &< \
        \infty \ , \label{later_main_ap_A_large}\\
        \sup_{0\leqslant t \leqslant T^*}
        \lp{F(t)}{\dot{H}^k} \ &< \ \infty \ ,
        &\sup_{0\leqslant t \leqslant T^*}
        \lp{\partial_t F\, (t)}{\dot{H}^{k-1}} \ &< \
        \infty \ , \label{later_main_ap_F_large}
\end{align}
the following set of stronger bounds holds:
\begin{align}
        \sup_{0\leqslant t \leqslant T^*}
        \lp{F(t)}{\dot{H}^\frac{n-4}{2}} \ &\leqslant \ N^{-1} C
        \td{\epsilon}_0 \ ,
        &\sup_{0\leqslant t \leqslant T^*}
        \lp{\partial_t F\, (t)}{\dot{H}^\frac{n-6}{2}} \ &\leqslant \
        N^{-1} C \td{\epsilon}_0 \ , \label{better_later_main_ap_F_small}\\
        \sup_{0\leqslant t \leqslant T^*}
        \lp{F(t)}{\dot{H}^k} \ &\leqslant \ N^{-1}_k C_k \td{M}_k \ ,
        &\sup_{0\leqslant t \leqslant T^*}
        \lp{\partial_t F\, (t)}{\dot{H}^{k-1}} \ &\leqslant \
        N^{-1}_k C_k \td{M}_k \ . \label{better_later_main_ap_F_large}
\end{align}
\end{thm}\ret

\begin{rem}
The bounds involving \eqref{later_main_ap_F_large} and
\eqref{better_later_main_ap_F_large} express the fact that the
control we provide here is at the critical level. That is, bounds on
the higher norms are completely irrelevant in the bootstrapping
procedure, except for the fact that they are finite. The only place
where we need higher norms to accomplish anything here is in the
local existence theorem \ref{local_coulomb_prop}. The way we will
prove Theorem \ref{apriori_thm} is by first establishing control at
the critical level through a bootstrapping argument. The control of
the higher norms will then be provided through an a-priori estimate
who's proof is essentially identical to that of the critical
bootstrapping bound, and will therefore be left to the reader.
\end{rem}\ret

\begin{rem}
The reader my find it useful to have a brief description of the
various constants appearing in Proposition \ref{local_coulomb_prop}
and Theorem \ref{apriori_thm}. The constants $L,L_k,N,N_k$ are input
into the a-priori machine, and these are meant to cover the
transition to and from estimates involving the connection and
curvature. The set $L,L_k$ is only needed to deal with the initial
data. This is necessary because we must have an account of bounds
involving the quantities $\partial_t F$. The other constants $N,N_k$
govern comparison type estimates similar to \eqref{AF_equiv}. The
constants $C,C_k$ are byproducts of the proof of the a-priori
estimate itself. These will very much depend on the $L,L_k,N,N_k$,
but are \emph{independent} of $\td{\epsilon}_0$ when it is small
enough. Finally, the main adjusting parameter $\td{\epsilon}_0$ has
two important roles. First and foremost, it is needed to prove the
a-priori estimate itself. However, it has a second purpose which is
also crucial, and that is to keep the dependence of $C,C_k$ on
$L,L_k,N,N_k$ from creating a feedback loop. Specifically, we need
our various comparison estimates to have constants which \emph{do
not} depend on the large constants $C,C_k$. Since the critical
energy of the curvature can grow by a factor of $C$, we will need
the extra influence of $\td{\epsilon}_0$ to make sure this does not
cycle back to $L,L_k,N,N_k$.
\end{rem}\ret

\begin{proof}[Proof that Theorem \ref{apriori_thm} and Proposition
\ref{local_coulomb_prop} together imply Theorem
\ref{main_coulomb_Th}] The proof here is more or less
straightforward and will be largely left to the reader. Everything
relies on two sets of estimates. The first has to do with showing
that the initial data bounds
\eqref{initial_coulomb_bounds1}--\eqref{initial_coulomb_bounds3}
imply the initial control assumed in
\eqref{main_ap_F_small}--\eqref{main_ap_F_large}. This is just a
matter of bounding the time derivatives $\partial_t F$, and is why
we have included the set of auxiliary constants $L, L_k$. Using now
the field equations \eqref{YM_eq1}--\eqref{YM_eq2} (we have not
included them in the system \eqref{2nd_order_YM}, but we may assume
they hold), we have the general schematic identity at time $t=0$:
\begin{equation}
        \partial_t F\, (0) \ = \ \nabla_x F\, (0) + [a,F(0)] \ ,
        \label{dt_F_iden}
\end{equation}
where we have generically set $a=(a_0,\{a_i\})$. Therefore, to
establish the control
\eqref{main_ap_F_small}--\eqref{main_ap_F_large}, we only need to
prove the estimates:
\begin{align}
        \lp{[a,F(0)]}{\dot{H}^\frac{n-6}{2}} \ &\lesssim \
        \td{\epsilon}_0 \ ,
        &\lp{[a,F(0)]}{\dot{H}^{k-1}} \ &\lesssim \ \td{M}_k \ ,
        \label{AF_initial_prod_est}
\end{align}
assuming that the bounds
\eqref{initial_coulomb_bounds1}--\eqref{initial_coulomb_bounds3}
hold. Notice that while these initial norms do not contain estimates
on the quantities $E_i = F_{0i}(0)$, we originally had bounds on
this from lines
\eqref{c_initial_cSob_bnds}--\eqref{h_initial_A_cSob_bnds} above.
Also, any estimates on $a_0$ which are needed in this process can be
provided, for instance, through the equation \eqref{initial_a0}.
Since the proof of estimate \eqref{AF_initial_prod_est} is a
straightforward paraproduct type bound, similar to what was done in
the proof of Lemma \ref{comp_lem} above, we leave it to the
interested reader (see below for more details).\\

The second set of estimates we need to prove here has to do with the
relationship between the later time norms
\eqref{later_main_ap_A_small}--\eqref{better_later_main_ap_F_large}
and the ones
\eqref{local_coulomb_bounds1}--\eqref{local_coulomb_bounds4}
contained in the proof of the local existence proposition. Since our
global regularity proof is by iteration of this latter result, we
need to first show that the weak control
\eqref{local_coulomb_bounds1}--\eqref{local_coulomb_bounds4} implies
the bootstrapping assumption
\eqref{later_main_ap_A_small}--\eqref{later_main_ap_F_large}. This
assertion is trivial for norms involving the potentials
\eqref{later_main_ap_A_small} and \eqref{later_main_ap_A_large}, as
well as the larger norms \eqref{later_main_ap_F_large} just by
applying the definition of curvature. Therefore, we only need to see
that \eqref{local_coulomb_bounds1}--\eqref{local_coulomb_bounds2}
implies the bounds \eqref{later_main_ap_F_small}. We first establish
the desired bounds for the undifferentiated term $F$. For the
spatial curvature and potentials $(\uF,\underline{A})$, this is just
the comparison principle form line \eqref{AF_equiv}, and we can
assume that the constants $N,N_k$ are large enough to cover that
case. To deal with potentials involving time derivatives of
$\underline{A}$ or the temporal potential $A_0$ we have the
following general calculation:
\begin{align}
        \lp{F}{\dot{H}^\frac{n-4}{2}} \ &\leqslant \
        \lp{dA}{\dot{H}^\frac{n-4}{2}} +
        \lp{[A,A]}{\dot{H}^\frac{n-4}{2}} \ , \notag\\
        &\lesssim \ \lp{A}{\dot{H}^\frac{n-2}{2}} +
        \lp{A}{\dot{H}^\frac{n-2}{2}}^2 \ , \notag
\end{align}
where the quadratic term follows from paraproduct decompositions,
H\"olders inequality, and Sobolev embeddings as in the proof of
\eqref{AF_equiv}. The desired result now follows from the smallness
of $\td{\epsilon}_0$ and the fact that we may assume the constant
$C$ in line \eqref{local_coulomb_bounds1} does not depend on it. To
establish the estimate for the quantity $\partial_t F$, we  use the
later time version of the identity \eqref{dt_F_iden}, as well as the
estimate which is responsible for the first estimate on line
\eqref{AF_initial_prod_est} above, which is:
\begin{equation}
        \lp{[A,F]}{\dot{H}^\frac{n-6}{2}} \ \lesssim \
        \lp{A}{\dot{H}^\frac{n-2}{2}}\cdot\lp{F}{\dot{H}^\frac{n-4}{2}}
        \ . \notag
\end{equation}
By again assume that the constant $\td{\epsilon}_0$ is sufficiently
small with respect to $C$ we have the desired bound.\\

The final thing we need to do here is to show that the improved
bounds
\eqref{better_later_main_ap_F_small}--\eqref{better_later_main_ap_F_large}
imply the assumed estimates of the local existence theorem
\eqref{initial_small_coulomb}--\eqref{initial_large_coulomb}. This
is again a comparison estimate either identical or similar to
\eqref{AF_equiv}. Note that we only need to bound the spatial
portion of the potentials $\{A_\alpha\}$ and their time derivatives.
The undifferentiated terms can be bounded directly by
\eqref{AF_equiv} because we may assume that the constant
$\td{\epsilon}_0$ on line \eqref{better_later_main_ap_F_small} is
small enough that the critical estimate \eqref{comp_lem_crit} holds.
To deal with the time differentiated potentials $\partial_t
\underline{A}$, one can simply differentiate the Hodge system
\eqref{2nd_order_YM2}--\eqref{2nd_order_YM3} with respect to time
and then apply essentially the same proof as was used to produce
\eqref{AF_equiv}. The details of this are left to the ambitious
reader.
\end{proof}\ret

\ret

\section{Proof of the Main Bootstrapping Estimate}

We are now ready to begin our proof of the (improved) 
main critical a-priori estimate
\eqref{better_later_main_ap_F_small}.
In order to do this, we will need to bootstrap in a function space
which is much stronger than the energy type spaces of Theorem
\ref{apriori_thm}. This will cost us another bootstrapping
procedure, but this will be easy to set up because it will be clear
the extra norms we create have good bounds on some very small
initial time interval due to the fact that we are assuming the
higher energy boundedness \eqref{later_main_ap_F_large} and that
these norms involve integration in time. All of the norms we
construct here will be of Strichartz type, with an $\ell^2$ Besov
structure in the spatial variable. It will also be necessary for us
to include an angular square sum structure in many of the
estimates we prove. This
may seem a bit odd at first because we will not need such bounds
directly in our proof of Theorem \ref{apriori_thm}. These extra bounds
will instead be used to give the fine control which is needed to
handle the linear part of the problem. At each fixed frequency, we
form the square-sum norms:
\begin{equation}
        \lp{P_\lambda A}{SL^P} \ = \ \sup_{\theta \lesssim
        1}  \ \left( \sum_{\substack{\phi\ :\\
     \omega_0 \in \Gamma_{\phi}} }
        \ \lp{{}^{\omega_0}\!\Pi_\theta  
	P_\lambda A}{L^p}^2 \right)^\frac{1}{2} \
        , \label{ang_sum_norms}
\end{equation}
where $\Gamma_\phi$ is taken to be a (uniformly) finitely
overlapping set of spherical caps such that $\mathbb{S}^{n-1}
= \cup_\phi \Gamma_\phi$, each of which has size $\sim \theta$
and constructed such a way that one has the bounds:
\begin{equation}
        \left(\sum_{\substack{\phi\ :\\
     \omega_0 \in \Gamma_{\phi}} }\
        \lp{{}^{\omega_0}\!\Pi_\theta P_\lambda A}
    {L^2}\right)^\frac{1}{2} \ \lesssim \
        \lp{P_\lambda A}{L^2} \ , \notag
\end{equation}
independent of the size of $\theta$. Here we take the condition
$\omega_0\in\Gamma_\phi$ to mean that the variable $\omega_0$ is
essentially in the center of that spherical cap $\Gamma_\phi$. The
exact placement is not essential. Notice that by construction, these
norms are contained in the usual $L^p$ spaces because we can assume
that one set of angular sectors we are summing over contains the
whole sphere.\\

Next, using the same prescription that defined the Besov spaces
\eqref{Besov_norm}, we define the angular square sum Besov spaces to
be:
\begin{equation}
        \lp{A}{ S \dot{B}_2^{p,(q,s)} } \ = \ \left( \sum_\lambda\
        \lambda^{2s-2n(\frac{1}{q} - \frac{1}{p})}\   \lp{P_\lambda
        A}{SL^p}^2\right)^\frac{1}{2} \ . \label{SB_norm}
\end{equation}
We now define the main dispersive component of the function spaces we
will be working with. These are $L^2_t$ based Strichartz spaces,
built on the norms \eqref{SB_norm} and \eqref{Besov_norm}. These are
all defined on a finite time interval $[0,T^*]$, which will for the
most part be left implicit:
\begin{align}
        \lp{A}{\dot{Z}^s} \ &= \
        \lp{A}{ L_t^2( \dot{B}_2^{\frac{2(n-1)}{n-3}
        ,(2,s+\frac{1}{2})})[0,T^*] } \ , \
        \label{Z_norm}\\
        \lp{A}{S\dot{Z}^s} \ &= \
        \lp{A}{ L_t^2(S \dot{B}_2^{\frac{2(n-1)}{n-3}
        ,(2,s+\frac{1}{2})})[0,T^*] } \ . \
        \label{BZ_norm}
\end{align}
To gain some intuition about these spaces, notice that they all
scale like $L^\infty(\dot{H}^s)$ under the change of variables
\eqref{spatial_scale_trans}. Therefore, they all scale like
solutions to the wave equations with $\dot{H}^s$ initial data.
Indeed, these spaces are consistent with the available range of
Strichartz estimates for the usual scalar wave equation, and it will
be our goal to show that one has bounds on the norm \eqref{BZ_norm}
for solutions of the covariant wave operator on the left hand side
of \eqref{2nd_order_YM}.\\

To form the overall spaces we will bootstrap in, we add the
above space-time norms to the energy type norms  used in the
statement of the main a-priori estimate \eqref{apriori_thm}:
\begin{align}
        \dot{X}^s \ &= \ L^\infty[0,T^*](\dot{H}^s)\, \cap\, S\dot{Z}^s
        \ , \label{X_norm}\\
        \dot{Y}^s \ &= \ L^\infty[0,T^*](\dot{H}^s)\, \cap\, \dot{Z}^s
        \ . \label{Y_norm}
\end{align}
It will also be necessary for us to estimate time derivatives in the
above spaces. Since differentiation will decrease the scaling by one
unit, we use the norms:
\begin{equation}
        \lp{A}{\dot{X}^s \times \partial_t^{-1}(\dot{X}^{s-1})}
        \ = \ \lp{A}{\dot{X}^s} + \lp{\partial_t A}{\dot{X}^{s-1}} \
        , \notag
\end{equation}
with an analogous definition for $\dot{Y}^s \times
\partial_t^{-1}(\dot{Y}^{s-1})$. \\

\subsection{Proof of the Critical Bootstrapping Estimate}
We are now ready to prove the critical component of
Theorem \eqref{apriori_thm} (we will now
change notation from $\td{\epsilon}_0$ back to $\epsilon_0$):\\

\begin{prop}[Critical bootstrapping estimate in the $\dot{X}^s$
spaces]\label{critical_bs_prop} Let the dimension be $6\leqslant n$.
Let the collection $(F,A)$ be a space-time connection curvature pair
which obeys the general smoothness conditions
\eqref{later_main_ap_A_large}--\eqref{later_main_ap_F_large}, and
which satisfies the system of equations \eqref{2nd_order_YM}. Let
$L,N$ be given constants such that one has the initial bounds:
\begin{equation}
        \lp{F(0)}{\dot{H}^\frac{n-4}{2}}
        + \lp{\partial_t F\, (0)}{\dot{H}^\frac{n-6}{2}} \
        \leqslant \ L\epsilon_0 \ . \label{initial_critical_bs_est}
\end{equation}
Then there exists a constant $C$ which depends only on $L,N$ and the
dimension such that if one has the bootstrapping bounds on a time
interval $[0,T^*]$:
\begin{align}
        \sup_{0\leqslant t \leqslant T^*}
        \lp{(\underline{A},\partial_t\underline{A})(t)}
        {\dot{H}^\frac{n-2}{2}\times\dot{H}^\frac{n-4}{2}} \ &\leqslant \
        2N C {\epsilon}_0 \ ,
        \label{critical_bs_assumed_est_A}\\
        \lp{F}{\dot{X}^\frac{n-4}{2}\times
        \partial_t^{-1}(\dot{X}^\frac{n-6}{2})}
        \ &\leqslant \ 2NC\epsilon_0 \ ,
        \label{critical_bs_assumed_est}
\end{align}
then for $\epsilon_0$ sufficiently small, we have that the following
improved bounds  on the same time interval $[0,T^*]$:
\begin{equation}
        \lp{F}{\dot{X}^\frac{n-4}{2}\times
        \partial_t^{-1}(\dot{X}^{\frac{n-6}{2}})}
        \ \leqslant \ N^{-1} C\epsilon_0 \ .
        \label{critical_bs_imp_est}
\end{equation}
\end{prop}\ret

The proof of Proposition \ref{critical_bs_prop} will be accomplished
through the standard use of Littlewood-Paley paraproduct
decompositions, and the application of space-time estimates.
All of the linear bounds we will need are
provided by the following, which is the main technical result of
this work:\\

\begin{thm}[Gauge covariant angular square-sum Strichartz estimates
for Yang-Mills connections]
\label{Str_thm} Let the number of dimensions be such that
$6\leqslant n$, and let $d+\td{\underline{A}}$ be a space-time
connection defined \emph{defined on all of Minkowski space}
$\mathcal{M}^{n+1}$ such that it satisfies the conditions:
\begin{subequations}\label{red_conctn_cond}
\begin{align}
        \td{\underline{A}}_0 \ &= \ 0 \
        &\hbox{(Temporal Gauge)} \ , \label{red_conctn_cond1}\\
        d^*\td{\underline{A}}\ &= \ 0 \
        &\hbox{(Coulomb Gauge)} \ , \label{red_conctn_cond2}\\
        P_{|\xi| \ll |\tau| }(
        \td{\underline{A}}) \ &= \ 0 \
        &\hbox{(Space-time frequency localization)} \ ,
        \label{red_conctn_cond4}\\
        \lp{\td{\underline{A}}}{\dot{X}^\frac{n-2}{2}} \
        &\leqslant \ \mathcal{E}\
        &\hbox{(Space-time estimate)} \ , \label{red_conctn_cond5}\\
        \Box \td{\underline{A}} \ &= \
        \td{\mathcal{P}}([B,H]) \  &\hbox{(Structure equation)} \ ,
        \label{red_conctn_cond6}\\
        \lp{(B,H)}{\dot{Y}^\frac{n-2}{2}\times
        \dot{Y}^\frac{n-4}{2}} \ &\leqslant \ \mathcal{E} \
        &\hbox{(Structure estimates)} \ , \label{red_conctn_cond7}
\end{align}
\end{subequations}
where $(B,H)$ is an auxiliary set of $\mathfrak{g}$ valued functions
defined on all of $\mathcal{M}^{n+1}$. The symbol $\td{\mathcal{P}}$
denotes a composition of the Leray projection $\mathcal{P}$ with some
frequency cutoff function which is bounded on all mixed Lebesgue-Besov
spaces of the type $L^p(\dot{B}_2^{p,(2,s)})$. We assume also that the
connection $d+\td{\underline{A}}$ satisfies the general smoothness
bounds:
\begin{align}
        \sup_{-T^* \leqslant t \leqslant T^*}\
        \lp{\td{\underline{A}}(t)}{\dot{H}^k} \ &< \
	\infty \ ,
        &\frac{n-2}{2} \ < \ k \ , \label{general_smoothness}
\end{align}
for each fixed time $T^*$. Let now $F$ be any other $\mathfrak{g}$
valued function which satisfies the inhomogeneous equation:
\begin{equation}
        \Box_{\td{\underline{A}}}  F \ = \
        G \ , \label{inhomog_gc_wave}
\end{equation}
with Cauchy data:
\begin{align}
        F\, (0) \ &= \ f \ ,
        &\partial_t  F\, (0) \ = \ \dot{f} \ .
        \label{general_cuachy_data}
\end{align}
Then if the constant $\mathcal{E}$ in lines \eqref{red_conctn_cond5}
and \eqref{red_conctn_cond7} above is sufficiently small, one has
the following family of space-time estimates:
\begin{equation}
        \lp{ F}{\dot{X}^\frac{n-4}{2} \times
        \partial_t^{-1}(\dot{X}^\frac{n-6}{2})}
        \ \ \lesssim \ \
        \lp{(f,\dot{f})}{\dot{H}^\frac{n-4}{2}\times
        \dot{H}^\frac{n-6}{2}} + \lp{G}{L^1(\dot{H}^\frac{n-6}{2})}
        \ . \label{the_strichartz_est}
\end{equation}
\end{thm}\ret

\begin{rem}
In the above Theorem,  the Strichartz estimates have a preferred
scaling. This is consistent with the application we have in mind. In
general, it is not possible to prove estimates of the type
\eqref{the_strichartz_est} for higher Sobolev indices without
assuming that the connection $\td{\underline{A}}$ itself has more
regularity. In the case where $\td{\underline{A}}$ does have better
regularity, a proof similar to that given after Proposition
\ref{freq1_parametrix_prop} below
can be used to show estimates for those higher norms.
\end{rem}\ret

\begin{proof}[Proof of Proposition \ref{critical_bs_prop}]
The proof requires another bootstrapping argument. This will be done
on subintervals $[0,T^{**}]\subseteq [0,T^*]$. Using the initial bounds
\eqref{initial_critical_bs_est} and the general smoothness
assumption \eqref{later_main_ap_F_large} we may assume that for
$T^{**}\ll 1$ we have the estimate \eqref{critical_bs_assumed_est}.
Therefore, it suffices to prove that \eqref{critical_bs_assumed_est}
implies \eqref{critical_bs_imp_est} on all subintervals $[0,T^{**}]$.
But this is just the same as proving Proposition
\ref{critical_bs_prop} itself since $T^*$ is arbitrary.\\

The proof will be accomplished in a series of steps. Our first goal
will be to derive $\dot{X}^s$ and $\dot{Z}^s$ type bounds for the
connection $d+A$. We will then split this connection into a sum of
two pieces $d+\td{A} + \td{\td{A}}$, where the potentials $\td{A}$
satisfy the criteria of Theorem \ref{Str_thm} and the remainder
term $\td{\td{A}}$ obeys the better $L^1(L^\infty)$
space-time estimate. This is enough to be able to write the equation
\eqref{2nd_order_YM1} schematically as:
\begin{equation}
        \Box_{\td{A}} F \ = \ [\nabla\td{\td{A}},F] + [\td{\td{A}},\nabla
        F] + \big[\td{A},[\td{\td{A}} ,F]\big] +
        \big[\td{\td{A}},[\td{\td{A}} ,F]\big] + [F,F] \ .
        \label{wave_schematic}
\end{equation}
One is then in a position where Theorem \ref{Str_thm} can be applied
directly, and we only need to choose our constant $C$ depending on
$L,N$ and the constant which appears on line
\eqref{the_strichartz_est}. The key thing is that the dangerous term
$[\td{\td{A}},\nabla F]$ can safely be put in
$L^1(\dot{H}^\frac{n-6}{2})$ using the improved space-time estimate
for $\td{\td{A}}$ and the energy estimate for $F$. Throughout
the proof we will use the usual splitting
$\{A_\alpha\}=(A_0,\underline{A})$ of $d+A$ into its temporal and
spatial components.\\

\subsection*{$\bullet$\ \  $\dot{X}^\frac{n-2}{2}$ estimates for
$\{\underline{A}_i\}$} Here we write $\uF$ for the spatial
components of  the field strength and use the Hodge system
\eqref{2nd_order_YM2}--\eqref{2nd_order_YM3} to  write
schematically:
\begin{equation}
        \underline{A} \ = \ \nabla_x \Delta^{-1} ( -\, \uF +
        [\underline{A},\underline{A}]) \ . \label{uA_int_eq}
\end{equation}
As a preliminary first step, we will show that the potentials
$\{\underline{A}_i\}$ can be estimated in $\dot{Y}^\frac{n-2}{2}$
with bounds comparable to $NC \epsilon_0$. Now, it is not too
difficult to see directly from the definition that:
\begin{equation}
        \nabla_x \Delta^{-1}\ : \ \dot{Y}^\frac{n-4}{2} \
        \hookrightarrow \ \dot{Y}^\frac{n-2}{2} \ . \notag
\end{equation}
Next, notice that we have the bilinear estimate:
\begin{equation}
        \nabla_x \Delta^{-1}\ : \ L^\infty(\dot{H}^\frac{n-2}{2})
        \cdot \dot{Y}^\frac{n-2}{2} \ \hookrightarrow \
        \dot{Y}^\frac{n-2}{2} \ , \label{bilin_to_Y}
\end{equation}
which follows integrating the bound \eqref{general_besov_embed}.
Note that in this case, the range restrictions
\eqref{sc_cond}--\eqref{Lb_cond} are easily satisfied. Therefore,
using the critical bounds \eqref{critical_bs_assumed_est_A} as well
as the general smoothness criteria \eqref{later_main_ap_A_large} (so
that in particular we may assume  the $\dot{Y}^\frac{n-2}{2}$ norm
of $\{\underline{A}_i\}$  is finite) we see we may absorb the
quadratic term on the right hand side of \eqref{uA_int_eq} onto the
left in the desired estimates.\\

Our task is now to show the more restrictive $\dot{X}^\frac{n-2}{2}$
estimates for the potentials $\{\underline{A}_i\}$. Again from the
definition, it is not hard to see that we have the embedding:
\begin{equation}
        \nabla_x\Delta^{-1}\ : \ \dot{X}^\frac{n-4}{2} \
        \hookrightarrow \ \dot{X}^\frac{n-2}{2} \ . \notag
\end{equation}
Therefore, keeping in mind the $\dot{Y}^\frac{n-2}{2}$ bounds just
proved, we see that is suffices to be able to show the bilinear
estimate:
\begin{equation}
        \nabla_x\Delta^{-1}\ : \ \dot{Y}^\frac{n-2}{2}
        \cdot \dot{Y}^\frac{n-2}{2} \
        \hookrightarrow \ \dot{X}^\frac{n-2}{2} \  .
        \label{bilin_X_embed}
\end{equation}
The main issue here is, of course, to be able to include the angular
square sum structure. This turns out to be very simple. Notice first
that by orthogonality and the general nesting \eqref{Besov_nesting}
we have the inclusion (on any finite time interval $[0,T^*]$):
\begin{equation}
        L^\infty(\dot{H}^\frac{n-2}{2})\cap L^2(\dot{H}^\frac{n-1}{2}) \
        \subseteq \ \dot{X}^\frac{n-2}{2} \ .
        \notag
\end{equation}
Therefore, to conclude \eqref{bilin_X_embed} we see that it suffices
to be able to show the set of bilinear estimates:
\begin{align}
        \nabla_x\Delta^{-1}\ : \ \dot{Y}^\frac{n-2}{2}
        \cdot \dot{Y}^\frac{n-2}{2} \
        \hookrightarrow \  L^\infty(\dot{H}^\frac{n-2}{2}) \  ,
        \label{bilin_X_embed1} \\
        \nabla_x\Delta^{-1}\ : \ \dot{Y}^\frac{n-2}{2}
        \cdot \dot{Y}^\frac{n-2}{2} \
        \hookrightarrow \  L^2(\dot{H}^\frac{n-1}{2}) \  .
        \label{bilin_X_embed2}
\end{align}
The first of these embedding follows easily from:
\begin{equation}
        \nabla_x\Delta^{-1}\ : \ L^\infty(\dot{H}^\frac{n-2}{2})\cdot
        L^\infty(\dot{H}^\frac{n-2}{2}) \
        \hookrightarrow \  L^\infty(\dot{H}^\frac{n-2}{2}) \  ,
        \notag
\end{equation}
which in turn follows directly from \eqref{general_besov_embed}. The
second estimate \eqref{bilin_X_embed2} above is more bilinear in
nature. It follows from applying a trichotomy and then summing the
following two fixed frequency bilinear inclusions:
\begin{align}
        \nabla_x\Delta^{-1}\ : \ P_{\bullet \ll \lambda }
        \big(L^2(\dot{B}_2^{\frac{2(n-1)}{n-3},(2,\frac{n-1}{2})})\big)
        \cdot P_\lambda \big(L^\infty(\dot{H}^\frac{n-2}{2})\big)
        \ &\hookrightarrow \  P_\lambda\big(
        L^2(\dot{H}^\frac{n-1}{2})\big) \  .
        \label{bilin_X_embed2a}\\
        \nabla_x\Delta^{-1}\ : \ P_{ \lambda }
        \big(L^2(\dot{B}_2^{\frac{2(n-1)}{n-3},(2,\frac{n-1}{2})})\big)
        \cdot P_\lambda \big(L^\infty(\dot{H}^\frac{n-2}{2})\big)
        \ &\hookrightarrow \  \left(\frac{\mu}
	{\lambda}\right)^\delta P_\mu\big(
        L^2(\dot{H}^\frac{n-1}{2})\big) \  ,
        \label{bilin_X_embed2b}
\end{align}
where we have set $\delta = n(\frac{n-2}{n-1}) -\frac{3}{2}$ to be
the ``gap'' constant. The estimates
\eqref{bilin_X_embed2a}--\eqref{bilin_X_embed2b} follow directly
from the frequency localized bounds
\eqref{freq_loc_general_besov_embed1}--\eqref{freq_loc_general_besov_embed2}.
Note that in this case, the various positivity conditions are
satisfied.\\

\subsection*{$\bullet$\ \ $\dot{Y}^\frac{n-2}{2}\times \dot{Y}^\frac{n-4}{2}$
bounds for the pair $(A_0,\partial_t A_0)$} Our first step here is
to deal with the variable $A_0$. We integrate equation
\eqref{2nd_order_YM5} and write it schematically as:
\begin{equation}
        A_0 \ = \ \Delta^{-1}( \nabla_x [A_0 , \underline{A}] +
        [\underline{A},F]) \ . \label{A0_int_schem}
\end{equation}
The desired estimate now follows by constructing $A_0$ from scratch
by iteration, using the already established estimates and bilinear
embedding \eqref{bilin_to_Y} and the following:
\begin{equation}
        \Delta^{-1} \ : \ \dot{Y}^\frac{n-2}{2}\cdot
        \dot{Y}^\frac{n-4}{2} \ \hookrightarrow \
        \dot{Y}^\frac{n-2}{2} \ . \label{YY_bilin}
\end{equation}
This last embedding follows in turn from the pair of estimates:
\begin{align}
        \Delta^{-1} \ : \
        L^\infty(\dot{H}^\frac{n-2}{2})\cdot
        L^\infty(\dot{H}^\frac{n-4}{2}) \ &\hookrightarrow \
        L^\infty(\dot{H}^\frac{n-2}{2}) \ , \notag\\
        \Delta^{-1} \ : \
        L^2(\dot{B}_2^{\frac{2(n-1)}{n-3},(2,\frac{n-1}{2})} )\cdot
        L^\infty(\dot{H}^\frac{n-4}{2}) \ &\hookrightarrow \
        L^2( \dot{B}_2^{\frac{2(n-1)}{n-3},(2,\frac{n-1}{2})} ) \ . \notag
\end{align}
Both of these are easy consequences of \eqref{general_besov_embed}
and we leave the numerology to the reader.\\

To establish the $\dot{Y}^\frac{n-4}{2}$ bound for $\partial_t A_0$,
we can use the equation \eqref{2nd_order_YM6} to treat it as a
separate variable. In that equation we have quantities of the form
$\partial_t \underline{A}$. We can use the curvature equation
\eqref{2nd_order_YM2} to swap this for spatial derivatives as
follows:
\begin{equation}
        \partial_t\underline{A} \ = \ \nabla_x A_0 -
    [A_0,\underline{A}] + F \ .
        \label{dt_A_trade}
\end{equation}
This allows us to write schematically:
\begin{equation}
        (\partial_t A_0) \ = \ \nabla_x \Delta^{-1}\big(
        [A,(\partial_t A_0)] + [A,\nabla_x A] +
    \big[A, [A,A]\big]
    + [A,F]\big) \ ,
        \label{A0_schematic}
\end{equation}
where $A$ now denotes any of the full set of potentials
$\{A_\alpha\}$ which we have estimated in the space
$\dot{Y}^\frac{n-2}{2}$. We may now iterate the equation
\eqref{A0_schematic} in the space $\dot{Y}^\frac{n-4}{2}$ to
constructively obtain the desired bounds using the bilinear
embedding:
\begin{equation}
        \nabla_x \Delta^{-1} \ : \ \dot{Y}^\frac{n-2}{2}\cdot
        \dot{Y}^\frac{n-4}{2} \ \hookrightarrow \
        \dot{Y}^\frac{n-4}{2} \ . \notag
\end{equation}
which follows from differentiating \eqref{YY_bilin} above. Notice that
the needed inclusion $[A,A]\hookrightarrow \dot{Y}^\frac{n-4}{2}$
follows, for instance, from differentiating the embedding
\eqref{bilin_to_Y}.\\

\subsection*{$\bullet$\ \ Splitting the spatial potentials}
Our next goal is to split the spatial potentials
$\{\underline{A}_i\}$ into a sum of two pieces which are each more
easily managed. This will be done using the ``structure'' equation
\eqref{2nd_order_YM4}. Using the formula \eqref{dt_A_trade} to get
rid of terms of the form $\partial_t \underline{A}$ on the right
hand side of this equation, and using the various $\dot{Y}^s$ space
embeddings we have just shown (on the time interval $[0,T^*]$), we
may write this equation in the schematic form:
\begin{equation}
        \Box \underline{A} \ = \ \mathcal{P}([B,H]) \ ,
        \label{uA_hyp_schematic}
\end{equation}
where the quantities $(B,H)$ obey the estimate:
\begin{equation}
        \lp{(B,H)}{\dot{Y}^\frac{n-2}{2}\times
        \dot{Y}^\frac{n-4}{2}} \ \lesssim \ NC \epsilon_0 \ ,\notag
\end{equation}
where the implicit constant in the above inequality comes from the
estimates just shown. Using Duhamel's principle and (sharp) time
cutoffs, we now extend \eqref{uA_hyp_schematic} to all possible
times. This is done simply by writing:
\begin{equation}
        \underline{A}(t) \ = \ \underline{A}^{(0)}(t) + \int_0^t\
        \frac{\sin((t-s)\sqrt{-\Delta})}{\sqrt{-\Delta}}\
        \mathcal{P}([B,H])(s)\cdot\chi_{[0,T^*]}(s)
        \ ds \ , \label{uA_hyp_schematic_int}
\end{equation}
where $\underline{A}^{(0)}$ denotes to propagation of
$\big(\underline{A}(0), \partial_t\underline{A}\, (0)\big)$ as a
solution to the free scalar wave equation.  Also, here
$\chi_{[0,T^*]}$ denotes the indicator function of the time interval
$[0,T^*]$. This implies that we have the condition:
\begin{align}
        \Box \underline{A}\ (t) \ &= \ 0 \ ,
        &t \ &< \ 0 \ ,
        &T^* \ &< \ t \ . \notag
\end{align}
Now, from the bootstrapping assumption
\eqref{critical_bs_assumed_est_A} we have the pair of bounds:
\begin{align}
        \lp{\big(\underline{A}(0), \partial_t
        \underline{A}\, (0)\big)}{\dot{H}^\frac{n-2}{2}
        \times\dot{H}^\frac{n-4}{2}} \ &\leqslant \ NC\epsilon_0 \ ,
        \notag\\
        \lp{\big(\underline{A}(T^*), \partial_t
        \underline{A}\, (T^*)\big)}{\dot{H}^\frac{n-2}{2}
        \times\dot{H}^\frac{n-4}{2}} \ &\leqslant \ NC\epsilon_0 \ .
        \notag
\end{align}
Therefore, using the bounds we have just shown in conjunction with
the usual Strichartz estimates for the wave equation, we have that
this extension of the potentials $\{\underline{A}_i\}$ satisfies the
bounds:
\begin{equation}
        \lp{\underline{A}}{\dot{X}^\frac{n-2}{2}} \ \lesssim \
        NC \epsilon_0 \ . \notag
\end{equation}
Notice that the angular square function structure inherent in the
$\dot{X}^s$ norms is provided automatically by the fact that the
usual wave equation commutes with the angular cutoffs
$\oPi_\theta$.\\

Our next step to introduce the space--time frequency cutoff
$S_{|\tau|\lesssim |\xi|}$, which cuts off smoothly on the region
$|\tau|\lesssim |\xi|$. That is, the compound multipliers $P_\lambda
S_{|\tau|\lesssim |\xi|}$ all have $L^1$ kernels with uniform
bounds. We denote by $S_{|\xi|\ll|\tau|} = I - S_{|\tau|\lesssim
|\xi|}$. Our decomposition of $\{\underline{A}_i\}$ is now given by
the formula:
\begin{align}
        \td{\underline{A}} \ &= \ S_{|\tau|\lesssim |\xi|}
        {\underline{A}} \ ,
        &\td{\td{\underline{A}}} \ &= \ S_{|\xi|\ll|\tau|}
        \underline{A} \ . \notag
\end{align}
We now need to  show that both the potential sets
$\{\td{\underline{A}_i}\}$ and $\{\td{\td{\underline{A}_i}}\}$ obey
good $\dot{X}^\frac{n-2}{2}$ estimates. Since the original
collection of extended potentials does, we only need to prove this
assertion for one of these sets. This is most easily shown for the
collection $\{\td{\underline{A}_i}\}$. As we have already mentioned,
the cutoffs $P_\lambda S_{|\tau|\lesssim |\xi|}$ are bounded on all
mixed Lebesgue spaces. Therefore, the entire multiplier
$S_{|\tau|\lesssim |\xi|}$ is bounded on any mixed Lebesgue-Besov
space of the type $L^q(\dot{B}_2^{p,(2,s)})$. This implies that this
multiplier is in fact bounded on the $\dot{X}^s$ spaces, which is
enough to support our claim.\\

Finally, we would like to prove two fixed frequency multiplier
estimates which will be useful in the sequel when dealing with the
two sets of potentials $\{\td{\underline{A}_i}\}$ and
$\{\td{\td{\underline{A}_i}}\}$. The first is:
\begin{align}
        \lp{\partial_t P_\lambda S_{|\tau|\lesssim |\xi|} A}{L^p} \
        \lesssim \ &\lambda \ \lp{A}{L^p} \
        &1\leqslant \ p \ \leqslant \ \infty \ . \label{dt_to_dx}
\end{align}
This is easily demonstrated by rescaling to frequency $\lambda = 1$
and using the $L^1$ bound on the convolution kernel of $\partial_t
P_\lambda S_{|\tau|\lesssim |\xi|}$. Combining this with the remarks
made above, we see that we have the estimate:
\begin{equation}
        \lp{\partial_t \td{\underline{A}}}{\dot{X}^\frac{n-4}{2}}
        \ \lesssim \ NC\epsilon_0 \ . \notag
\end{equation}
In particular, from everything we have shown, the potential set
$\{\td{\underline{A}_i}\}$ satisfies all of the requirements
\eqref{red_conctn_cond} of Theorem \ref{Str_thm} when $\epsilon_0$
is sufficiently small.\\

The second fixed frequency multiplier bound that will be of use
shortly is the space--time estimate:
\begin{equation}
        \lp{\varXi^{-1} P_\lambda S_{|\xi|\ll|\tau|} A}{L^q(L^p)}
        \ \lesssim \ \lambda^{-2}\ \lp{A}{L^q(L^p)} \ .
        \label{wave_inv_st_bound}
\end{equation}
Here $\varXi$ is the multiplier with symbol $\varXi(\tau,\xi)=\tau^2
- |\xi|^2$. To prove this, we employ a family of Littlewood-Paley
space-time cutoffs which we denote by $S_\mu$. By this we mean that
the space-time frequency support of these is supported where $|\tau|
+ |\xi| \sim \mu$. As usual, these are all chosen so as to have
uniform $L^1$ bounds on their convolution kernels. Using the support
restrictions of the $S_{|\xi|\ll|\tau|}$ multiplier, we have the
formula:
\begin{equation}
        P_\lambda S_{|\xi|\ll|\tau|} A \ = \ \sum_{\substack{\mu \ : \\
        \lambda \lesssim \mu}}\ P_\lambda S_\mu S_{|\xi|\ll|\tau|} A
        \ . \notag
\end{equation}
Therefore, by dyadic summing and the boundedness of the multiplier
$P_\lambda $, to prove \eqref{wave_inv_st_bound} it suffices to be
able to show that:
\begin{equation}
        \lp{\varXi^{-1} S_{|\xi|\ll|\tau|} S_\mu A}{L^q(L^p)}
        \lesssim \ \mu^2\ \lp{A}{L^q(L^p)} \ . \notag
\end{equation}
This last bound follows easily from rescaling to frequency $\mu=1$
and the appropriate differential bounds on the symbol of
$\varXi^{-1} S_{|\xi|\ll|\tau|}$ which we leave to the reader.\\

\subsection*{$\bullet$\ \ $L^1(L^\infty)$ bounds for the potentials
$\{\td{\td{A}}_\alpha\} = (A_0 , \{\td{\td{\underline{A}}}\})$} Our
goal here is to show the $\ell^1$ type Besov estimate:
\begin{equation}
        \lp{(A_0 , \{\td{\td{\underline{A}}}\})}
        {  L^1( \dot{B}_1^{\infty,(2,\frac{n}{2})} ) } \ \lesssim \
        NC\epsilon_0 \ . \label{imp_L1_est}
\end{equation}
By repeatedly using the estimate \eqref{wave_inv_st_bound}, we have
that the multiplier $\varXi^{-1} \Delta S_{|\xi|\ll|\tau|}$ is
bounded on the space $L^1( \dot{B}_1^{\infty,(2,\frac{n}{2})} )$.
Furthermore, from all of the estimates we have shown above, and by
distributing the derivative in the first term on the right hand side
of \eqref{A0_int_schem}, we see that the right hand side of the
schematics \eqref{A0_int_schem} and \eqref{uA_hyp_schematic} are
equivalent. Therefore, we have the following heuristic schematic for
the potentials $\{\td{\td{A}}_\alpha\}$:
\begin{equation}
        \td{\td{A}} \ = \ \Delta^{-1} ([B,H]) \ , \notag
\end{equation}
where the pair $(B,H)$ enjoys the bounds:
\begin{equation}
        \lp{(B,H)}{\dot{Y}^\frac{n-2}{2}\times
        \dot{Y}^\frac{n-4}{2}} \ \lesssim \
        NC\epsilon_0 \ . \notag
\end{equation}
The bound \eqref{imp_L1_est} now follows from the bilinear estimate:
\begin{equation}
        \Delta^{-1} \ : \
        \dot{Y}^\frac{n-2}{2}\cdot
        \dot{Y}^\frac{n-4}{2} \ \hookrightarrow \
        L^1( \dot{B}_1^{\infty,(2,\frac{n}{2})} ) \ . \notag
\end{equation}
This in turn follows from the product estimate:
\begin{equation}
        \Delta^{-1} \ : \
        L^2(\dot{B}_2^{\frac{2(n-1)}{n-3},
	(2,\frac{n-1}{2})  })\cdot
        L^2( \dot{B}_2^{\frac{2(n-1)}{n-3},(2,\frac{n-3}{2})  
	}  ) \ \hookrightarrow \
        L^1( \dot{B}_1^{\infty,(2,\frac{n}{2})} ) \ . \notag
\end{equation}
This last estimate follows at once from \eqref{general_besov_embed}.
The check on the conditions \eqref{sc_cond}--\eqref{Lb_cond} is left
to the reader.\\

\subsection*{$\bullet$\ \ Improving the curvature}
This is the final part of the proof of Proposition
\ref{critical_bs_prop}. Recalling the schematic
\eqref{wave_schematic} and using the Strichartz estimates
\eqref{the_strichartz_est}, our goal here is to show the following
four bounds:
\begin{align}
        \lp{[\nabla\td{\td{A}},F]}{L^1(\dot{H}^\frac{n-6}{2})} \
        &\lesssim \ N^2C^2 \epsilon_0^2 \ , \label{F_source_bound1}\\
        \lp{  [\td{\td{A}},\nabla F]}{L^1(\dot{H}^\frac{n-6}{2})} \
        &\lesssim \ N^2C^2 \epsilon_0^2 \ , \label{F_source_bound2}\\
        \lp{  \big[\td{A},[\td{\td{A}} ,F]\big] }{L^1(\dot{H}^\frac{n-6}{2})} \
        &\lesssim \ N^2C^2 \epsilon_0^2 \ , \label{F_source_bound3}\\
        \lp{  \big[\td{\td{A}},[\td{\td{A}}
        ,F]\big] }{L^1(\dot{H}^\frac{n-6}{2})} \
        &\lesssim \ N^2C^2 \epsilon_0^2 \ , \label{F_source_bound3.5}\\
        \lp{[F,F]}{L^1(\dot{H}^\frac{n-6}{2})} \
        &\lesssim \ N^2C^2 \epsilon_0^2 \ . \label{F_source_bound4}
\end{align}
For $\epsilon_0$ sufficiently small, this will be enough for us to
conclude the improved bootstrapping estimates
\eqref{critical_bs_imp_est} by choosing $C$ to be such that
$\frac{1}{2} (LN)^{-1} C$ is equal to the constant appearing on the
right hand side of estimate \eqref{the_strichartz_est}. This works
because the implicit constants which appear in
\eqref{F_source_bound1}--\eqref{F_source_bound4} above have only
been manufactured in the estimates of this proof, and can all be
chosen to be independent of $N$ and $C$ if $\epsilon_0$ is chosen
small enough.\\

To prove these bounds, first notice that the estimates
\eqref{F_source_bound1} and
\eqref{F_source_bound3}--\eqref{F_source_bound4} are essentially
identical. This follows from the equivalence (in terms of
$\dot{Y}^s$ spaces) $\nabla \td{\td{A}} \approx F$. We also have
the equivalences $[\td{A},\td{\td{A}}] \approx F$ and
$[\td{\td{A}},\td{\td{A}}] \approx F$. These are given by
the inclusion:
\begin{equation}
        \dot{Y}^\frac{n-2}{2}\cdot \dot{Y}^\frac{n-2}{2} \ \subseteq
        \ \dot{Y}^\frac{n-4}{2} \ . \label{YY_incl}
\end{equation}
This is easily demonstrated, as we have already mentioned,
by differentiating the
inclusion \eqref{bilin_to_Y} and using the boundedness of
$\nabla^2_x\Delta^{-1}$ on the various $\dot{Y}^s$ component spaces.
Therefore, to prove \eqref{F_source_bound1} and
\eqref{F_source_bound3}--\eqref{F_source_bound4} we only need to
know that:
\begin{equation}
        L^2( \dot{B}_2^{\frac{2(n-1)}{n-3},(2,\frac{n-3}{2})} )
        \cdot L^2( \dot{B}_2^{\frac{2(n-1)}{n-3},(2,\frac{n-3}{2})}  )
        \ \hookrightarrow \ L^1(\dot{H}^\frac{n-6}{2}) \ .
    \label{basic_FF_embed}
\end{equation}
This is yet again a consequence of our general Besov calculus
\eqref{general_besov_embed}, and we leave the various additions to
the reader.\\

Our final task here is to prove the estimate
\eqref{F_source_bound2}. This needs to be frequency decomposed using
a trichotomy. Specifically, we have the following set of fixed
frequency estimates in the three cases (note that in the first two
estimates below the square summing needs to be done \emph{inside}
the time integral):
\begin{align}
        P_{\bullet\ll\lambda}\big(L^1(\dot{B}_1^{\infty,(n,\frac{n}{2})} )\big)
        \cdot P_\lambda \big( L^\infty(\dot{H}^\frac{n-6}{2}) \big)
        \ &\hookrightarrow \ P_\lambda \big( L^1(\dot{H}^\frac{n-6}{2}) \big)
        \ , \label{last_AF_est1}\\
        P_\lambda \big( L^2( \dot{B}_2^{\frac{2(n-1)}{n-3},(2,\frac{n-1}{2})}
        \big)\cdot
        P_{\bullet\ll\lambda}\big(L^2 (  
	\dot{B}_2^{\frac{2(n-1)}{n-3},(2,\frac{n-3}{2})} )
        \big) \ &\hookrightarrow \ P_\lambda 
	\big( L^1(\dot{H}^\frac{n-6}{2}) \big)
        \ , \label{last_AF_est2}\\
        P_\lambda \big( L^2( \dot{B}_2^{\frac{2(n-1)}{n-3},(2,\frac{n-1}{2})}
        \big)\cdot
        P_{\lambda}\big(L^2 (  \dot{B}_2^{\frac{2(n-1)}
	{n-3},(2,\frac{n-3}{2})} )
        \big) \ &\hookrightarrow \ \left(\frac{\mu}{\lambda}
        \right)^\delta P_\mu \big( L^1(\dot{H}^\frac{n-6}{2}) \big)
        \ , \label{last_AF_est3}
\end{align}
where the quantity $\delta$ in the last estimate
\eqref{last_AF_est3} above can be computed to be $\delta =
n(\frac{n-3}{n-1}) -3$. The estimate \eqref{last_AF_est1} follows
from inspection. The latter two estimates
\eqref{last_AF_est2}--\eqref{last_AF_est3} follow from
\eqref{freq_loc_general_besov_embed1}--\eqref{freq_loc_general_besov_embed1}
of Remark \ref{freq_loc_rem}. This completes the proof of
Proposition \ref{critical_bs_prop}.
\end{proof}\ret

\ret

\section{Reduction to Approximate Half-Wave Operators}

This is a  preliminary technical section where we reduce the
proof of the Strichartz estimates \eqref{the_strichartz_est} to a
more easily managed form. This material more or less
standard, and we again follow closely what was done
in \cite{RT_MKG}. Our first step here is
to reduce the proof of Theorem \ref{Str_thm} to the following:\\

\begin{prop}[Existence of a fixed frequency parametrix]\label{freq1_parametrix_prop}
Let the number of dimensions be $6\leqslant n$, and let $d+\uAl$ be a
connection which satisfies the conditions
\eqref{red_conctn_cond}. In addition assume that we have the frequency
localization condition:
\begin{equation}
    P_{\lambda\lesssim \bullet}( \uAl) \ = \ 0 \ , \label{freq_loc_cond}
\end{equation}
where $P_{\lambda\lesssim \bullet}$ is a  frequency cutoff on the
region where $2^{-10a}\lambda\leqslant |\xi|$, where $1 \leqslant a$
is some fixed parameter. Then if the constant $\mathcal{E}$ on lines
\eqref{red_conctn_cond5} and \eqref{red_conctn_cond7} is
sufficiently small, there exists a family of  approximate
propagation operators $W^\lambda_{\uAl}(s)$ (or just $W^\lambda_s$
for short) such that if $(f_\lambda,g_\lambda)$ is any set of
$\lambda$--frequency initial data with Fourier support in the region $
2^{-a}\lambda \leqslant |\xi| \leqslant 2^a \lambda$, the following
estimates hold:
\begin{subequations}\label{parametrix_estimates}
\begin{align}
        \lp{W^\lambda_s
        (f_\lambda,g_\lambda)}{\dot{X}^0\times
        \partial^{-1}_t(\dot{X}^{-1})} \
        &\lesssim \  E^\frac{1}{2}(f_\lambda,g_\lambda)
        \ , \label{parametrix_estimate1}\\
        \lp{W^\lambda_s(f_\lambda,g_\lambda)(s) -
        f_\lambda}{L^2} \ &\lesssim \
        \mathcal{E}^\frac{1}{2}\, E^\frac{1}{2}
        (f_\lambda,g_\lambda)  \ ,
        \label{parametrix_estimate2}\\
        \lp{\partial_t W^\lambda_s
        (f_\lambda,g_\lambda)(s) - g_\lambda}{L^2} \ &\lesssim \ \lambda\,
        \mathcal{E}^\frac{1}{2}\, E^\frac{1}{2}(f_\lambda,g_\lambda) \ ,
        \label{parametrix_estimate3}\\
        \lp{ \Box_{\uAl}W^\lambda_s
        (f_\lambda,g_\lambda) }{L^1(L^2)} \ &\lesssim
        \ \lambda\, \mathcal{E} \, E^\frac{1}{2}(f_\lambda,g_\lambda)  \ . \label{parametrix_estimate4}
\end{align}
\end{subequations}
Here we have set $E(f_\lambda,g_\lambda)$ to the $L^2$ normalized
energy:
\begin{equation}
        E(f_\lambda,g_\lambda) \ = \
        \lp{f_\lambda}{L^2}^2 + \lambda^{-2}\, \lp{g_\lambda}{L^2}^2 \ .
        \notag
\end{equation}
Finally, we have that the frequency support of the parametrix
is contained in the set $2^{-2a}\lambda \leqslant |\xi| \leqslant
2^{2a}\lambda$, where $a$ is as above.
\end{prop}\ret

\begin{proof}[Proof that Proposition \ref{freq1_parametrix_prop}
implies Theorem \ref{Str_thm}] The first step here is to reduce the
estimate \eqref{the_strichartz_est} to the case where $G\equiv 0$.
This is done in the usual way via Duhamel's principle. We define the
true propagation operator $U_s(t)$ via the formulas:
\begin{align}
        U_s(s)(f,g) \ &= \ f \ ,
        &\partial_t U_s(s)(f,g) \ &= \ g \ , \notag
\end{align}
and:
\begin{equation}
        \Box_{\underline{A}} U_s(f,g) \ = \ 0 \ , \notag
\end{equation}
We then have that:
\begin{equation}
        F(t) \ = \ U_0(t)(f,\dot{f}) \ + \ \int_0^t\ U_s(t)(0,G(s))\
        ds \ , \label{Duhamel_formula}
\end{equation}
solves the problem \eqref{inhomog_gc_wave}--\eqref{general_cuachy_data}.
In particular, by Minkowski's triangle inequality we easily have that:
\begin{equation}
        \lp{\int_0^t\ U_s(t)(0,G(s))\
        ds}{\dot{X}^\frac{n-4}{2}\times
        \partial_t^{-1}(\dot{X}^\frac{n-6}{2})} \ \leqslant \
        \int_0^\infty \
        \lp{U_s(0,G(s))}{\dot{X}^\frac{n-4}{2}\times\partial_t^{-1}
        (\dot{X}^\frac{n-6}{2})}
        \ ds \ . \notag
\end{equation}
Therefore, we  are trying to show:
\begin{equation}
        \lp{U_s(f,g)}{\dot{X}^\frac{n-4}{2}
    \times\partial_t^{-1}(\dot{X}^\frac{n-6}{2})}
        \ \leqslant \ C\,  \lp{(f,g)}{\dot{H}^\frac{n-4}{2}
    \times\dot{H}^\frac{n-6}{2}} \ ,
        \label{strichartz_to_prove}
\end{equation}
for any pair of functions $(f,g)$ and any initial time $s$. Since it
is easy to see that the conditions \eqref{red_conctn_cond} are
translation invariant, it suffices to show this estimates for
$s=0$.\\

The estimate \eqref{strichartz_to_prove} will be shown using a
bootstrapping procedure. This will be done inside of the compact
intervals $[0,T^*]$. What we will do is to first assume that
\eqref{strichartz_to_prove} is true  for all $0\leqslant s \leqslant
T^*$ on all time intervals of the form $[0,s]$ and $[s,T^*]$, where
the constant on the left hand side of \eqref{strichartz_to_prove} is
replaced by $2C$.  Our goal is then to improve the constant by
proving the desired bound \eqref{strichartz_to_prove} on the time
subintervals of $[0,T^*]$. Once this is accomplished, we can easily
extend the bound $\eqref{strichartz_to_prove}$ to all subintervals
of a slightly larger time interval $[0,T^* + \gamma]$, where the
constant $0 < \gamma \ll 1$ is determined by the bound
\eqref{general_smoothness}. This is provided by the usual local
existence theory based on energy and $L^\infty$ estimates. Once this
is done, the bootstrapping closes. Notice again that, by using the
local existence theory and the bound \eqref{general_smoothness}, we
may begin the argument for some very small time interval
$[0,\gamma]$.\\

We are now assuming that \eqref{strichartz_to_prove} holds on our
time interval $[0,T^*]$ with constant $2C$ which we will decide
on in a moment. We are working with a solution:
\begin{equation}
    \Box_{{\underline{A}}} F \ = \ 0 \ , \label{homog_cov_wave}
\end{equation}
where the connection $d+\underline{A}$ satisfies
\eqref{red_conctn_cond},
and where we have the initial data:
\begin{align}
    F(0) \ &= \ f ,
    &\partial_t F\, (0) \ = \ g \ . \label{homog_cov_wave_data}
\end{align}
We now split this initial data into a sum frequency localized pieces:
\begin{align}
    f \ &= \ \sum_\lambda\ P_\lambda(f) \ = \ \sum_\lambda\
    f_\lambda \ , \notag\\
    g \ &= \ \sum_\lambda\ P_\lambda(g) \ = \ \sum_\lambda\
    g_\lambda \ , \notag
\end{align}
and then repeatedly use Proposition \ref{freq1_parametrix_prop}
to construct an approximate solution to
\eqref{homog_cov_wave}--\eqref{homog_cov_wave_data} as follows:
\begin{equation}
    \td{F} \ = \ \sum_\lambda \ \td{F}_\lambda
    \ = \ \sum_\lambda\ W^\lambda_0(f_\lambda,g_\lambda)
    \ . \notag
\end{equation}
By summing over the parametrix estimate \eqref{parametrix_estimate1} we
automatically have that:
\begin{equation}
    \lp{\td{F}}{\dot{X}^\frac{n-4}{2}\times\partial_t^{-1}(\dot{X}^\frac{n-6}{2})}
    \ \leqslant \ \frac{1}{2}C \
    \lp{(f,g)}{\dot{H}^\frac{n-4}{2}
    \times\dot{H}^\frac{n-6}{2}} \ , \notag
\end{equation}
where $C$ is some fixed constant. We choose this to be our
definition of the constant on the right hand side of
\eqref{strichartz_to_prove}. Thus, our goal is to conclude that:
\begin{equation}
    \lp{F - \td{F}}{\dot{X}^\frac{n-4}{2}\times\partial_t^{-1}(\dot{X}^\frac{n-6}{2})}
    \ \leqslant \ \frac{1}{2}C \
    \lp{(f,g)}{\dot{H}^\frac{n-4}{2}
    \times\dot{H}^\frac{n-6}{2}} \ . \label{difference_est_to_show}
\end{equation}
To do this, we use the Duhamel formula \eqref{Duhamel_formula} to
express everything in terms of the operators $U_s(t)$:
\begin{equation}
        F(t) - \td{F}(t) \ = \
        U_0(t)\big(f - \td{F}(0), g - \partial_t
        \td{F}\, (0)\big) -
        \int_0^t\ U_s(t)\big(0 , \Box_{\underline{A}}
    \td{F}\, (s) \big) \ ds \ . \notag
\end{equation}
By combining the assumed estimate \eqref{strichartz_to_prove} and
the approximation bounds
\eqref{parametrix_estimate2}--\eqref{parametrix_estimate3}, we have
that:
\begin{equation}
    \lp{U_0\big(f - \td{F}(0), g - \partial_t
        \td{F}\, (0)\big)
    }{\dot{X}^\frac{n-4}{2}\times\partial_t^{-1}(\dot{X}^\frac{n-6}{2})}
    \ \lesssim \ C\mathcal{E}^\frac{1}{2} \,
    \lp{(f,g)}{\dot{H}^\frac{n-4}{2}
    \times\dot{H}^\frac{n-6}{2}} \ . \notag
\end{equation}
Therefore, by using Minkowski's triangle inequality and again using
the bootstrapping assumption \eqref{strichartz_to_prove}, we see
that in order to conclude \eqref{difference_est_to_show} we only
need to show the following remainder estimate on the time interval
$[0,T^*]$:
\begin{equation}
    \lp{\Box_{\underline{A}} \td{F}}{L^1(\dot{H}^\frac{n-6}{2})}
    \ \lesssim \ C\mathcal{E}\, \lp{(f,g)}{\dot{H}^\frac{n-4}{2}
    \times\dot{H}^\frac{n-6}{2}} \ . \label{L1L2_last_piece}
\end{equation}\ret

To show the estimate \eqref{L1L2_last_piece}, we use a family of
frequency cutoffs:
\begin{equation}
    I \ = \ P_{\bullet \ll \lambda} +
    P_{\lambda\lesssim \bullet} \ , \notag
\end{equation}
for each scale $\lambda$ such that they all have $L^1$ kernels with
uniform bounds, and such that the cutoff $P_{\bullet \ll \lambda}$
is consistent with the definition of $d+\uAl$ in the statement of
Proposition \ref{freq1_parametrix_prop}. This allows us to
schematically write:
\begin{multline}
    \Box_{\underline{A}}\td{F}\ = \
    \sum_\lambda \Big( \Box_{\uAl}\td{F}_\lambda
    + [\nabla_x \uAlg , \td{F}_\lambda] +
     [\uAlg , \nabla_x \td{F}_\lambda] \\
     + \ \big[[\uAl,\uAlg], \td{F}_\lambda\big] +
    \big[[\uAlg,\uAlg], \td{F}_\lambda\big]\Big)
    \ . \label{tdF_big_sum}
\end{multline}
The bound \eqref{L1L2_last_piece} for the term $\sum_\lambda
\Box_{\uAl}\td{F}_\lambda$ is a direct consequence of repeatedly
applying the estimate \eqref{parametrix_estimate4} while using the
fact that each term in this sum is supported in frequency where
$|\xi|\sim \lambda$ to gain the orthogonality needed to obtain
bounds in terms of the pair $(f,g)$. Therefore, we are reduced to
showing the following family of error estimates:
\begin{align}
    \sum_\lambda\ \lp{[\nabla_x \uAlg , \td{F}_\lambda]}
    {L^1(\dot{H}^\frac{n-6}{2})} \ \lesssim \
    \mathcal{E}\, \lp{(f,g)}{\dot{H}^\frac{n-4}{2}
    \times\dot{H}^\frac{n-6}{2}} \ , \label{L1L2_error1}\\
    \sum_\lambda\ \lp{[\uAlg , \nabla_x \td{F}_\lambda] }
    {L^1(\dot{H}^\frac{n-6}{2})} \ \lesssim \
    \mathcal{E}\, \lp{(f,g)}{\dot{H}^\frac{n-4}{2}
    \times\dot{H}^\frac{n-6}{2}} \ , \label{L1L2_error2}\\
     \lp{\sum_\lambda\ \ \big[[\uAl,\uAlg],
    \td{F}_\lambda\big] }
    {L^1(\dot{H}^\frac{n-6}{2})} \ \lesssim \
    \mathcal{E}\, \lp{(f,g)}{\dot{H}^\frac{n-4}{2}
    \times\dot{H}^\frac{n-6}{2}} \ , \label{L1L2_error3}\\
    \lp{\sum_\lambda\ \big[[\uAlg,\uAlg], \td{F}_\lambda\big] }
    {L^1(\dot{H}^\frac{n-6}{2})} \ \lesssim \
    \mathcal{E}\, \lp{(f,g)}{\dot{H}^\frac{n-4}{2}
    \times\dot{H}^\frac{n-6}{2}} \ . \label{L1L2_error4}
\end{align}
These estimates are all very similar to each other, and to estimates
we have already proved in the last section, in particular
\eqref{F_source_bound1}--\eqref{F_source_bound4}. To prove the first
estimate \eqref{L1L2_error1} above, we further decompose the left
hand side into frequencies and use the triangle inequality to bound:
\begin{equation}
    \hbox{(L.H.S.)}\eqref{L1L2_error1} \ \leqslant \
    \sum_{\substack{\lambda , \mu \ : \\ \lambda\lesssim \mu}}\
    \lp{[\nabla_x P_\mu(\underline{A}) ,  \td{F}_\lambda] }
    {L^1(\dot{H}^\frac{n-6}{2})} \ . \notag
\end{equation}
Thus, by Young's inequality, it suffices to show the following
family of fixed frequency estimates:
\begin{equation}
    \lp{[\nabla_x P_\mu(\underline{A}) , \td{F}_\lambda] }
    {L^1(\dot{H}^\frac{n-6}{2})} \ \lesssim \
    \left(\frac{\lambda}{\mu}\right)^\delta\
    \lp{P_\mu(\underline{A})}{\dot{Z}^\frac{n-2}{2}}\cdot
    \lp{(f_\lambda,g_\lambda)}{\dot{H}^\frac{n-4}{2}
    \times\dot{H}^\frac{n-6}{2}} \ , \notag
\end{equation}
where we have set $\delta = \frac{3}{2} - \frac{n}{n-1}$. Notice
that we have used the $\dot{Z}^\frac{n-2}{2}$ norm for the
$\{\underline{A}_i\}$ on the right hand side. This allows us to
reconstruct norms through square-summing. For $\lambda\sim\mu$ this
estimate is nothing but a fixed frequency version of the estimate
\eqref{basic_FF_embed} above, so it suffices to consider case
$\lambda \ll \mu$. Using the simple inclusion
$\nabla_x\dot{X}^\frac{n-2}{2}\subseteq \dot{X}^\frac{n-4}{2}$, this
is a consequence of the fixed frequency embedding:
\begin{equation}
    P_{\mu}\big(L^2(\dot{B}_2^{\frac{2(n-1)}{n-3}
    ,(2,\frac{n-3}{2})}) \big)\cdot
     P_\lambda\big(L^2(
    \dot{B}_2^{\frac{2(n-1)}{n-3},(2,\frac{n-3}{2})} )\big)
    \ \hookrightarrow \ \left(\frac{\lambda}{\mu}
    \right)^\delta
    \ L^1(\dot{H}^\frac{n-6}{2}) \ , \label{fixed_freq_bilin_mu}
\end{equation}
which follows at once from the fixed frequency estimate
\eqref{gen_besov_fixed_freq} which helps to generate the general
estimate \eqref{general_besov_embed}. Notice that the proof of the
second estimate \eqref{L1L2_error2} above is very similar to what we
have just done. In fact, there is more room because the derivative
is on the low frequency term. We leave the details to the reader.\\

It remains to prove the two estimates
\eqref{L1L2_error3}--\eqref{L1L2_error4}. Since these follow from
essentially identical reasoning, we concentrate on proving the
second of these estimates. This one in fact requires a bit more work
than the fist because it has more frequency overlap. Applying a
trichotomy to the product, we see that it suffices to be able to
show the following three estimates:
\begin{align}
     \begin{split}
    &\int_0^{T^*}\ \big( \sum_\lambda\
    \big( \sum_{\substack{\mu \ :\\ \mu\ll \lambda}}\
    \lp{\big[P_\mu([\uAlg,\uAlg]), \td{F}_\lambda\big](s) }
    {\dot{H}^\frac{n-6}{2}}\big)^2\big)^\frac{1}{2}\ ds \\
    &\ \ \ \ \ \ \ \ \ \ \ \ \ \ \
    \ \ \ \ \ \ \ \ \ \ \ \ \ \ \ \ \ \ \ \ \ \ \ \ \ \
    \lesssim \ \lp{\underline{A}}{\dot{X}^\frac{n-2}{2}}^2\cdot
    \lp{(f,g)}{\dot{H}^\frac{n-4}{2}
    \times\dot{H}^\frac{n-6}{2}} \ ,
     \end{split}\label{cubic_tric1}\\
     \begin{split}
    &\int_0^{T^*}\ \big( \sum_\mu\
    \big( \sum_{\substack{\lambda \ :\\ \lambda\ll \mu}}\
    \lp{\big[P_\mu([\uAlg,\uAlg]), \td{F}_\lambda\big](s) }
    {L^1(\dot{H}^\frac{n-6}{2})}\big)^2\big)^\frac{1}{2}\ ds \\
    &\ \ \ \ \ \ \ \ \ \ \ \ \ \ \
    \ \ \ \ \ \ \ \ \ \ \ \ \ \ \   \ \ \ \ \ \ \ \ \ \ \
    \lesssim \ \lp{\underline{A}}{\dot{X}^\frac{n-2}{2}}^2\cdot
    \lp{(f,g)}{\dot{H}^\frac{n-4}{2}
    \times\dot{H}^\frac{n-6}{2}} \ ,
     \end{split}\label{cubic_tric2}\\
     \begin{split}
        &\sum_{\substack{\lambda,\mu\ : \\ \lambda\sim \mu}}\
    \lp{\big[P_\mu([\uAlg,\uAlg]), \td{F}_\lambda\big] }
    {L^1(\dot{H}^\frac{n-6}{2})} \\
    &\ \ \ \ \ \ \ \ \ \ \ \ \ \ \
    \ \ \ \ \ \ \ \ \ \ \ \ \ \ \   \ \ \ \ \ \ \ \ \ \ \
    \lesssim \
    \lp{\underline{A}}{\dot{X}^\frac{n-2}{2}}^2\cdot
    \lp{(f,g)}{\dot{H}^\frac{n-4}{2}
    \times\dot{H}^\frac{n-6}{2}} \ .
     \end{split}\label{cubic_tric3}
\end{align}
The first two estimates \eqref{cubic_tric1}--\eqref{cubic_tric2}
follow from first fixing time and then proving the fixed frequency
estimate:
\begin{multline}
    \lp{\big[P_\mu([\uAlg,\uAlg]), \td{F}_\lambda\big](s) }
    {\dot{H}^\frac{n-6}{2}} \\
    \lesssim \ \
    \min_\pm \left(\frac{\lambda}{\mu}\right)^{\pm\delta}\
    \lp{  P_\mu([\uAlg,\uAlg])(s) }
    {\dot{B}^{\frac{2(n-1)}{n-3},(2,\frac{n-3}{2})}}
    \cdot
    \lp{\td{F}_\lambda(s)}
    {\dot{B}^{\frac{2(n-1)}{n-3},(2,\frac{n-3}{2})}
    } \  , \notag
\end{multline}
where $\delta$ is the same constant from estimate
\eqref{fixed_freq_bilin_mu}. Indeed, this last line follows from the
non-time integrated version of that estimate. Applying Young's
inequality to this, integrating in time and applying
Cauchy-Schwartz, using the parametrix bound
\eqref{parametrix_estimate1}, the product embedding \eqref{YY_incl},
and the fact that for each fixed value of $\lambda$ the multipliers
$P_{\bullet\ll\lambda}$ and $P_{\lambda \lesssim \bullet}$ are
bounded on the $\dot{X}^s$ spaces we arrive at the
desired pair of estimates.\\

It remains for us to prove the last estimate \eqref{cubic_tric3}
above. After another application of the embedding
\eqref{basic_FF_embed} and a Cauchy-Schwartz, followed by the
parametrix estimate \eqref{parametrix_estimate1}, we are left with
showing the bound:
\begin{equation}
    \big( \sum_{\substack{\lambda,\mu\ : \\ \lambda\sim \mu}}\
    \lp{\big[P_\mu([\uAlg,\uAlg])}{L^2(
    \dot{B}^{\frac{2(n-1)}{n-3},(2,\frac{n-3}{2})})
    }^2 \big)^\frac{1}{2}
    \ \lesssim \ \lp{\underline{A}}
    {\dot{X}^\frac{n-2}{2}}^2 \ . \notag
\end{equation}
This last estimate follows from applying a further trichotomy, and
then using Young's inequality after reduction to the various fixed
frequency versions of the product estimate \eqref{YY_incl} which are
provided by the general fixed frequency estimates
\eqref{freq_loc_general_besov_embed1}--\eqref{freq_loc_general_besov_embed1}.
We leave the details to the diligent reader. This completes the
proof of our reduction of Theorem \ref{Str_thm} to Proposition
\ref{freq1_parametrix_prop}.
\end{proof}\ret

The final thing we will do in this section is to make one further
reduction of the Strichartz estimates \eqref{the_strichartz_est}.
This involves the following proposition:\\

\begin{prop}[Existence of approximate half-wave
parametrices]\label{half_wave_prop} Let the number of dimensions be
$6\leqslant n$, and let $d+\uA$ be a connection which satisfies the
conditions \eqref{red_conctn_cond} as well as the frequency
localization condition \eqref{freq_loc_cond} for $\lambda=1$. Then
there exists pair of evolution operators $\Phi^\pm(\widehat{f})(t)$
from $L^2(\mathbb{R}^n_\xi)$ to $L^2(\mathbb{R}^n_x)$
such that the fixed time adjoints
$(\Phi^\pm(t))^*$ are always supported in the region
$2^{-a}\leqslant  |\xi|\leqslant 2^a$ for some fixed $1 \leqslant
a$, and such that they obey the following estimates:
\begin{subequations}\label{half_wave_est}
\begin{align}
    \lp{\Big( P_1 \Phi^\pm(\widehat{f}) ,
    \Phi^\pm(\widehat{f})
    \Big)}{\dot{X}^0\times L^2_x }
    \ &\lesssim \ \lp{\widehat{f}}{L^2_\xi} \ ,
    \label{half_wave_est1}\\
    \lp{\nabla_x \Phi^\pm(\widehat{f}) }{L^2_t(L^\frac{2(n-1)}{n-3}_x)}
    \ &\lesssim \ \lp{\widehat{f}}{L^2_\xi} \ ,
    \label{half_wave_est1.5}\\
    \lp{\partial_t P_1 \Phi^\pm(\widehat{f}) \mp
     P_1 \Phi^\pm(2\pi i|\xi| \widehat{f})}{\dot{X}^0}
    \ &\lesssim \ \mathcal{E}\ \lp{\widehat{f}}{L^2_\xi} \ ,
    \label{half_wave_est2}\\
    \lp{\Phi^\pm(0)\big((2\pi  |\xi|)^\alpha (\Phi^\pm(0))^*\big)
    g - (-\Delta)^\frac{\alpha}{2} P_1(g)}{L^2_x}
    \ &\lesssim \ \mathcal{E}^\frac{1}{2} \ \lp{g}{L^2_x} \ ,
    \label{half_wave_est3}\\
    \lp{\Box_{\uA} \Phi^\pm(\widehat{f}) }{L^1_t(L^2_x)}
    \ &\lesssim \ \mathcal{E}\ \lp{\widehat{f}}{L^2_\xi}
    \ . \label{half_wave_est4}
\end{align}
\end{subequations}
\end{prop}\ret

\begin{proof}[Proof that Proposition \ref{half_wave_prop} implies
Proposition \ref{freq1_parametrix_prop}] This is a simple matter,
and we explain it briefly. Notice first that it suffices to prove
Proposition \ref{freq1_parametrix_prop} on the scale $\lambda=1$
because everything in sight is scale invariant. We now let
$(f_1,g_1)$ be any pair of unit frequency initial data, and we
define the approximate unit frequency wave propagator:
\begin{multline}
    W_0^1(f_1,g_1)(t) \ = \ P_1\Big(\ \frac{1}{2}\Phi^+(t) (\Phi^+(0))^*f_1
    + \frac{1}{2}\Phi^-(t) (\Phi^-(0))^*f_1 \\
    + \ \ \Phi^+(t) \big(\frac{1}{4\pi i|\xi|}(\Phi^+(0))^*\big)g_1 -
    \Phi^-(t) \big(\frac{1}{4\pi i|\xi|}(\Phi^-(0))^*\big)g_1
    \ \Big) \ . \label{half_wave_to_param}
\end{multline}
Here $P_1$ is defined to be the cutoff on line
\eqref{half_wave_est3} which is also chosen large enough such that
$P_1(f_1,g_1)=(f_1,g_1)$. From the boundedness of the $P_1$
multiplier, the estimates
\eqref{half_wave_est1} and \eqref{half_wave_est2}, the frequency
support of the adjoints, and the dualized $L_x^2\to L^2_\xi$
estimate contained in \eqref{half_wave_est1}, we easily have that
the operator \eqref{half_wave_to_param} obeys the estimate
\eqref{parametrix_estimate1}. Next, notice that by applying
\eqref{half_wave_est3}  with $\alpha=0$ and $\alpha=-1$, and using
the unit frequency condition which implies the boundedness of
$(-\Delta)^{-\frac{1}{2}}$, we have the estimate
\eqref{parametrix_estimate2}. Furthermore, by using estimate
\eqref{half_wave_est2} in conjunction with \eqref{half_wave_est3},
where this time we use the indices $\alpha=0$ and  $\alpha=1$, and
using the boundedness of $(-\Delta)^{\frac{1}{2}}$ at unit
frequency, we have the second accuracy estimate
\eqref{parametrix_estimate3}. Therefore, it remains to show that we
have the error estimate \eqref{parametrix_estimate4}. By  the
estimate \eqref{half_wave_est4} and by again making use of the dual
$L^2_x\to L^2_\xi$ adjoint bound, we are reduced to proving
(operator) commutator bounds of the type:
\begin{equation}
    \lp{[\Box_{\uA} , P_1] \Phi^\pm (\widehat{h})}{L_t^1(L_x^2)}
    \ \lesssim  \ \mathcal{E}\ \lp{\widehat{h}}{L^2_\xi} \ . \notag
\end{equation}
Using the commutator estimate \eqref{Taos_est} in conjunction with
the parametrix bounds \eqref{half_wave_est1}--\eqref{half_wave_est1.5}
(this is where the extra bound on the gradient comes in),
this reduces to showing the two bounds:
\begin{align}
    \lp{\nabla_x \uA}{L^2_t(L^{n-1}_x)} \ &\lesssim \
    \lp{\uA}{\dot{X}^\frac{n-2}{2}} \ , \label{der_on_low_est}\\
    \lp{\nabla_x [\uA,\uA]}{L^1_t(L^\infty_x)} \ &\lesssim \
    \lp{\uA}{\dot{X}^\frac{n-2}{2}}^2 \ . \label{easy_L1Linfty_est}
\end{align}
The first estimate follows easily from integrating
the following Besov and low frequency Besov nestings:
\begin{equation}
    P_{\bullet\lesssim 1}
    (\dot{B}_2^{\frac{2(n-1)}{n-3},(2,\frac{n-3}{2})})
    \ \subseteq \
    \dot{B}_2^{n-1,(2,n(\frac{n-3}{2(n-1)}))}
    \ \subseteq \ L^{n-1}
    \ . \notag
\end{equation}
The second estimate follows as easily from first distributing the
derivative and then integrating the two low frequency nestings:
\begin{equation}
    P_{\bullet\lesssim 1}(\dot{B}_2^{\frac{2(n-1)}{n-3}
    ,(2,\frac{n-3}{2})}) \ , \
    P_{\bullet\lesssim 1}(\dot{B}_2^{\frac{2(n-1)}{n-3}
    ,(2,\frac{n-1}{2})})
    \ \ \subseteq \ \ B_1^{\infty,(2,\frac{n}{2})} \
    \subseteq \ L^\infty \ . \notag
\end{equation}
This completes the proof that Proposition \ref{half_wave_prop}
implies Proposition \ref{freq1_parametrix_prop}.
\end{proof}

\ret

\section{Construction of the half wave operators}

We now begin construction of our approximate solutions $\Phi^\pm$ to
the reduced covariant wave equation $\Box_{\uA}$. This will be
accomplished by integrating over a collection of gauge
transformations designed to eliminate the highest order effect of
troublesome term $\uA^\alpha \nabla_\alpha$. In order to understand
what such a gauge transformation should be, we begin with a simple
calculation. We consider the covariant wave equation $\Box_{\oA}$,
where the connection ${}^\omega\! D =d+\oA$ will be determined in a
moment, acting on a vector valued plane wave $e^{2\pi i \lambda\,
\ou_\pm} \widehat{f}$. Here $\widehat{f}$ is a constant complex
valued matrix in $\Co$, and the $\ou^\pm$ are the standard plane
wave optical functions:
\begin{align}
        \ou^+ \ &= \  t + \omega\cdot x \ ,
        &\ou^- \ &= \ - t + \omega\cdot x \ . \notag
\end{align}
In particular, $\nabla^\alpha (\ou^\pm) = (\oL^\mp)^\alpha$, where the
$\oL^\pm$ are the associated null hyper-surface generators:
\begin{align}
        \oL^+ \ &= \  \nabla_t + \omega\cdot \nabla_x \ ,
        &\oL^- \ &= \ - \nabla_t + \omega\cdot \nabla_x \ . \notag
\end{align}
With these identifications, we easily have the calculation:
\begin{equation}
        \Box_{\oA}( e^{2\pi i \lambda\,  \ou_\pm} \widehat{f}) \ = \
        e^{2\pi i \lambda\,  \ou_\pm}\cdot
        \left(  4\pi i \lambda \big[ \oA(\oL^\mp) , \widehat{f}\big]
        + D^{\oA}_\alpha \big[\oA^\alpha , \widehat{f} \big]
        \right) \ . \label{plane_wave_error}
\end{equation}
Using the heuristic\footnote{For those who are familiar with this
kind of problem, this is precisely a reduction to the famous
$Low\times High$ frequency interaction $\uA^\alpha \nabla_\alpha
\Phi_1$.} that terms of the form $ \nabla (\oA)$ and $[\oA,\oA]$ are
lower order, and splitting the potentials $\{\oA_\alpha\}$ into the
sets $\{\oA^\pm_\alpha\}$ associated with the optical functions
$\ou_\pm$ (resp.), we see that in order eliminate the highest order
term on the right hand side of \eqref{plane_wave_error} would need
to assume this connection is in the backward (resp. forward)
\emph{$\omega$-null-gauge}:
\begin{align}
        \oA^+(\oL^-) \ &= \ 0 \ , &\oA^-(\oL^+) \ &= \ 0\  .
        \label{null_gauge_cond}
\end{align}
Of course, it is not possible to assume that a given \emph{fixed
connection} will simultaneously be in the null-gauge for every
direction $\omega$. However, it is more or less clear that since
these gauges are  of Cr\"onstrom type, it is always possible to
transform a given connection so that it is in the null-gauge for a
\emph{fixed direction}.  This motivates the following form of an
approximate solution to $\Box_{\uA}$:
\begin{equation}
        \Phi^\pm(\widehat{f}) \ = \ \int_{\mathbb{R}^n} \
        e^{2\pi i \lambda \ou^\pm} \ \og_\pm^{-1} \,
        \widehat{f}(\lambda \omega)
        \, \og_\pm \
    \chi_{(\frac{1}{2},2)}(\lambda)\
    \lambda^{n-1} d\lambda d\omega  \ , \label{parametrix}
\end{equation}
where $\chi_{(\frac{1}{2},2)}$ is a smooth bump function such that
$\chi_{(\frac{1}{2},2)}\equiv 1$ on the interval $[2^{-1},2]$ and
such that $\chi_{(\frac{1}{2},2)}\equiv 0$ outside of $[4^{-1},4]$
(the variable width assumption of Proposition \ref{half_wave_prop}
can be achieved with similar bump functions). Here, the gauge
transformation:
\begin{equation}
        {\oB^\pm} \ = \ \og_\pm \underline{A}_{\, \bullet \ll 1}
    (\og_\pm^{-1})    +  \og_\pm\, d(\og_\pm^{-1})
    \ , \label{approx_null_trans}
\end{equation}
will be chosen so that ${\oB^\pm}$ approximately satisfies
\eqref{null_gauge_cond}. It seems that there are in fact many
choices of how to do this, although the naive choice of letting
${\oB^\pm}$ satisfy \eqref{null_gauge_cond} directly by solving the
appropriate transport equations\footnote{This would  end up being
the usual  a frequency based Hadamard parametrix for the operator
$\Box_{\uA}$.} leads to group elements with poor regularity
properties. Therefore, the procedure for arriving at the correct
choice deserves some motivation.\\

The heart of the matter is two-fold. First and foremost, we need to
come up with a construction that gives us explicit formulas so that
we may perform certain standard calculations on the integral
\eqref{parametrix}. In particular, we will need to perform
integration by parts with respect to the variable $\omega$. Since
$G$ is assumed to be non-abelian, and since we will not be able to
localize things to a neighborhood of any fixed point on the
group\footnote{This is an artifact of the critical nature of the
problem. Specifically, the group elements have the heuristic form
$\og =\hbox{exp}(\nabla^{-1}\oA)$. Since we do not have $L^\infty$
control on $\nabla^{-1}\oA$ we cannot localize its image.}, this is
actually a non-trivial matter. For example, it is not possible to do
this directly through a use of the exponential map because we would
run into trouble with conjugate points.\\

Secondly, we will need to replace the transport equation which
defines the naive pure null-gauge transformation, with something
that has more ``elliptic'' features. That such a choice is possible
is, strangely enough, determined by the fact that the connection
$\{\underline{A}_{\, \bullet \ll 1} \}$ is not arbitrary, but
instead evolves according to a hyperbolic equation. This is taken
into account by condition \eqref{red_conctn_cond6}. This kind of
structure seems to be ubiquitous in geometric wave equations, both
semi and quasi-linear, and the observation that it makes the crucial
difference goes back to work of Klainerman-Rodnianski on
quasi-linear wave equations \cite{KR_quasi}. The particular form we
will use it in here is almost identical to that of \cite{RT_MKG},
but since everything we do is non-abelian, the derivation will seem
a bit different at first.\\

The first observation we use is that just like the Cr\"onstrom
gauge, the null-gauge allows one to recover the potentials directly
from the curvature. However, since we aim to derive an
(sub)-elliptic equation, we do not do this by simply integrating
along null directions. Instead, we write:
\begin{equation}
        \oL^\mp \oB^\pm_\alpha \ =  \ F^{\oB^\pm}( \oL^\mp , \partial_\alpha)
        \ . \label{null_curvature_form}
\end{equation}
Making now the approximate assumption that the $\{\oB^\pm\}$ are
simply a solution to the scalar wave equation $\Box =
\nabla_\alpha\nabla^\alpha$, which we write as:
\begin{equation}
       \Box \ = \  \oL^\pm \, \oL^\mp \ + \ \Delta_{\omega^\perp} \ ,
\end{equation}
the identity \eqref{null_curvature_form} can be written in the
integral form:
\begin{equation}
        \oB^\pm_\alpha \ = \ -\, \oL^\pm \Delta_{\omega^\perp}^{-1}\
    F^{\oB^\pm}( \oL^\mp , \partial_\alpha)
    \ . \label{null_curvature_form_int}
\end{equation}
Here $\Delta_{\omega^\perp} = \Delta - \nabla_\omega^2$ is 
the Laplacean on the plane perpendicular to the $\omega$ direction
in $\mathbb{R}^n$. We would now like to make
\eqref{null_curvature_form_int} our ``choice'' for the gauge
transformed connection on the right hand side of
\eqref{approx_null_trans}. For example, even though it was based on
the approximate assumption the $\{\oB^\pm\}$ satisfy the scalar wave
equation, it still respects the null-gauge \eqref{null_gauge_cond}
simply by the skew-symmetry property of the curvature.
Unfortunately, \eqref{null_curvature_form_int} has several
undesirable features. Firstly, we would like an expression which
involves the curvature of $\{\uA\}$, not the curvature
$F^{\oB^\pm}$. Secondly, the sub-Laplacean on the right hand side of
this expression needs to be smoothed out in some way so that its
dependence on the angular variable $\omega$ is not so rough.\\

To get around the first of these problems, we simply pretend that
the various differential operators on the right hand side of
\eqref{null_curvature_form_int} are gauge covariant. Assuming this
and then conjugating both sides of that expression by $\og_\pm$, 
moving these group elements past the differential operators on the
right, and throwing away quadratic terms from the curvature while
assuming that the reduced connection $\underline{A}_{\, \bullet \ll
1}$ satisfies the usual homogeneous wave equation, we are left with
the approximate identities:
\begin{align}
        \og_\pm^{-1} \oB^\pm_\alpha\og_\pm \ &\approx
        \ -\, \oL^\pm \Delta_{\omega^\perp}^{-1}\
        F^{\underline{A}_{\, \bullet \ll 1}}
        ( \oL^\mp , \partial_\alpha) \ , \notag \\
        &\approx \ (\underline{A}_{\, \bullet \ll 1})_\alpha +
        \nabla_\alpha \oL^\pm \Delta_{\omega^\perp}^{-1}
        \underline{A}_{\, \bullet \ll 1}(\oL^\mp) \ . \notag
\end{align}
To get around the second problem, we  mollify  the angular variable
of the second term on the right hand side of this last expression.
Doing this  and looking back  on the definition
\eqref{approx_null_trans}, we see that we would like our group
elements to be such that:
\begin{equation}
        \og_\pm^{-1} d(\og_\pm) \ \approx \
        -\, \ooPi^{(\frac{1}{2} - \delta)}
        \nabla_x \oL^\pm \Delta_{\omega^\perp}^{-1}
        \underline{A}_{\, \bullet \ll 1}(\partial_\omega) \ .
        \label{approx_group_elem_def}
\end{equation}
Here we have set:
\begin{equation}
    0 \ < \  \gamma \ \ll \ \delta \ \ll \ 1
    \ , \label{dimensional_constant}
\end{equation}
where $\gamma$ is our small all purpose constant from line
\eqref{constants_line}
above. Now the problem is, of course, that right hand side of the
above formula does not in general represent a flat connection.
However, as one can see immediately, \emph{its curvature is small}
in some sense because it is a quadratic expression. At this point,
the problem now looks essentially like what happens for
wave-maps\footnote{It is very much our philosophy here that this
problem is essentially equivalent to wave-maps after a
microlocalization. Of course, as the reader will see, this
microlocalization is quite costly and introduces many objects that
are not present in the original wave-maps problem.} (see e.g. \cite{T_wm1} and
\cite{SS_wm}). In particular, it is clear that the right way to define
the group elements $\og_\pm$ so that the approximate formula
\eqref{approx_group_elem_def} holds is to flatten out the right hand
side of that expression as much as possible by using the potential
version \ref{A_Uhl_lem} of the Uhlenbeck lemma. Therefore, what we
need to do is to show the fixed time estimate:
\begin{equation}
        \lp{\ooPi^{(\frac{1}{2} - \delta)}
        \nabla_x \oL^\pm
        \Delta_{\omega^\perp}^{-1}\,  \uA(\partial_\omega)
        }{L^n} \ \lesssim \ \mathcal{E} \ ,
        \label{degen_A_tdA_Ln_smallness}
\end{equation}
and then assume that $\mathcal{E}$ is chosen small enough to that we
may use it as the constant in \eqref{A_tdA_Ln_smallness}. Because of
its utility in the sequel, we will in fact prove the more general
estimate:
\begin{equation}
        \lp{\ooPi^{(\frac{1}{2} - \delta)}
        \nabla_x \oL^\pm
        \Delta_{\omega^\perp}^{-1}\,  \uA(\partial_\omega)
        }{\dot{B}_{2,10n}^{p_\gamma , (2,\frac{n-2}{2})}}
    \ \lesssim \ \mathcal{E} \ ,
        \label{Besov_degen_tdA_Ln_smallness}
\end{equation}
where $p_\gamma$ is a dimension dependent Lebesgue index which
we set to:
\begin{equation}
    p_\gamma \ = \ \frac{2(n-1)}{n-3 -2\gamma}  \ . \label{p_gamma_line}
\end{equation}
Here $0 < \gamma \ll 1$ is again the all-purpose constant which we
have fixed in section \ref{basic_not_sect} to be small enough so
that it is compatible with its use here. Notice that
\eqref{Besov_degen_tdA_Ln_smallness} implies the estimate
\eqref{degen_A_tdA_Ln_smallness} thanks to the embedding
\eqref{Lebesgue_besov_incl} and the fact that for $\gamma$
sufficiently small there is plenty of room in
the inequality $p_\gamma < n$.\\

Now, because the norm $\dot{B}_{2,10n}^{p_\gamma ,
(2,\frac{n-2}{2})}$ is $\ell^2$ based, by orthogonality and the
$L^\infty(\dot{H}^\frac{n-2}{2})$ 
estimate contained in the bootstrapping assumption
\eqref{red_conctn_cond5}, we see that in order to conclude
\eqref{Besov_degen_tdA_Ln_smallness} it is enough to show the fixed
frequency estimate (note that there are no high frequencies here):
\begin{equation}
        \lp{ \nabla_x \oL^\pm
        \Delta_{\omega^\perp}^{-1}\,  (\uA)_\mu (\partial_\omega)
        }{L^{p_\gamma}} \ \lesssim  \  \mu^{
        n(\frac{1}{2}-\frac{1}{p_\gamma})}
        \lp{(\uA)_\mu }{L^2} \ .
        \notag
\end{equation}
Decomposing the spatial frequency variable into fixed dyadic angular
sectors spread from the direction $\omega$:\ $P_\mu = \sum_\theta
\oPi_\theta P_\mu$, this estimate further reduces (after dyadic
summing) to being able being able to prove that:
\begin{equation}
        \lp{\oPi_\theta
        \nabla_x \oL^\pm
        \Delta_{\omega^\perp}^{-1}\,  (\uA)_\mu (\partial_\omega)
        }{L^{p_\gamma}} \ \lesssim  \  \theta^\gamma
    \mu^{n(\frac{1}{2} - \frac{1}{p_\gamma})}
        \lp{(\uA)_\mu }{L^2} \ .
        \label{trunc_A_Ln_bound}
\end{equation}
We are now almost at the point
where we can apply the angular Bernstein inequality
\eqref{Bernstein2} directly, because in the
current localized setting we have the symbol bounds:
\begin{equation}
        \oPi_\theta
        \nabla_x \oL^\pm
        \Delta_{\omega^\perp}^{-1}\, S_{|\tau|\lesssim |\xi|}
	P_\mu  \ \approx \
        \theta^{-2} P_\mu \ , \label{main_degen_symbol_bound}
\end{equation}
where we are enforcing the heuristic notation introduced on line
\eqref{heuristic_op_bnd}. However, since Bernstein only nets us a
savings of:
\begin{equation}
        \theta^{(n-1) ( \frac{1}{2} - \frac{1}{p_\gamma})} \ = \
        \theta^{1+\gamma} \ , \notag
\end{equation}
in this context, we need to be a bit more careful in order to gain
an extra power of $\theta$. This is provided by the fact that the
potentials $\{\uA\}$ are in the Coulomb gauge. Notice that if say,
$\frac{1}{10} < \theta$ there is nothing to worry about and we have
estimate \eqref{trunc_A_Ln_bound} without any problem. On the
other-hand, if it is the case that $\theta < \frac{1}{10}$, then we
can use the fact that $\oPi_\theta\nabla_\omega^{-1}$ is elliptic
(in terms of symbol bounds) in conjunction with the gauge condition
$d^* \uA = 0 $ to write:
\begin{equation}
        \oPi_\theta \uA(\partial_\omega)
    \ = \ \nabla_\omega^{-1} \oPi_\theta \sd^* \suA \
    \approx \ \theta \ \oPi_\theta \suA \ . \label{main_coulomb_savings}
\end{equation}
Here $\{\suA\}$ the induced connection (angular portion) on the
hyperplane $\mathcal{H}_{\omega^\perp}$ perpendicular to $\omega$,
and $\sd^*$ is the associated divergence. We note here that this
identity will turn out to be very useful and will be used many times
throughout the sequel. With these extra savings in mind, an
application of Bernstein now directly yields the desired estimate
\eqref{trunc_A_Ln_bound}.\\

We have now constructed the infinitesimal group elements $\og_\pm$
in equations \eqref{approx_null_trans}, which is explicitly defined
by the formulas \eqref{C_dc} in Lemma \ref{A_Uhl_lem} applied to the
connection:
\begin{equation}
        {\uoA}^\pm \ = \ -\,  \ooPi^{(\frac{1}{2} - \delta)}
        \nabla_{x} \oL^\pm \Delta_{\omega^\perp}^{-1}
        \underline{A}_{\, \bullet \ll 1}(\partial_\omega) \ .
        \label{A_omega_def}
\end{equation}
This has the pleasant effect that we will never need to explicitly
refer to the connection $\{\oB^\pm\}$ in line
\eqref{approx_null_trans}. We can calculate the conjugated right
hand side of that expression to be:
\begin{equation}
        \og_\pm^{-1}{\oB^\pm}\og_\pm \ = \  \underline{A}_{\, \bullet \ll 1}
        -  \oC^\pm
        \ , \label{B_conj_eq}
\end{equation}
where we have set:
\begin{equation}
        \og_\pm^{-1} d(\og_\pm) \ = \ \oC^\pm
        \ . \label{C_def}
\end{equation}
Using the formulas \eqref{C_dc}, we have the following expressions
for the spatial components $\{\uoC^\pm\}$:
\begin{subequations}\label{uCpm_system}
\begin{align}
        (\uoC^\pm)^{df} \ &= \ d^* \, \Delta^{-1} [\uoC^\pm , \uoC^\pm ]
        \ , \label{uCpm_system1}\\
        (\uoC^\pm)^{cf} \ &= \  \uoA^\pm
        \ - \ \nabla_x \Delta^{-1} [\uoA^\pm , \uoC^\pm ]
        \ . \label{uCpm_system2}
\end{align}
\end{subequations}
In order to compute a formula for the temporal potential $\oC^\pm_0$, we
simply use the fact that $F^{\oC^\pm} = 0$ and the formula
\eqref{uCpm_system2} which together imply (by computing $d^* E^{\oC^\pm}$):
\begin{equation}
        \oC^\pm_0 \ = \ \oA^\pm_0
        \ - \ \nabla_t \Delta^{-1} [  {\uoA^\pm} , \uoC^\pm]
        \ - \ d^* \Delta^{-1} [\oC^\pm_0 , \uoC^\pm ] \ , \label{C0_system}
\end{equation}
where we have:
\begin{equation}
        \oA^\pm_0 \ = \ -\,  \ooPi^{(\frac{1}{2} - \delta)}
        \nabla_{t} \oL^\pm \Delta_{\omega^\perp}^{-1}
        \underline{A}_{\, \bullet \ll 1}(\partial_\omega) \ . \notag
\end{equation}\ret

We remark here that the importance of the system of equations
\eqref{uCpm_system1}--\eqref{C0_system} is that they give the
following decomposition of the infinitesimal gauge transformation
$\{\oC^\pm \}$:
\begin{equation}
        \oC^\pm \ = \ -\, \nabla_{t,x} \ \ooPi^{(\frac{1}{2} - \delta)}
    \oL^\pm \Delta_{\omega^\perp}^{-1}\,
    \uA(\partial_\omega) \ + \ \hbox{ \{Quadratic Error\} } \ .
    \label{rough_decomp}
\end{equation}
The linear term in the above expression is enough to kill off the
worst error term when differentiating the parametrix
\eqref{parametrix}. It should be noted that this linear term is
precisely what one gets more directly in the abelian case studied in
\cite{RT_MKG}. We should also point out here that the quadratic
error on the right hand side of \eqref{rough_decomp} above is
\emph{much} more delicate than the quadratic error resulting form
the cancellation involving the linear term in this expression. In
order to control this, we will need the full force of the
orthogonality properties of our parametrix, which are contained in
the bootstrapping assumption \eqref{red_conctn_cond5}, as well as
some rather technical function spaces and multilinear estimates
which we will develop in Section \ref{Decomp_sect}.\\

To close out this section, we apply the truncated covariant wave
operator $\Box_{\uA}$ to the parametrix \eqref{parametrix} and record
the various error terms which result. We gather this together in the
following proposition:\\

\begin{prop}[Error terms for the differentiated parametrix]
Consider the parametrix $\Phi^\pm(\widehat{f})$ defined by the
formula \eqref{parametrix}, with infinitesimal gauge transformations
given by equations \eqref{uCpm_system1}--\eqref{C0_system}. Then one
has the identity:
\begin{align}
        &\Box_{\uA} \Phi^\pm(\widehat{f}) \label{error_terms} \\
    = \ &4\pi i \ \int_{\mathbb{R}^n} \
        e^{2\pi i \lambda \ou^\pm} \ \big[
        \uA(\oL^\mp) - \oC^\pm(\oL^\mp) \ , \
        \og_\pm^{-1} \,
        \widehat{f}(\lambda \omega)
        \, \og_\pm \big] \
        \chi_{(2^{-1},2)}(\lambda)\
        \lambda^{n} d\lambda d\omega \notag\\
        &- \ \ \int_{\mathbb{R}^n} \
        e^{2\pi i \lambda \ou^\pm} \ \big[
        D^{\uA}_\alpha \big({\oC^\pm}\big)^\alpha \ , \
        \og_\pm^{-1} \,
        \widehat{f}(\lambda \omega)
        \, \og_\pm \big] \ \chi_{(2^{-1},2)}(\lambda)\
        \lambda^{n-1} d\lambda d\omega \notag\\
        &+\ \ \int_{\mathbb{R}^n} \
        e^{2\pi i \lambda \ou^\pm} \ \Big[ \uA^\alpha -
        (\oC^\pm)^\alpha \ , \ \big[
        (\uA)_\alpha - \oC^\pm_\alpha \  , \
        \og_\pm^{-1} \,
        \widehat{f}(\lambda \omega)
        \, \og_\pm \big] \Big] \
        \chi_{(2^{-1},2)}(\lambda)\
        \lambda^{n-1} d\lambda d\omega \ . \notag
\end{align}
\end{prop}\ret

\begin{rem}
The worst error term in the expression \eqref{error_terms} is of
course the ``derivative fall on high'' term which is the first on
the right hand side. However, using the structure equation
\eqref{red_conctn_cond6}, this takes the form:
\begin{align}
       &\uA(\oL^\mp) - \oC^\pm(\oL^\mp) \ , \label{rough_main_error}\\
       = \ \ &\uA(\partial_\omega) \ + \
       \oPi^{(\frac{1}{2} - \delta)} \oL^\mp \oL^\pm
       \Delta_{\omega^\perp}^{-1}\,  \uA(\partial_\omega) \ + \
    \hbox{\{Quadratic Error\}} \ , \notag\\
    = \ \ &(I - \oPi^{(\frac{1}{2} - \delta)} )\uA(\partial_\omega) \ + \
    \hbox{\{Quadratic Error\}} \ . \notag
\end{align}
The key observation now is  that since the operator $(I -
\oPi^{(\frac{1}{2} - \delta)} )$ cuts off on such a small angular
sector with respect to the spatial frequency, an application of
Bernstein's inequality gains enough extra \emph{spatial} derivatives
to put this term in the  mixed Lebesgue space $L^2(L^{n-1})$.
Furthermore, the quadratic error term which is left over involves
enough bilinear interactions to go in $L^1(L^\infty)$. So in this
sense, as we have mentions before, the problem reduces to something
which is reminiscent of wave-maps. Of course, there is a somewhat
heavy price to pay for this ``renormalization'', which is that it
must take place under an integral sign. Finally, it is worth
pointing out that this top order cancellation is completely analogous
to what happens in the abelian case \cite{RT_MKG}.
\end{rem}\ret

\begin{proof}[Proof of the error identity \eqref{error_terms}]
The proof is a simple consequence of using gauge transformations in
conjunction with the identity \eqref{plane_wave_error}. Applying the
truncated covariant wave operator, and differentiating under the
integral sign, we see that:
\begin{align}
        &\Box_{\uA} \Phi^\pm_f \ , \notag\\
    = \
        &\int_{\mathbb{R}^n} \ \Box_{\uA} \Big(
        e^{2\pi i \lambda \ou^\pm} \
        \og_\pm^{-1} \,
        \widehat{f}(\lambda \omega)
        \, \og_\pm \Big)  \
    \chi_{(2^{-1},2)}(\lambda)\
        \lambda^{n-1} d\lambda d\omega \ ,
        \notag\\
        = \ &\int_{\mathbb{R}^n} \ \og_\pm^{-1}\ \Box_{{\oB^\pm} } \Big(
        e^{2\pi i \lambda \ou^\pm} \
        \widehat{f}(\lambda \omega) \Big)
        \ \og_\pm   \
    \chi_{(2^{-1},2)}(\lambda)\
        \lambda^{n-1} d\lambda d\omega \ , \notag\\
        = \ &\int_{\mathbb{R}^n} \ e^{2\pi i \lambda\,  \ou_\pm}\ \og_\pm^{-1}
        \left(  4\pi i \lambda \big[ {\oB^\pm}(\oL^\mp) , \widehat{f}\big]
        + D^{\oB^\pm}_\alpha \big[{\oB^\pm}^\alpha , \widehat{f} \big]
        \right)\og_\pm \
    \chi_{(2^{-1},2)}(\lambda)\
        \lambda^{n-1} \ d\lambda d\omega \ ,
        \notag\\
   \begin{split}
        = \ &4\pi i \ \int_{\mathbb{R}^n} \ e^{2\pi i \lambda\,  \ou_\pm}\
        \big[ \og_\pm^{-1}
        {\oB^\pm}(\oL^\mp) \, \og_\pm \ , \ \og_\pm^{-1}
        \widehat{f}\ \og_\pm \big]\ \lambda^{n-1} \
    \chi_{(2^{-1},2)}(\lambda)\
        d\lambda d\omega \\
        &\ \ \ \ \ \ \ + \ \int_{\mathbb{R}^n} \ e^{2\pi i \lambda\,  \ou_\pm}\
        D^{\uA}_\alpha \big[ \og_\pm^{-1}{\oB^\pm}^\alpha \og_\pm \
        , \ \og_\pm^{-1} \widehat{f} \ \og_\pm \big]
        \ \lambda^{n-1} \
    \chi_{(2^{-1},2)}(\lambda)\
            d\lambda d\omega \ ,
   \end{split}\notag\\
   \begin{split}
        = \ &4\pi i \ \int_{\mathbb{R}^n} \ e^{2\pi i \lambda\,  \ou_\pm}\
        \big[ \uA(\oL^\mp) - \oC^\pm(\oL^\mp)
        \ , \ \og_\pm^{-1}
        \widehat{f}\ \og_\pm \big]\ \lambda^{n-1} \
    \chi_{(2^{-1},2)}(\lambda)\
    d\lambda d\omega \\
        &\ \ \ \ \ \ \ - \ \int_{\mathbb{R}^n} \ e^{2\pi i \lambda\,  \ou_\pm}\
        \big[ \nabla_\alpha (\oC^\pm)^\alpha \
        , \ \og_\pm^{-1} \widehat{f} \ \og_\pm \big]
        \ \lambda^{n-1} \
    \chi_{(2^{-1},2)}(\lambda)\
    d\lambda d\omega\\
        &\ \ \ \ \ \ \ + \ \int_{\mathbb{R}^n} \ e^{2\pi i \lambda\,  \ou_\pm}\
        \big[ (\uA)_\alpha -\oC^\pm_\alpha \
        , \ \nabla^\alpha \ \big( \og_\pm^{-1} \widehat{f} \ \og_\pm\big) \big]
        \ \lambda^{n-1} \
    \chi_{(2^{-1},2)}(\lambda)\
    d\lambda d\omega\\
        &\ \ \ \ \ \ \ + \ \int_{\mathbb{R}^n} \ e^{2\pi i \lambda\,  \ou_\pm}\
        \big[\ (\uA)^\alpha \
        , \ [ (\uA)_\alpha - \oC^\pm_\alpha \
        , \ \og_\pm^{-1} \widehat{f} \ \og_\pm] \big]
        \ \lambda^{n-1} \
    \chi_{(2^{-1},2)}(\lambda)\
    d\lambda d\omega
   \end{split}\notag\\
        = \ &\hbox{(L.H.S.)}\eqref{error_terms} \ . \notag
\end{align}
Notice that the equality on the last line follows from:
\begin{equation}
        \nabla_\alpha \ \big( \og_\pm^{-1} \widehat{f} \ \og_\pm\big)
        \ = \ \big[ \ \og_\pm^{-1} \widehat{f} \ \og_\pm
        \ , \ \oC^\pm_\alpha
        \ \big] \ , \notag
\end{equation}
which is a consequence of line \eqref{C_def} above, followed by the
Jacobi identity:
\begin{align}
        &\big[ \ (\uA)^\alpha - (\oC^\pm)^\alpha \ , \
        [\ \og_\pm^{-1} \widehat{f} \ \og_\pm
        \ , \ \oC^\pm_\alpha
        \ ] \ \big] \ , \notag \\
        = \ -\, &\big[ \ \oC^\pm_\alpha \ , \
        [\ (\uA)^\alpha - (\oC^\pm)^\alpha
        \ , \ \og_\pm^{-1} \widehat{f} \ \og_\pm
        \ ] \ \big] \notag \\
        &\ \ \ \ \ \ \ - \ \big[ \ \og_\pm^{-1} \widehat{f} \ \og_\pm \ , \
        [\ \oC^\pm_\alpha \ , \
        (\uA)^\alpha - (\oC^\pm)^\alpha
        \ ] \ \big] \ , \notag\\
        = \ -\, &\big[ \ \oC^\pm_\alpha \ , \
        [\ (\uA)^\alpha - (\oC^\pm)^\alpha
        \ , \ \og_\pm^{-1} \widehat{f} \ \og_\pm
        \ ] \ \big] \notag \\
        &\ \ \ \ \ \ \ - \ \big[ \
        [\ (\uA)_\alpha \ ,
        \ (\oC^\pm)^\alpha \ ] \ , \
        \og_\pm^{-1} \widehat{f} \ \og_\pm
        \ ] \ \big] \ . \notag
\end{align}
This completes the proof of \eqref{error_terms}.
\end{proof}

\ret

\section{Fixed Time $L^2$ Estimates for the Parametrix}

We now begin our proof of the estimates \eqref{half_wave_est} for
the integral operator \eqref{parametrix} introduced in the last
section. Here we cover bounds which are of non-differentiated energy
type. Specifically, we will show the undifferentiated 
$L^\infty(L^2)$ estimate contained in
\eqref{half_wave_est1}, as well as the multiplier-approximation
bound \eqref{half_wave_est3}. Both of these will follow from the
same set of estimates. At a heuristic level, they are not much more
involved that a standard $TT^*$ argument followed by some
integration by parts, although the details turn out to be a bit
involved. Things will be computed more or less directly by an appeal
to the explicit equations \eqref{uCpm_system1}--\eqref{C0_system},
taking a little bit of care to use them properly. This will be done
by considering them as ``path lifting'' formulas from Minkowski
space $\mathcal{M}^n$ to the compact group $G$. This allows us to
employ an integral form of the intermediate value theorem from
elementary calculus which is valid in the context of Lie groups. It
turns out that this identity can be differentiated as many times as
necessary with respect to the angular frequency 
variable, although this fact
is provided through a surprisingly delicate bootstrapping argument.
Here the unitarity of the group is needed in a crucial way to keep
everything from collapsing. Once the bootstrapping is complete, the
estimates themselves will be proved using a ``trace-Bernstein'' type
inequality that we construct by hand using various multipliers. Once
the integration by parts portion of things is taken care of, we will
close the $L^2$ estimate by showing that a ``non-smooth''
remainder kernel has small amplitudes after integration in the
angular frequency variable. This involves some fairly technical
bilinear estimates because the necessary othogonality arguments are
difficult to pass through Hodge systems. The details of these
procedures are as follows.\\

Throughout this section we will replace the specific cutoff
function $\chi_{(\frac{1}{2},2)}$ appearing in the definition of
parametrix \eqref{parametrix} with an arbitrary smooth scalar bump
function $\chi(\xi)$ that
we may assume to be supported in the frequency annulus
$\{4^{-1} < |\xi|  < 4\}$. At fixed time $t_0$, we define the
operator $T(\widehat{f}) = \Phi(\widehat{f})(t_0)$, where we have
suppressed the $\pm$ notation because it will be irrelevant for what
we do here. Our first goal is the prove the bound:
\begin{equation}
        \lp{T(\widehat{f})}{L^2} \ \lesssim \ \lp{\widehat{f}}{L^2} \
        . \label{first_fixed_time_bound}
\end{equation}
Squaring this, it suffices to show that (here $f$ has no relation to
$\widehat{f}$ and simply represents a function of the physical-space
variables):
\begin{equation}
        \lp{TT^*(f)}{L^2} \ \lesssim \ \lp{f}{L^2} \
        , \label{sq_first_fixed_time_bound}
\end{equation}
where the adjoint $T^*$ is taken with  respect to the Killing form
\eqref{explicit_killing}. A quick calculation of the kernel of this
operator shows that:
\begin{equation}
        K^{TT^*}(x,y) \ = \
    \int_{\mathbb{R}^n} \
    e^{2\pi i ( x-y)\cdot \xi}\
    \og^{-1}(x) \og(y) \big[\bullet \big]\
     \og^{-1}(y) \og(x) \ \chi(\xi)\ d\xi \ , \label{L2_TT*_kernel}
\end{equation}
where we use the $[\bullet]$ notation to emphasize the fact that this
operator acts via conjugation. Our task is now to show the estimates:
\begin{equation}
        \lp{K^{TT^*}}{L^\infty_y(L^1_x)} \ , \
        \lp{K^{TT^*}}{L^\infty_x(L^1_y)}
        \ \lesssim \ 1\ . \label{TT*_L2_est}
\end{equation}
Since $K^{TT^*}$ is essentially symmetric in $(x,y)$, we may
concentrate on the first such estimate.\\

To proceed, we first decompose the product
physical space $\RR\times \RR$ into the dyadic regions:
\begin{equation}
        \mathcal{D}_\sigma \ = \ \{|x-y| \sim \sigma \ \big| \ \sigma
        = 2^i \ , \ i\in \mathbb{N}\} \ . \label{dyadic_physical_scales}
\end{equation}
We then decompose the kernel $TT^*$ kernel into the dyadic sum:
\begin{equation}
        K^{TT^*} \ = \ \sum_\sigma \ \chi_{\mathcal{D}_\sigma}
        K^{TT^*} \ = \ \sum_\sigma \ K_\sigma^{TT^*} \ . \notag
\end{equation}
By dyadic summing, to show \eqref{TT*_L2_est} it suffices to be able
to show the single estimate:
\begin{equation}
        \lp{K_\sigma^{TT^*}}{L^\infty_y(L^1_x)} \ \lesssim \
        \sigma^{-\gamma} \ , \label{TT*_L2_est_dyadic}
\end{equation}
where $0 < \gamma\ll 1$  now
represents a small savings in physical space decay.
Now \eqref{TT*_L2_est_dyadic} would be easy to show if we had the
absolute decay estimate:
\begin{equation}
        |K_\sigma^{TT^*}(x,y)| \ \lesssim \ |x - y |^{-(n + \gamma)} \
         , \notag
\end{equation}
and this is almost true. Unfortunately, there is a regularity
problem due to the degeneracy of the sub-Laplacean
$\Delta_{\omega^\perp}$ used in the connection \eqref{uCpm_system}
which provides the group elements $\og$. This forces us to write the
kernel $K_\sigma^{TT^*}$ as a sum of two terms:
\begin{equation}
        K_\sigma^{TT^*} \ = \ \td{K}_\sigma^{TT^*} \ + \
    \mathcal{R}^{TT^*}_\sigma \ . \label{L2_TT*_splitting}
\end{equation}
We will then prove that both:
\begin{align}
        |\td{K}_\sigma^{TT^*}(x,y)| \ &\lesssim \ |x - y |^{-(n +
         \gamma)} \ , \label{abs_TT*_bound}\\
    \lp{\mathcal{R}^{TT^*}_\sigma}{L_y^\infty(L^1_x)} \ &\lesssim \
    \sigma^{-\gamma} \ . \label{Lp_TT*_bound}
\end{align}\\

To define the splitting \eqref{L2_TT*_splitting}, we factor the
group elements $\og$ into a product of \emph{smooth} and
\emph{small} parts. This is completely analogous to the procedure
used in \cite{RT_MKG}, but since things are non-abelian (and hence
non-linear) here, the estimates required are quite a bit more
involved. What we will do is construct another gauge transformation
$\td{\og}$, which is based on a further smoothing of the connection
\eqref{A_omega_def}. This will produce a group element which can be
treated as a standard symbol. To this end, we define the scale
mollified connection:
\begin{equation}
        \td{\uoA^{(\sigma)}} \ = \ -\,  \ooPi_{\sigma^{-1 + \gamma}
        < \bullet }\ooPi^{(\frac{1}{2} - \delta)}
        \nabla_x \oL\
        \Delta_{\omega^\perp}^{-1}\,  \uA(\partial_\omega) \ ,
        \label{approx_null_generator_smoothed}
\end{equation}
where  $\gamma$ is, again, the small dimensional constant from line
\eqref{dimensional_constant}. Again, we have dropped the $\pm$
notation because it is irrelevant. Following the proof of
\eqref{degen_A_tdA_Ln_smallness}, and using the fact that the
multipliers  $\ooPi_{\sigma^{-1 + \gamma} < \bullet}$ are bounded on
frequency localized Lebesgue spaces, we may apply Lemma
\ref{A_Uhl_lem} to the connection $\{\td{\uoA^{(\sigma)}}\}$. This
produces a group element $\td{\og}$, which is defined by the
infinitesimal generator:
\begin{equation}
        \td{\og}^{-1} d (\td{\og}) \  = \ \td{\uoC} \ . \label{tdC_ODE}
\end{equation}
Furthermore, this generator is itself defined via the Hodge system:
\begin{subequations}\label{td_uCpm_system}
\begin{align}
        (\td{\uoC})^{df} \ &= \ d^* \, \Delta^{-1} [\td{\uoC} , \td{\uoC} ]
        \ , \label{td_uCpm_system1}\\
        (\td{\uoC})^{cf} \ &= \ \widetilde{\uoA^{(\sigma)}}
        \ - \ \nabla_x \Delta^{-1} [\widetilde{\uoA^{(\sigma)}} , \td{\uoC} ]
        \ . \label{td_uCpm_system2}
\end{align}
\end{subequations}
Using this new group element $\td{\og}$, we define the remainder
group element ${\oh}$ via the product:
\begin{equation}
        \og \ = \ {\oh}\, \td{\og} \ . \label{gh_prod_def}
\end{equation}
To compute the infinitesimal generator of ${\oh}$, we first use the
identity:
\begin{align}
        d ({\oh}) \ &= \ d(\og) \td{\og}^{-1} + g\, d(\td{\og}^{-1}) \
        , \notag\\
    &= \ {\oh} \ \td{\og} \,
    \big( \uoC - \td{\uoC}\big)\,  \td{\og}^{-1} \ . \label{dh_C_eq}
\end{align}
This leads us to define the difference connection:
\begin{equation}
        \td{\td{\uoC}} \ = \ \uoC - \td{\uoC} \ . \label{Cdt_def}
\end{equation}
A quick calculation using the systems \eqref{uCpm_system} and
\eqref{td_uCpm_system} shows that this new connection can be pinned
down via the Hodge system:
\begin{subequations}\label{td_td_uDpm_system}
\begin{align}
        (\td{\td{\uoC}})^{df} \ &= \ d^* \, \Delta^{-1} \big(
        [\td{\uoC} , \td{\td{\uoC}} ] \ + [\td{\td{\uoC}} , \td{\uoC}]
        \ \big) \ , \label{td_td_uDpm_system1}\\
        (\td{\td{\uoC}})^{cf} \ &= \ \widetilde{\uoA} -
        \widetilde{\uoA^{(\sigma)}} \
        - \  \nabla_x \Delta^{-1} \Big([\widetilde{\uoA} -
        \widetilde{\uoA^{(\sigma)}}, \td{\uoC} ]
        \ + \ [ \widetilde{\uoA^{(\sigma)}} , \td{\td{\uoC}}]
        \Big)\ , \label{td_td_uDpm_system2}
\end{align}
\end{subequations}
where a simple computation shows that:
\begin{equation}
        \uoA - \td{\uoA^{(\sigma)}} \ = \ -\,  \ooPi_{\bullet\leqslant
        \sigma^{ -1 + \gamma} }\ooPi^{(\frac{1}{2} - \delta)}
        \nabla_x \oL\
        \Delta_{\omega^\perp}^{-1}\,  \uA(\partial_\omega) \ ,
        \label{uoA_difference}
\end{equation}\ret

We now define the decomposition \eqref{L2_TT*_splitting} along the
following decompositions of the group element products in the kernel
\eqref{L2_TT*_kernel}:
\begin{align}
        \og^{-1}(x)\og(y) \ &= \ \td{\og}^{-1}(x)\td{\og}(y) +
    \td{\og}^{-1}(x)\left(  {\oh}^{-1}(x){\oh}(y)
    - I \right) \td{\og}(y) \ , \label{gh_left}\\
    \og^{-1}(y)\og(x) \ &= \ \td{\og}^{-1}(y)\td{\og}(x) +
     \td{\og}^{-1}(y)\left(  {\oh}^{-1}(y){\oh}(x)
    - I \right)\td{\og}(x) \ . \label{gh_right}
\end{align}
Accordingly, we define:
\begin{equation}
        \td{K}^{TT^*}(x,y) \ = \
    \int_{\mathbb{R}^n} \
    e^{2\pi i ( x-y)\cdot \xi}\
    \td{\og}^{-1}(x) \td{\og}(y) \big[\bullet \big]\
     \td{\og}^{-1}(y) \td{\og}(x) \ \chi(\xi)\ d\xi
     \ , \label{smoothed_L2_TT*_kernel}
\end{equation}
and then define $\mathcal{R}_\sigma^{TT^*}$ according to the formula
\eqref{L2_TT*_splitting}. The idea now is that while one can only
perform integration by parts in the kernel
\eqref{smoothed_L2_TT*_kernel} above, the group element
${\oh}^{-1}(x){\oh}(y)$ and its inverse, which must be contained as
at least one factor in the remainder, are so close to the identity
matrix that the resulting difference expression can be estimated
without use of the oscillations which take place under the integral
sign.\\

We now begin our proof of the estimate \eqref{abs_TT*_bound}. To do
this, we simply integrate by parts as may times as necessary with
respect to the variable $\xi$ in order to pick up the needed point-wise
decay. Doing this, we see that in order to
draw our conclusion, it suffices to show the
following symbol bounds for $1\leqslant k$:
\begin{align}
        \chi_{\mathcal{D}_\sigma}\
        \lp{ \nabla_\xi^k \big( \td{\og}^{-1}(x) \td{\og}(y)\big) }
        \ \lesssim \ \mathcal{E}\cdot
        \sigma^{k( 1 - \gamma)} \ ,
        \label{L2_symbol_bds1}\\
        \chi_{\mathcal{D}_\sigma}\
        \lp{ \nabla_\xi^k \big(  \td{\og}^{-1}(y) \td{\og}(x)    \big) }
        \ \lesssim \ \mathcal{E}\cdot \sigma^{k( 1 - \gamma)} \ .
        \label{L2_symbol_bds2}
\end{align}
In fact, we shall prove the following more general bounds, which contain
\eqref{L2_symbol_bds1}--\eqref{L2_symbol_bds2} as a special case, and
which will be useful in the sequel:\\

\begin{prop}[Symbol bounds for the smoothed amplitudes
 $\td{\og}^{-1}(t,x) \td{\og}(s,y)$ and $\td{\og}^{-1}(s,y) \td{\og}(t,x)$]
 \label{symbol_bound_prop}
Let the group elements $\td{\og}$ be defined infinitesimally by the
Hodge system \eqref{td_uCpm_system}, where the parameter
$\sigma^{-1 + \gamma}$ is replaced by $M^{-1}$, where $M$
lies in the range:
\begin{equation}
    (|t-s| + |x-y|)^\frac{1}{2} \ \leqslant \ M \ \leqslant
    |t-s| + |x-y| . \label{M_range}
\end{equation}
Then for any integer $1\leqslant k$, one has the following symbol
bounds assuming that the bootstrapping constant $\mathcal{E}$ from
line \eqref{red_conctn_cond5} is chosen sufficiently small (with
respect to each fixed $k$):
\begin{align}
        \lp{ \nabla_\xi^k \big( \td{\og}^{-1}(t,x) \td{\og}(s,y)\big) }
        \ \lesssim \ \mathcal{E}
    \cdot M^{k} \ ,
        \label{st_symbol_bds1}\\
        \lp{ \nabla_\xi^k \big(  \td{\og}^{-1}(s,y) \td{\og}(t,x)    \big) }
        \ \lesssim \ \mathcal{E}\cdot M^{k} \ .
        \label{st_symbol_bds2}
\end{align}
Here the $\nabla_\xi^k$ notation is shorthand for all $k^{th}$ order
partial derivatives involving the variable $\xi$, and $\lp{\cdot}{}$
is the standard matrix vector-norm from line \eqref{matrix_inner}.
The implicit constants on the right hand side depend on $k$, but are
uniform in the parameter $M$ for each fixed $k$.
\end{prop}\ret

\begin{proof}[Proof of the estimates
\eqref{st_symbol_bds1}--\eqref{st_symbol_bds2}] It suffices for us
to prove the first bound \eqref{st_symbol_bds1}, as the second
follows from virtually identical reasoning. The goal is to reduce
this via an ODE bootstrapping type argument to an associated
estimate involving the connection $\{\td{\oC}\}$. This associated
estimate will then be proved by another  bootstrapping argument in
certain mixed Lebesgue-Besov spaces naturally associated with the
ODE problem from the first step. The goal of the second
bootstrapping will be to reduce things to proving the Besov
estimates for the connection $\{\uA\}$ which appears as the linear
term
on the right hand side of the Hodge system \eqref{td_uCpm_system1}.\\

Before proceeding, we first make a preliminary reduction on the product
$\td{\og}^{-1}(t,x) \td{\og}(s,y)$. We would like be set up as to only have to
handle products which involve the same space or same time variables. This
is easily accomplished via the product decomposition:
\begin{equation}
    \td{\og}^{-1}(t,x) \td{\og}(s,y) \ = \
    \td{\og}^{-1}(t,x) \td{\og}(t,y)\cdot \td{\og}^{-1}(t,y)\td{\og}(s,y)
    \ . \label{simple_product}
\end{equation}
It is clear that if we can produce the bounds
\eqref{st_symbol_bds1} for each of the terms on the right
hand side of \eqref{simple_product} separately, then by the
product rule for derivatives
we have the estimate \eqref{st_symbol_bds1} for the full term.
Since they require
slightly different arguments, we will proceed separately for each of these
two factors.\\

Our first task is to prove the bound \eqref{st_symbol_bds1} for the spatial
product $\td{\og}^{-1}(t,x) \td{\og}(t,y)$. This will be done inductively with
respect to the value of $k$. Since we will proceed via a bootstrapping type
procedure, we first assume that we can prove the desired bounds over small
intervals and then try to use this knowledge to extend things to longer
intervals. To do this, we differentiate the product
$\td{\og}^{-1}(t,\ell) \td{\og}(t,y)$, where $[y,\ell]$ is some shorter line segment
inside of $[y,x]$, with respect to the
operators $(M^{-1} \nabla_\xi)^k$. This yields the equation:
\begin{multline}
        (M^{-1}\nabla_\xi)^k\big(\td{\og}^{-1}(\ell) \td{\og}(y)\big)
        \ = \\
        \sum_{i=0}^{k-1} \
        (M^{-1}\nabla_\xi)^{k-i}\big(
        \td{\og}^{-1}(\ell) \td{\og}(x_1)
        \big)
        \cdot (M^{-1} \nabla_\xi)^i\big(
        \td{\og}^{-1}(x_1)\td{\og}(y)\big) \  \\
        + \ \ \big(
        \td{\og}^{-1}(\ell) \td{\og}(x_1)
        \big)
    \cdot(M^{-1}\nabla_\xi)^k\big(
        \td{\og}^{-1}(x_1)\td{\og}(y)
        \big) \ . \label{large_G_k_diff}
\end{multline}
In the above identity, we have dropped the dependence on time as it
no longer has any bearing on how we proceed. Also $[x_1,\ell]$
denotes an even  smaller interval embedded in the overall
bootstrapping line segment $[y,\ell]$. We will let this smaller
segment go to zero. Before doing this, we collect the last term on
the right hand side of \eqref{large_G_k_diff} onto the left, apply
the matrix norm \eqref{matrix_inner} and the reverse triangle
inequality, and use the isometric identity \eqref{matrix_isom} to
arrive at the bound:
\begin{multline}
        \Big|
    \llp{(M^{-1}\nabla_\xi)^k\big(\td{\og}^{-1}(\ell) \td{\og}(y)\big)}
    - \llp{(M^{-1}\nabla_\xi)^k\big(
        \td{\og}^{-1}(x_1)\td{\og}(y)
        \big)}\Big|
        \ \leqslant \\
        \sum_{i=0}^{k-1} \
        \llp{(M^{-1}\nabla_\xi)^{k-i}\big(
        \td{\og}^{-1}(\ell) \td{\og}(x_1)
        \big)}
        \cdot\llp{ (M^{-1} \nabla_\xi)^i\big(
        \td{\og}^{-1}(x_1)\td{\og}(y)\big)} \ . \label{large_G_k_diff_abs}
\end{multline}
We now divide both sides of this last expression by the small
interval length $|x_1-\ell|$ and let the resulting expression go the
limit $x_1\to\ell$. To compute this, we only need to handle the
expressions:
\begin{equation}
    \lim_{x_1 \to \ell}\ \ |x_1-\ell|^{-1}\cdot\llp{(M^{-1}\nabla_\xi)^{k-i}\big(
        \td{\og}^{-1}(\ell) \td{\og}(x_1)
        \big)} \ , \label{uC_limit}
\end{equation}
where we have the important restriction $1\leqslant k-i$. We do this
by using the fact that the gauge equation \eqref{tdC_ODE} gives us
an explicit realization of the product $\td{\og}^{-1}(\ell)
\td{\og}(x_1)$ as an integral over the interval $[x_1,\ell]$:
\begin{equation}
        \td{\og}^{-1}(\ell)\td{\og}(x_1) \ = \
        \int_{x_1}^{\ell} \ \td{\og}^{-1}(x_1)\td{\og}(s)\
        \td{\uoC}_{\alpha(\ell)}(s)\ ds \  + \ I  \ . \label{sht_path_int}
\end{equation}
Here the $\alpha(\ell)$ index denotes the component of the connection
$\{\uoC\}$ in the direction of the line segment $[y,x]$. Plugging this
last expression into the limit \eqref{uC_limit} and using the fundamental
theorem of calculus on the resulting identity we arrive at the simple equation:
\begin{equation}
    \eqref{uC_limit} \ = \ \llp{(M^{-1}\nabla_\xi)^{k-i}\Big(
         \td{\uoC}_{\alpha(\ell)}(\ell)\Big)} \ . \label{uC_limit_iden}
\end{equation}
Notice that the identity matrix on line \eqref{sht_path_int} drops
out because of the condition $1\leqslant k-i$, and that all terms
where the derivatives fall on the group elements are zero because
when $x_1=\ell$ these are again just derivatives of the identity
matrix $I$. Now, substituting \eqref{uC_limit_iden} into the
limiting version of \eqref{large_G_k_diff_abs} we have the
differential inequality:
\begin{multline}
    \Big|
    \llp{ (M^{-1}\nabla_\xi)^k\big(\td{\og}^{-1}(\ell) \td{\og}(y)\big)}'\Big|
        \ \leqslant \\
        \sum_{i=0}^{k-1} \
    \llp{(M^{-1}\nabla_\xi)^{k-i}\Big(
        \td{\uoC}_{\alpha(\ell)}(\ell)\Big)}
        \cdot\llp{ (M^{-1} \nabla_\xi)^i\big(
        \td{\og}^{-1}(\ell)\td{\og}(y)\big)} \ .  \label{large_G_k_diff_ineq}
\end{multline}
Assuming now that we have proved the inductive bound:
\begin{equation}
    \sup_{ 0\leqslant  i\leqslant {k-1} } \ \llp{ (M^{-1} \nabla_\xi)^i\big(
        \td{\og}^{-1}(\ell)\td{\og}(y)\big)} \ \lesssim \ 1
    \ , \notag
\end{equation}
which is easy when $k-1 =0$ on account of the compactness of the
group $O(m)$, we see that by integrating the expression $\llp{
(M^{-1}\nabla_\xi)^k\big(\td{\og}^{-1}(\ell) \td{\og}(y)\big)}'$ the
proof of \eqref{L2_symbol_bds1} at the $k^{th}$ step boils down to
being able to establish  the line integral estimate:
\begin{equation}
    \sum_{i=0}^{k-1} \ \ \int_{y}^x\
    \llp{(M^{-1}\nabla_\xi)^{k-i}\Big(
        \td{\uoC}_{\alpha(\ell)}(\ell)\Big)}\ d\ell \ \lesssim \
    \mathcal{E} \ . \label{L1_uC_int_xy}
\end{equation}
The reason this bound will be possible is that we have taken care to
make sure that there is always at least one copy of the operator
$(M^{-1}\nabla_\xi)$ in each of the above integrals, and it is the
presence of the extra factor $M^{-1}$ in conjunction with the range
restriction \eqref{M_range} that will be enough to provide the
needed integrability. In fact, using the condition that $M^{-1}
\leqslant \ |x-y|^{-\frac{1}{2}}$ and the Cauchy-Schwartz
inequality, we see that it suffices to be able to prove the bound:
\begin{equation}
    \sum_{i=0}^{k-1} \ \ \int_{y}^x\
    \llp{(M^{-1}\nabla_\xi)^{i}\ \nabla_\xi \Big(
        \td{\uoC}_{\alpha(\ell)}(\ell)\Big)}^2\ d\ell \ \lesssim \
    \mathcal{E}^2 \ . \label{L2_uC_int_xy}
\end{equation}
This last integral can now be bounded in terms of energy type
estimates once one applies the $L^\infty \to L^2$ trace theorem to
it. However, because of the various angular degeneracies involved in
the potentials $\{\nabla_\xi \td{\uoC}\}$, it will be necessary for
us to use a more refined ``trace-Bernstein'' type inequality.
Furthermore, since the connection $\{ \td{\uoC}\}$ is only defined
implicitly via the Hodge system \eqref{td_uCpm_system}, it will be
necessary for us to prove estimate \eqref{L2_uC_int_xy} via a
bootstrapping argument in mixed Lebesgue spaces. What we will do is
to show the following somewhat more restrictive estimate which
yields \eqref{L2_uC_int_xy}
as a consequence:\\

\begin{lem}\label{non_iso_spat_lem}
Let the connection $\{\td{\uoC}\}$ be defined via the Hodge system
\eqref{td_uCpm_system}:
\begin{subequations}\label{lemma_td_uC}
\begin{align}
        (\td{\uoC})^{df} \ &= \ d^* \, \Delta^{-1} [\td{\uoC} , \td{\uoC} ]
        \ , \label{lemma_td_uC1}\\
        (\td{\uoC})^{cf} \ &= \ \widetilde{\uoA^{(M)}}
        \ - \ \nabla_x \Delta^{-1} [\widetilde{\uoA^{(M)}} , \td{\uoC} ]
        \ . \label{lemma_td_uC2}
\end{align}
\end{subequations}
where we have set:
\begin{equation}
        \widetilde{\uoA^{(M)}} \ = \ -\,
        \nabla_{x}\ \ooPi_{M^{-1}
        < \bullet} \ooPi^{(\frac{1}{2} - \delta)}
        \ \oL\ \Delta_{\omega^\perp}^{-1}\,  \uA(\partial_\omega) \ .
        \label{ucalA_shorthand}
\end{equation}
Furthermore, the parameter $M^{-1}$ which lies in the range
\eqref{M_range} (although this is not essential). Then the following
mixed Lebesgue space estimates of Besov type hold:
\begin{equation}
        \sum_{i=0}^{k-1} \ \sum_\mu \
        \lp{(M^{-1}\nabla_\xi)^{i}\ \nabla_\xi \
        P_\mu\big( \td{\uoC}\big) }{L^\infty_{\ell^\perp}(L^2_\ell)} \  \lesssim \
        \mathcal{E} \ . \label{LinftyL2_uC_int}
\end{equation}
\end{lem}\ret

\begin{proof}[Proof of estimate \eqref{LinftyL2_uC_int}]
Things will be a bit easier if we prove the following more
restrictive estimate:
\begin{equation}
        \sum_{i=0}^{k-1} \ \sum_\mu \ \mu^{-\gamma}( 1 +\mu)^{n}\
        \lp{(M^{-1}\nabla_\xi)^{i}\ \nabla_\xi \
        P_\mu\big( \td{\uoC}\big) }{L^2_\ell(L^\infty_{\ell^\perp})} \  \lesssim \
        \mathcal{E} \ . \label{L2Linfty_uC_int}
\end{equation}
That \eqref{LinftyL2_uC_int} is a consequence of
\eqref{L2Linfty_uC_int} is a simple matter applying the Minkowski
inequality for mixed Lebesgue spaces and the fact that the weights
in \eqref{L2Linfty_uC_int} are clearly more restrictive. Now, the
proof of this second estimate is essentially no more complicated
than using the Bernstein inequality in  the hyperplane plane
$\mathbb{R}_{\ell^\perp}^{n-1}$ to turn things into the energy
estimate contained in the bootstrapping norm
\eqref{red_conctn_cond5}. To see this, we begin our proof of
\eqref{L2Linfty_uC_int} by first establishing this bound for the
reduced Coulomb potentials $\{\widetilde{\uoA^{(M)}}\}$. \\

We are now trying to prove that:
\begin{equation}
        \sum_{j=0,1}\ \sum_{i=0}^{k-1} \ \sum_\mu \ \mu^{-\gamma}( 1 +\mu)^{n}\
        \lp{(M^{-1}\nabla_\xi)^{i}\ \nabla_\xi^j \
        P_\mu\big( \widetilde{\uoA^{(M)}} \big) }
    {L^2_\ell(L^\infty_{\ell^\perp})} \  \lesssim \
        \mathcal{E} \ . \label{L2Linfty_ucalA_int}
\end{equation}
For each fixed frequency in the above sum, we decompose things into
all frequencies corresponding to the $\mathbb{R}_{\ell^\perp}^{n-1}$
plane, as well as all possible dyadic angular sectors spread from
the $\omega$ (fixed) direction:
\begin{equation}
        P_\mu \ = \ \sum_{\substack{\theta,\lambda\ :
    \\ \lambda\leqslant \mu }}\
        \oPi_\theta Q_\lambda P_\mu \ , \notag
\end{equation}
where $Q_\lambda$ is an $(n-1)$ dimensional fixed frequency
multiplier which is defined in analogy with $P_\lambda$. Freezing
all frequencies, our goal will be to show the following estimate:
\begin{equation}
        \lp{(M^{-1}\nabla_\xi)^{i}\ \nabla_\xi^j \
        \oPi_\theta Q_\lambda
        P_\mu\big( \widetilde{\uoA^{(M)}}
    \big) }{L^2_\ell(L^\infty_{\ell^\perp})} \ \ \lesssim
        \ \
        \theta^\gamma\left(\frac{\lambda}{\mu}\right)^\gamma\mu^{2\gamma}
        \cdot \mathcal{E} \ . \label{L2Linfty_ucalA_int_dyadic}
\end{equation}
By adding in the weights $\mu^{-\gamma}( 1 +\mu)^{n}$, using the
fact that the potentials $\{\widetilde{\uoA^{(M)}}\}$ are truncated
to frequencies $\mu\ll 1$, and dyadic summing, the fixed frequency
estimate \eqref{L2Linfty_ucalA_int_dyadic} implies
\eqref{L2Linfty_ucalA_int} with room to spare. To deal with all of
the $\xi$ derivatives, notice that we have the following heuristic
multipliers bounds:
\begin{equation}
        (M^{-1}\nabla_\xi)^{i}\ \nabla_\xi^j \
        \oPi_\theta Q_\lambda
        P_\mu\big( \widetilde{\uoA^{(M)}}\big) \ \approx \
        \theta^{-2}\ \oPi_\theta \ooPi^{(\frac{1}{2} - \delta)}
        Q_\lambda
        P_\mu\big(\uA\big) \ , \label{Mxi_heuristic}
\end{equation}
where we are enforcing the notation introduced on line
\eqref{heuristic_op_bnd}. That is, the left hand side of the above
identity satisfies \emph{all} mixed Lebesgue space bounds as the
right hand side with the same constants. Notice that this bound uses
the extra Coulomb savings introduced on line
\eqref{main_coulomb_savings} above to kill off one power of
$\theta^{-1}$ from the degenerate Laplacean $\Delta_{\omega^\perp}$.
The other power of $\theta^{-1}$ on the right hand side of
\eqref{Mxi_heuristic} comes from the operator $\nabla_\xi$ which has
no smoothing factor of $M^{-1}$. This is precisely what one pays for
passing from the $L^1$ integral  \eqref{L1_uC_int_xy} to the more
manageable $L^2$ integral \eqref{L2_uC_int_xy}. Finally, it is
important to point out that although we have not emphasized it, the
multipliers $Q_\lambda$ depend on $\omega$, but the fact that
$\lambda \ll \theta$ implies that the multiplier product on the left
hand side of \eqref{Mxi_heuristic} is zero prevents the derivatives
of $Q_\lambda$ with respect to $\xi$ from costing more than
derivatives of $\oPi_\theta$ (alternatively, we could have applied
the $Q_\lambda$ multipliers on the outside of the $\nabla_\xi^k$
operators, because differentiation will not change the support of
the various multipliers).\\

Now, to use the Bernstein inequality on the
$\mathbb{R}_{\ell^\perp}^{n-1}$ plane, we simply note that one has
the multiplier identity:
\begin{equation}
        \oPi_\theta Q_\lambda P_\mu \ = \ {}^{\omega || \ell^\perp}\!
        B_{(\mu\theta)}
        \oPi_\theta Q_\lambda P_\mu \ , \label{parallel_muliplier}
\end{equation}
where ${}^{\omega || \ell^\perp}\! B_{(\mu\theta)}$ is a 
(smooth symbol) block type
cutoff in the $\mathbb{R}_{\ell^\perp}^{n-1}$ frequency plane of
dimensions $1\times(\mu\theta)\times\ldots\times(\mu\theta)$ which
has its long side centered along the projection\footnote{In the case
that $\omega$ lies perfectly in the $\ell$ direction, we will just
be wasteful and choose any direction.} of the unit vector $\omega$
onto the $\mathbb{R}_{\ell^\perp}^{n-1}$ (frequency) plane. The
crucial fact about the geometry of the multiplier
\eqref{parallel_muliplier} is that is has support contained in a box
of size $\lambda\times(\mu\theta)\times\ldots\times(\mu\theta)$ in
the  $\mathbb{R}^{n-1}_\xi$ (frequency) plane. Using now the
identities \eqref{Mxi_heuristic} and \eqref{parallel_muliplier}, as
well as the $n-1$ dimensional Bernstein inequality, we see that we
may estimate:
\begin{multline}
        \lp{(M^{-1}\nabla_\xi)^{i}\ \nabla_\xi^j \
        \oPi_\theta Q_\lambda
        P_\mu\big( \widetilde{\uoA^{(M)}} \big) }{L^2_\ell(L^\infty_{\ell^\perp})}
    \  \lesssim\\
        \  \theta^{-2}\cdot
        \lambda^\frac{1}{2}(\mu\theta)^\frac{n-2}{2} \
    \lp{P_\mu\big( \uA \big)}{L^2} \ .
        \label{dyadic_ucalA_Bern_bound}
\end{multline}
To deal with the weights on the right hand side, we use the
truncation condition that $\mu^{ \frac{1}{2} - \delta}\leqslant
\theta$, as well as the fact that $\lambda\leqslant \mu$ to conclude
the bound:
\begin{equation}
        \theta^{-2}\cdot\lambda^\frac{1}{2}(\mu\theta)^\frac{n-2}{2} \ \leqslant \
        \mu^\frac{n-2}{2}\cdot\theta^{\gamma}\left(\frac{\lambda}{\mu}\right)^\gamma
        \mu^{2\gamma} \ . \notag
\end{equation}
Substituting this into the right hand side of estimate
\eqref{dyadic_ucalA_Bern_bound} and using the $L^\infty
(\dot{H}^\frac{n-2}{2})$ bound
contained in the bootstrapping estimate \eqref{red_conctn_cond5}, we
have achieved the desired result
\eqref{L2Linfty_ucalA_int_dyadic}.\\

It is now our task to use \eqref{L2Linfty_ucalA_int} and the
Hodge system \eqref{lemma_td_uC} to pass to the more general
estimate \eqref{LinftyL2_uC_int}. In order to do this, it will be
necessary for us to first prove some \emph{critical} estimates for
the potentials $\{\td{\uoC}\}$. These will then be used as a
reference point in certain bilinear estimates involving the space
used to define estimate \eqref{LinftyL2_uC_int}. While we're at it,
this will also give us a chance to prove some estimates which will
be used many times in the sequel. What we will show is that:
\begin{align}
    \lp{(M^{-1}\nabla_\xi)^k\
    \uoC}{ \dot{B}_{2,10n}^{p_\gamma ,(2,\frac{n-2}{2})} } \ &\lesssim \
    \mathcal{E} \ , \label{d_omega_C_bnd1}\\
    \lp{(M^{-1}\nabla_\xi)^k
    \nabla_t \uoC }{\dot{B}_{2,10n}^{p_\gamma ,(2,\frac{n-4}{2})}}
    \ &\lesssim \ \mathcal{E} \ , \label{d_omega_C_bnd2}
\end{align}
where $p_\gamma$ is exponent defined on line \eqref{p_gamma_line}
above.  Both of the bounds
\eqref{d_omega_C_bnd1}--\eqref{d_omega_C_bnd2} will easily follow
via our general Besov embedding \eqref{general_besov_embed}
once we have established them for the linear term on the right hand
side of the Hodge system \eqref{lemma_td_uC}. That is, we fist
establish that:
\begin{align}
    \lp{(M^{-1}\nabla_\xi)^k\
    \widetilde{\uoA^{(M)}}
    }{ \dot{B}_{2,10n}^{p_\gamma ,(2,\frac{n-2}{2})}} \ &\lesssim \
    \mathcal{E} \ , \label{d_omega_calA_bnd1}\\
    \lp{(M^{-1}\nabla_\xi)^k
    \nabla_t\ \widetilde{\uoA^{(M)}} }{\dot{B}_{2,10n}
    ^{p_\gamma ,(2,\frac{n-4}{2})}}
    \ &\lesssim \ \mathcal{E} \ , \label{d_omega_calA_bnd2}
\end{align}
These follow from immediately from the steps used to prove
\eqref{Besov_degen_tdA_Ln_smallness} above, and the following
heuristic identity which follows our convention established on line
\eqref{heuristic_op_bnd}:
\begin{align}
        \nabla^k_\xi \left(\oPi_\theta P_\mu\, \widetilde{\uoA^{(M)}}\right) \
        &\approx \ \theta^{-k}\
        \oPi_\theta \ooPi_{M^{-1}< \bullet}\ P_\mu \,
        \nabla_x\ \oL\ \Delta_{\omega^\perp}^{-1}\,
        \uA(\partial_\omega)\ , \label{heuristic_dxi_bnds1}\\
        \nabla^k_\xi \left(\oPi_\theta P_\mu\, \nabla_t\
        \widetilde{\uoA^{(M)}}\right) \
        &\approx \ \mu \theta^{-k}\
        \oPi_\theta \ooPi_{M^{-1} < \bullet}\ P_\mu \,
        \nabla_x\ \oL\ \Delta_{\omega^\perp}^{-1}\,
    \uA(\partial_\omega) \ , \label{heuristic_dxi_bnds2}
\end{align}
Notice that the space-time frequency localization \eqref{red_conctn_cond4}
allows us to trade the $\nabla_t$ with the factor of $\mu$
on the second line above.\\

We now prove the estimates
\eqref{d_omega_C_bnd1}--\eqref{d_omega_C_bnd2} by proceeding
inductively on the value of $k$. If $k=0$ the first estimate
\eqref{d_omega_C_bnd1} holds because one can solve the system
\eqref{lemma_td_uC} via Picard iteration in the space
$\dot{B}_{2,10n}^{p_\gamma ,(2,\frac{n-2}{2})}$ thanks to the
bilinear embedding \eqref{general_besov_embed} which furnishes the
embedding:
\begin{equation}
    \nabla_x\Delta^{-1} \ : \
    \dot{B}_{2,10n}^{p_\gamma ,(2,\frac{n-2}{2})}\cdot
    \dot{B}_{2,10n}^{p_\gamma ,(2,\frac{n-2}{2})} \ \hookrightarrow\
    \dot{B}_{2,10n}^{p_\gamma ,(2,\frac{n-2}{2})} \ .
    \label{reisz_-1_critical}
\end{equation}
The key thing to point out here is that for $\gamma$ sufficiently small,
and in dimensions $6\leqslant n$
we have the bound $p_\gamma < n$, which is
all that is needed to satisfy the gap condition \eqref{gap_cond} in
this case. The other conditions of Lemma \ref{Besov_lem} are also easily
seen to be satisfied for this set of indices.\\

To establish \eqref{d_omega_C_bnd1} for $0 < k$, we simply differentiate the
system \eqref{lemma_td_uC} $k$ times with respect to the operator
$(M^{-1}\nabla_\xi)$. Doing this yields the linearized set of
equations:
\begin{align}
        (M^{-1}\nabla_\xi)^k(\td{\uoC})^{df} \
        &= \ \sum_{j=0}^k\
        d^* \, \Delta^{-1} \Big[(M^{-1}\nabla_\xi)^{k-j}\ \td{\uoC}
        \ , \
        (M^{-1}\nabla_\xi)^j \ \td{\uoC} \Big]
        \ , \label{d_omega_td_uCpm_system1}\\
    \begin{split}
        (M^{-1}\nabla_\xi)^k(\td{\uoC})^{cf} \
        &= \  (M^{-1}\nabla_\xi)^k \ \widetilde{\uoA^{(M)}}
        \ - \\
        &\ \   \ \sum_{j=0}^k\ \nabla_x \Delta^{-1}
        \Big[(M^{-1}\nabla_\xi)^{k-j}\ \widetilde{\uoA^{(M)}}
        \ , \
        (M^{-1}\nabla_\xi)^j\ \td{\uoC} \Big] \ ,
    \end{split} \label{d_omega_td_uCpm_system2}
\end{align}
which can again be solved in the Besov space $\dot{B}_{2,10n}^{p_\gamma
,(2,\frac{n-2}{2})}$ by using the already
established estimate \eqref{d_omega_calA_bnd1} for the linear term,
in conjunction with the (inductive) hypothesis that estimate
\eqref{d_omega_C_bnd1} holds for $k-1$, and absorbing the highest
derivative (involving $(M^{-1}\nabla_\xi)$ falling on
$\td{\uoC}$) term to the left hand side. All of this is permissible
by referring to the embedding \eqref{reisz_-1_critical}.\\

To prove the second estimate \eqref{d_omega_C_bnd2} above, we first
apply the time derivative $\nabla_t$ to both sides of the system
\eqref{d_omega_td_uCpm_system1}--\eqref{d_omega_td_uCpm_system2}
above. The resulting system of equations, which we will not write
down, can easily be solved in the derivative critical Besov space
$\dot{B}_{2,10n}^{p_\gamma ,(2,\frac{n-4}{2})}$ by again using an
induction on $k$, the already established estimate
\eqref{d_omega_calA_bnd2} for the linear term, and the following
bilinear Besov estimate which is again a special case of
\eqref{general_besov_embed}:
\begin{equation}
    \nabla_x\Delta^{-1} \ : \
    \dot{B}_{2,10n}^{p_\gamma ,(2,\frac{n-2}{2})}\cdot
    \dot{B}_{2,10n}^{p_\gamma ,(2,\frac{n-4}{2})} \ \hookrightarrow\
    \dot{B}_{2,10n}^{p_\gamma ,(2,\frac{n-4}{2})} \ .
    \label{reisz_-1_dcritical}
\end{equation}
Notice that \eqref{reisz_-1_dcritical}
is permissible because for $\gamma$ sufficiently small, we have
the condition $p_\gamma < \frac{2n}{3}$ in dimensions $6\leqslant n$
which is necessary to get around the gap condition \eqref{gap_cond}.
The other conditions of \eqref{general_besov_embed} are easily satisfied
for this choice of indices.\\

Armed with estimates \eqref{L2Linfty_ucalA_int} and
\eqref{d_omega_C_bnd1}, we now move back to the proof of estimate
\eqref{LinftyL2_uC_int}. We set the norm in that latter bound equal
to:
\begin{equation}
        \lp{A}{\mathcal{N}_1^{-\gamma,2,\infty}} \ = \
        \sum_\mu\ \mu^{-\gamma}(1+\mu)^{n}\
        \lp{P_\mu(A)}{L^2_{\ell}(L^\infty_{\ell^\perp})}
        \ . \notag
\end{equation}
By differentiating the system \eqref{lemma_td_uC}
with respect to the operators $(M^{-1}\nabla_\xi)^k\nabla_\xi$, 
we see that the
claim will now follow once we can prove the bilinear Riesz operator
bound:
\begin{equation}
        \nabla_x\Delta^{-1} \ : \
    \dot{B}_{2,10n}^{p_\gamma ,(2,\frac{n-2}{2})}
        \cdot\mathcal{N}_1^{-\gamma,2,\infty} \
        \hookrightarrow \ \mathcal{N}_1^{-\gamma,2,\infty} \ .
        \label{bilin_reisz_N_bnd}
\end{equation}
We now let $A$ and $C$ be any two elements of the two spaces on the
left hand side of \eqref{bilin_reisz_N_bnd}. By applying the
trichotomy, we see that it suffices to be able to prove the three
estimates:
\begin{align}
    \begin{split}
        \sum_{\substack{\lambda , \mu_i \ :\\
        \mu_1 \ll \mu_2 \\
        \lambda\sim \mu_2 }}\
        \lambda^{-\gamma}(1+\lambda)^n \
        \lp{\nabla_x\Delta^{-1} P_\lambda\left(P_{\mu_1}A\cdot P_{\mu_2}C
        \right)}
        {L^2_{\ell}(L^\infty_{\ell^\perp})} \ &\lesssim \\
        &\hspace{-1.5in} \lp{A}{\dot{B}_{2,10n}^{p_\gamma ,
        (2,\frac{n-2}{2})}
        }\cdot\lp{C}{\mathcal{N}_1^{-\gamma,2,\infty}}
        \ ,
    \end{split} \label{tric_N1}\\
    \begin{split}
        \sum_{\substack{\lambda , \mu_i \ :\\
        \mu_2 \ll \mu_1 \\
        \lambda\sim \mu_1 }}\
        \lambda^{-\gamma}(1+\lambda)^n \
        \lp{\nabla_x\Delta^{-1} P_\lambda\left(P_{\mu_1}A\cdot P_{\mu_2}C
        \right)}
        {L^2_{\ell}(L^\infty_{\ell^\perp})} \ &\lesssim \\
        &\hspace{-1.5in} \lp{A}{\dot{B}_{2,10n}^{p_\gamma ,
        (2,\frac{n-2}{2})}
        }\cdot\lp{C}{\mathcal{N}_1^{-\gamma,2,\infty}}
        \ ,
    \end{split} \label{tric_N2}\\
    \begin{split}
        \sum_{\substack{\lambda , \mu_i \ :\\
        \mu_1 \sim \mu_2 \\
        \lambda\lesssim \mu_1,\mu_2 }}\
        \lambda^{-\gamma}(1+\lambda)^n \
        \lp{\nabla_x\Delta^{-1} P_\lambda\left(P_{\mu_1}A\cdot P_{\mu_2}C
        \right)}
        {L^2_{\ell}(L^\infty_{\ell^\perp})} \ &\lesssim \\
        &\hspace{-1.5in} \lp{A}{\dot{B}_{2,10n}^{p_\gamma ,
        (2,\frac{n-2}{2})}
        }\cdot\lp{C}{\mathcal{N}_1^{-\gamma,2,\infty}}
        \ .
    \end{split} \label{tric_N3}
\end{align}
The proofs of \eqref{tric_N1}--\eqref{tric_N2} are very simple, and
follow from essentially the same principle. First of all, we use the
fact that the kernel of the fixed frequency operator $\lambda\cdot
\nabla_x\Delta^{-1}P_\lambda$ is in $L^1$ with norm independent of
$\lambda$. Thus, it is bounded on all mixed Lebesgue spaces. This,
used in conjunction with the estimate:
\begin{equation}
        \sum_{\substack{ \mu_1 \ :\\ \mu_1 \lesssim \lambda}} \
        \lambda^{-1}\ \lp{P_{\mu_1} A}{L^\infty}\ \lesssim \
        \lp{A}{\dot{B}_{2,10n}^{p_\gamma ,
        (2,\frac{n-2}{2})}} \ , \label{A_linfty_besov}
\end{equation}
which follows easily from Bernstein's inequality and dyadic summing,
is enough for us to conclude the first estimate \eqref{tric_N1}.
To conclude
the second estimate, \eqref{tric_N2}, we simply employ the fixed
frequency version of \eqref{A_linfty_besov} and then estimate:
\begin{align}
        \sum_{\substack{\lambda , \mu_2 \ :\\
        \mu_2\lesssim \lambda}}\
        \lambda^{-\gamma}\
        \lp{P_{\mu_2}C}{L^2_{\ell}(L^\infty_{\ell^\perp})}\ &= \
        \sum_{\substack{\lambda , \mu_2 \ :\\
        \mu_2\lesssim \lambda}}\
         \left(\frac{\mu_2}{\lambda}\right)^\gamma
        \mu_2^{-\gamma}
        \lp{P_{\mu_2}C}{L^2_{\ell}(L^\infty_{\ell^\perp})} \ ,
        \notag\\
        &= \ \sum_{\mu_2}\ \mu_2^{-\gamma}\
        \lp{P_{\mu_2}C}{L^2_{\ell}(L^\infty_{\ell^\perp})} \ \cdot \
        \sum_{\substack{\lambda \ :\\ \mu_2
        \lesssim \lambda}}\
        \left(\frac{\mu_2}{\lambda}\right)^\gamma \ , \notag\\
        &\lesssim \ \lp{C}{\mathcal{N}_1^{-\gamma,2,\infty}}
        \ . \notag
\end{align}
We now move to the last estimate \eqref{tric_N3}. This is only
slightly more complicated than what we have already done. Here we
will only bother to estimate things for $\lambda\lesssim 1$. The
case where $1 \ll \lambda$ is much easier due to the extra smoothing
in the norms we are working with and is left to the reader. As a
first step, we freeze all frequencies and decompose $P_\lambda =
\sum_{\substack{\sigma\ :\\ \sigma\leqslant \lambda}}\
Q_{\sigma}P_\lambda$, where $Q_{\sigma}$ is again the fixed
frequency cutoff in the $\mathbb{R}^{n-1}_{\ell^\perp}$ frequency
plane. Using Bernstein, this allows us to estimate:
\begin{align}
        &\lambda^{-\gamma} \
        \lp{\nabla_x\Delta^{-1} P_\lambda\left(P_{\mu_1}A\cdot P_{\mu_2}C
        \right)}
        {L^2_{\ell}(L^\infty_{\ell^\perp})} \ , \label{tric3_N_start}\\
        \lesssim \ &\sum_{\substack{ \sigma\ : \\ \sigma\leqslant
        \lambda}}\ \lambda^{-1-\gamma}\
        \lp{Q_\sigma \left(P_{\mu_1}A\cdot P_{\mu_2}C
        \right)}
        {L^2_{\ell}(L^\infty_{\ell^\perp})}\ , \notag\\
        \lesssim \ &\sum_{\substack{ \sigma\ : \\ \sigma\leqslant
        \lambda}}\ \lambda^{-1-\gamma}\sigma^{\frac{n-1}{p_\gamma}}\
        \lp{P_{\mu_1}A\cdot P_{\mu_2}C}
        {L^2_{\ell}(L^{p_\gamma}_{\ell^\perp})}\ , \notag\\
        \lesssim \ &\ \lambda^{\frac{n-1}{p_\gamma}-1-\gamma}\
        \lp{P_{\mu_1}A\cdot P_{\mu_2}C}
        {L^2_{\ell}(L^{p_\gamma}_{\ell^\perp})}\ . \notag
\end{align}
By using H\"olders inequality and rearranging weights, and using the
index bound:
\begin{equation}
        1 + 2\gamma < \frac{n-1}{p_\gamma} \ , \notag
\end{equation}
which follows on account of the fact that $\gamma\ll 1$  and we are
in $n\leqslant 6$ dimensions, we see that we have the fixed
frequency estimate:
\begin{equation}
        \eqref{tric3_N_start} \ \lesssim \
        \left(\frac{\lambda}{\mu_1}\right)^\gamma\
        \mu_1^{\frac{n-1}{p_\gamma} -1}
        \lp{P_{\mu_1}A}{L_\ell^\infty(L_{\ell^\perp}^{p_\gamma})}
        \cdot \mu_2^{-\gamma}\lp{P_{\mu_2}C}{L^2_\ell(L^\infty_{\ell^\perp})}
        \ . \label{tric_N_finish}
\end{equation}
We are done once we deal with the first factor on the right hand
side of the above inequality. To do this, we run a multiplier
decomposition $P_{\mu_1} = \sum_{\substack{\sigma\ :\\
\sigma\leqslant \mu_1}} \td{Q}_\sigma P_{\mu_1}$, where this time
$\td{Q}_\sigma$ is a fixed frequency cutoff on the $\mathbb{R}_\ell$
line. Using Minkowski's integral inequality and Bernstein on the
real line, we  estimate:
\begin{align}
        \mu_1^{\frac{n-1}{p_\gamma} -1}
        \lp{P_{\mu_1}A}{L_\ell^\infty(L_{\ell^\perp}^{p_\gamma})} \
        &\lesssim \ \sum_{\substack{\sigma\ :\\ \sigma\leqslant
        \mu_1}}\
        \mu_1^{\frac{n-1}{p_\gamma} -1}
        \lp{\td{Q}_\sigma
        P_{\mu_1}A}{L_{\ell^\perp}^{p_\gamma}(L^\infty_\ell)} \ ,
        \notag\\
        &\lesssim \ \mu_1^{\frac{n}{p_\gamma} -1}
        \lp{P_{\mu_1}A}{L^{p_\gamma}} \ , \notag\\
        &= \ \lp{P_{\mu_1}A}{\dot{B}_2^{p_\gamma,(2,\frac{n-2}{2})}} \ . \notag
\end{align}
Plugging this last bound into the right hand side of
\eqref{tric_N_finish} and summing over all $\lambda \lesssim
\mu_1,\mu_2$ yields the bound \eqref{tric_N3} as was to be shown.
This completes the proof of the bilinear estimate
\eqref{bilin_reisz_N_bnd} and hence the proof of
\eqref{LinftyL2_uC_int}.
\end{proof}\ret

To wrap things up here, we need to prove the symbol bounds
\eqref{st_symbol_bds1}--\eqref{st_symbol_bds2} for the second factor
in the product \eqref{simple_product}. By repeating the steps which
started on line \eqref{large_G_k_diff} and culminated in the
differential inequality \eqref{large_G_k_diff_ineq} for this term,
we arrive at the temporal differential inequality:
\begin{multline}
    \Big|
    \llp{ (M^{-1}\nabla_\xi)^k\big(\td{\og}^{-1}(\ell) \td{\og}(s)\big)}'\Big|
        \ \leqslant \\
        \sum_{i=0}^{k-1} \
    \llp{(M^{-1}\nabla_\xi)^{k-i}\Big(
        \td{\oC}_{0}(\ell)\Big)}
        \cdot\llp{ (M^{-1} \nabla_\xi)^i\big(
        \td{\og}^{-1}(\ell)\td{\og}(s)\big)} \ ,  \label{temporal_diff_ineq}
\end{multline}
where this time $\ell$ denotes a single variable which lies in the
range $s\leqslant \ell \leqslant t$, and $\td{\oC}_0$ is the
temporal potential which is defined via the equation:
\begin{equation}
        \td{\og}^{-1}\nabla_t(\td{\og}) \ = \ \td{\oC}_0 \ . \notag
\end{equation}
Via integration in time of the quantity $\llp{
(M^{-1}\nabla_\xi)^k\big(\td{\og}^{-1}(\ell) \td{\og}(s)\big)}'$,
and keeping in mind the derivation of the temporal potential
equation \eqref{C0_system} above, we see that to prove the estimates
\eqref{st_symbol_bds1} for the product $\td{\og}^{-1}(t)\td{\og}(s)$
it suffices to show (the same estimate works to establish
\eqref{st_symbol_bds2}):\\

\begin{lem}
Let the temporal potential $\td{C}_0$ be defined via the elliptic
equation:
\begin{equation}
        \td{\oC}_0 \ = \   \td{\oA^{(M)}_0}
        \ - \ \nabla_t \Delta^{-1} [\widetilde{\uoA^{(M)}} , \td{\uoC} ]
        \ - \ d^* \Delta^{-1} [\td{\oC}_0 , \td{\uoC} ]
        \ , \label{td_C0_system}
\end{equation}
where the spatial connection $\{\td{\uoC}\}$ is defined via the
Hodge system \eqref{lemma_td_uC}, and where we have set:
\begin{equation}
        \td{\oA^{(M)}_0} \ = \ -\,
        \nabla_{t}\ \ooPi_{M^{-1}
        < \bullet} \ooPi^{(\frac{1}{2} - \delta)}
        \ \oL\ \Delta_{\omega^\perp}^{-1}\,  \uA(\partial_\omega) \ ,
        \label{calA0_shorthand}
\end{equation}
The parameter $M^{-1}$ which lies in the range \eqref{M_range} (this 
is essential). 
Then the following mixed Lebesgue
space estimates of Besov type hold:
\begin{equation}
        \sum_{i=1}^{k} \ \sum_\mu \
        \lp{(M^{-1}\nabla_\xi)^{i} \
        P_\mu\big( \td{\oC}_0\big) }{L^\infty_x(L^1_t[s,t])} \  \lesssim \
        \mathcal{E} \ . \label{LinftyL1_C0_int}
\end{equation}
\end{lem}\ret

\begin{proof}[Proof of the estimate \eqref{L1Linfty_C0_int}]
As with the proof of \eqref{LinftyL2_uC_int} above, it will be
convenient to prove the somewhat more restrictive estimate:
\begin{equation}
        \sum_{i=1}^{k} \ \sum_\mu \ \mu^{-\gamma}(1 + \mu)^n\
        \lp{(M^{-1}\nabla_\xi)^{i} \
        P_\mu\big( \td{\oC}_0\big) }{L^1_t[s,t](L^\infty_x)} \  \lesssim \
        \mathcal{E} \ . \label{L1Linfty_C0_int}
\end{equation}
This will again be done by essentially proving that this estimate is
true for the potentials $\{\td{\oA^{(M)}}\}$ contained in the right hand
side of \eqref{td_C0_system}, and then transferring that knowledge
to $\td{\oC}_0$ through that elliptic equation. A little care needs
to be taken in this regard due to the effect of bad $High\times
High\Rightarrow Low$ frequency interactions coming from the
$\Delta^{-1}$ in the second term on the right hand side of
\eqref{td_C0_system} which sits by itself because the time
derivative must be distributed. This all needs to be tempered
against the fact that we need to recover enough in the low
frequencies to apply $L^2(L^q)$ Strichartz estimates to integrate
over the line segment $[s,t]$. The Lebesgue exponent which is
important in this regard is the following:
\begin{equation}
    q_\gamma \ = \ \frac{2(n-1)}{n-5-2\gamma} \ . \label{q_gamma_def}
\end{equation}
The significance of $q_\gamma$ is that it is the smallest Lebesgue
exponent such that one can recover an extra angular weight of $\theta^{-1}$
via running Bernstein from the Strichartz endpoint and still have an
extra factor of $\theta^\gamma$ to spare for dyadic summing.\\

We proceed with our proof of \eqref{L1Linfty_C0_int} by first establishing
a fixed time estimate. Notice that the Besov norm
$\dot{B}_{1,10n}^{\infty,(2,\frac{n}{2}-\gamma)}$ embeds into the spatial
norm on the left hand side of \eqref{L1Linfty_C0_int}. Thus, our first step
will be to establish the following fixed time estimate for $1\leqslant k$:
\begin{equation}
        \lp{(M^{-1}\nabla_\xi)^k\
        \td{\oC}_0 (t_0)}{\dot{B}_{1,10n}
    ^{\infty,(2,\frac{n}{2}-\gamma)}}
        \ \lesssim \
        \sum_{i=1}^k\ \lp{(M^{-1}\nabla_\xi)^i\
        \td{\oA^{(M)}} (t_0)}{\dot{B}_{2,10n}
    ^{q_\gamma,(2,\frac{n}{2}-2\gamma)}} \ ,
        \label{diff_C_Linfty_A}
\end{equation}
where $\{\td{\oA^{(M)}}\}$ is the full space-time connection defined on lines
\eqref{ucalA_shorthand} and \eqref{calA0_shorthand}. The important thing
about estimate \eqref{diff_C_Linfty_A} is that we retain at least one
copy of the operator $(M^{-1}\nabla_\xi)$ on the right hand side so that
we may pass to an $L^2_t$ integral via Cauchy-Schwartz.
Our proof of \eqref{diff_C_Linfty_A} require that the
bootstrapping constant $\mathcal{E}$
from line \eqref{red_conctn_cond5} is sufficiently small. Based on previous work,
our task here is largely finished.
Our first step here will be to
differentiate the equation \eqref{L1Linfty_C0_int} as many times as
necessary with respect to the operators $(M^{-1}\nabla_\xi)$. Doing this
and distributing the time derivative in the second term on the right hand
side yields the equation:
\begin{align}
    (M^{-1}\nabla_\xi)^k\ \td{\oC}_0 \ &= \
    (M^{-1}\nabla_\xi)^k\  \td{\oA^{(M)}_0} \ \label{diff_td_C0_system} \\
    &- \ \ \sum_{i=0}^k\ \Big(
    \Delta^{-1} [(M^{-1}\nabla_\xi)^{k-i}\nabla_t\widetilde{\uoA^{(M)}}
    , (M^{-1}\nabla_\xi)^i \td{\uoC} ]\notag\\
    &+ \ \
    \Delta^{-1} [(M^{-1}\nabla_\xi)^{k-i}\widetilde{\uoA^{(M)}}
    , (M^{-1}\nabla_\xi)^i \nabla_t \td{\uoC} ]\Big) \notag\\
    &- \ \ \sum_{i=0}^k \  d^* \Delta
    [(M^{-1}\nabla_\xi)^{k-i}\td{\oC}_0 ,
    (M^{-1}\nabla_\xi)^i \td{\uoC} ]  \ , \notag\\
    &= \ T_1 \ + \ T_2 \ + \ T_3 \ . \notag
\end{align}
Our second step is to prove the intermediate estimate:
\begin{multline}
        \lp{(M^{-1}\nabla_\xi)^k\
        \td{\oC}_0 (t_0)}{\dot{B}_{1,10n}
        ^{\infty,(2,\frac{n}{2}-\gamma)}}
        \ \lesssim \\
        \lp{(M^{-1}\nabla_\xi)^k\
        \td{\oA^{(M)}} (t_0)}{\dot{B}_{2,10n}
        ^{q_\gamma,(2,\frac{n}{2}-2\gamma)}} \ +
        \ \sum_{i=1}^k\ \lp{(M^{-1}\nabla_\xi)^i\
        \td{\uoC} (t_0)}{\dot{B}_{2,10n}
        ^{q_\gamma,(2,\frac{n}{2}-\gamma)}}\\
        + \ \sum_{i=1}^k\ \lp{(M^{-1}\nabla_\xi)^i\
        \nabla_t \td{\uoC} (t_0)}{\dot{B}_{2,10n}
        ^{q_\gamma,(2,\frac{n-2}{2}-\gamma)}} \ .
        \label{prelim_diff_C_Linfty_A}
\end{multline}
This in turn is a consequence of the three estimates:
\begin{align}
        \lp{T_1}{\dot{B}_{1,10n}
        ^{\infty,(2,\frac{n}{2}-\gamma)}}
        \ &\lesssim \ \lp{(M^{-1}\nabla_\xi)^k\
        \td{\oA^{(M)}} (t_0)}{\dot{B}_{2,10n}
        ^{q_\gamma,(2,\frac{n}{2}-2\gamma)}} \ , \label{p_diff_C_Linfty_A1}\\
   \begin{split}
        \lp{T_2}{\dot{B}_{1,10n}
        ^{\infty,(2,\frac{n}{2}-\gamma)}}
        \ &\lesssim \ \lp{(M^{-1}\nabla_\xi)^k\
        \td{\oA^{(M)}} (t_0)}{\dot{B}_{2,10n}
        ^{q_\gamma,(2,\frac{n}{2}-\gamma)}} \\
        &\ \ \ \ \ \ + \
        \sum_{i=1}^k\ \lp{(M^{-1}\nabla_\xi)^i\
        \td{\uoC} (t_0)}{\dot{B}_{2,10n}
        ^{q_\gamma,(2,\frac{n}{2}-\gamma)}}\\
        &\ \ \ \ \  \ + \ \sum_{i=1}^k\ \lp{(M^{-1}\nabla_\xi)^i\
        \nabla_t \td{\uoC} (t_0)}{\dot{B}_{2,10n}
        ^{q_\gamma,(2,\frac{n-2}{2}-\gamma)}} \ ,
   \end{split}\label{p_diff_C_Linfty_A2}\\
   \begin{split}
        \lp{T_3}{\dot{B}_{1,10n}^{\infty,(2,\frac{n}{2}-\gamma)}}
        \ &\lesssim \
        \lp{(M^{-1}\nabla_\xi)^k\
        \td{\oC}_0 (t_0)}{\dot{B}_{1,10n}
        ^{\infty,(2,\frac{n}{2}-\gamma)}}\cdot
        \lp{ \td{\uoC} (t_0)}{\dot{B}_{2,10n}
        ^{p_\gamma,(2,\frac{n-2}{2})}}\\
        &+ \ \sum_{i=1}^k\ \lp{(M^{-1}\nabla_\xi)^i\
        \td{\uoC} (t_0)}{\dot{B}_{2,10n}
        ^{q_\gamma,(2,\frac{n}{2}-\gamma)}} \ .
   \end{split}\label{p_diff_C_Linfty_A3}
\end{align}
Notice that these all combine to give \eqref{prelim_diff_C_Linfty_A} because
the estimate \eqref{d_omega_C_bnd1}
in conjunction with the assumption that $\mathcal{E}$
is sufficiently small allows one to absorb the first term on the right hand
side of \eqref{p_diff_C_Linfty_A3} into the left hand side. The proof of the
first estimate, \eqref{p_diff_C_Linfty_A1} above is a trivial consequence
of the Besov nesting \eqref{Besov_nesting},
and the fact that we allow for an extra
power of $\mu^\gamma$ to sum over the low frequencies to turn the $\ell^2$
sum into and $\ell^1$ sum.\\

The proof of \eqref{p_diff_C_Linfty_A2} is the most involved, and is why
we have been forced to work with the exponent $q_\gamma$. There are several
cases to consider, depending on whether or not the time derivative falls
on the term containing at least one copy of the operator
$(M^{-1}\nabla_\xi)$. An inspection of the structure of the $T_2$ term
shows that these can all be taken into account through an application
of the already established estimates
\eqref{d_omega_C_bnd1}--\eqref{d_omega_C_bnd2}
and \eqref{d_omega_calA_bnd1}--\eqref{d_omega_calA_bnd2},
and application of the bilinear embeddings:
\begin{align}
    \Delta^{-1} \ : \ \dot{B}_{2,10n}
    ^{p_\gamma,(2,\frac{n-2}{2})}\cdot
    \dot{B}_{2,10n}
    ^{q_\gamma,(2,\frac{n-2}{2}-\gamma)} \ &\hookrightarrow\
    \dot{B}_{1,10n}
    ^{\infty,(2,\frac{n}{2}-\gamma)} \ , \label{p_diff_C_Linfty_A2_besov1}\\
    \Delta^{-1} \ : \ \dot{B}_{2,10n}
    ^{p_\gamma,(2,\frac{n-4}{2})}\cdot
    \dot{B}_{2,10n}
    ^{q_\gamma,(2,\frac{n}{2}-\gamma)} \ &\hookrightarrow\
    \dot{B}_{1,10n}
    ^{\infty,(2,\frac{n}{2}-\gamma)} \ . \label{p_diff_C_Linfty_A2_besov2}
\end{align}
A quick calculation shows that one has the needed gap bound
(by condition \eqref{gap_cond}):
\begin{equation}
    2 + \gamma \ < \ n(\frac{1}{p_\gamma} + \frac{1}{q_\gamma}) \ ,
    \notag
\end{equation}
for $\gamma$ sufficiently small when the dimension satisfies the
bound $6\leqslant n$. For example, when $n=6$ we have that $p_\gamma
= \frac{10}{3} + \epsilon$ and $q_\gamma = 10 + \epsilon$ where
$\epsilon \to 0$ as $\gamma \to 0$. Notice that there is not a whole
lot of room in this. The other condition in
\eqref{sc_cond}--\eqref{Lb_cond} are easily verified in the above
estimates. Notice that for the term in $T_2$ where all the
$(M^{-1}\nabla_\xi)$ derivatives, as well as the time derivative
$\nabla_t$ fall on the linear potentials
$\{\widetilde{\uoA^{(M)}}\}$, we use \eqref{d_omega_C_bnd1} and the
first embedding \eqref{p_diff_C_Linfty_A2_besov1} above, together
with the easy bound (which follows from the truncation
\eqref{red_conctn_cond4}):
\begin{equation}
    \lp{(M^{-1}\nabla_\xi)^k\nabla_t \widetilde{\uoA^{(M)}}}
    {\dot{B}_{2,10n}
    ^{q_\gamma,(2,\frac{n-2}{2}-\gamma)}} \ \lesssim \
    \lp{(M^{-1}\nabla_\xi)^k\ \widetilde{\uoA^{(M)}}}
    {\dot{B}_{2,10n}
    ^{q_\gamma,(2,\frac{n}{2}-\gamma)}} \ , \label{A_dt_to_dx}
\end{equation}
to conclude \eqref{p_diff_C_Linfty_A2} for that portion of things.\\

To conclude our second main step in the proof of \eqref{diff_C_Linfty_A}, we
need to establish the estimate \eqref{p_diff_C_Linfty_A3}. To tie things
down, we first need to know that $\td{\oC}_0$ satisfies a critical estimate
similar to \eqref{d_omega_C_bnd1}. This is:
\begin{equation}
    \lp{(M^{-1}\nabla_\xi)^i \td{\oC}_0}
    {\dot{B}_{2,10n}^{p_\gamma,(2,\frac{n-2}{2})}} \ \lesssim \
    \mathcal{E} \ . \label{critical_C0_bound}
\end{equation}
This in turn is provided through applying the already established estimates
\eqref{d_omega_C_bnd1}--\eqref{d_omega_C_bnd2}
and \eqref{d_omega_calA_bnd2}
to the equation \eqref{td_C0_system} with the help of the bilinear estimate
\eqref{reisz_-1_critical} and the following embedding
which follows as yet
another  special case of our general bound
\eqref{general_besov_embed}:
\begin{equation}
    \Delta^{-1} \ : \ \dot{B}_{2,10n}
    ^{p_\gamma,(2,\frac{n-4}{2})}\cdot
    \dot{B}_{2,10n}
    ^{p_\gamma,(2,\frac{n-2}{2})} \ \hookrightarrow\
    \dot{B}_{1,10n}
    ^{p_\gamma,(2,\frac{n-2}{2})} \ . \label{p_diff_C_Linfty_A3_besov1}
\end{equation}
Armed with the estimate \eqref{critical_C0_bound}, we can proceed to prove
\eqref{p_diff_C_Linfty_A3} by applying the following bilinear
estimates to the various terms contained in $T_3$:
\begin{align}
    \nabla_x \Delta \ : \ \dot{B}_{2,10n}
    ^{\infty,(2,\frac{n}{2}-\gamma)}\cdot
    \dot{B}_{2,10n}
    ^{p_\gamma,(2,\frac{n-2}{2})} \ &\hookrightarrow\
    \dot{B}_{1,10n}
    ^{\infty,(2,\frac{n}{2}-\gamma)} \ , \label{p_diff_C_Linfty_A3_besov2}\\
    \nabla_x \Delta \ : \ \dot{B}_{2,10n}
    ^{p_\gamma,(2,\frac{n-2}{2})}\cdot
    \dot{B}_{2,10n}
    ^{q_\gamma,(2,\frac{n}{2}-\gamma)} \ &\hookrightarrow\
    \dot{B}_{1,10n}
    ^{q_\gamma,(2,\frac{n}{2}-\gamma)} \ . \label{p_diff_C_Linfty_A3_besov3}
\end{align}
We use the second embedding \eqref{p_diff_C_Linfty_A3_besov3} in conjunction
with the nesting \eqref{Besov_nesting} to derive the second term on the
right hand side of \eqref{p_diff_C_Linfty_A3}.
We note that in estimates
\eqref{p_diff_C_Linfty_A3_besov2}--\eqref{p_diff_C_Linfty_A3_besov3},
it is a simple matter to check the validity of the conditions
\eqref{sc_cond}--\eqref{Lb_cond}.
We leave this as an exercise for the reader.\\

To complete this portion of the proof, we need to establish the
implication
\eqref{prelim_diff_C_Linfty_A}$\Rightarrow$\eqref{diff_C_Linfty_A}.
This will be done once we can show that (keeping in mind the bounds
of the form \eqref{A_dt_to_dx}):
\begin{align}
        \sum_{i=1}^k\ \lp{(M^{-1}\nabla_\xi)^i\
        \td{\uoC} (t_0)}{\dot{B}_{2,10n}
        ^{q_\gamma,(2,\frac{n}{2}-\gamma)}}
        \ &\lesssim \ \sum_{i=1}^k\ \lp{(M^{-1}\nabla_\xi)^i\
        \td{\oA^{(M)}} (t_0)}{\dot{B}_{2,10n}
        ^{q_\gamma,(2,\frac{n}{2}-\gamma)}}\ , \label{to_diff_C_Linfty_A1}\\
   \begin{split}
        \sum_{i=1}^k\ \lp{(M^{-1}\nabla_\xi)^i\
        \nabla_t \td{\uoC} (t_0)}{\dot{B}_{2,10n}
        ^{q_\gamma,(2,\frac{n-2}{2}-\gamma)}} \
        &\lesssim \ \sum_{i=1}^k\ \lp{(M^{-1}\nabla_\xi)^i\
        \nabla_t \td{\oA^{(M)}} (t_0)}{\dot{B}_{2,10n}
        ^{q_\gamma,(2,\frac{n-2}{2}-\gamma)}}\\
        &\hspace{-2.5in} \ + \
        \sum_{i=1}^k\ \lp{(M^{-1}\nabla_\xi)^i\
        \td{\oA^{(M)}} (t_0)}{\dot{B}_{2,10n}
        ^{q_\gamma,(2,\frac{n}{2}-\gamma)}}\ + \
        \sum_{i=1}^k\ \lp{(M^{-1}\nabla_\xi)^i\
        \td{\uoC} (t_0)}{\dot{B}_{2,10n}
        ^{q_\gamma,(2,\frac{n}{2}-\gamma)}}\ .
   \end{split} \label{to_diff_C_Linfty_A2}
\end{align}
The estimate \eqref{to_diff_C_Linfty_A1} is a simple consequence of
applying the embedding \eqref{p_diff_C_Linfty_A3_besov3} to the differentiated
Hodge system \eqref{d_omega_td_uCpm_system1}--\eqref{d_omega_td_uCpm_system2},
while using the already established critical estimates
\eqref{d_omega_C_bnd1} and \eqref{d_omega_calA_bnd1} to tie things down.
To prove the second estimate \eqref{to_diff_C_Linfty_A2} above, we apply
the time derivative $\nabla_t$ to the system
\eqref{d_omega_td_uCpm_system1}--\eqref{d_omega_td_uCpm_system2}, and
then employ the embeddings:
\begin{align}
    \nabla_x \Delta \ : \ \dot{B}_{2,10n}
    ^{p_\gamma,(2,\frac{n-4}{2})}\cdot
    \dot{B}_{2,10n}
    ^{q_\gamma,(2,\frac{n}{2}-\gamma)} \ &\hookrightarrow\
    \dot{B}_{1,10n}
    ^{q_\gamma,(2,\frac{n-2}{2}-\gamma)} \ , \label{to_diff_C_besov1}\\
    \nabla_x \Delta \ : \ \dot{B}_{2,10n}
    ^{p_\gamma,(2,\frac{n-2}{2})}\cdot
    \dot{B}_{2,10n}
    ^{q_\gamma,(2,\frac{n-2}{2}-\gamma)} \ &\hookrightarrow\
    \dot{B}_{1,10n}
    ^{q_\gamma,(2,\frac{n-2}{2}-\gamma)} \ . \label{to_diff_C_besov2}
\end{align}
Notice that these estimates have the same small amount of
room as \eqref{p_diff_C_Linfty_A2_besov1}--\eqref{p_diff_C_Linfty_A2_besov2}
above when measuring the gap condition \eqref{gap_cond}
for this set of exponents. Using
\eqref{to_diff_C_besov1}--\eqref{to_diff_C_besov2} in conjunction
with the already established estimates
\eqref{d_omega_C_bnd1}--\eqref{d_omega_C_bnd2}
and \eqref{d_omega_calA_bnd1}--\eqref{d_omega_calA_bnd2}, we may conclude
\eqref{to_diff_C_Linfty_A2} when the bootstrapping constant $\mathcal{E}$
is sufficiently small.\\

We have now established the estimate \eqref{diff_C_Linfty_A}. Integrating
this in time, and applying a Cauchy-Schwartz with respect to the
time integration and using the condition $|t-s|^\frac{1}{2} \leqslant M$, we
have the estimate:
\begin{equation}
        \hbox{(L.H.S.)}\eqref{L1Linfty_C0_int} \ \lesssim \
        \sum_{i=0}^{k-1}\ \lp{(M^{-1}\nabla_\xi)^i\
        \nabla_\xi\ \td{\oA^{(M)}}}{L^2_t\big(\dot{B}_{2,10n}
        ^{q_\gamma,(2,\frac{n}{2}-2\gamma)}\big)} \ . \notag
\end{equation}
Therefore, to conclude the estimate \eqref{LinftyL1_C0_int}, we simply need
to prove the bound:
\begin{equation}
        \sum_{i=0}^{k-1}\ \lp{(M^{-1}\nabla_\xi)^i\
        \nabla_\xi\ \td{\oA^{(M)}}}{L^2_t\big(\dot{B}_{2,10n}
        ^{q_\gamma,(2,\frac{n}{2}-2\gamma)}\big)}
        \ \lesssim \ \mathcal{E} \ . \label{temp_int_last_step}
\end{equation}
At a heuristic level, this estimate is true because there is enough
room in the norm $L^2_t\big(\dot{B}_{2,10n}^
{q_\gamma,(2,\frac{n}{2}-2\gamma)}\big)$ verses the bootstrapping
norm \eqref{red_conctn_cond5} to save precisely $\frac{1}{2} -
2\gamma$ derivatives. This, used in conjunction with the truncation
condition coming form the operator $\ooPi^{(\frac{1}{2}-\delta)}$,
is enough to absorb the extra angular factor $\theta^{-1}$ produced
by the unsmoothed derivative $\nabla_\xi$. All throughout the
calculation, the exponent $q_\gamma$ is high enough that the
intrinsic angular singularity contained in the potentials
$\{\td{\oA^{(M)}}\}$ can be recovered by an application of Bernstein's
inequality to the endpoint Strichartz spatial exponent
$L^\frac{2(n-1)}{n-3}$.
We spell out briefly the details of this procedure as follows:\\

Freezing now the frequency and the number of $(M^{-1}\nabla_\xi)$
derivatives on the right hand side of \eqref{temp_int_last_step},
and using the bootstrapping condition \eqref{red_conctn_cond5},
we see that it suffices to show the bound
(note that \eqref{temp_int_last_step}
is already in square function form):
\begin{equation}
        \lp{(M^{-1}\nabla_\xi)^i\
        \nabla_\xi\ P_\mu( \td{\oA^{(M)}})}{L^2_t\big(\dot{B}_{2,10n}
        ^{q_\gamma,(2,\frac{n}{2}-2\gamma)}\big)}
        \ \lesssim \  \lp{P_\mu (\uA)}
        {L^2\big(\dot{B}_{2,10n}^{\frac{2(n-1)}{n-3} ,
        (2,\frac{n-1}{2})}\big)} \ . \notag
\end{equation}
In fact, after a further localization in the angle, we will show that:
\begin{equation}
        \lp{(M^{-1}\nabla_\xi)^i\
        \nabla_\xi\ \oPi_\theta\
        P_\mu( \td{\oA^{(M)}})}{L^2_t\big(\dot{B}_{2,10n}
        ^{q_\gamma,(2,\frac{n}{2}-2\gamma)}\big)}
        \ \lesssim \  \theta^\gamma\ \lp{P_\mu (\uA)}
        {L^2\big(\dot{B}_{2,10n}^{\frac{2(n-1)}{n-3} ,
        (2,\frac{n-1}{2})}\big)} \ . \notag
\end{equation}
By an application of the Bernstein inequality, and recalling the definition
\eqref{q_gamma_def} of the exponent $q_\gamma$, this last estimate is a
consequence of being able to show that:
\begin{multline}
        \lp{(M^{-1}\nabla_\xi)^i\
        \nabla_\xi\ \oPi_\theta\
        P_\mu( \td{\oA^{(M)}})}{L^2_t\big(\dot{B}_{2,10n}
        ^{\frac{2(n-1)}{n-3},(2,\frac{n}{2}-2\gamma)}\big)}
        \ \lesssim \\
        \theta^{-1}\ \lp{P_\mu (\uA)}
        {L^2\big(\dot{B}_{2,10n}^{\frac{2(n-1)}{n-3} ,
        (2,\frac{n-1}{2})}\big)} \ . \label{temp_int_last_step_dyadic}
\end{multline}
Using now the heuristic operator bound \eqref{main_degen_symbol_bound}
in conjunction with the Coulomb savings \eqref{main_coulomb_savings} and
the heuristic symbol type bounds
\eqref{heuristic_dxi_bnds1}, we have the following
heuristic identity which follows our strict convention
\eqref{heuristic_op_bnd}:
\begin{equation}
        (M^{-1}\nabla_\xi)^i\
        \nabla_\xi\ \oPi_\theta\
        P_\mu( \td{\oA^{(M)}}) \ \approx \
        \theta^{-2}\ P_\mu (\uA) \ . \notag
\end{equation}
Plugging this last bound into the left hand side of
\eqref{temp_int_last_step_dyadic}, and
rearranging the Besov weights, we have the estimate:
\begin{multline}
        \lp{(M^{-1}\nabla_\xi)^i\
        \nabla_\xi\ \oPi_\theta\
        P_\mu( \td{\oA^{(M)}})}{L^2_t\big(\dot{B}_{2,10n}
        ^{\frac{2(n-1)}{n-3},(2,\frac{n}{2}-2\gamma)}\big)}
        \ \lesssim \\
        \left(\frac{\mu^{\frac{1}{2} - 2\gamma}}
        {\theta}\right)\theta^{-1}\
        \lp{P_\mu (\uA)}
        {L^2\big(\dot{B}_{2,10n}^{\frac{2(n-1)}{n-3} ,
        (2,\frac{n-1}{2})}\big)} \ . \notag
\end{multline}
The truncation condition that $\mu^{\frac{1}{2} - \delta} \leqslant \theta$
(note that $\mu\lesssim 1$) now guarantees that we have the bound:
\begin{equation}
    \left(\frac{\mu^{\frac{1}{2} - 2\gamma}}
    {\theta}\right) \ \lesssim \ 1 \ , \notag
\end{equation}
when $\gamma$ is sufficiently small compared to  $\delta$. This
completes the proof of \eqref{temp_int_last_step_dyadic}, and hence
\eqref{temp_int_last_step}, which in turn finishes our proof of
estimate \eqref{L1Linfty_C0_int}.
\end{proof}\ret

Having now established the symbol bounds
\eqref{st_symbol_bds1}--\eqref{st_symbol_bds2} separately for each
of the two terms in the product \eqref{simple_product}. By using the
Leibniz rule for derivatives, these together imply the bounds
\eqref{st_symbol_bds1}--\eqref{st_symbol_bds2} for the full product
on the left hand side of \eqref{simple_product}. This completes our
proof of Proposition \ref{symbol_bound_prop}.
\end{proof}\ret

We now proceed to prove the second main estimate \eqref{Lp_TT*_bound} for
remainder kernel in the splitting \eqref{L2_TT*_splitting}. This involves
a sum of kernels, each of which according to the identities
\eqref{gh_left}--\eqref{gh_right} has at
least one copy of the terms ${\oh}^{-1}(x){\oh}(y)
- I$ and ${\oh}(x){\oh}^{-1}(y) - I $. Therefore, without loss
of generality, we may assume that we are trying to prove the estimate:
\begin{equation}
    	\int_{\mathbb{R}_x^n}\ \chi_{\mathcal{D}_\sigma}(x)\ \llp{
    	\int_{\mathbb{R}_\xi^n} \
        e^{2\pi i ( x-y)\cdot \xi}\
        {}^\omega G(x,y) \chi(\xi)\ d\xi}\ dx
    	\ \lesssim \ \sigma^{-\gamma}
    	\ . \label{small_L2_est}
\end{equation}
where we have set:
\begin{equation}
    {}^\omega G(x,y) \ = \
    \td{\og}^{-1}(x)\left({\oh}^{-1}(x){\oh}(y) - I
    \right) \td{\og}(y) \big[\bullet \big]\
        \og^{-1}(y) \og(x) \ . \notag
\end{equation}
We note here that the corresponding estimates
for the other terms in $\mathcal{R}^{TT^*}_\sigma$
are similar and are left to the reader.\\

To prove \eqref{small_L2_est}, we use following angular cutoff functions
to split:
\begin{equation}
    	{}^\omega G \ = \
    	\chi_{|\cos (\theta_{\xi,x-y} )| 
	\geqslant |x-y|^{{-1}+\gamma}}\, {}^\omega G \ + \
	\chi_{|\cos (\theta_{\xi,x-y} )| <  
	|x-y|^{-1+\gamma} }\, {}^\omega G \ . \notag
\end{equation}
Therefore, using the triangle and Minkowski inequalities, we see that it suffices
to prove the pair of bounds:
\begin{align}
    \int_{\mathbb{R}_x^n}\ \chi_{\mathcal{D}_\sigma}(x)\ \llp{
    \int_{\mathbb{R}_\xi^n} \
        e^{2\pi i ( x-y)\cdot \xi}\
    \chi_{|\cos (\theta_{\xi,x-y} )| \geqslant |x-y|^{{-1}+\gamma}}
        {}^\omega G(x,y) \chi(\xi)\ d\xi}\ dx
    \ &\lesssim \ \sigma^{-\gamma}
    \ , \label{small_L2_estA}\\
    \int_{\mathbb{R}_x^n}\ \chi_{\mathcal{D}_\sigma}(x)\ \llp{
    \int_{\mathbb{R}_\xi^n} \
        e^{2\pi i ( x-y)\cdot \xi}\
    \chi_{|\cos (\theta_{\xi,x-y} )| < |x-y|^{{-1}+\gamma}}
        {}^\omega G(x,y) \chi(\xi)\ d\xi}\ dx
    \ &\lesssim \ \sigma^{-\gamma}
    \ .
    \label{small_L2_estB}
\end{align}
The proof of the first estimate, \eqref{small_L2_estA}, is a simple matter of
integrating by parts as many times as necessary with respect to the
weighted radial derivative $\frac{1}{2\pi i |x-y| \cos(\theta_{\xi,x-y})
} \partial_{|\xi|}$, taking account of the fact that ${}^\omega G$ is
independent of the variable $|\xi|$. Assuming that $|x-y|\sim\sigma$ is
sufficiently large, we will eventually have that:
\begin{equation}
    \Big| (\frac{1}{2\pi i |x-y| \cos(\theta_{\xi,x-y})
    } \partial_{|\xi|})^k \chi(\xi)\Big| \
    \lesssim \ \sigma^{-n -\gamma} \ , \label{diff_chi_bnd}
\end{equation}
at which point we may stop the integration by parts and put absolute value signs
around the remaining integral. The right hand side of
\eqref{small_L2_estA} will then
follow as a direct consequence of \eqref{diff_chi_bnd}
and the simple bounds:
\begin{align}
    \int_{\mathbb{R}_x^n}\ \chi_{\mathcal{D}_\sigma}(x)\ dx \ &\lesssim \
    \sigma^n \ , \notag\\
    \sup_{x,\omega} \ \llp{{}^\omega G(x,y)} \ &\lesssim \ 1 \ . \notag
\end{align}
To conclude the proof of estimate \eqref{small_L2_est}, we need to show the
second estimate \eqref{small_L2_estB} above. At this point, we have stripped
things down to where oscillations under the integral sign are no longer
of any use, so we simply strive to estimate the absolute value of the
integrand. Here the smallness of the function ${}^\omega G(x,y)$ is essential.
To make use of this, we rearrange the order in the absolute integral
and use H\"olders inequality to bound:
\begin{multline}
    \hbox{(L.H.S.)}\eqref{small_L2_estB} \ \lesssim \
    \int_{\mathbb{S}^{n-1}}\ \sup_{x\in \mathcal{D}_\sigma}
    \llp{{}^\omega G(x,y) }\ d\omega \\
    \cdot \ \ \sup_\omega\  \int_{\mathbb{R}^n_x} \
    \chi_{|\cos (\theta_{\xi,x-y} )| < |x-y|^{{-1}+\gamma}}(x)\
    \chi_{\mathcal{D}_\sigma}(x)\ dx
    \ . \label{G_chi_prelim_est}
\end{multline}
To bound the second integral on the right hand side of the above product, we
translate by the vector $y$ and then apply a rotation to reduce the
bound we wish to show to the following:
\begin{equation}
    \int_{|x| \sim \sigma} \
    \chi_{|\cos (\theta_{(1,0),x} )| < |x|^{{-1}+\gamma}}(x)\
    dx \ \lesssim \ \sigma^{n-1 + \gamma}
    \ . \label{truncates_x_int}
\end{equation}
The validity of \eqref{truncates_x_int} follows trivially from the fact that
if we split $x=(x_1,x')$, we have the bounds $|x_1|\lesssim \sigma^\gamma$ and
$|x'|\lesssim \sigma$ over the range of integration thanks to the angular
cutoff and the identity:
\begin{equation}
    \cos (\theta_{(1,0),x} ) \ = \ \frac{x_1}{ |x|} \ . \notag
\end{equation}
Thus, keeping in mind the bound \eqref{truncates_x_int}, we see from
estimate \eqref{G_chi_prelim_est} that the proof of
\eqref{small_L2_estB} follows from a Cauchy-Schwartz on the sphere
$\mathbb{S}^{n-1}$ and the following integrated bounds:
\begin{equation}
    \left(\int_{\mathbb{S}^{n-1}}\ \sup_{x\in \mathcal{D}_\sigma}
    \llp{{}^\omega G(x,y) }^2\ d\omega
    \right)^\frac{1}{2} \ \lesssim \
    \sigma^{1-n - 2\gamma} \ . \label{int_G_bound}
\end{equation}
Due to its use in the next section, we will in fact show the following
more general set of estimates which includes \eqref{int_G_bound} as
a special case:\\

\begin{prop}[Estimates for integrated remainder group elements ]\label{h_int_rem_est}
Let the group elements ${\oh}$ be defined infinitesimally via the equations
\eqref{dh_C_eq}--\eqref{Cdt_def} and the Hodge system
\eqref{td_td_uDpm_system}, where the parameter $\sigma^{-1+\gamma}$ is replaced by
$M^{-1}$. Then upon integration, one has the following bounds:
\begin{align}
    \left( \int_{\mathbb{S}^{n-1}}\ \sup_{|x-y|\sim N}\
    \llp{ {\oh}^{-1}(t,x) {\oh}(s,y) - I }^2\
    d\omega \right)^\frac{1}{2} \ &\lesssim \ \mathcal{E}\,
    (1 + |t-s| + N )\cdot M^{-n - \delta} \ , \label{h_int_bnd1}\\
    \left(\int_{\mathbb{S}^{n-1}}\ \sup_{|x-y|\sim N}\
    \llp{ {\oh}(t,x) {\oh}^{-1}(s,y) - I }^2\
    d\omega\right)^\frac{1}{2}  \ &\lesssim \ \mathcal{E}\,
    (1 + |t-s| + N )\cdot M^{-n - \delta} \ , \label{h_int_bnd2}
\end{align}
where $\mathcal{E}$ is the bootstrapping constant from line
\eqref{red_conctn_cond5}. The above estimates
are uniform in the value of $M$ when it is sufficiently large.
\end{prop}\ret

\begin{proof}[Proof of the estimates \eqref{h_int_bnd1}--\eqref{h_int_bnd2}]
As will become apparent to the reader, it suffices to show the first bound
\eqref{h_int_bnd1}, as the second follows from essentially identical reasoning.
Our first step here is to disentangle the products, and to work exclusively with
either spatially separated or temporally separated products.
This is accomplished via the following simple algebraic identity:
\begin{multline}
    {\oh}^{-1}(t,x) {\oh}(s,y) - I \ = \
    {\oh}^{-1}(t,x) {\oh}(s,x) -I \\ + \ \
    {\oh}^{-1}(t,x) {\oh}(s,x)\cdot
    \left({\oh}^{-1}(s,x) {\oh}(s,y) - I\right)
    \ . \label{temp_spat_h_split}
\end{multline}
Working now, for the moment, with the second term in
this last expression substituted
into the estimate \eqref{h_int_bnd1} we expand:
\begin{equation}
    {\oh}^{-1}(s,x) {\oh}(s,y) - I \ =
    \ \int^y_x \ {\oh}^{-1}(s,x)
    \partial_\ell({\oh}(s,\ell))\ d\ell \ . \notag
\end{equation}
Integrating this last expression in $L^2_\omega$ we
are reduced to proving the following:\\

\begin{lem}\label{small_spat_int_lem}
Let the (spatial) connection $\{\td{\td{\uoC}}\}$ be defined via the
Hodge system:
\begin{subequations}\label{lemma_td_td_uDpm_system}
\begin{align}
    (\td{\td{\uoC}})^{df} \ &= \ d^* \, \Delta^{-1} \big(
        [\td{\uoC} , \td{\td{\uoC}} ] \ + [\td{\td{\uoC}} , \td{\uoC}]
        \ \big) \ , \label{lemma_td_td_uDpm_system1}\\
        (\td{\td{\uoC}})^{cf} \ &= \
    \td{\td{\uoA^{(M)}}} \
        -  \ \nabla_x \Delta^{-1} \Big([ \td{\td{\uoA^{(M)}}}
    , \td{\uoC} ]
        \ + \ [ \td{\uoA^{(M)}} , \td{\td{\uoC}}]
        \Big)\ , \label{lemma_td_td_uDpm_system2}
\end{align}
\end{subequations}
where we have set:
\begin{equation}
    \td{\td{\uoA^{(M)}}} \ = \ -\, \nabla_{x}\
    \ooPi_{\bullet\leqslant M^{-1}}
    \ooPi^{(\frac{1}{2} - \delta)}
        \ \oL\ \Delta_{\omega^\perp}^{-1}\,  \uA(\partial_\omega) \ ,
    \label{lemma2_tdtdA_shorthand}
\end{equation}
and where the spatial connections $\{\td{\uoA^{(M)}}\}$ and
$\{\td{\uoC}\}$ are defined on
the lines \eqref{lemma_td_uC} and \eqref{ucalA_shorthand} above.
Then one has the following integrated estimate uniform in the parameter
$M$:
\begin{equation}
     \left(\int_{\mathbb{S}^{n-1}}\ \sup_x
        \llp{ \td{\td{\uoC}}(x)
     }^2\
        d\omega \right)^\frac{1}{2} \ \lesssim \ \mathcal{E}\,
        \cdot M^{-n - \delta} \ . \label{tdtdC_int_bnd}
\end{equation}
\end{lem}\ret

\begin{proof}[Proof of the estimate \eqref{tdtdC_int_bnd}]
Our strategy here is similar to the previous lemmas. We first prove
things for the linear term in \eqref{lemma_td_td_uDpm_system2}, and
then use the critical embeddings \eqref{d_omega_calA_bnd1} and
\eqref{d_omega_C_bnd1} to transfer things to the
connection $\{\td{\td{\uoC}}\}$ via the Hodge system
\eqref{lemma_td_td_uDpm_system}.\\

Our first step then is to show that:
\begin{equation}
    \left(\int_{\mathbb{S}^{n-1}}\ \sup_x
        \llp{ \td{\td{  \uoA^{(M)}   }}(x)
     }^2\
        d\omega \right)^\frac{1}{2} \ \lesssim \ \mathcal{E}\,
        \cdot M^{-n - \delta} \ . \label{tdtdA_int_bnd}
\end{equation}
In fact, we will show the following somewhat stronger estimate which
will easily imply \eqref{tdtdA_int_bnd}, and which is more robust with
respect to Hodge systems:
\begin{equation}
        \left(\int_{\mathbb{S}^{n-1}}\
        \lp{ \td{\td{  \uoA^{(M)}   }}
        }{\dot{B}_{2,10n}^{\infty,(2,\frac{n}{2}-\gamma)}}^2\
        d\omega \right)^\frac{1}{2} \ \lesssim \ \mathcal{E}\,
        \cdot M^{-n - \delta} \ . \label{better_tdtdA_int_bnd}
\end{equation}
This last estimate is a simple matter of using Bernstein's inequality and
orthogonality which will net us the factor $M^{1-n}$, followed by
the condition that $\mu^{\frac{1}{2}-\delta} \lesssim M^{-1 -
\delta}$ at each fixed frequency thanks to the
$\ooPi_{\bullet\leqslant M^{-1}} \ooPi^{(\frac{1}{2} - \delta)}$
multiplier which nets us the remaining powers of $M^{-1}$. The
implementation is as follows: We first decompose things into the sum
over all frequencies $\mu\lesssim 1$ and angles $\theta \lesssim 1$:
\begin{equation}
    \td{\td{\uoA^{(M)}}} \ = \ \sum_{\substack{\theta,\mu\ :\\
    \mu\lesssim 1}}  \
        \oPi_\theta P_\mu \td{\td{\uoA^{(M)}}} \ . \notag
\end{equation}
Keeping in mind the spatial frequency truncation of $\td{\td{\uoA^{(M)}}}$,
and by the square sum definition of the Besov norms, the
triangle inequality, and dyadic summing, we see that it
suffices to show the following fixed frequency estimate:
\begin{equation}
        \left(\int_{\mathbb{S}^{n-1}}\ \sup_x
        \llp{  \oPi_\theta P_\mu \td{\td{\uoA^{(M)}}}(x)
         }^2\
        d\omega \right)^\frac{1}{2} \ \lesssim \
        \mu^{\gamma} \theta^\gamma
        \cdot\mathcal{E}\,
        \cdot M^{-n - \delta} \ . \label{tdtdA_int_bnd_fixed}
\end{equation}
For each fixed $\omega$, we use Bernstein's inequality and the
equivalence:
\begin{equation}
    \oPi_\theta P_\mu \td{\td{\uoA^{(M)}}} \ \approx \
    \theta^{-1} \oPi_\theta\ooPi_{\bullet\leqslant M^{-1}}
    \ooPi^{(\frac{1}{2} - \delta)}
    P_\mu\big(\uA\big) \ , \notag
\end{equation}
to compute that:
\begin{align}
        &\sup_x\ \llp{  \oPi_\theta P_\mu \td{\td{\uoA^{(M)}}}(x)
        } \ , \notag\\
        \lesssim \ \ &\theta^{-1}\mu\cdot
        \theta^{\frac{n-1}{2}}
        \cdot\, \lp{\oPi_\theta
        \ooPi_{\bullet\leqslant M^{-1}}
        \ooPi^{(\frac{1}{2} - \delta)}
        P_\mu\big(\uA\big)
        }{\dot H^\frac{n-2}{2}_x}\ , \notag\\
        \lesssim \ \ &\mu^{\gamma} \theta^\gamma \cdot
        M^{-\frac{n+1}{2}- \delta}\, \lp{
        \ooPi_{\bullet\leqslant M^{-1}}
        P_\mu\big(\uA\big)
        }{\dot H^\frac{n-2}{2}_x} \ . \notag
\end{align}
Notice that this last line follows from the truncation condition
$\mu^{1-2\delta} \lesssim \theta^2$ as well as the small constant
bounds $\gamma\ll \delta$.
The proof of \eqref{tdtdA_int_bnd_fixed} is now a result of the
following simple calculation involving Plancherel:
\begin{align}
        &\left(\int_{\mathbb{S}^{n-1}}\
        \lp{
        \ooPi_{\bullet\leqslant M^{-1}}
        P_\mu\big(\uA\big)
        }{\dot H^\frac{n-2}{2}_x}^2\
        d\omega \right)^\frac{1}{2} \ , \label{begin_orth_comp}\\
        \lesssim \
        &\left(\int_{\RR_\xi} \int_{\mathbb{S}^{n-1}}\
        \llp{
        (b^\omega_{\bullet\leqslant M^{-1}}
    + b^{-\omega}_{\bullet\leqslant M^{-1}})
        \ |\xi|^\frac{n-2}{2}p_\mu \widehat{\uA}(\xi)
        }^2\
        d\omega d\xi \right)^\frac{1}{2} \ , \notag\\
        \lesssim \ &M^{-\frac{n-1}{2}}\,
        \lp{
        P_\mu\big(\uA\big)
        }{\dot H^\frac{n-2}{2}_x}
        \ , \notag \\
        \lesssim \ &\mathcal{E}\cdot M^{-\frac{n-1}{2}} \ . \notag
\end{align}\ret

To finish the proof of \eqref{tdtdC_int_bnd}, we simply need to pass
the estimate \eqref{better_tdtdA_int_bnd} onto the set of spatial
potentials $\{\td{\td{\uoC}}\}$. To do this, we set up auxiliary
spaces $L^\infty_\omega (\dot
B_{2,10n}^{p_\gamma,(2,\frac{n-2}{2})})$ and $L^2_\omega (\dot
B_{2,10n}^{\infty,(2,\frac{n}{2}-\gamma)})$. From the estimates
\eqref{d_omega_calA_bnd1} and \eqref{d_omega_C_bnd1} we immediately
have that:
\begin{align}
        \lp{\td{\uoA^{(M)}}}{L^\infty_\omega (\dot
    B_{2,10n}^{p_\gamma,(2,\frac{n-2}{2})})}
        \ &\lesssim \ \mathcal{E} \ , \label{Linfty_omega_A_bnd}\\
        \lp{\td{\uoC}}{L^\infty_\omega (\dot B_{2,10n}^{p_\gamma,(2,\frac{n-2}{2})})}
        \ &\lesssim \ \mathcal{E} \ , \label{Linfty_omega_uC_bnd}
\end{align}
where the index $p_\gamma$ is the exponent from the line
\eqref{p_gamma_line} above.
The desired result now follows from the bilinear estimate:
\begin{equation}
        \nabla_x\Delta^{-1} \ : \
        L^2_\omega(\dot{B}_{2,10n}^{\infty ,(2,\frac{n}{2}-\gamma)})
    \cdot L^\infty_\omega(\dot{B}_{2,10n}^{p_\gamma
    ,(2,\frac{n-2}{2})}) \ \hookrightarrow\
        L^2_\omega(\dot{B}_{2,10n}^{\infty ,(2,\frac{n}{2}-\gamma)}) \ .
       \label{L2_omega_besov_bound}
\end{equation}
This is a simple consequence of the condition $p_\gamma < n$
which allows us to fulfill the condition \eqref{gap_cond} of
the general embedding \eqref{general_besov_embed}.  The result follows from
integrating this bound  in $L^2_\omega$.
\end{proof}\ret

We now turn our attention to proving the bound \eqref{h_int_bnd1}
for the temporally separated product which is the first term on the
right hand side of equation \eqref{temp_spat_h_split} above. Expand
the integrand here an the derivative of another integral over time
line, we have that:
\begin{equation}
    {\oh}^{-1}(t,x) {\oh}(s,x) -I \ = \ \int_s^t \
    {\oh}^{-1}(t,x)\partial_t( {\oh}(\ell,x) ) \ d\ell \
    . \notag
\end{equation}
After integrating in $L^2_\omega$ the right hand side of this last
expression, we see that we are reduced to proving that:\\

\begin{lem}
Let the quantity $\td{\td{\oC}}_0$ be defined implicitly
via the elliptic equation:
\begin{multline}
    \td{\td{\oC}}_0 \ = \ \td{\td{\oA_0^{(M)}}}
    \ - \ \nabla_t \Delta^{-1}\Big(
     [\td{\td{ {\uoA^{(M)}} }} , \td{\uoC} ] +
     [\td{ {\uoA^{(M)}} } , \td{\td{\uoC}} ] \Big)
        \\ - \ d^* \Delta^{-1}\Big( [\td{\td{\oC}}_0 , \td{\uoC} ] +
    \ [\td{\oC}_0 , \td{\td{\uoC}} ] \Big)
        \ , \label{lemma2_td_C0_system}
\end{multline}
where we have set:
\begin{equation}
        \td{\td{\oA_0^{(M)}}} \ = \ -\, \nabla_{t}\
        \ooPi_{\bullet\leqslant M^{-1}}
        \ooPi^{(\frac{1}{2} - \delta)}
        \ \oL\ \Delta_{\omega^\perp}^{-1}\,  \uA(\partial_\omega) \ ,
    \label{lemma2_tdtdA0_shorthand}
\end{equation}
and where the connections $\{\td{ \td{\uoA^{(M)}} }\}$ and
$\{\td{ \td{\uoC} }\}$ are as in Lemma \ref{small_spat_int_lem}, and where the
quantity $\td{\oA^{(M)}_0}$ is defined on line \eqref{calA0_shorthand}
above. Then one has the following integrated estimate
uniform in the parameter $M$:
\begin{equation}
     \left(\int_{\mathbb{S}^{n-1}}\ \sup_x
        \llp{ \td{\td{\oC}}_0(x) }^2\
        d\omega \right)^\frac{1}{2} \ \lesssim \ \mathcal{E}\,
        \cdot M^{-n - \delta} \ . \label{tdtdC0_int_bnd}
\end{equation}
\end{lem}\ret

\begin{proof}[Proof of the estimate \eqref{tdtdC0_int_bnd}]
As in the proof of the previous Lemma, our goal here is to first
prove the
$L^2_\omega(\dot{B}_{2,10n}^{\infty,(n,\frac{n}{2}-\gamma)})$
improvement of this claim for the terms on the right hand side of
the equation \eqref{lemma2_td_C0_system} which do not involve the
variable $\td{\td{\oC}}_0$. The desired bound can then be achieved
via iteration or bootstrapping using the bilinear estimate
\eqref{L2_omega_besov_bound} and the estimate
\eqref{Linfty_omega_uC_bnd} to deal with the term involving
$\td{\td{\oC}}_0$ on the right hand side of
\eqref{lemma2_td_C0_system}. Therefore, we are trying to show the following
three estimates:
\begin{align}
        \left(\int_{\mathbb{S}^{n-1}}\
        \lp{ \td{\td{  \oA_0^{(M)}   }}
        }{\dot{B}_{2,10n}^{\infty,(2,\frac{n}{2}-\gamma)}}^2\
        d\omega \right)^\frac{1}{2} \ &\lesssim \ \mathcal{E}\,
        \cdot M^{-n - \delta} \ , \label{better_tdtdA0_int_bnd}\\
        \left(\int_{\mathbb{S}^{n-1}}\
        \lp{  d^* \Delta^{-1}
        \ [\td{\oC}_0 , \td{\td{\uoC}} ]  }
        {\dot{B}_{2,10n}^{\infty,(2,\frac{n}{2}-\gamma)}}^2\
        d\omega \right)^\frac{1}{2} \ &\lesssim \ \mathcal{E}\,
        \cdot M^{-n - \delta} \ , \label{better_easydx_int_bnd}\\
        \left(\int_{\mathbb{S}^{n-1}}\
        \lp{  \frac{\nabla_t}{ \Delta}\Big(
        [\td{\td{ {\uoA^{(M)}} }} , \td{\uoC} ] +
        [\td{ {\uoA^{(M)}} } , \td{\td{\uoC}} ] \Big)   }
        {\dot{B}_{2,10n}^{\infty,(2,\frac{n}{2}-\gamma)}}^2\
        d\omega \right)^\frac{1}{2} \ &\lesssim \ \mathcal{E}\,
        \cdot M^{-n - \delta} \ . \label{better_baddt_int_bnd}
\end{align}
The proof of the first estimate \eqref{better_tdtdA0_int_bnd} is
essentially identical to that of \eqref{better_tdtdA_int_bnd} above,
once one takes into account the truncation condition
\eqref{red_conctn_cond4}. The proof of the second estimate
\eqref{better_easydx_int_bnd} follows from the
$L^2_\omega(\dot{B}_{2,10n}^{\infty,(2,\frac{n}{2}-\gamma)})$ bound
proved for the potentials $\{\td{\td{\uoC}}\}$ proved in the
previous lemma, the bilinear estimate \eqref{L2_omega_besov_bound},
and the following:
\begin{equation}
        \lp{\td{\oC_0}}{L^\infty_\omega(\dot{B}_{2,10n}^{p_\gamma
        ,(2,\frac{n-2}{2})})} \ \lesssim \ \mathcal{E} \ , \notag
\end{equation}
which is a direct consequence of \eqref{critical_C0_bound} above.\\

Therefore, it remains for us to prove the last bound
\eqref{better_baddt_int_bnd}. Unfortunately, this does not follow
directly from the procedure we have been using so far. The trouble
is that the time derivative $\nabla_t$ will in general not cancel
with the Laplacean, and it is not possible to prove a bilinear
estimate which is morally of the form
$\dot{B}^{p_\gamma,(2,\frac{n-4}{2})}_2\cdot\dot{B}_2^{\infty,(2,\frac{n}{2})}
\subseteq \Delta \dot{B}_2^{\infty,(2,\frac{n}{2})}$ due to bad
$High\times High$ frequency interactions in dimension $n=6$. The
only way around this seems to be to  do something which is quite a
bit more involved. The way we will prove
\eqref{better_baddt_int_bnd} is in a series of steps designed to
reduce things to a term which, in some sense, represents the central
difficulty. This last term will be dealt with using a scale of
non-isotropic spaces which are similar to the ones employed in the
proof of Lemma \ref{non_iso_spat_lem} above. The argument we will
present here is largely ad-hoc, and there are many variations. Furthermore,
we will proceed by proving certain estimates which may be cut out at
this stage of the overall paper but will turn out to be useful in
the sequel.\\

The first step we make here is to recall that, although we have been
suppressing it, there is additional polarity information in the
definition of the connections $d+ \oA $ (see \eqref{A_omega_def}).
This comes from the choice of null vector-field $\oL^\pm$. For
convenience, we will use here an implicitly defined notation which
we call $\uL$, to denote the \emph{opposite} vector-field for any
given choice of polarization. That is, we always have the formula:
\begin{equation}
        \Box \ = \ \uL\oL + \Delta_{\omega^\perp} \ .
        \label{LLb_wave}
\end{equation}
Now, for a given choice of polarization, we can always write $\pm
\nabla_t = \oL^\pm - \omega\cdot\nabla_x$. Therefore, modulo
proving estimates which are identical to those Lemma
\ref{small_spat_int_lem} above, and distributing the $\uL$
derivative, we that the proof of \eqref{better_baddt_int_bnd} can be
reduced to the proof of the following three bilinear estimates:
\begin{align}
        \lp{\Delta^{-1}[\uL\td{ {\uoA^{(M)}} } , \td{\td{\uoC}}]}
        {L^2_\omega(\dot{B}_{2,10n}^{\infty,(2,\frac{n}{2}-\gamma)})}
        \ &\lesssim \ \mathcal{E}\cdot M^{-n - \delta} \ ,
        \label{better_baddt_int_bnd1}\\
        \lp{\Delta^{-1}[\td{\td{ {\uoA^{(M)}} }} , \uL\td{\uoC}]}
        {L^2_\omega(\dot{B}_{2,10n}^{\infty,(2,\frac{n}{2}-\gamma)})}
        \ &\lesssim \ \mathcal{E}\cdot M^{-n - \delta} \ ,
        \label{better_baddt_int_bnd2}\\
        \lp{\Delta^{-1}\Big([\uL\td{\td{ {\uoA^{(M)}} }} , \td{\uoC}]+
        [\td{ {\uoA^{(M)}} } , \uL\td{\td{\uoC}}]\Big)}
        {L^2_\omega(\dot{B}_{2,10n}^{\infty,(2,\frac{n}{2}-\gamma)})}
        \ &\lesssim \ \mathcal{E}\cdot M^{-n - \delta} \ .
        \label{better_baddt_int_bnd3}
\end{align}
Our first step is to prove the estimates
\eqref{better_baddt_int_bnd1}--\eqref{better_baddt_int_bnd2}. To do
this, we introduce the auxiliary index:
\begin{equation}
        r_\gamma \ = \ \frac{2n(n-1)}{(n-2)(n+1) - 3\gamma n} \ .
        \label{r_gamma_def}
\end{equation}
We  now show that one has the following improvements over the
estimates \eqref{d_omega_calA_bnd2}, \eqref{d_omega_C_bnd2}:
\begin{align}
        \lp{\uL\td{ {\uoA^{(M)}} }}{ L_\omega^\infty(\dot{B}_{2,10n}^{r_\gamma,
        (2,\frac{n-4}{2})}) } \ &\lesssim \ \mathcal{E} \ ,
        \label{bL_tdA_imp}\\
        \lp{\uL\td{\uoC}}{ L^\infty_\omega(\dot{B}_{2,10n}^{r_\gamma,
        (2,\frac{n-4}{2})})} \ &\lesssim \ \mathcal{E} \ .
        \label{bL_tdC_imp}
\end{align}
With the help of \eqref{bL_tdA_imp}--\eqref{bL_tdC_imp}, the proof of
the estimates
\eqref{better_baddt_int_bnd1}--\eqref{better_baddt_int_bnd2} follows
from the
$L^2_\omega(\dot{B}_{2,10n}^{\infty,(2,\frac{n}{2}-\gamma)})$ 
bound shown
in the previous Lemma, and the following bilinear embedding:
\begin{align}
        \Delta^{-1} \ : \ L^2_\omega(\dot{B}_{2,10n}^{\infty,
        (2,\frac{n}{2}-\gamma)}) \cdot L_\omega^\infty(
        \dot{B}_{2,10n}^{r_\gamma, (2,\frac{n-4}{2})})
        \ &\hookrightarrow \
        L^2_\omega(\dot{B}_{2,10n}^{\infty,
        (2,\frac{n}{2}-\gamma)}) \ . \label{biline_better_index1}
\end{align}
Notice that the validity of this last estimate follows  from the
condition \eqref{gap_cond}, because for $0 < \gamma\ll 1$ we have
the index bounds:
\begin{equation}
        2 + \gamma \ < \  \frac{n}{r_\gamma}  \ , \notag
\end{equation}
which follows easily from the definitions of $r_\gamma$. To prove
\eqref{bL_tdA_imp}--\eqref{bL_tdC_imp} we proceed by first
showing \eqref{bL_tdA_imp}, and then  use the Hodge system
\eqref{lemma_td_uC} to show that
$\eqref{bL_tdA_imp}\Rightarrow\eqref{bL_tdC_imp}$. The details 
follow.\\

We are now trying to show \eqref{bL_tdA_imp}. We use the identity
\eqref{LLb_wave} and the definition \eqref{ucalA_shorthand} and the
structure equation \eqref{red_conctn_cond6} to compute that:
\begin{multline}
        \uL\td{ {\uoA^{(M)}} } \ = \ \nabla_{x}\ \ooPi_{M^{-1}
        < \bullet} \ooPi^{(\frac{1}{2} - \delta)}
        \   \uA(\partial_\omega) \\
        - \ \ \nabla_{x}\ \ooPi_{M^{-1}
        < \bullet} \ooPi^{(\frac{1}{2} - \delta)}
        \  \Delta_{\omega^\perp}^{-1}\, \td{\mathcal{P}}
        \big([B,H]\big)(\partial_\omega) \ . \label{bLA_identity}
\end{multline}
The estimate \eqref{bL_tdA_imp} for the first term on the right hand
side of this last expression is a trivial consequence the
$\dot{H}^\frac{n-4}{2}$ bound (at fixed time) for that term which is
provided through the energy type norm contained in the bootstrapping
assumption \eqref{red_conctn_cond5}, and the Besov nesting
\eqref{Besov_nesting}. Therefore, we strive to bound the second term
on the right hand side of \eqref{bLA_identity} above. To do this, we
first decompose things into a sum over all possible  angles spread
from the $\omega$ direction and write:
\begin{multline}
        \nabla_{x}\ \ooPi_{M^{-1}
        < \bullet} \ooPi^{(\frac{1}{2} - \delta)}
        \  \Delta_{\omega^\perp}^{-1}\,  \td{\mathcal{P}}
        \big([B,H]\big)(\partial_\omega) \\
        = \ \
        \sum_{\theta} \ \nabla_{x}\ \oPi_\theta \ \ooPi_{M^{-1}
        < \bullet} \ooPi^{(\frac{1}{2} - \delta)}
        \  \Delta_{\omega^\perp}^{-1}\,  \td{\mathcal{P}}
        \big([B,H]\big)(\partial_\omega) \ . \notag
\end{multline}
For each angularly localized piece in this last expression, we may
make use of the Coulomb savings \eqref{main_coulomb_savings} to show
the following heuristic multiplier bound (again making use of our
convention explained below \eqref{heuristic_op_bnd} above):
\begin{equation}
        \nabla_{x}\ \oPi_\theta \ \ooPi_{M^{-1}
        < \bullet} \ooPi^{(\frac{1}{2} - \delta)}
        \  \Delta_{\omega^\perp}^{-1}\, P_\mu
        \td{\mathcal{P}} \big([B,H]\big)(\partial_\omega)
        \ \approx \ (\mu\theta)^{-1} \
        \oPi_\theta \  P_\mu P_{\bullet \ll 1} \big([B,H]\big) \ . \notag
\end{equation}
Therefore, dropping the small frequency multiplier, our goal is to
show the following fixed angle estimate:
\begin{equation}
        \lp{
        \oPi_\theta |D_x|^{-1}\big([B,H]\big)}
        {  \dot{B}_{2}^{r_\gamma,
        (2,\frac{n-4}{2})}  } \
        \lesssim \ \theta^{1 + \gamma}\cdot \mathcal{E} \ . \notag
\end{equation}
Using Bernstein's inequality on each fixed dyadic block in the Besov
nesting \eqref{Besov_nesting}, and making use of a small numerical
calculation which we leave to the reader, one finds that this last
estimate is a consequence of the following non-localized Besov space
estimate:
\begin{equation}
        \lp{ |D_x|^{-1}\big([B,H]\big)}
        {  \dot{B}_{2}^{\frac{2n}{n+2-\gamma},
        (2,\frac{n-4}{2})}  } \
        \lesssim \  \mathcal{E} \ . \notag
\end{equation}
This last bound is now a consequence of the bootstrapping structure
estimate \eqref{red_conctn_cond7} and the following bilinear
embedding:
\begin{equation}
       |D_x|^{-1} \ : \  \dot{B}_2^{2,(2,\frac{n-2}{2})}\cdot
        \dot{B}_2^{2,(2,\frac{n-4}{2})} \ \hookrightarrow  \
        \dot{B}_2^{\frac{2n}{n+2-\gamma},(2,\frac{n-4}{2})}  \ . \notag
\end{equation}
Notice that the reason we are forced to work with the relatively
high space $L^\frac{2n}{n+2-\gamma}$ is because of $Low\times High$
frequency interactions. This is why we are forced to work in the
less aesthetic space $L^{r_\gamma}$ above instead of $L^2$. This
completes the proof of \eqref{bL_tdA_imp}.\\

Our next step is to establish the implication
$\eqref{bL_tdA_imp}\Rightarrow\eqref{bL_tdC_imp}$. This follows
immediately from differentiation of the Hodge system
\eqref{lemma_td_uC} with respect to the $\uL$ vector-field, and then
using the following bilinear estimate to bootstrap:
\begin{equation}
        \nabla_x\Delta^{-1} \ : \ \dot{B}_{2,10n}^{r_\gamma,(2,\frac{n-4}{2})}
        \cdot\dot{B}_{2,10n}^{p_\gamma,(2,\frac{n-2}{2})} \
        \hookrightarrow \ 
	\dot{B}_{2,10n}^{r_\gamma,(2,\frac{n-4}{2})} \ . \notag
\end{equation}
We leave it to the reader to check that the various conditions of
estimate \eqref{general_besov_embed} are satisfied in this case.\\

It remains for us to show the bound \eqref{better_baddt_int_bnd3}.
To make this a bit easier, we employ the skew symmetry of the Lie
brackets in that expression to write it as:
\begin{multline}
        \Delta^{-1}\Big([\uL\td{\td{ {\uoA^{(M)}} }} , \td{\uoC}]+
        [\td{ {\uoA^{(M)}} } , \uL\td{\td{\uoC}}]\Big) \\
        = \ \
        \Delta^{-1}\Big([\uL\td{\td{ {\uoA^{(M)}} }} , \td{\uoC} -
        \td{ {\uoA^{(M)}} }]+
        [\td{ {\uoA^{(M)}} } , \uL( \td{\td{\uoC}} - \td{\td{ {\uoA^{(M)}} }})]\Big)
        \ . \notag
\end{multline}
From this we see that the proof of \eqref{better_baddt_int_bnd3}
will follow once we can establish the three separate estimates:
\begin{align}
        \lp{\Delta^{-1}\Big([\uL\td{\td{ {\uoA^{(M)}} }} , \td{\uoC} -
        \td{ {\uoA^{(M)}} }]\Big)}{L^2_\omega(\dot{B}_{2,10n}^{\infty,
        (2,\frac{n}{2}-\gamma)})} \ \lesssim \ \mathcal{E}\cdot
        M^{-n-\delta} \ , \label{baddt_term_a}\\
        \lp{\Delta^{-1}\Big(
        [\td{ {\uoA^{(M)}} } , \uL \td{\td{\uoC}} ]
        \Big)}{L^2_\omega(\dot{B}_{2,10n}^{\infty,
        (2,\frac{n}{2}-\gamma)})} \ \lesssim \ \mathcal{E}\cdot
        M^{-n-\delta} \ , \label{baddt_term_b}\\
        \lp{\Delta^{-1}\Big(
        [\td{ {\uoA^{(M)}} } , \uL\td{\td{ {\uoA^{(M)}} }}]
        \Big)}{L^2_\omega(\dot{B}_{2,10n}^{\infty,
        (2,\frac{n}{2}-\gamma)})} \ \lesssim \ \mathcal{E}\cdot
        M^{-n-\delta} \ . \label{baddt_term_c}
\end{align}
To prove the first estimate \eqref{baddt_term_a} above, we make use
of the fact that $\td{\uoC} - \td{ {\uoA^{(M)}} }$ obeys a better
bound than either term in that expression does individually:
\begin{equation}
        \lp{\td{\uoC} - \td{ {\uoA^{(M)}} }}{\dot{B}_{2,10n}^{s_\gamma,(2,\frac{n-2}{2})}}
        \ \lesssim \ \mathcal{E} \ , \label{better_AC_diff_bound}
\end{equation}
where we have set the index $s_\gamma$ to be:
\begin{equation}
        s_\gamma \ = \ \frac{n p_\gamma }{ n+p_\gamma} + \gamma
        \ . \notag
\end{equation}
The proof of \eqref{better_AC_diff_bound} follows immediately from
the entirely quadratic structure of the terms in the expression
$\td{\uoC} - \td{ {\uoA^{(M)}} }$, in conjunction with following
bilinear estimate whose proof is a simple consequence of the
definition of the $p_\gamma$ indices, \eqref{p_gamma_line}, and the
general embedding \eqref{general_besov_embed}:
\begin{equation}
        \nabla_x\Delta^{-1} \ : \ \dot{B}_{2,10n}^{p_\gamma,(2,\frac{n-2}{2})}\cdot
        \dot{B}_{2,10n}^{p_\gamma,(2,\frac{n-2}{2})} \ \hookrightarrow  \
        \dot{B}_{2,10n}^{s_\gamma,(2,\frac{n-2}{2})} \ . \notag
\end{equation}
Furthermore, by taking the $\uL$ derivative of the potentials in the
estimate \eqref{tdtdA_int_bnd_fixed}, and making use of the
truncation condition \eqref{red_conctn_cond4}, we easily have the
following:
\begin{equation}
        \lp{\uL\td{\td{ {\uoA^{(M)}} }} }{\dot{B}_{2,10n}^{\infty,(2,\frac{n-2}{2}-\gamma)}}
        \ \lesssim \ \mathcal{E}\cdot M^{-n-\delta} \ .
        \label{LbA_inv_der_bnd}
\end{equation}
The proof of \eqref{baddt_term_a} now follows from combining
estimates \eqref{better_AC_diff_bound} and \eqref{LbA_inv_der_bnd}
in to the following bilinear embedding whose validity follows easily
from \eqref{general_besov_embed} and the condition $2 + \gamma <
\frac{n}{s_\gamma}$ (say for $n=6$ or higher):
\begin{equation}
        \Delta^{-1} \ : \ L^2_\omega(\dot{B}_{2,10n}^{\infty,
        (2,\frac{n-2}{2}-\gamma)}) \cdot L_\omega^\infty(
        \dot{B}_{2,10n}^{s_\gamma, (2,\frac{n-2}{2})})
        \ \hookrightarrow \
        L^2_\omega(\dot{B}_{2,10n}^{\infty,
        (2,\frac{n}{2}-\gamma)}) \ . \label{biline_better_index2}
\end{equation}\ret

We have now come to the point where the current techniques reach an
impasse. Notice that while the terms $\uL \td{\td{\uoC}}$ and
$\uL\td{\td{ {\uoA^{(M)}} }}$ do seem to have a better structure at
first glance via the equation \eqref{bLA_identity}, it is
surprisingly difficult to pass this into integrated estimates of the
form \eqref{tdtdA_int_bnd_fixed}. This is because while the linear
term on the right hand side of \eqref{bLA_identity} is quite nice,
the only saving grace of the quadratic term in that expression is
that it can go in lower spatial $L^p$ space, which is not
particularly useful when half of the needed savings in the estimate
\eqref{tdtdA_int_bnd_fixed} comes form orthogonality (meaning that
anything below $L^2$ gets wasted). A way to get rid of this problem
is to employ non-isotropic spaces. Specifically, we define the norm:
\begin{equation}
        \lp{A}{{}^\omega\!\mathcal{N}_{1,10n}^
        {-\frac{1}{2}-\gamma,2,\infty}} \ = \ \sum_\mu\
        \mu^{-\frac{1}{2} - \gamma}(1 + \mu)^{10n}\ \lp{P_\mu (A)}
        {L^2_{\omega^{||}}(L^\infty_{\omega^\perp})} \ . \notag
\end{equation}
Our goal is now to show the following  estimate which represents a
more manageable form of the differentiated version of
\eqref{better_tdtdA_int_bnd}:
\begin{equation}
        \lp{\uL\td{\td{ {\uoA^{(M)}} }}}{L^2_\omega({}^\omega\!\mathcal{N}_{1,10n}^
        {-\frac{1}{2}-\gamma,2,\infty})} \ \lesssim \
        \mathcal{E}\cdot M^{-n-\delta} \ . \label{non-iso_tdtdA}
\end{equation}
Having done this, our next goal will be to pass on estimates of this
form on to the non-linear potential $\td{\td{\uoC}}$. For reasons
which will become apparent in a moment, it is more convenient to
state this estimate for the following sum of spaces:
\begin{equation}
        \lp{  \uL \td{\td{\uoC}}
        }{L^2_\omega({}^\omega\!\mathcal{N}_{1,10n}^
        {-\frac{1}{2}-\gamma,2,\infty})
        + L^2_\omega(\dot{B}_{2,10n}^{n,(2,\frac{n-2}{2}-\gamma)})
        } \ \lesssim \
        \mathcal{E}\cdot M^{-n-\delta} \ . \label{non-iso_tdtdC}
\end{equation}
Once this is accomplished, the proof of
\eqref{baddt_term_b}--\eqref{baddt_term_b} will follow from the two
bilinear estimates:
\begin{align}
        \Delta^{-1} \ : \
        L^\infty_\omega\big(P_{\bullet\ll 1}
        \Delta_{\omega^\perp}^{-\frac{1}{2}} \dot{H}^\frac{n-4}{2}
        \big)
        \cdot L^2_\omega({}^\omega\!\mathcal{N}_{1,10n}^
        {-\frac{1}{2}-\gamma,2,\infty}) \ &\hookrightarrow \
        L^2_\omega(\dot{B}_{2,10n}^{\infty,
        (2,\frac{n}{2}-\gamma)}) \ ,
        \label{main_noniso_savings_est}\\
        \Delta^{-1} \ : \
        L^\infty_\omega(\dot{B}_{2,10n}^{p_\gamma,(2,\frac{n-2}{2})})\cdot
        L^2_\omega(\dot{B}_{2,10n}^{n,(2,\frac{n-4}{2})} ) \
        &\hookrightarrow \ L^2_\omega(\dot{B}_{2,10n}^{\infty,
        (2,\frac{n}{2}-\gamma)}) \ .
        \label{main_iso_easy_est}
\end{align}
Here the space in the first term in the product on the left hand
side of \eqref{main_noniso_savings_est} above is given by the norm:
\begin{equation}
        \lp{A}{\Delta_{\omega^\perp}^{-\frac{1}{2}} \dot{H}^\frac{n-4}{2}}
        \ = \ \lp{\Delta_{\omega^\perp}^{\frac{1}{2}} A
        }{\dot{H}^\frac{n-4}{2}} \ . \notag
\end{equation}
That the set of potentials $\{\td{ {\uoA^{(M)}} }\}$ is in this
space with norm $\lesssim \mathcal{E}$ follows from the explicit
formula \eqref{ucalA_shorthand} and the Coulomb gauge savings
\eqref{main_coulomb_savings}. Having now outlined the general
strategy, we move to the proofs of the individual estimates.\\

To prove \eqref{non-iso_tdtdA} we use the spatial truncation
condition \eqref{freq_loc_cond},  the triangle inequality,
and dyadic summing  to
reduce things to the following single frequency estimate:
\begin{equation}
        \left( \int_{\mathbb{S}^{n-1}} \
        \lp{P_\mu (\uL\td{\td{ {\uoA^{(M)}} }})}{L^2
        _{\omega^{||}}(L^\infty_{\omega^\perp})}^2\ d\omega
        \right)^\frac{1}{2} \ \lesssim \ \mu^{\frac{1}{2}+ 2\gamma} \
        \mathcal{E}\cdot M^{-n-\delta} \ .
        \label{frozen_N_int_bound}
\end{equation}
Now freeze $\omega$ and run a Littlewood-Paley decomposition in the
$\mathbb{R}_{\omega^\perp}^{n-1}$ frequency plane:
\begin{align}
        P_\mu (\uL\td{\td{ {\uoA^{(M)}} }}) \ &= \ \sum_{\substack{\lambda \ : \\
        \lambda \lesssim \ M^{-1}\mu}}\
        \nabla_x \uL \oL \Delta_{\omega^\perp}^{-1}
        \ooPi_{\bullet\leqslant M^{-1}}
        \ooPi^{(\frac{1}{2} - \delta)}
        Q_\lambda P_\mu\big(\uA\big)(\partial_\omega) \ , \notag\\
        &\approx \ M^{-1} \ \sum_{\substack{\lambda \ : \\
        \lambda \lesssim \ M^{-1}\mu}}\
        \mu^3 \lambda^{-2}\
        \ooPi_{\bullet\leqslant M^{-1}}
        \ooPi^{(\frac{1}{2} - \delta)}
        Q_\lambda P_\mu\big(\uA\big) \ , \label{bL_A_heur_mult}
\end{align}
where the last line follows from the truncation condition
\eqref{red_conctn_cond4} and our heuristic multiplier convention.
Notice that the sum restriction in these formulas comes because of
the presence of the cutoff $\ooPi_{\bullet\leqslant M^{-1}}$. The
extra $M^{-1}$ factor comes from this same angular cutoff and the
Coulomb gauge savings \eqref{main_coulomb_savings}. Working now with
the right hand side of \eqref{bL_A_heur_mult}, we use Bernstein's
inequality and dyadic summing to compute that:
\begin{align}
        &\lp{P_\mu (\uL\td{\td{ {\uoA^{(M)}} }})}{L^2
        _{\omega^{||}}(L^\infty_{\omega^\perp})} \ , \notag\\
        \lesssim \
        &M^{-1} \ \sum_{\substack{\lambda \ : \\
        \lambda \lesssim \ M^{-1}\mu}}\
        \mu^3 \lambda^{-2}\ \lp{Q_\lambda
        \ooPi_{\bullet\leqslant M^{-1}}
        \ooPi^{(\frac{1}{2} - \delta)}
         P_\mu\big(\uA\big) }{L^2
        _{\omega^{||}}(L^\infty_{\omega^\perp})} \ , \notag\\
        \lesssim \ &M^{-\frac{n-3}{2}} \mu^\frac{3}{2} \
        \lp{\ooPi_{\bullet\leqslant M^{-1}}
        \ooPi^{(\frac{1}{2} - \delta)}\ P_\mu\big( \uA \big) }{\dot{H}_x^\frac{n-2}{2}}
        \ , \notag\\
        \lesssim \ &M^{-\frac{n+1}{2}-\delta } \mu^{\frac{1}{2}+2\gamma} \
        \lp{\ooPi_{\bullet\leqslant M^{-1}}
        \ooPi^{(\frac{1}{2} - \delta)}\  \uA  }{\dot{H}_x^\frac{n-2}{2}}
        \ . \notag
\end{align}
Integrating now this last line in $L^2_\omega$, and using the
orthogonality computation which began on line
\eqref{begin_orth_comp} above we have achieved
\eqref{frozen_N_int_bound} as was to be shown.\\

Our goal is now to pass the estimate \eqref{non-iso_tdtdA} on to
potentials $\{\uL \td{\td{\uoC}}\}$ modulo terms which are in the
more regular space
$L^2_\omega(\dot{B}_{2,10n}^{n,(2,\frac{n-2}{2}-\gamma)})$. To do
this, we differentiate the system \eqref{lemma_td_td_uDpm_system}
with respect to the vector-field $\uL$, and write it heuristically
as:
\begin{multline}
        \uL \td{\td{\uoC}} \ = \  \uL \td{\td{\uoA^{(M)}}} \
        +  \ \nabla_x \Delta^{-1} \Big([\uL \td{\td{\uoA^{(M)}}}
        , \td{\uoC} ] + [ \td{\td{\uoA^{(M)}}}
        ,\uL \td{\uoC} ]\\
         + \ [\uL {\td{\uoA^{(M)}}}
        , \td{\td{\uoC}} ] + [ {\td{\uoA^{(M)}}}
        ,\uL \td{\td{\uoC}} ] \ + \
         [ \uL \td{\uoC} , \td{\td{\uoC}}] \ + \
        [ \td{\uoC} ,\uL \td{\td{\uoC}}] \Big) \ . \notag
\end{multline}
Therefore, the desired bound will follow from a bootstrapping
argument the estimates \eqref{non-iso_tdtdA},
\eqref{Linfty_omega_A_bnd}--\eqref{Linfty_omega_uC_bnd}, and
\eqref{bL_tdA_imp}--\eqref{bL_tdC_imp} with the help of the
following three bilinear estimates:
\begin{align}
        \nabla_x\Delta^{-1} \ : \
     L^2_\omega(\dot{B}_{2,10n}^{\infty,(2,\frac{n}{2}-\gamma)})\cdot
        L^\infty_\omega(\dot{B}_{2,10n}^{r_\gamma,
    (2,\frac{n-4}{2})}) \ &\hookrightarrow  \
        L^2_\omega(\dot{B}_{2,10n}^{n,(2,\frac{n-2}{2}-\gamma)}) \ ,
        \label{bLtdtdC_good_est1}\\
        \nabla_x\Delta^{-1} \ : \
     L^2_\omega(\dot{B}_{2,10n}^{n,(2,\frac{n-2}{2}-\gamma)})\cdot
        L^\infty_\omega (\dot{B}_{2,10n}^{p_\gamma,(2,\frac{n-2}{2})}) \
     &\hookrightarrow  \
        L^2_\omega( \dot{B}_{2,10n}^{n,(2,\frac{n-2}{2}-\gamma)}) \ ,
        \label{bLtdtdC_good_est2}\\
        \nabla_x\Delta^{-1} \ : \ L^2_\omega({}^\omega\!\mathcal{N}_{1,10n}^
        {-\frac{1}{2}-\gamma,2,\infty})\cdot L^\infty_\omega(
        \dot{B}_{2,10n}^{p_\gamma,(2,\frac{n-2}{2})})
        \ &\hookrightarrow  \  L^2_\omega({}^\omega\!\mathcal{N}_{1,10n}^
        {-\frac{1}{2}-\gamma,2,\infty}) \ . \label{bLtdtdC_bad_est}
\end{align}
The estimates \eqref{bLtdtdC_good_est1}--\eqref{bLtdtdC_good_est2}
are again an integrated form of the general Besov embedding
\eqref{general_besov_embed}, and we leave it to the reader to check
that the indices $p_\gamma,r_\gamma$ are in the right range to
satisfy the conditions \eqref{gap_cond}--\eqref{Lb_cond}. It remains
for us to prove the inclusion \eqref{bLtdtdC_bad_est}. We do this
now. Let $A$ and $C$ be two test matrices. By performing a
trichotomy, we see that it suffices to prove the following three
frequency localized summation estimates for fixed values of
$\omega$:
\begin{align}
    \begin{split}
        \sum_{\substack{\lambda , \mu_i \ :\\
        \mu_1 \ll \mu_2 \\
        \lambda\sim \mu_2 }}\
        \lambda^{-\frac{1}{2} -\gamma}(1+\lambda)^{10n} \
        \lp{\nabla_x\Delta^{-1} P_\lambda\left(P_{\mu_1}A\cdot P_{\mu_2}C
        \right)}
        {L^2_{\omega^{||}}(L^\infty_{\omega^\perp})} \ &\lesssim \\
        &\hspace{-1.5in}
        \lp{A}{{}^\omega\! \mathcal{N}_{1,10n}^{-\frac{1}{2} -\gamma,2,\infty}}\cdot
        \lp{C}{\dot{B}_{2,10n}^{p_\gamma ,
        (2,\frac{n-2}{2})}}
        \ ,
    \end{split} \label{second_tric_N1}\\
    \begin{split}
        \sum_{\substack{\lambda , \mu_i \ :\\
        \mu_2 \ll \mu_1 \\
        \lambda\sim \mu_1 }}\
        \lambda^{-\frac{1}{2}-\gamma}(1+\lambda)^{10n} \
        \lp{\nabla_x\Delta^{-1} P_\lambda\left(P_{\mu_1}A\cdot P_{\mu_2}C
        \right)}
        {L^2_{\omega^{||}}(L^\infty_{\omega^\perp})} \ &\lesssim \\
        &\hspace{-1.5in} \lp{A}{{}^\omega\!
        \mathcal{N}_{1,10n}^{-\frac{1}{2} -\gamma,2,\infty}}\cdot
        \lp{C}{\dot{B}_{2,10n}^{p_\gamma ,
        (2,\frac{n-2}{2})}}
        \ ,
    \end{split} \label{second_tric_N2}\\
    \begin{split}
        \sum_{\substack{\lambda , \mu_i \ :\\
        \mu_1 \sim \mu_2 \\
        \lambda\lesssim \mu_1,\mu_2 }}\
        \lambda^{-\frac{1}{2}-\gamma}(1+\lambda)^{10n} \
        \lp{\nabla_x\Delta^{-1} P_\lambda\left(P_{\mu_1}A\cdot P_{\mu_2}C
        \right)}
        {L^2_{\omega^{||}}(L^\infty_{\omega^\perp})} \ &\lesssim \\
        &\hspace{-1.5in} \lp{A}{{}^\omega\! \mathcal{N}_{1,10n}^
        {-\frac{1}{2} -\gamma,2,\infty}}\cdot
        \lp{C}{\dot{B}_{2,10n}^{p_\gamma ,
        (2,\frac{n-2}{2})}}
        \ .
    \end{split} \label{second_tric_N3}
\end{align}
The proof of \eqref{second_tric_N1}--\eqref{second_tric_N3} is
essentially identical to the proof of the three estimates
\eqref{tric_N1}--\eqref{tric_N3} we have shown earlier, although the
proof of the last estimate \eqref{second_tric_N3} requires a
slightly more delicate argument due to the presence of additional
low frequency weights. We leave
\eqref{second_tric_N1}--\eqref{second_tric_N2} to the reader. To
show the last estimate above, we follow the proof of \eqref{tric_N3}
which begins on line \eqref{tric3_N_start}, although we do so
without throwing away the $P_\lambda$ multiplier so soon. This
leaves us with the fixed frequency estimate, which we expand out
into all frequencies in the $\omega^{||}$ variable, calling the
appropriate multipliers $\td{Q}_\sigma$:
\begin{align}
        &\lambda^{-\frac{1}{2}-\gamma} \
        \lp{\nabla_x\Delta^{-1} P_\lambda\left(P_{\mu_1}A\cdot P_{\mu_2}C
        \right)}
        {L^2_{\omega^{||}}(L^\infty_{\omega^\perp})} \ , \notag\\
        \lesssim \
        &\lambda^{\frac{n-1}{p_\gamma}-\frac{3}{2}-\gamma}\
        \lp{P_\lambda( P_{\mu_1}A\cdot P_{\mu_2}C)}
        {L^2_{\omega^{||}}(L^{p_\gamma}_{\omega^\perp})}\ , \notag\\
        \lesssim \ &\lambda^{\frac{n-1}{p_\gamma}-\frac{3}{2}-\gamma}\
        \sum_{\sigma \lesssim \lambda}\
    \lp{\td{Q}_\sigma( P_{\mu_1}A\cdot P_{\mu_2}C)}
        {L^2_{\omega^{||}}(L^{p_\gamma}_{\omega^\perp})} \ .
        \label{second_N_tric_end}
\end{align}
Now we use the fact that the multiplier $\td{Q}_\sigma$ only acts in
the $\omega^{||}$ variable. In that variable its action can be
written in terms of a kernel $K^{\td{Q}_\sigma}$ which 
has uniformly bounded
$L^1_{\omega^{||}}$ norm (in terms of the parameter
$\sigma$) and has amplitude
$\sim \sigma$. Therefore, via Young's and then H\"older's
inequality, and a little dyadic summing, this allows us to bound:
\begin{align}
        \hbox{(L.H.S.)}\eqref{second_N_tric_end} \ &\lesssim \
        \lambda^{\frac{n-1}{p_\gamma}-\frac{3}{2}-\gamma}\
        \sum_{\sigma \lesssim \lambda}\ \lp{\, \Big(|K^{\td{Q}_\sigma}|*
        \lp{( P_{\mu_1}A\cdot P_{\mu_2}C) }{L^{p_\gamma}_{\omega^\perp}
        }\Big) }{L^2_{\omega^{||}}} \ , \notag\\
        &\lesssim \
        \lambda^{\frac{n-1}{p_\gamma}-\frac{3}{2}-\gamma}\
        \sum_{\sigma \lesssim \lambda}\ \sigma^{\frac{1}{p_\gamma}}\ \lp{
         P_{\mu_1}A\cdot P_{\mu_2}C
        }{L^\frac{2 p_\gamma}{2 + p\gamma}
        _{\omega^{||}} (L^{p_\gamma}_{\omega^\perp})} \ , \notag\\
        &\lesssim \
        \lambda^{\frac{n}{p_\gamma}-\frac{3}{2}-\gamma}\
        \lp{P_{\mu_1}A}{L^2_{\omega^{||}}(L^\infty_{\omega^\perp})}
        \cdot\lp{P_{\mu_2}C}{L^{p_\gamma}} \ , \notag\\
        &\lesssim \ \left(\frac{\lambda}{\mu_1}\right)
        ^{\frac{n}{p_\gamma}-\frac{3}{2}-\gamma}\
        \mu_1^{-\frac{1}{2} - \gamma}
        \lp{P_{\mu_1}A}{L^2_{\omega^{||}}(L^\infty_{\omega^\perp})} \cdot
        \mu_2^{\frac{n}{p_\gamma} - 1} \lp{P_{\mu_2}C}{L^{p_\gamma}}
        \ . \notag
\end{align}
This last line and the condition $0 <
\frac{n}{p_\gamma}-\frac{3}{2}-\gamma$ allows us to safely make the
sum on the left hand side of  \eqref{second_tric_N3} and then
proceed via Cauchy-Schwartz to arrive at the desired bound. This
completes our proof of the bilinear estimate
\eqref{bLtdtdC_bad_est}.\\

The last thing we need to do here is to prove the two final
estimates \eqref{main_noniso_savings_est} and
\eqref{main_iso_easy_est}. The second of these is of course simply
an integrated version of the general estimate
\eqref{general_besov_embed}. Therefore we concentrate on proving the
first. To do this, we proceed as we did in the proof of estimate
\eqref{bLtdtdC_bad_est} and run a trichotomy on a product of test
matrices $A\cdot C$. This leaves us with establishing the three
estimates (forgetting about the extra high frequency weights which
are not central):
\begin{align}
        \sum_{\substack{\lambda , \mu_i \ :\\
        \mu_1 \ll \mu_2 \\
        \lambda\sim \mu_2 }}\
        \lambda^{ -\gamma} \
        \lp{\Delta^{-1} P_\lambda\left(P_{\mu_1}A\cdot P_{\mu_2}C
        \right)}
        {L^\infty} \ &\lesssim \
        \lp{A}{\Delta_{\omega^\perp}^{-\frac{1}{2}} \dot{H}^\frac{n-4}{2} }\cdot
        \lp{C}{ {}^\omega\!
        \mathcal{N}_{1,10n}^{-\frac{1}{2} -\gamma,2,\infty}}
        \ ,  \label{last_tric_N1}\\
        \sum_{\substack{\lambda , \mu_i \ :\\
        \mu_2 \ll \mu_1 \\
        \lambda\sim \mu_1 }}\
        \lambda^{-\gamma}\
        \lp{\Delta^{-1} P_\lambda\left(P_{\mu_1}A\cdot P_{\mu_2}C
        \right)}
        {L^\infty} \ &\lesssim \
        \lp{A}{\Delta_{\omega^\perp}^{-\frac{1}{2}} \dot{H}^\frac{n-4}{2} }\cdot
        \lp{C}{ {}^\omega\!
        \mathcal{N}_{1,10n}^{-\frac{1}{2} -\gamma,2,\infty}}
        \ ,  \label{last_tric_N2}\\
        \sum_{\substack{\lambda , \mu_i \ :\\
        \mu_1 \sim \mu_2 \\
        \lambda\lesssim \mu_1,\mu_2 }}\
        \lambda^{-\gamma} \
        \lp{\Delta^{-1} P_\lambda\left(P_{\mu_1}A\cdot P_{\mu_2}C
        \right)}
        {L^\infty} \ &\lesssim \
        \lp{A}{\Delta_{\omega^\perp}^{-\frac{1}{2}} \dot{H}^\frac{n-4}{2} }\cdot
        \lp{C}{ {}^\omega\!
        \mathcal{N}_{1,10n}^{-\frac{1}{2} -\gamma,2,\infty}}
        \ . \label{last_tric_N3}
\end{align}
The proofs of the two $Low\times High$ interaction estimates,
\eqref{last_tric_N1}--\eqref{last_tric_N2}, are both similar and
very simple. They follow from the pair of $L^\infty$ estimates:
\begin{align}
        \lp{P_{\mu_1}( A)}{L^\infty} \ &\lesssim \
        \mu_1 \ \lp{P_{\mu_1}( A)}
        {\Delta_{\omega^\perp}^{-\frac{1}{2}} \dot{H}^\frac{n-4}{2}}
        \ , \label{non_iso_Linfty_A}\\
        \lp{P_{\mu_2}( C)}{L^\infty} \ &\lesssim \
        \mu_2^{1+\gamma} \ \lp{P_{\mu_2}(C)}{{}^\omega\!
        \mathcal{N}_{1,10n}^{-\frac{1}{2} -\gamma,2,\infty}} \ .
        \label{non_iso_Linfty_C}
\end{align}
The proof of \eqref{non_iso_Linfty_A} follows easily from the kind
of angular decomposition and Bernstein inequality tricks used to prove
estimate \eqref{Besov_degen_tdA_Ln_smallness} above. To prove the
second estimate \eqref{non_iso_Linfty_C}, we let $\td{Q}_\sigma$
again denote a family of frequency cutoffs in the
$\mathbb{R}_{\omega^{||}}$ variable and we compute via Bernstein:
\begin{align}
        \lp{P_{\mu_2}( C)}{L^\infty} \ &\lesssim \ \sum_{\sigma
        \lesssim \mu_1} \ \lp{\td{Q}_\sigma P_{\mu_2}( C)}
        {L_{\omega^\perp}^\infty(L^\infty_{\omega^{||}})} \ ,
        \notag\\
        &\lesssim \ \sum_{\sigma
        \lesssim \mu_1} \ \sigma^\frac{1}{2} \ \lp{ P_{\mu_2}( C)}
        {L_{\omega^\perp}^\infty(L^2_{\omega^{||}})} \ , \notag\\
        &\lesssim \ \mu_2^{1+\gamma} \ \lp{P_{\mu_2}(C)}{{}^\omega\!
        \mathcal{N}_{1,10n}^{-\frac{1}{2} -\gamma,2,\infty}} \ .
        \notag
\end{align}
Using now \eqref{non_iso_Linfty_A}--\eqref{non_iso_Linfty_C} we have
the pair of fixed frequency bounds:
\begin{align}
        &\lambda^{ -\gamma} \
        \lp{\Delta^{-1} P_\lambda\left(P_{\mu_1}A\cdot P_{\mu_2}C
        \right)}{L^\infty} &\notag \\
        \lesssim \
        &\left(\frac{\mu_1}{\mu_2}\right)\ \lp{P_{\mu_1}(A)}{\Delta_{\omega^\perp}
        ^{-\frac{1}{2}} \dot{H}^\frac{n-4}{2} }\cdot
        \lp{P_{\mu_2}(C)}{ {}^\omega\!
        \mathcal{N}_{1,10n}^{-\frac{1}{2} -\gamma,2,\infty}}
        &\mu_1 \ &\ll \ \mu_2 \ , \notag
\end{align}
and:
\begin{align}
        &\lambda^{ -\gamma} \
        \lp{\Delta^{-1} P_\lambda\left(P_{\mu_1}A\cdot P_{\mu_2}C
        \right)}{L^\infty} &\notag \\
        \lesssim \
        &\left(\frac{\mu_2}{\mu_1}\right)^{1+\gamma} \
        \lp{P_{\mu_1}(A)}{\Delta_{\omega^\perp}
        ^{-\frac{1}{2}} \dot{H}^\frac{n-4}{2} }\cdot
        \lp{P_{\mu_2}(C)}{ {}^\omega\!
        \mathcal{N}_{1,10n}^{-\frac{1}{2} -\gamma,2,\infty}}
        &\mu_2 \ &\ll \ \mu_1 \ . \notag
\end{align}
These may easily be summed over the respective ranges on the left
hand side of \eqref{last_tric_N1}--\eqref{last_tric_N2} to yield the
desired bounds.\\

Our final task here is to establish the $High \times High$ frequency
interaction estimate \eqref{last_tric_N3}. This is where the
non-isotropic spaces really shine. In what follows, we let
$Q_{\sigma_1}$ denote a frequency cutoff in the
$\mathbb{R}_{\omega^\perp}^{n-1}$ frequency plane, and
$\td{Q}_{\sigma_2}$ a cutoff in the orthogonal direction. We compute
that:
\begin{align}
        &\lambda^{-\gamma} \
        \lp{\Delta^{-1} P_\lambda\left(P_{\mu_1}A\cdot P_{\mu_2}C
        \right)}
        {L^\infty} \ , \notag\\
        \lesssim \ &\lambda^{-2-\gamma}\ \sum_{\substack{\sigma_i \
        :\\ \sigma_i \lesssim \lambda}}\
        \lp{\td{Q}_{\sigma_2} Q_{\sigma_1}\left(P_{\mu_1}A\cdot P_{\mu_2}C
        \right)}{L^\infty_{\omega^\perp}(L^\infty_{\omega^{||}})  }\
        , \notag\\
        \lesssim \ &\lambda^{-1-\gamma}\ \sum_{\substack{\sigma_1 \
        :\\ \sigma_1 \lesssim \lambda}}\
        \lp{ Q_{\sigma_1}\left(P_{\mu_1}A\cdot P_{\mu_2}C
        \right)}{L^\infty_{\omega^\perp}(L^1_{\omega^{||}})  }\
        , \notag\\
        \lesssim \ &\lambda^{-1-\gamma}\ \sum_{\substack{\sigma_1 \
        :\\ \sigma_1 \lesssim \lambda}}\
        \lp{ Q_{\sigma_1}\left(P_{\mu_1}A\cdot P_{\mu_2}C
        \right)}{L^1_{\omega^{||}}(L^\infty_{\omega^\perp})  }\
        , \notag\\
        \lesssim \ &\lambda^{-1-\gamma}\ \sum_{\substack{\sigma_1 \
        :\\ \sigma_1 \lesssim \lambda}}\ \sigma_1^\frac{n-3}{2}
        \lp{ Q_{\sigma_1}\left(P_{\mu_1}A\cdot P_{\mu_2}C
        \right)}{L^1_{\omega^{||}}(L^\frac{2(n-1)}{n-3}_{\omega^\perp})  }\
        , \notag\\
        \lesssim \ &\lambda^{\frac{n-5}{2} -\gamma}\
        \lp{P_{\mu_1}(A)}{L^2_{\omega^{||}}(L^\frac{2(n-1)}{n-3}_{\omega^\perp})}
        \cdot\lp{P_{\mu_2}(C)}
        {L^2_{\omega^{||}}(L^\infty_{\omega^\perp}) } \ , \notag\\
        \lesssim \ &\left(\frac{\lambda}{\mu_1}\right)^{\frac{n-5}{2} -\gamma}
        \ \lp{\Delta^\frac{1}{2}_{\omega^\perp}
        P_{\mu_1}(A)}{\dot{H}_x^\frac{n-4}{2}}\cdot
        \lp{P_{\mu_2}(C)}{ {}^\omega\!
        \mathcal{N}_{1,10n}^{-\frac{1}{2} -\gamma,2,\infty}} \ .
        \notag
\end{align}
Notice that the last line above follows from the $\dot{H}^1$ Sobolev
embedding in the $\mathbb{R}^{n-1}_{\omega^\perp}$ plane. This
estimate can now be safely summed using the condition that
$6\leqslant n$ to sum over the lower
dyadics, and then using Cauchy-Schwartz to
sum over the frequency localized pieces. This completes our proof of
the bilinear estimate \eqref{bLtdtdC_bad_est}, and hence our proof
of the integrated bound \eqref{tdtdC0_int_bnd}.
\end{proof}\ret

\noindent Having now established the proof of both the integrated
bounds \eqref{tdtdC_int_bnd}--\eqref{tdtdC0_int_bnd}, we have proved
the integrated group element bounds
\eqref{h_int_bnd1}--\eqref{h_int_bnd2}. This ends our proof of
Proposition \ref{h_int_rem_est}.
\end{proof}\ret

\subsection{Proof of the Accuracy estimate \eqref{half_wave_est3}}
We will now give a short proof of the multiplier equivalence bound
\eqref{half_wave_est3}. This will follow almost directly from the
estimates we have already shown. We compute the kernel of the
operator $\Phi(0)\big((2\pi  |\xi|)^\alpha (\Phi(0))^*\big)
 - (-\Delta)^\frac{\alpha}{2} P_1$ to be (again suppressing $\pm$ notations):
\begin{multline}
        K^{\alpha}(x,y) \ = \
        \int_{\mathbb{R}^n} \
        e^{2\pi i ( x-y)\cdot \xi}\
        \og^{-1}(x) \og(y) \big[\bullet \big]\
        \og^{-1}(y) \og(x) \ \chi^\alpha(\xi)\ d\xi \\
    - \ \ \int_{\mathbb{R}^n} \
        e^{2\pi i ( x-y)\cdot \xi}\
         \big[\bullet \big]\
         \ \chi^\alpha(\xi)\ d\xi \ , \label{Kalpha_kernel_sum}
\end{multline}
where $\chi^\alpha(\xi) = (2\pi|\xi|)^\alpha
\chi_{(-\frac{1}{2},2)}(\xi)$. Notice that this cutoff function
satisfies the general requirements of the generic bump function
${\chi}$ used throughout this section. In particular, there exist
constants $C_k$ which depend only on $\alpha$ and the original
$\chi_{(-\frac{1}{2},2)}$ such that:
\begin{equation}
    \int_{\mathbb{R}^n}\ |\nabla_\xi^k \chi^\alpha(\xi) |
    \ d\xi \ \leqslant\ C_k \ . \label{D_chi_alpha_bound}
\end{equation}
We now decompose the kernel $K^{\alpha} = \sum_\sigma K^{\alpha}_\sigma$
according to the dyadic physical space decomposition
\eqref{dyadic_physical_scales}. For each fixed value of the small
constant $\mathcal{E}$ on line \eqref{half_wave_est3} we write this
sum in terms of two pieces, a ``close'' part and a ``far'' part:
\begin{align}
    K^{\alpha} \ &= \ K^\alpha_{\bullet \leqslant
    \mathcal{E}^{-\frac{1}{2(n+1)}}}
    + K^\alpha_{ \mathcal{E}^{-\frac{1}{2(n+1)}} < \bullet}
    \ , \label{near_far_decomp}\\
    &= \ \sum_{\substack{\sigma\ :\\ \sigma
    \leqslant \mathcal{E}^{-\frac{1}{2(n+1)}}}}\ K^{\alpha}_\sigma
    \ \ + \ \
    \sum_{\substack{\sigma\ :\\
    \mathcal{E}^{-\frac{1}{2(n+1)}}
    < \sigma}}\ K^{\alpha}_\sigma \ . \notag
\end{align}
To estimate the near portion of things, we do a little algebraic
manipulation and write the kernel as:
\begin{multline}
    K^\alpha_{\bullet \leqslant
    \mathcal{E}^{-\frac{1}{2(n+1)}}} \ = \
    \chi_{\mathcal{D}_{ \bullet \leqslant
    \mathcal{E}^{-\frac{1}{2(n+1)} }}}\ \Big(
    \int_{\RR}\ e^{2\pi i ( x-y)\cdot \xi}\
        \big(\og^{-1}(x) \og(y)-I\big) \big[\bullet \big]\
        \og^{-1}(y) \og(x) \ \chi^\alpha(\xi)\ d\xi \\
    + \ \int_{\RR}\ e^{2\pi i ( x-y)\cdot \xi}\
         \big[\bullet \big]\
        \big( \og^{-1}(y) \og(x) -I\big)
    \ \chi^\alpha(\xi)\ d\xi \Big) \ . \notag
\end{multline}
By a direct application of the pair of integrated bounds
\eqref{h_int_bnd1}--\eqref{h_int_bnd2} (with $M\sim 1$)
this last expression gives us the absolute kernel bound:
\begin{equation}
    | K^\alpha_{\bullet \leqslant
    \mathcal{E}^{-\frac{1}{2(n+1)}}}(x,y) | \
    \lesssim \ \mathcal{E}\cdot (1 + |x-y|)
    \chi_{\mathcal{D}_{ \bullet \leqslant
    \mathcal{E}^{-\frac{1}{2(n+1)} }}}(|x-y|) \ . \notag
\end{equation}
By integrating the right hand side of this last inequality
we easily arrive at the pair of Schur-test bounds:
\begin{equation}
    \lp{K^\alpha_{\bullet \leqslant
    \mathcal{E}^{-\frac{1}{2(n+1)}}}}{L^\infty_y(L^1_x)} \  , \
    \lp{K^\alpha_{\bullet \leqslant
    \mathcal{E}^{-\frac{1}{2(n+1)}}}}{L^\infty_x(L^1_y)}
     \ \ \lesssim \ \ \mathcal{E}^\frac{1}{2} \ .
     \label{small_schur_test_bounds}
\end{equation}
To estimate the second kernel on the right hand side of
\eqref{near_far_decomp}, we do things separately for each term in
the sum \eqref{Kalpha_kernel_sum}. For the second term, which does
not contain the group elements, a simple application of the estimate
\eqref{D_chi_alpha_bound} and integration by parts shows that one
has the absolute bounds:
\begin{align}
    &\big|\chi_{\mathcal{D}_{
    \mathcal{E}^{-\frac{1}{2(n+1)} }
    < \bullet }}(|x-y|)
    \int_{\mathbb{R}^n}\ e^{2\pi i ( x-y)\cdot \xi}\
         \big[\bullet \big]\
         \ \chi^\alpha(\xi)
    \ d\xi\big| \ , \label{far_int_bounds_begin}\\
    \lesssim  \
    &\chi_{\mathcal{D}_{
    \mathcal{E}^{-\frac{1}{2(n+1)} }
    < \bullet }}(|x-y|)\cdot(1 + |x-y|)^{-2(n+1)} \ , \notag\\
    \lesssim \ &\mathcal{E}^\frac{1}{2}\cdot
    (1 + |x-y|)^{-(n+1)} \ . \notag
\end{align}
This easily yields Schur-test bounds of the form
\eqref{small_schur_test_bounds}. Therefore, it remains to prove
these bounds for the first integral expression on the right hand
side \eqref{Kalpha_kernel_sum} after it has been cut off in the far
region $\mathcal{E}^{-\frac{1}{2(n+1)}} < |x-y|$. This follows at
once from writing this kernel as a sum over various dyadic regions,
and using the symbol bounds
\eqref{L2_symbol_bds1}--\eqref{L2_symbol_bds2} as well as the
reduction to the integrated estimates
\eqref{h_int_bnd1}--\eqref{h_int_bnd2}. The key thing to notice 
is that there are only two places where we do not pick up the  factor
of $\mathcal{E}$ in the resulting estimates. The first is in the main
integration by parts argument when the derivatives $\nabla^k_\xi$
\emph{all} fall on the cutoff function $\chi^\alpha$. In that case
we can simply use the compactness of the group elements and proceed
in a way that is analogous to the computation which started on line
\eqref{far_int_bounds_begin} above. The second place is where we 
estimate the integral \eqref{small_L2_estA}. 
In that case we can easily upgrade 
the bound \eqref{diff_chi_bnd} to have the factor $\sigma^{-2(n+1)}$ on the
right hand side. We are then essentially in the same
situation as was reached starting on line 
\eqref{far_int_bounds_begin} above. This completes our proof of the
general multiplier approximation estimate \eqref{half_wave_est3}.

\ret

\section{The Dispersive Estimate}

In this section, we complete our proof of the non-microlocalized
version of the Strichartz estimates contained in \eqref{half_wave_est1}.
Using the abstract
machinery of \cite{KT_str}, these will follow once we can show that
the parametrix \eqref{parametrix} satisfies a dispersive estimate.
If at fixed time $t$ we write that operator as:
\begin{equation}
        T(t)(\widehat{f}) \ = \ \Phi(t)(\widehat{f}) \ , \notag
\end{equation}
where we have suppressed the $\pm$ notation, then we seek to prove
the bound (where $f$ has nothing to do with the original $\widehat{f}$,
but just represents a function of the physical space variables):
\begin{equation}
        \lp{T(t)T^*(s) f}{L^\infty_x} \ \lesssim \ (1 +
        |t-s|)^{-\frac{n-1}{2}} \ \lp{f}{L^1_x} \ . \label{disp_est}
\end{equation}
Now, a calculation similar that used to produce
\eqref{L2_TT*_kernel} shows that the kernel of the above operator
can be computed to be:
\begin{multline}
        K^{TT^*}(t,s;x,y) \ = \\
        \int_{\mathbb{R}^n} \
        e^{2\pi i\big( (t-s)|\xi| + ( x-y)\cdot \xi\big)}\
        \og^{-1}(t,x) \og(s,y) \big[\bullet \big]\
        \og^{-1}(s,y) \og(t,x) \ \chi(\xi)\ d\xi \ . \label{L2_TT*_disp_kernel}
\end{multline}
Therefore, as is usually the case, we see that it suffices to  show
the fixed time uniform bound:
\begin{equation}
        \lp{K^{TT^*}(t,s;\cdot,\cdot)}{L^\infty_{x,y}} \ \lesssim \
        (1 + |t-s|)^{-\frac{n-1}{2}} \ . \label{kernal_disp_est}
\end{equation}
The proof of \eqref{kernal_disp_est} turns out to be a
straightforward consequence of the bounds established in the
previous section. The strategy we follow here is almost identical.
We first decompose the $K^{TT*}$ kernel into a sum of two pieces:
\begin{equation}
        K_\sigma^{TT^*} \ = \ \td{K}^{TT^*} \ + \
        \mathcal{R}^{TT^*} \ , \notag
\end{equation}
for which we'll show the bound \eqref{kernal_disp_est} individually.
The $\td{K}^{TT^*}$ kernel will be smooth enough that we can use a
standard stationary phase computation on it. The remainder kernel
$\mathcal{R}^{TT^*}$ will be small in absolute value without using
any sophisticated integration by parts (although, as in the previous
section, there will be some use for oscillations in this term also).
As in the previous section, the definition of $\td{K}^{TT^*}$ will
depend on a physical space scale, in this case the value of $(1 +
|t-s| + |x-y|)$. This will again be effected by the choice of an
auxiliary gauge transformation $\td{\og}$. This time we define
$\td{\og}$ to be the transformation into the Coulomb gauge of the
smoothed out potential:
\begin{equation}
        \td{\uoA^{(M)}} \ = \ -\,  \ooPi_{M^{-1}
        < \bullet }\ooPi^{(\frac{1}{2} - \delta)}
        \nabla_x \oL\
        \Delta_{\omega^\perp}^{-1}\,  \uA(\partial_\omega) \ ,
        \label{approx_null_generator_disp}
\end{equation}
where we define the scale $M$ to be such that:
\begin{equation}
        M \ = \ (1 + |t-s| + |x-y| )^\frac{1}{2} \ . \notag
\end{equation}
As before, we use the splitting \eqref{gh_left}--\eqref{gh_right} to
compute:
\begin{multline}
        \td{K}^{TT^*}(t,s;x,y) \ = \\
        \int_{\mathbb{R}^n} \
        e^{2\pi i\big( (t-s)|\xi| + ( x-y)\cdot \xi\big)}\
        \td{\og}^{-1}(t,x) \td{\og}(s,y) \big[\bullet \big]\
        \td{\og}^{-1}(s,y) \td{\og}(t,x) \ \chi(\xi)\ d\xi \ . 
	\label{L2_TT*_disp_tdkernel}
\end{multline}
Our first step here is to notice that it suffices to show
\eqref{kernal_disp_est} for the kernel \eqref{L2_TT*_disp_tdkernel}
under the condition that $|x-y| > \frac{1}{2}(1 + |t-s|)$, for if
this were not the case then we could simply integrate by parts as
many times as necessary with respect to the variable $\lambda =
|\xi|$ in the expression \eqref{L2_TT*_disp_tdkernel} and easily
achieve \eqref{kernal_disp_est}. Therefore, we will now show that:
\begin{equation}
        \llp{\td{K}^{TT^*}(t,s;x,y)} \ \lesssim \
         |x-y|^{-\frac{n-1}{2}} \ . \label{td_kernel_disp_est}
\end{equation}
We now factor the phase in \eqref{L2_TT*_disp_tdkernel} as:
\begin{equation}
        e^{2\pi i\big( (t-s)|\xi| + ( x-y)\cdot \xi\big)}
        \ = \ e^{2\pi i (t-s)\lambda}\, e^{2\pi i \lambda\, |x-y| \,
        \cos(\Theta_{x-y,\omega})} \ , \notag
\end{equation}
where we are using the frequency polar coordinates $\xi = \lambda
\omega$. Integrating first on the sphere $\mathbb{S}^{n-1}$, we see
that to conclude \eqref{td_kernel_disp_est} it is enough to show
that:
\begin{multline}
        \llp{ \int_{\mathbb{S}^{n-1}} \
        e^{2\pi i  \lambda\, |x-y| \,
        \cos(\Theta_{x-y,\omega})}\
        \td{\og}^{-1}(t,x) \td{\og}(s,y) \big[\bullet \big]\
        \td{\og}^{-1}(s,y) \td{\og}(t,x) \ d\omega } \\
        \lesssim \ \  |x-y|^{-\frac{n-1}{2}} \ .
        \label{sph_int_disp}
\end{multline}
This last estimate will follow easily from the Morse lemma and the
already established symbol bounds
\eqref{st_symbol_bds1}--\eqref{st_symbol_bds2}. To implement this,
we first cut off the above integral into small neighborhoods of
stationary points of the phase and a remainder. We do this with the
smooth partition of unity:
\begin{equation}
        1 \ =\ \chi_{| 1 - \cos(\Theta_{x-y,\omega})| < \frac{1}{8}}
        + \chi_{| 1 + \cos(\Theta_{x-y,\omega})| < \frac{1}{8}} +
        \td{\chi} \ . \notag
\end{equation}
The cutoff $\td{\chi}$ cuts off on the region where
$\cos(\Theta_{x-y,\omega})$ is bounded away from $\pm 1$, and there
we have the gradient estimate:
\begin{equation}
        c \ < \ |\nabla_\omega \cos(\Theta_{x-y,\omega})| \ , \notag
\end{equation}
for a sufficiently small constant $c$. Using this, and integrating
by parts $n-1$ times while using the symbol bounds
\eqref{st_symbol_bds1}--\eqref{st_symbol_bds2}, we easily have that:
\begin{multline}
        \llp{ \int_{\mathbb{S}^{n-1}} \
        e^{2\pi i  \lambda\, |x-y| \,
        \cos(\Theta_{x-y,\omega})}\
        \td{\og}^{-1}(t,x) \td{\og}(s,y) \big[\bullet \big]\
        \td{\og}^{-1}(s,y) \td{\og}(t,x) \ \td{\chi}(\omega) \ d\omega } \\
        \lesssim \ \  \lambda^{1-n}\cdot
	|x-y|^{-\frac{n-1}{2}} \ . \notag
\end{multline}
This proves \eqref{sph_int_disp} because we may assume that
$\frac{1}{4} < \lambda$. Our goal is now to prove the localized
estimate:
\begin{multline}
        \llp{ \int_{\mathbb{S}^{n-1}} \
        e^{2\pi i  \lambda\, |x-y| \,
        \cos(\Theta_{x-y,\omega})}\
        \td{\og}^{-1}(t,x) \td{\og}(s,y) \big[\bullet \big]\\
        \td{\og}^{-1}(s,y) \td{\og}(t,x) \
        \td{\chi}_{| 1 - \cos(\Theta_{x-y,\omega})| <
        \frac{1}{8}}(\omega) \ d\omega } \ \ \
        \lesssim \ \ \  |x-y|^{-\frac{n-1}{2}} \ . \notag
\end{multline}
It will become clear that the corresponding estimate for the region
where $| 1 + \cos(\Theta_{x-y,\omega})| < \frac{1}{8}$ follows from
identical calculations.\\

Now, the angular function $\cos(\Theta_{x-y,\omega})$ has a single
non-degenerate critical point in a neighborhood of the unit vector
$(x-y)/|x-y|$ with index $n-1$. Therefore, by the Morse lemma there
exists a diffeomorphism $\theta = \varphi(\omega)$ in a neighborhood
of this point such that:
\begin{equation}
        1 - \cos(\Theta_{x-y,\omega}) \ = \ \theta_1^2 + \ldots +
        \theta_{n-1}^2 \ . \notag
\end{equation}
By making this change of variables, we see that we are trying to
prove that:
\begin{multline}
        ||| \int_{\mathbb{R}^{n-1}} \
        e^{2\pi i  \lambda\, |x-y| \,
        |\theta|^2 }\
        {{}^{\varphi^{-1}(\theta)}\! \td{g}}^{-1}(t,x)
        {{}^{\varphi^{-1}(\theta)}\! \td{g}}(s,y) \big[\bullet \big]\\
        {{}^{\varphi^{-1}(\theta)}\! \td{g}}^{-1}(s,y)
        {{}^{\varphi^{-1}(\theta)}\! \td{g}}(t,x) \
        \chi(\theta) \ J_{\varphi^{-1}}(\theta)\
        d\theta ||| \ \ \
        \lesssim \ \ \ |x-y|^{-\frac{n-1}{2}} \ . \notag
\end{multline}
Here $J_{\varphi^{-1}}$ denotes the Jacobian matrix of
$\varphi^{-1}$, and $\chi$ is some smooth function which is
supported where $|\theta| \leqslant 1$. Making now the simple change
of variables $\sqrt{\lambda |x-y|} \theta = \theta'$, it suffices to
be able to show that:
\begin{multline}
        ||| \int_{\mathbb{R}^{n-1}} \
        e^{2\pi i   \, |\theta'|^2 }\
        {{}^{\td{\varphi} (\theta')}\! \td{g}}^{-1}(t,x)
        {{}^{\td{\varphi} (\theta')}\! \td{g}}(s,y) \big[\bullet \big]\
        {{}^{\td{\varphi} (\theta')}\! \td{g}}^{-1}(s,y)
        {{}^{\td{\varphi} (\theta')}\! \td{g}}(t,x) \
        \td{J}(\theta')\
        d\theta' ||| \\ \ \
        \lesssim \ \ \ 1 \ . \label{frensel_int}
\end{multline}
Here $\td{J}(\theta')$ denotes a smooth function with (large)
compact support and uniform gradient bounds:
\begin{equation}
        |\nabla^k_{\theta'} \td{J}| \ \lesssim \ 1 \ . \notag
\end{equation}
Furthermore, the function $\td{\varphi} (\theta')$ obeys the
gradient bounds:
\begin{equation}
        |\nabla_{\theta'}^k \td{\varphi}| \ \lesssim \
        |x-y|^{-\frac{k}{2}} \ . \notag
\end{equation}
Combining this last estimate with the symbol bounds
\eqref{st_symbol_bds1}--\eqref{st_symbol_bds2} and the truncation
condition $M = |x-y|^\frac{1}{2}$, we have the uniform gradient
estimates:
\begin{align}
        \llp{ \nabla_{\theta'}^k \big( {}^{\td{\varphi}}\!
        \td{g}^{-1}(t,x) {}^{\td{\varphi}}\!
        \td{g}(s,y) \big)}
        \ \lesssim \ 1 \ , \notag \\
        \llp{ \nabla_{\theta'}^k \big( {}^{\td{\varphi}}\!
        \td{g}^{-1}(s,y) {}^{\td{\varphi}}\!
        \td{g}(t,x) \big)} \ \lesssim \ 1\ . \notag
\end{align}
Using these bounds, we can prove the bound \eqref{frensel_int} by
treating the quantity on the left hand side as a Fresnel-type integral
and performing $n$ integrations by parts in the region where $1 <
|\theta'|$.\\

To complete our proof of \eqref{kernal_disp_est} we need to show
that:
\begin{equation}
        \lp{\mathcal{R}^{TT^*}(t,s;\cdot,\cdot)}{L^\infty_{x,y}} \ \lesssim \
        (1 + |t-s|)^{-\frac{n-1}{2}} \ , \label{remainder_disp}
\end{equation}
where $\mathcal{R}^{TT^*}$ is the kernel which is defined by
subtracting \eqref{L2_TT*_disp_tdkernel} from
\eqref{L2_TT*_disp_kernel}. Using the splitting
\eqref{gh_left}--\eqref{gh_right} we see that this has at least one
factor involving the  expressions ${\oh}^{-1}(x){\oh}(y) - I$ or
${\oh}(x){\oh}^{-1}(y) - I $ under the integral sign. There are
several such combination, but we will choose to estimate only one
such term and leave the others to reader as they can be treated
analogously. Therefore, we may without loss of generality assume
that we are trying to prove the bound:
\begin{multline}
       \llp{ \int_0^\infty \int_{\mathbb{S}^{n-1}} \
        e^{2\pi i\lambda \big( (t-s) + ( x-y)\cdot\omega\big)}\
        {}^\omega G(t,x;s,y) \ \chi(\lambda)\
        \lambda^{n-1} \ d\lambda d\omega }\\
        \lesssim \ \ (1 + |t-s|)^{-\frac{n-1}{2}} \ ,
        \label{small_disp_case}
\end{multline}
where we have set:
\begin{equation}
        {}^\omega G(t,x;s,y) \ = \
        \td{\og}^{-1}(t,x)\left({\oh}^{-1}(t,x){\oh}(s,y) - I
        \right) \td{\og}(s,y) \big[\bullet \big]\
        \og^{-1}(s,y) \og(t,x) \ . \notag
\end{equation}
As in the proof of \eqref{kernal_disp_est} above for the smoothed
out kernel $\td{K}^{TT^*}$, we may without loss of generality assume
that we trying to prove \eqref{small_disp_case} in the region where
$|x-y| > \frac{1}{2}( 1 +|t-s|)$ because otherwise we may integrate
as many times as necessary with respect to the radial frequency
variable to pick up the desired decay.\\

To proceed further, we will first decompose the range of frequency
integration into a small set and a remainder where we can again
integrate by parts with respect to $\lambda$. This is accomplished
by using the  angular partition of unity:
\begin{equation}
        1 \ = \ \chi_{| \frac{t-s}{|x-y|} + \cos(\Theta_{x-y,\omega}) |
        > |x-y|^{\gamma-1}} + \chi_{| \frac{t-s}{|x-y|} 
	+ \cos(\Theta_{x-y,\omega}) |
        \leqslant  |x-y|^{\gamma-1}} \ . \notag
\end{equation}
To deal with the bound \eqref{small_disp_case} for the first cutoff
function above, we need to show that:
\begin{multline}
        \llp{ \int_{\mathbb{S}^{n-1}} \
        {}^\omega G(t,x;s,y) \ d\omega \\
        \cdot \ \int_0^\infty
        e^{2\pi i\lambda \big( (t-s) + ( x-y)\cdot\omega\big)}\
        \chi_{| \frac{t-s}{|x-y|} + \cos(\Theta_{x-y,\omega}) |
        > |x-y|^{\gamma-1}}\ \chi(\lambda)\
        \lambda^{n-1} \ d\lambda  }\\
        \lesssim \ \ |x-y|^{-\frac{n-1}{2}} \ . \notag
\end{multline}
This bound follows easily from radial integration by parts
in the inner integral, followed by the simple compactness estimate:
\begin{equation}
        \int_{\mathbb{S}^{n-1}} \
        \llp{ {}^\omega G(t,x;s,y)} \ d\omega \
    \lesssim \ 1 \ , \notag
\end{equation}
which is of course uniform in the variables $(t,x;s,y)$.\\

To wrap things up here, we need to show the absolute estimate:
\begin{equation}
        \int_{\RR} \ \llp{{}^\omega G(t,x;s,y) }\
        \chi_{| \frac{t-s}{|x-y|} + \cos(\Theta_{x-y,\omega}) |
        \leqslant |x-y|^{\gamma-1}}\ \chi(\xi)\
        d\xi  \
        \lesssim \ \ |x-y|^{-\frac{n-1}{2}} \ . \notag
\end{equation}
After a Cauchy-Schwartz, this will follow once we can establish that
both:
\begin{align}
        \left(\int_{\RR}\ \chi_{| \frac{t-s}{|x-y|} + \cos(\Theta_{x-y,\omega}) |
        \leqslant |x-y|^{\gamma-1}}\ \chi(\xi)\
        d\xi\right)^\frac{1}{2} \ &\lesssim \
        |x-y|^{\frac{1}{2}(\gamma-1)} \ , \label{small_disp_int1}\\
        \left(\int_{\mathbb{S}^{n-1}}\
        \llp{{}^\omega G(t,x;s,y) }^2 \ d\omega \right)^\frac{1}{2} \
        &\lesssim \
        |x-y|^{-\frac{1}{2}({n-2} + \gamma)} \ . \label{small_disp_int2}
\end{align}
The first estimate, \eqref{small_disp_int1} follows from elementary
bounds. Notice first that after a rotation, it suffices to assume
that the vector $x-y$ lies along the $(1,0)$ direction. Then the
cutoff function is supported in the region where:
\begin{equation}
        \frac{\xi_1}{|\xi|} \  = \ -\frac{t-s}{|x-y|} + O( |x-y|^{\gamma-1}) \ ,
        \notag
\end{equation}
which is a conical set about the $\xi_1$-axis of volume
no greater than a constant times $|x-y|^{\gamma-1}$ in the region where
$|\xi| \lesssim 1$. The second estimate \eqref{small_disp_int2}
above we have already shown. It is a special case of the bound
\eqref{h_int_bnd1} which was proved in the previous section. This
completes our proof of \eqref{remainder_disp}, and hence our
demonstration of the dispersive estimate \eqref{kernal_disp_est}.

\ret

\section{The Decomposable Function Spaces: Proof of the Square-Sum
and Differentiated Strichartz Estimates}\label{Decomp_sect}

We now introduce a piece of machinery which will be of central
importance for the remainder of the paper. This is a suitable
reinterpretation of the important ``decomposable function''
criterion from the work \cite{RT_MKG}. In our context, we set
up the general situation as follows:  Suppose we are given
an $M(m\times m)$ valued Fourier integral operator:
\begin{equation}
        \Phi(\widehat{f})(t,x) \ = \
    \int_{\RR} \ e^{2\pi i \psi(t,x;\xi)} e^{2\pi i
        x\cdot\xi}\ g_1(t,x;\xi) \ \widehat{f}(\xi) \ g_2(t,x;\xi)
    \ d\xi \ , \label{gen_FIO}
\end{equation}
where the $g_i$ are arbitrary matrix valued functions, such that
this operator satisfies certain mixed Lebesgue space mapping
properties (uniform in $y_0$):
\begin{equation}
        \lp{\Phi_{y_0}(\widehat{f})}{L^{q_1}(L^{r_1})} \ \lesssim \
        \lp{\widehat{f}}{L^2} \ , \label{Lp_FIO_est}
\end{equation}
where $\Phi_{y_0}$ is the same operator as \eqref{gen_FIO} but with
phase $\psi$ replaced by $\psi_{y_0}= \psi(t,x-y_0;\xi)$. Suppose
now that we are given a matrix valued function $C(t,x;\omega)$ which
only depends on the angular variable $\omega=\xi/|\xi|$ in
frequency. We would like to prove estimates for the coupled operator
(we only discuss left multiplication here, the case of right
multiplication is analogous):
\begin{equation}
        \widetilde{\Phi}(\widehat{f})(t,x) \ = \ \int_{\RR} \
    e^{2\pi i \psi(t,x;\xi)} e^{2\pi i
        x\cdot\xi}\   C(t,x;\omega)\cdot g_1(t,x;\xi)\
    \widehat{f}(\xi) \ g_2(t,x;\xi)
    \ d\xi \ .\label{coupled_gen_FIO}
\end{equation}
These should be done in a way that the decay properties of the
function $C(t,x;\omega)$ can be used to improve the range of the
estimates \eqref{Lp_FIO_est}. A robust way for doing this has 
been worked out in the paper of Rodnianski--Tao
\cite{RT_MKG}. The answer is to fix an angular scale, say $\theta$, and
then to form the norm (``classical'' decomposable norm):
\begin{equation}
        \lp{C}{D^{cl}_\theta\big(L_t^{q_2}(L_x^{r_2})\big)}^2 \ = \
        \sum_{k=0}^{10n}\ \theta^{-n+1}\ \int_{\mathbb{S}_\omega^{n-1}}\
        \lp{(\theta \nabla_\xi )^k \ C }{L_t^{q_2}(L_x^{r_2})}^2\
        d\omega \ . \label{RT_decomp_norm}
\end{equation}
By decomposing the frequency variable in \eqref{coupled_gen_FIO}
into angular sectors of size $\sim \theta$, a straightforward
computation then shows that one has the estimate:
\begin{equation}
        \lp{\widetilde{\Phi}(\widehat{f})}{L^q(L^r)} \ \lesssim \
        \lp{C}{D^{cl}_\theta\big(L_t^{q_2}(L_x^{r_2})\big)}
        \cdot\lp{\widehat{f}}{L^2} \ , \label{Lp_FIO_est_decomp}
\end{equation}
whenever estimate \eqref{Lp_FIO_est} holds with $\frac{1}{q} =
\frac{1}{q_1} + \frac{1}{q_2}$ and $\frac{1}{r} = \frac{1}{r_1} +
\frac{1}{r_2}$.\\

There are two problems which occur when trying to apply
\eqref{RT_decomp_norm} in the present context. The first is that
this norm is for a single scale, which causes problems in products
where many different scales interact with each other. The other
problem, which is conceptually much more serious, is that the
estimate \eqref{RT_decomp_norm} contains the highly singular factor
of $\theta^{-\frac{n-1}{2}}$, which needs to be eliminated with a
delicate orthogonality argument, the kind which is not preserved in
this problem for a variety of reasons (non-linear Hodge systems, a
covariant wave equation that does not commute with angular cutoffs,
etc). However, with only a slight reworking the basic idea behind
\eqref{RT_decomp_norm} can be shown to be surprisingly robust. First
of all, for a fixed scale we replace \eqref{RT_decomp_norm} with a
square function norm which has the same effect, and which will be
very easy to verify in the present context. Since we will be using
multiple scales in a moment, we introduce the \emph{solid} angular
cutoff functions $\overline{b^{\phi}}_\theta(\omega)$ (not to be
confused with the \emph{hollow} multipliers $b^\omega_\theta(\xi)$
introduced in Section \ref{basic_anal_sect}), such that:
\begin{equation}
        \overline{b^{\phi}}_\theta(\omega) \ \equiv \ 1 \ ,
\end{equation}
when $\omega\in\Gamma_\phi$, for the angular sector $\Gamma_\phi$
which we interpret as a cap in a finitely overlapping
collection on the sphere $\mathbb{S}_\omega^{n-1}=
\cup_\phi \Gamma_\phi$. Here the scale is determined by the condition
$|\Gamma_\phi|\sim \theta$. On this
scale, we replace \eqref{RT_decomp_norm} with the norm:
\begin{equation}
        \lp{C}{D_\theta \big(L_t^{q_2}(L_x^{r_2})\big)} \ = \
    \big\Vert \Big( \sum_{k=0}^{10n} \
        \sum_\phi \ \sup_\omega\ \lp{\overline{b^{\phi}}_\theta\
    (\theta \nabla_\xi )^k \ C }
        {L_x^{r_2}}^2\Big)^\frac{1}{2}\big\Vert_{L_t^{q_2}}
    \ . \label{theta_decomp_norm}
\end{equation}
It is not difficult to see that by decomposing the integral on the
right hand side of \eqref{RT_decomp_norm} into fine and course
scales, and applying H\"older's on the fine (continuous) scales,
that the Rodnianski-Tao norm \eqref{RT_decomp_norm} with the time
integral on the outside is bounded by the norm
\eqref{theta_decomp_norm}. Furthermore, it is easy to see from the
proof given in \cite{RT_MKG} that having the time integral on the
outside does not effect the bound \eqref{Lp_FIO_est_decomp} so long
as the index $q_1$ implicitly appearing in this bound is such that
$2\leqslant q_1$. This allows one to use Minkowski's inequality to
pull the square sum on the parametrix through the time integral. For
us this index condition will always hold because we are working with
Strichartz type norms. We leave it to the reader to work out the
details of these claims.\\

We now form an $\ell^1$ Banach space
based on incorporating the norms \eqref{theta_decomp_norm} over all
dyadic angular scales $\theta \lesssim 1$. The elements of this space we
denote by $\{C\}= \{C^{(\theta)}\}$, and we define its norm
$\ell^1(D_\theta)$ norm as:
\begin{equation}
        \lp{\{C\}}{ \ell^1 \big( D_\theta (L_t^{q_2}(L_x^{r_2}))\big) }
        \ = \ \sum_\theta\ \lp{C^{(\theta)}}{D_\theta
    \big(L_t^{q_2}(L_x^{r_2})\big)}
        \ . \label{sum_theta_decomp_norm}
\end{equation}
There is also the forgetful map from the space $\ell^1(D_\theta)$ to
functions which define as:
\begin{equation}
        \{C\} \ \rightsquigarrow \ C \ = \ 
	\sum_\theta \ C^{(\theta)} \ ,
        \label{forget_map}
\end{equation}
and we will in practice abusively identify $\{C\}$ with $C$ via the
map \eqref{forget_map}. The main point is that given any function
$C$, there may be a \emph{variety} of ways which we embed $C$ in the
space $\ell^1(D_\theta)$, and it is up to the structure of the
application to decide how this should be done. Of course, given the
square function norms \eqref{ang_sum_norms} we are working with,
our choice here is somewhat canonical.\\

Now, if we consider the $C$ in \eqref{forget_map} as embedded in the
integral \eqref{coupled_gen_FIO}, we easily have the estimate:
\begin{equation}
        \lp{\widetilde{\Phi}(\widehat{f})}{L^q(L^r)} \ \lesssim \
        \lp{\{C\}}{  \ell^1 \big( D_\theta (L_t^{q_2}(L_x^{r_2}))\big)   }
        \cdot\lp{\widehat{f}}{L^2} \ . \label{l1_decomp_norm}
\end{equation}
We also form spatial Besov versions of the norm
\eqref{l1_decomp_norm}, which we denote as $\ell^1
D_\theta(L^{q}(\dot{B}_2^{r,(2,s)}))$. This leads us to the basic
notation of this section:\\

\begin{defn}\label{decomp_norm_def}
For a given matrix valued function, we say it is in the \emph{decomposable
space} $D(L^{q}(\dot{B}_2^{r,(2,s)})) $ if the following norm
is finite:
\begin{equation}
    \lp{C}{D(L^{q}(\dot{B}_2^{r,(2,s)}))} \ \ = \ \
    \inf_{C=\sum_\theta C^{(\theta)} }\ \ \Big\{
    \sum_\theta \lp{C^{(\theta)}}{D_\theta(L^{q}(\dot{B}_2^{r,(2,s)}))}
    \Big\} \ . \label{inf_decomp_norm}
\end{equation}
We also define the low frequency analog of these spaces, which we
denote by $D(L^{q}(\dot{B}_{2,10n}^{r,(2,s)}))$, similarly.
\end{defn}\ret

\noindent We remark here that it is easy to see that the norm
\eqref{inf_decomp_norm} leads to a Banach space. This will be
important in a moment. Also, it is easy to show that the various
Besov-Lebesgue space inclusions
\eqref{Besov_nesting}--\eqref{Linfty_besov_incl} hold for these
spaces if we define $D(L^p)$ analogously to \eqref{inf_decomp_norm}.
This is a simple consequence of the fact that the Littlewood-Paley
theory commutes with the derivatives $\nabla_\xi^k$. We now show
that this space satisfies the expected range of bilinear Riesz
operator estimates:

\begin{lem}[A decomposable Besov calculus]\label{decomp_besov_calc}
Let the indices $0\leqslant \sigma$, $1\leqslant q_i,r_i \leqslant \infty$,
and $s_i$ be given. Then one has the following family of bilinear estimates:
\begin{equation}
        |D_x|^{-\sigma} \ : \
        D\big(L_t^{q_1} (\dot{B}_2^{r_1,(2,s_1)})\big)\cdot
         D\big(L_t^{q_2}
        (\dot{B}_2^{r_2,(2,s_2)})\big)
        \ \hookrightarrow \
        D\big(L_t^{q_3}(\dot{B}_1^{r_3,(2,s_3)})\big) \ ,
        \label{general_decomp_besov_embed}
\end{equation}
where the various indices satisfy the conditions:
\begin{align}
        s_3 \ &= \ s_1 + s_2 + \sigma -\frac{n}{2}
    \ , \label{Dsc_cond} \\
        \sigma + \frac{n}{2} - s_3 \ &< \ n(\frac{1}{r_1} + \frac{1}{r_2}) \
        ,  \label{Dgap_cond} \\
        s_{1} \ &< \ \frac{n}{2} + \min \{n( \frac{1}{r_2}
    - \frac{1}{r_3} )\ , \ 0\} \ ,
         \label{Dpos_cond1} \\
        s_{2} \ &< \ \frac{n}{2} + \min \{( \frac{1}{r_1} -
    \frac{1}{r_3})\ , \ 0 \} \ ,
         \label{Dpos_cond2}\\
        \frac{1}{q_3} \ &= \ \frac{1}{q_1} + \frac{1}{q_2}
        \ , \label{DLb_cond1} \\
        \frac{1}{r_3} \ &\leqslant \ \frac{1}{r_1} + \frac{1}{r_2}
         \ . \label{DLb_cond2}
\end{align}
\end{lem}\ret

\begin{proof}[Proof of estimate \eqref{general_decomp_besov_embed}]
The proof of \eqref{general_decomp_besov_embed} is largely a
triviality given that it is true for the norms
$L_t^{q_1}(\dot{B}_2^{r_1,(2,s_1)})$ without the decomposable
structure. First of all, notice that from the definition
\ref{decomp_norm_def} it suffices to establish things with the norms
$D\big(L_t^{q_i} (\dot{B}_2^{r_i,(2,s_i)})\big)$ replaced by their
vector generalizations $\ell^1 D_\theta\big(L_t^{q_i}
(\dot{B}_2^{r_i,(2,s_i)})\big)$. This follows at once from working
with two test matrices $A$ and $C$ and decomposing them into sums:
\begin{align}
    A \ &= \ \sum_\theta \ A^{(\theta)}\ ,
    &C \ &= \ \sum_\theta \ C^{(\theta)} \ , \notag
\end{align}
where the $\{A\}$ and $\{C\}$ collections have norms no greater than
twice that of $A$ and $C$ respectively.\\

Suppose now that we are given two test elements $\{A^{(\theta)}\}$ and
$\{C^{(\theta)}\}$. The we write their product  under the map
\eqref{forget_map} as:
\begin{align}
        A\cdot C \ &= \ \sum_{\theta_1,\theta_2} \ A^{(\theta_1)}\cdot
        C^{(\theta_2)} \ , \label{AC_prod_decomp}\\
        &= \ \sum_\theta \ \Big(\sum_{\substack{\theta_1 \
    :\\ \theta < \theta_1}}\
        A^{(\theta_1)}\cdot C^{(\theta)} \ \ + \ \
        \sum_{\substack{\theta_2 \ :\\ \theta \leqslant \theta_2}}\
        A^{(\theta)}\cdot C^{(\theta_2)}
        \Big) \ , \notag\\
        &= \ \sum_\theta \ T_1^{(\theta)} \ + \ T_2^{(\theta)} \ .
        \notag
\end{align}
Freezing the scale $\theta$, we will prove the following
two estimates:
\begin{align}
        \lp{|D_x|^{-\sigma}\
    T_1^{(\theta)}}{D_\theta\big(L_t^{q_3}
    (\dot{B}_1^{r_3,(2,s_3)})\big)}
        \ \lesssim \ \lp{\{A\}}{\ell^1 D_\theta\big(L_t^{q_1}
    \dot{B}_2^{r_1,(2,s_1)}\big)}\cdot
        \lp{C^{(\theta)}}{D_\theta\big(L_t^{q_2}
    \dot{B}_2^{r_2,(2,s_2)}\big)}
        \ , \label{Decomp_to_prove1}\\
        \lp{|D_x|^{-\sigma}\
     T_2^{(\theta)}}{D_\theta\big(L_t^{q_3}
     (\dot{B}_1^{r_3,(2,s_3)}  )\big)}
        \ \lesssim \ \lp{A^{(\theta)}}{D_\theta\big(L_t^{q_1}
    \dot{B}_2^{r_1,(2,s_1)}\big)}\cdot
        \lp{\{C\}}{\ell^1 D_\theta\big(L_t^{q_2}
    \dot{B}_2^{r_2,(2,s_2)}\big)}
        \ . \label{Decomp_to_prove2}
\end{align}
We will only concentrate on \eqref{Decomp_to_prove1}, as the second
estimate above follows from virtually identical reasoning. Expanding
out the sum in $T_1^{(\theta)}$, it suffices to show:
\begin{multline}
        \lp{|D_x|^{-\sigma}\ (A^{(\theta_1)}\cdot C^{(\theta)})
        }{D_\theta\big(L_t^{q_3}(\dot{B}_1^{r_3,(2,s_3)}  )\big)}
        \ \lesssim \\
        \lp{A^{(\theta_1)}}{D_{\theta_1}\big(L_t^{q_1}
    \dot{B}_2^{r_1,(2,s_1)}\big)}\cdot
        \lp{C^{(\theta)}}{D_\theta\big(L_t^{q_2}
    \dot{B}_2^{r_2,(2,s_2)}\big)}
        \ , \label{last_step}
\end{multline}
where we have the condition $\theta \leqslant \theta_1$. We now
compute the norm on the right hand side of this last equation.
For the remainder of the proof we fix the time variable. This can
then be dealt with at the end by integrating in time and using
H\"olders inequality because all of the action
in the norms \eqref{theta_decomp_norm} takes place under the time
integral. To proceed, we first fix the angular sector $\Gamma_\phi$
and the number of $(\theta \nabla_\xi)$
derivatives to compute that:
\begin{align}
        &\sup_\omega \ \lp{\overline{b^\phi}_\theta \ (\theta\nabla_\xi)^k
        |D_x|^{-\sigma}\ (A^{(\theta_1)}
    \cdot C^{(\theta)})}{ \dot{B}_1^{r_3,(2,s_3)} }
        \ , \notag\\
        \lesssim \ &\sum_{i=0}^k \ \sum_{\substack{\phi_1 \ :\\
        \Gamma_{\phi_1} \subseteq 10\Gamma_\phi}}\
        \sup_\omega \ \lp{
        |D_x|^{-\sigma}\ \Big(\overline{b^{\phi_1}}_{\theta}
        (\theta\nabla_\xi)^{k-i} A^{(\theta_1)}\ \cdot\
    \overline{b^\phi}_\theta
        (\theta\nabla_\xi)^i
        C^{(\theta)}\Big)}{ \dot{B}_2^{r_3,(2,s_3)} }
        \ , \notag\\
        \lesssim \ &\sum_{i=0}^k \ \sum_{\substack{ \phi_1 \ :\\
        \Gamma_{\phi_1} \subseteq 10\Gamma_\phi} }\
        \sup_\omega \ \lp{\overline{b^{\phi_1}}_{\theta}
        (\theta\nabla_\xi)^{k-i}
        A^{(\theta_1)}}{\dot{B}_2^{r_1,(2,s_1)}}\ \cdot\
    \sup_\omega\
        \lp{\overline{b^\phi}_\theta
    (\theta\nabla_\xi)^i C^{(\theta)}}
        {\dot{B}_2^{r_2,(2,s_2)}} \ . \notag
\end{align}
Square summing this last expression over angular sectors, and adding
over all $0\leqslant k \leqslant 10n$ we arrive at the estimate:
\begin{multline}
        \lp{|D_x|^{-\sigma}\ (A^{(\theta_1)}\cdot C^{(\theta)})
        }{D_\theta\big(\dot{B}_1^{r_3,(2,s_3)}\big)} \
        \lesssim \\
        \sup_\phi \
        \sum_{k=0}^{10n} \
        \sup_\omega \ \lp{\overline{b^\phi}_\theta
        (\theta\nabla_\xi)^{k}
        A^{(\theta_1)}}{\dot{B}_2^{r_1,(2,s_1)}} \cdot
        \lp{C^{(\theta)}}{D_\theta\big(
    \dot{B}_2^{r_2,(2,s_2)}\big)} \ . \notag
\end{multline}
We can now conclude \eqref{last_step} on account of the condition
$\theta \leqslant \theta_1$ which implies the trivial bound:
\begin{align}
        &\sup_\phi \
        \sum_{k=0}^{10n} \
        \sup_\omega \ \lp{\overline{b^\phi}_\theta
        (\theta\nabla_\xi)^{k}
        A^{(\theta_1)}}{\dot{B}_2^{r_1,(2,s_1)}} \ , \notag\\
        \lesssim \ &\sup_{\phi'} \
        \sum_{k=0}^{10n} \
        \sup_\omega \ \lp{\overline{b^{\phi'}}_{\theta_1}
        (\theta_1\nabla_\xi)^{k}
        A^{(\theta_1)}}{\dot{B}_2^{r_1,(2,s_1)}} \ ,
        \notag\\
        \lesssim \
     &\lp{A^{(\theta_1)}}{D_{\theta_1}
     \big(\dot{B}_2^{r_1,(2,s_1)}\big)} \ . \notag
\end{align}
This completes our proof of the estimate \eqref{general_decomp_besov_embed}.
\end{proof}\ret

We now establish the link which relates the norms \eqref{inf_decomp_norm}
to the $\dot{X}^s$ norms we have proved for the parametrix $\Phi$:\\

\begin{lem}[Core decomposable estimates for the potentials
$\{\oA^\pm\}$ and $\{\oC^\pm\}$]\label{core_decomp_besov_lem} Let
the sets of potentials $\{\oA^\pm\}$ and $\{\oC^\pm\}$ be defined as
on lines \eqref{A_omega_def}, \eqref{uCpm_system}, and
\eqref{C0_system} above. Then one has the following family of
decomposable bounds:
\begin{align}
        \lp{\oA^\pm}{
        D\big(L^\infty_t(\dot{B}_{2,10n}^{p_\gamma,(2,\frac{n-2}{2})})\big)
        } \ &\lesssim \ \mathcal{E} \ ,
        &\lp{\oA^\pm}{
        D\big(L^2_t(\dot{B}_{2,10n}^{q_\gamma,(2,\frac{n-1}{2})})\big)
        } \ &\lesssim \ \mathcal{E} \ , \ \label{main_A_decomp_ests1}\\
        \lp{\nabla_t \uoA^\pm}{
        D\big(L^\infty_t(\dot{B}_{2,10n}^{p_\gamma,(2,\frac{n-4}{2})})\big)
        } \ &\lesssim \ \mathcal{E} \ ,
        &\lp{\nabla_t \uoA^\pm}{
        D\big(L^2_t(\dot{B}_{2,10n}^{q_\gamma,(2,\frac{n-3}{2})})\big)
        } \ &\lesssim \ \mathcal{E} \ , \ \label{main_A_decomp_ests2}\\
        \lp{\oC^\pm}{
        D\big(L^\infty_t(\dot{B}_{2,10n}^{p_\gamma,(2,\frac{n-2}{2})})\big)
        } \ &\lesssim \ \mathcal{E} \ ,
        &\lp{\oC^\pm}{
        D\big(L^2_t(\dot{B}_{2,10n}^{q_\gamma,(2,\frac{n-1}{2})})\big)
        } \ &\lesssim \ \mathcal{E} \ , \
        \label{main_C_decomp_ests1}\\
        \lp{\nabla_t \uoC^\pm}{
        D\big(L^\infty_t(\dot{B}_{2,10n}^{p_\gamma,(2,\frac{n-4}{2})})\big)
        } \ &\lesssim \ \mathcal{E} \ ,
        &\lp{\nabla_t \uoC^\pm}{
        D\big(L^2_t(\dot{B}_{2,10n}^{q_\gamma,(2,\frac{n-3}{2})})\big)
        } \ &\lesssim \ \mathcal{E} \ , \ \label{main_C_decomp_ests2}
\end{align}
where $p_\gamma$ and $q_\gamma$ are the dimensional constants from
lines \eqref{p_gamma_line} and \eqref{q_gamma_def} above.
Furthermore, one has the following improved null-differentiated
space-time bounds:
\begin{align}
    \lp{ \big(\oL^\mp \uoA^\pm \ , \nabla_t\Delta^{-\frac{1}{2}}
    \oL^\mp \uoA^\pm\ \Big) }{
    D\big(L^2_t(\dot{B}_{2,10n}^{p_\gamma,(2,\frac{n-3}{2})})\big)
    } \ &\lesssim \ \mathcal{E} \ , \ \label{imp_main_A_decomp_ests2}\\
    \lp{  \big(\oL^\mp \uoC^\pm \ , \nabla_t\Delta^{-\frac{1}{2}}
    \oL^\mp \uoC^\pm\ \Big) }{
    D\big(L^2_t(\dot{B}_{2,10n}^{p_\gamma,(2,\frac{n-3}{2})})\big)
    } \ &\lesssim \ \mathcal{E} \ . \ \label{imp_main_C_decomp_ests2}
\end{align}
In all of these estimates, the small constant $\mathcal{E}$ is the
same as on lines \eqref{red_conctn_cond5} and
\eqref{red_conctn_cond7} above.
\end{lem}\ret

\begin{proof}[Proof of the estimates
\eqref{main_A_decomp_ests1}--\eqref{imp_main_C_decomp_ests2}] With
the current setup, the proof of these bounds is very simple and repeats
many of things we have already done in previous sections. Starting
with the estimates
\eqref{main_A_decomp_ests1}--\eqref{main_A_decomp_ests2}, we see
that using the truncation condition \eqref{red_conctn_cond4} it
suffices to prove the first collection, as the time differentiated
versions will follow from these with little fuss. We now follow
essentially the same steps used to prove estimates
\eqref{Besov_degen_tdA_Ln_smallness} and \eqref{temp_int_last_step}.
The only difference here is that we incorporate the square function
norms contained in the $\dot{X}^\frac{n-2}{2}$ spaces. In what
follows, we will in fact only prove the space-time estimate which is
the second bound on the right hand side of
\eqref{main_A_decomp_ests1} above. The first bound on this line
follows from similar reasoning and is left to the reader. The first
step is to define the scale decomposition (we now ignore $\pm$
notation for the remainder of the proof):
\begin{equation}
        \oA \ = \ \sum_\theta \ \oA^{(\theta)} \ = \
        \sum_\theta \ \oPi_\theta \ \oA \ . \notag
\end{equation}
Our goal is now to prove the following fixed time bounds which can
easily be summed over and then integrated to achieve the desired
goal:
\begin{equation}
        \lp{\oPi_\theta  \oA\, (t)}{D_\theta\big(
        \dot{B}_{2,10n}^{q_\gamma,(2,\frac{n-1}{2})}\big)} \
        \lesssim \ \theta^\gamma\
        \lp{ \uA\, (t)}{S
        \dot{B}_{2}^{\frac{2(n-1)}{n-3},(2,\frac{n-1}{2})}} \ .
        \notag
\end{equation}
By using the square function structure contained in the definition
of the various Besov and decomposable Besov norms and taking into
account the low frequency truncation of the potentials $\{\oA\}$ and
$\{\uA\}$, the proof of this last estimate reduces to the fixed
frequency bound:
\begin{equation}
        \lp{\oPi_\theta P_\mu \oA\, (t)}{D_\theta\big(
        \dot{B}_{2}^{q_\gamma,(2,\frac{n-1}{2})}\big)} \
        \lesssim \ \theta^\gamma\
        \lp{P_\mu \uA\, (t)}{S
        \dot{B}_{2}^{\frac{2(n-1)}{n-3},(2,\frac{n-1}{2})}} \ .
        \notag
\end{equation}
Expanding now the decomposable norm on the left hand side of this
last inequality, we see that the proof reduces to showing the square
function bounds:
\begin{multline}
        \sum_{k=0}^{10n} \
        \sum_\phi \ \sup_\omega\ \lp{\overline{b^{\phi}}_\theta(\omega)\
        (\theta \nabla_\xi )^k \  \oPi_\theta P_\mu \oA\, (t)}
        {  \dot{B}_{2}^{q_\gamma,(2,\frac{n-1}{2})}}^2 \\
        \lesssim \ \ \theta^{2\gamma}
        \sum_{\substack{\phi\ :\\
     \omega_0 \in \Gamma_{\phi}} }
    \ \lp{\td{{}^{\omega_0}\!\Pi
        }_\theta P_\mu \uA\, (t)}{
        \dot{B}_{2}^{\frac{2(n-1)}{n-3},(2,\frac{n-1}{2})}}^2 \ ,
        \label{sq_sum_diff_to_prove}
\end{multline}
where $\td{\oPi}_\theta$ is a fixed thickening of the multiplier
$\oPi_\theta$ such that one has the general quasi-idempotence bound:
\begin{align}
        \sup_\omega \lp{\overline{b^{\phi}}_\theta(\omega)\,
        \td{\td{\oPi}}_\theta A}
        {  L^q } \ &\lesssim \ \lp{\td{{}^{\omega_0}\!\Pi}_\theta A}{L^q}
        \ , &\omega_0\in\Gamma_\phi \ , \label{general_idem}
\end{align}
where $\td{\td{\oPi}}_\theta$ is any multiplier with frequency
support contained in the frequency support of $\oPi_\theta$ whose
convolution kernel satisfies comparable $L^1$ bounds. Here the
statement that $\omega_0\in\Gamma_\phi$ is taken to mean that
$\omega_0$ is in the center of the cap $\Gamma_\phi$, the very same
notion we used in the definition of the square function norms
\eqref{ang_sum_norms} above. Using now the general bound
\eqref{general_idem} as well as the heuristic multiplier identity:
\begin{equation}
        (\theta \nabla_\xi )^k \  \oPi_\theta P_\mu
        \ooPi^{(\frac{1}{2} - \delta)}
        \nabla_{t,x} \oL\,  \Delta_{\omega^\perp}^{-1}
        \underline{A}_{\, \bullet \ll 1}(\partial_\omega)\, (t)
        \ \approx \ \theta^{-1} \ \oPi_\theta P_\mu
        \uA\, (t) \ , \notag
\end{equation}
we have the bound:
\begin{equation}
        \hbox{(L.H.S)}\eqref{sq_sum_diff_to_prove} \
        \lesssim \ \sum_{\substack{\phi\ :\\
     \omega_0 \in \Gamma_{\phi}} }
        \ \theta^{-2}\ \lp{\td{{}^{\omega_0}\!\Pi
        }_\theta P_\mu \uA\, (t)}{
        \dot{B}_{2}^{q_\gamma,(2,\frac{n-1}{2})}}^2 \ . \notag
\end{equation}
The estimate \eqref{sq_sum_diff_to_prove} now follows from the
Bernstein nested-Besov inclusion:
\begin{equation}
        \td{{}^{\omega_0}\!\Pi}_\theta 
	(\dot{B}_{2}^{\frac{2(n-1)}{n-3} ,(2,\frac{n-1}{2})})
        \ \subseteq \ \theta^{1+\gamma}\
        \td{{}^{\omega_0}\!\Pi}_\theta
        (\dot{B}_{2}^{q_\gamma,(2,\frac{n-1}{2})}) \ .
        \notag
\end{equation}\ret

Our next goal is to pass the estimates
\eqref{main_A_decomp_ests1}--\eqref{main_A_decomp_ests2} on to the
non-linear set of potentials $(\oC_0,\{\uoC\})$. Since it is
a-priori not clear that these functions have finite
$D(L^{q}(\dot{B}_2^{r,(2,s)}))$ norms, we construct the bounds from
scratch by running a contraction mapping argument in these spaces on
the Picard iterates of the systems \eqref{uCpm_system} and
\eqref{C0_system}. To guarantee convergence of the resulting
sequences, we make use specific instances of the general embedding
\eqref{general_decomp_besov_embed}. Our general strategy here is the
following. We first establish the non-time differentiated estimates
\eqref{main_C_decomp_ests1} for the spatial potentials $\{\uoC\}$.
Then, assuming the non-time differentiated versions of the improved
estimates
\eqref{imp_main_A_decomp_ests2}--\eqref{imp_main_C_decomp_ests2}
(whose proof relies only on the previously established bounds) we
prove the time-differentiated estimates \eqref{main_C_decomp_ests2}.
Having established these, we then prove the estimates
\eqref{main_C_decomp_ests1} for the temporal potential $\oC_0$. Our
next order of business is to prove the non-time differentiated
versions of the improved null-differentiated bounds
\eqref{imp_main_A_decomp_ests2}--\eqref{imp_main_C_decomp_ests2}.
Finally, armed with all of this, we show the version of the
estimates
\eqref{imp_main_A_decomp_ests2}--\eqref{imp_main_C_decomp_ests2}
which contain the extra time derivatives. In what follows, we will
only list out the various bilinear estimates which yield the desired
bounds. Since these are almost identical to many of the estimates we
have dealt with in the past sections, we leave the verification of
the numerology to the reader.\\

To prove the non-time differentiated versions of
\eqref{main_C_decomp_ests1} for the collection $\{\uoC\}$ we use the
pair of bounds:
\begin{align}
        \nabla_x\Delta^{-1}\ : \
        D\big(L^\infty_t (\dot{B}_{2,10n}^{p_\gamma,(2,\frac{n-2}{2})})\big)\cdot
         D\big(L_t^\infty
        (\dot{B}_{2,10n}^{p_\gamma,(2,\frac{n-2}{2})})\big)
        \ &\hookrightarrow \
        D\big(L_t^{\infty}
        (\dot{B}_{2,10n}^{p_\gamma,(2,\frac{n-2}{2})})\big) \ ,
        \label{decomp_Lpgamma_est}\\
        \nabla_x\Delta^{-1}\ : \
        D\big(L^2_t (\dot{B}_{2,10n}^{q_\gamma,(2,\frac{n-1}{2})})\big)\cdot
        D\big(L_t^\infty
        (\dot{B}_{2,10n}^{p_\gamma,(2,\frac{n-2}{2})})\big)
        \ &\hookrightarrow \
        D\big(L_t^{2}
        (\dot{B}_{2,10n}^{q_\gamma,(2,\frac{n-1}{2})})\big) \ .
        \label{decomp_Lqgamma_est}
\end{align}
To establish the first bound on line \eqref{main_C_decomp_ests2} we
first differentiate the Hodge system \eqref{uCpm_system} with
respect to time and then apply the embedding:
\begin{equation}
        \nabla_x\Delta^{-1}\ : \
        D\big(L^\infty_t (\dot{B}_{2,10n}^{p_\gamma,(2,\frac{n-4}{2})})\big)\cdot
         D\big(L_t^\infty
        (\dot{B}_{2,10n}^{p_\gamma,(2,\frac{n-2}{2})})\big)
        \ \hookrightarrow \
        D\big(L_t^{\infty}
        (\dot{B}_{2,10n}^{p_\gamma,(2,\frac{n-4}{2})})\big) \ .
        \label{dt_decomp_Lpgamma_est}
\end{equation}
To prove the time integrated bound which is the second on line
\eqref{main_C_decomp_ests2} we decompose the vector-field $\nabla_t$
into $\pm \uL \mp\omega\cdot\nabla_x$ just as we did starting on
line \eqref{LLb_wave} above. Then, modulo estimates of the form
\eqref{decomp_Lqgamma_est}, and assuming that we have shown the
non-time differentiated versions of
\eqref{imp_main_A_decomp_ests2}--\eqref{imp_main_C_decomp_ests2} we
may reduce things to the embedding:
\begin{equation}
        \nabla_x\Delta^{-1}\ : \
        D\big(L^2_t (\dot{B}_{2,10n}^{p_\gamma,(2,\frac{n-3}{2})})\big)\cdot
         D\big(L_t^\infty
        (\dot{B}_{2,10n}^{p_\gamma,(2,\frac{n-2}{2})})\big)
        \ \hookrightarrow \
        D\big(L_t^{2}
        (\dot{B}_{2,10n}^{p_\gamma,(2,\frac{n-3}{2})})\big) \ .
        \label{dt_decomp_Lqgamma_est}
\end{equation}\ret

Our next step is to prove the estimates \eqref{main_C_decomp_ests1}
for the temporal potential $\oC_0$.  By an inspection of the
elliptic equation \eqref{C0_system}, we see that modulo  embeddings
of the form \eqref{decomp_Lpgamma_est}--\eqref{decomp_Lqgamma_est}
and the bounds we have already shown, we only need to establish
things for the term $\nabla_t \Delta^{-1}([\uoA,\uoC])$. Again
expanding the time derivative as $\pm \uL \mp\omega\cdot\nabla_x$
and distributing the $\uL$ derivative via the Leibniz rule, we are
reduced knowing the following (which just represent another form of
the embeddings
\eqref{dt_decomp_Lpgamma_est}--\eqref{dt_decomp_Lqgamma_est}):
\begin{align}
        \Delta^{-1}\ : \
        D\big(L^\infty_t (\dot{B}_{2,10n}^{p_\gamma,(2,\frac{n-4}{2})})\big)\cdot
         D\big(L_t^\infty
        (\dot{B}_{2,10n}^{p_\gamma,(2,\frac{n-2}{2})})\big)
        \ &\hookrightarrow \
        D\big(L_t^{\infty}
        (\dot{B}_{2,10n}^{p_\gamma,(2,\frac{n-2}{2})})\big) \ .
        \label{dinv_dt_decomp_Lpgamma_est}\\
        \Delta^{-1}\ : \
        D\big(L^2_t (\dot{B}_{2,10n}^{p_\gamma,(2,\frac{n-3}{2})})\big)\cdot
         D\big(L_t^\infty
        (\dot{B}_{2,10n}^{p_\gamma,(2,\frac{n-2}{2})})\big)
        \ &\hookrightarrow \
        D\big(L_t^{2}
        (\dot{B}_{2,10n}^{p_\gamma,(2,\frac{n-1}{2})})\big) \ .
        \label{dinv_dt_decomp_Lqgamma_est}
\end{align}\ret

Finally, we wish to show the improved bounds
\eqref{imp_main_A_decomp_ests2}--\eqref{imp_main_C_decomp_ests2}. We
work recursively here. First, we assume that the non-time
differentiated versions of these estimates are valid. By the
truncation condition \eqref{red_conctn_cond4}, we see that the proof
of the estimate \eqref{imp_main_A_decomp_ests2} with the extra
operator $\partial_t\Delta^{-\frac{1}{2}}$ follows from the proof of
this estimate without that operator. Thus, our aim is to establish
the estimate \eqref{imp_main_C_decomp_ests2} in the presence of the
extra $\partial_t\Delta^{-\frac{1}{2}}$ derivatives. Applying this
operator to the $\uL$ differentiated Hodge system
\eqref{uCpm_system}, we see that things can be handled with the help
of the two bilinear inclusions:
\begin{align}
        \nabla_x \Delta^{-\frac{3}{2}}\ : \
        D\big(L^2_t (\dot{B}_{2,10n}^{p_\gamma,(2,\frac{n-5}{2})})\big)\cdot
         D\big(L_t^\infty
        (\dot{B}_{2,10n}^{p_\gamma,(2,\frac{n-2}{2})})\big)
        \ &\hookrightarrow \
        D\big(L_t^{2}
        (\dot{B}_{2,10n}^{p_\gamma,(2,\frac{n-3}{2})})\big) \ .
        \label{singular_decomp_L2Lpgamma_est1}\\
        \nabla_x \Delta^{-\frac{3}{2}}\ : \
        D\big(L^2_t (\dot{B}_{2,10n}^{p_\gamma,(2,\frac{n-3}{2})})\big)\cdot
         D\big(L_t^\infty
        (\dot{B}_{2,10n}^{p_\gamma,(2,\frac{n-4}{2})})\big)
        \ &\hookrightarrow \
        D\big(L_t^{2}
        (\dot{B}_{2,10n}^{p_\gamma,(2,\frac{n-3}{2})})\big) \ .
        \label{singular_decomp_L2Lqgamma_est2}
\end{align}
Notice that the numerology in these last two estimates is a bit
tight in  $High\times High$ frequency regime. In particular, the
condition \eqref{Dgap_cond} only has room of about $1/10$ in $n=6$
dimensions. The next item on the stack for us is the estimates
\eqref{imp_main_A_decomp_ests2}--\eqref{imp_main_C_decomp_ests2}
without the extra derivative $\nabla_t \Delta^{-1}$. Assuming for
the moment that this is true for \eqref{imp_main_A_decomp_ests2}, we
see that the proof of \eqref{imp_main_C_decomp_ests2} in this case
follows easily from $\uL$ differentiating the Hodge system
\eqref{uCpm_system} and applying a less singular version of the
estimate \eqref{singular_decomp_L2Lpgamma_est1}. Therefore, we are
now at the point where everything has been reduced to the proof of
the first estimate \eqref{imp_main_A_decomp_ests2}. To do this we
apply the following instance of the identity \eqref{bLA_identity}:
\begin{equation}
        \uL \uoA \ = \ \nabla_{x}\
        \ooPi^{(\frac{1}{2} - \delta)}
        \   \uA(\partial_\omega) \
        -  \ \nabla_{x}\
        \ooPi^{(\frac{1}{2} - \delta)}
        \  \Delta_{\omega^\perp}^{-1}\, \td{\mathcal{P}}
        \big([B,H]\big)(\partial_\omega) \ .
        \label{bLA_iden_clean}
\end{equation}
The estimate \eqref{imp_main_A_decomp_ests2} for the first term on
the right hand side of \eqref{bLA_iden_clean} is very simple and
left to the reader. It follows from steps similar to the proof we
gave above of the estimates \eqref{main_A_decomp_ests1}. Notice that
there are no singular angular factors here so there is a lot of room
in this estimate if one takes into account the extra Coulomb savings
\eqref{main_coulomb_savings}.\\

We are now trying to prove \eqref{imp_main_A_decomp_ests2} for the
second term on the right hand side of \eqref{bLA_iden_clean} which
we decompose into angular scales as:
\begin{equation}
        \nabla_{x}\
        \ooPi^{(\frac{1}{2} - \delta)}
        \  \Delta_{\omega^\perp}^{-1}\, \td{\mathcal{P}}
        \big([B,H]\big)(\partial_\omega) \ = \
        \sum_\theta \ \nabla_x\ \oPi_\theta \ooPi^{(\frac{1}{2} - \delta)}
        \  \Delta_{\omega^\perp}^{-1}\, \td{\mathcal{P}}
        \big([B,H]\big)(\partial_\omega) \ . \label{LA_quad_ang_sum}
\end{equation}
By the definition of the norm
\eqref{inf_decomp_norm} and using some dyadic summing, we see
that it suffices to bound and then sum over 
the fixed time-fixed frequency expressions:
\begin{multline}
        \lp{\nabla_x\  \oPi_\theta \ooPi^{(\frac{1}{2} - \delta)}
        \  \Delta_{\omega^\perp}^{-1}\, \td{\mathcal{P}}
        \big([B,H]\big)(\partial_\omega)\, (t)
        }{D\big( \dot{B}_{2,10n}^{p_\gamma,(2,\frac{n-3}{2})}\big) } \ , \\
        = \ \Big(
        \sum_{k=0}^{10n}\ \sum_\phi\
        \sup_\omega \ \lp{\overline{b^\phi}_\theta
        (\theta \nabla_\xi )^k \
        \nabla_x\ \oPi_\theta \ooPi^{(\frac{1}{2} - \delta)}
        \  \Delta_{\omega^\perp}^{-1}\, \td{\mathcal{P}}
        \big([B,H]\big)(\partial_\omega)\, (t)
        }{\dot{B}_{2,10n}^{p_\gamma,(2,\frac{n-3}{2})}}^2\
         \Big)^\frac{1}{2} \ . \label{fixed_time_decomp_BH_to_prove}
\end{multline}
For each fixed value of $\theta$, and for each fixed spatial
frequency $\mu$ we have the following heuristic multiplier bound one
the Coulomb savings are taken into account:
\begin{equation}
        (\theta \nabla_\xi )^k \
        \nabla_x\ \oPi_\theta \ooPi^{(\frac{1}{2} - \delta)}
        \  \Delta_{\omega^\perp}^{-1}\, P_\mu\ \td{\mathcal{P}}
        \big([B,H]\big)(\partial_\omega)\, (t) \ \approx \
        (\mu\theta)^{-1}\ \oPi_\theta  P_\mu P_{\bullet \ll 1}
        \big([B,H]\big)\, (t) \ . \notag
\end{equation}
Taking this into account, and using the same multiplier reductions
used to prove \eqref{sq_sum_diff_to_prove} above, we have the
inequality:
\begin{align}
        \hbox{(L.H.S.)}\eqref{fixed_time_decomp_BH_to_prove}
        \ &\lesssim \  \theta^{-1}\ \Big(
        \sum_{\substack{\phi\ :\\
     \omega_0 \in \Gamma_{\phi}}}
        \lp{\td{ {}^{\omega_0}\!\Pi}_\theta P_{\bullet\ll 1}
        \big([B,H]\big)\, (t)
        }{\dot{B}_{2,10n}^{p_\gamma,(2,\frac{n-5}{2})}}^2\
        \Big)^\frac{1}{2} \ , \notag\\
        &\lesssim \ \theta^\gamma \ \Big(
        \sum_{\substack{\phi\ :\\
     \omega_0 \in \Gamma_{\phi}}}
        \lp{\td{{}^{\omega_0}\!\Pi}_\theta
        \big([B,H]\big)\, (t)
        }{\dot{B}_{2}^{2,(2,\frac{n-5}{2})}}^2\
        \Big)^\frac{1}{2} \ , \notag\\
        &\lesssim \ \theta^\gamma\
        \lp{ \big([B,H]\big)\, (t)
        }{\dot{B}_{2}^{2,(2,\frac{n-5}{2})}} \ .
        \label{arrived_fixed_time_bound}
\end{align}
This last set of inequalities results from the localized Besov
inclusion:
\begin{equation}
        \td{{}^{\omega_0}\!\Pi}_\theta
        (\dot{B}_{2}^{2,(2,\frac{n-5}{2})})
        \ \subseteq \ \theta^{1+\gamma}\
        \td{{}^{\omega_0}\!\Pi}_\theta
        (\dot{B}_{2}^{p_\gamma,(2,\frac{n-5}{2})}) \ ,
        \notag
\end{equation}
an orthogonality argument. Integrating the bound
\eqref{arrived_fixed_time_bound} $L^2$ in time, and using some
dyadic summing, we see that our
proof of \eqref{imp_main_A_decomp_ests2} is reduced to showing the
following:
\begin{equation}
        \lp{[B,H]}{L^2_t(\dot{H}^\frac{n-5}{2})}
        \ \lesssim \ \mathcal{E} \ . \notag
\end{equation}
Keeping in mind the bootstrapping estimates
\eqref{red_conctn_cond7}, we see that this last line is simply a
more singular version of the embedding \eqref{bilin_X_embed2} shown
above. In the $Low\times High$ case the proof follows from
\eqref{bilin_X_embed2a}. In the $High\times Low$ case there is even
more room and one can again use something similar to
\eqref{bilin_X_embed2a}. In the $High \times High$ case we use the
embedding:
\begin{equation}
        P_{ \lambda }
        \big(L^2(\dot{B}^{\frac{2(n-1)}{n-3},(2,\frac{n-1}{2})})\big)
        \cdot P_\lambda \big(L^\infty(\dot{H}^\frac{n-4}{2})\big)
        \ \hookrightarrow \  \left(\frac{\mu}{\lambda}\right)^\sigma P_\mu\big(
        L^2(\dot{H}^\frac{n-5}{2})\big) \  ,
       \notag
\end{equation}
where $\sigma = n(\frac{n-2}{n-1}) - \frac{5}{2}$. This last bound
follows from the general frequency localized embedding
\eqref{freq_loc_general_besov_embed2}. Note that in dimensions
$6\leqslant n$ we have the necessary condition $0 < \sigma$. This
completes our proof of the estimate \eqref{imp_main_A_decomp_ests2}
and therefore our demonstration of Lemma
\ref{core_decomp_besov_lem}.
\end{proof}\ret

\subsection{Proof of the Square Sum Strichartz Estimates}
We now come to what is perhaps the linchpin of our argument. These
are the square sum structure estimates contained in
\eqref{half_wave_est1}. With the current machinery in hand, these
will be quite easy to establish. At the heart of things
is whether the angular multipliers $\oPi_\theta$ ``commute with the
dynamics'' of the covariant wave operator $\Box_{\uA}$. At a quick
first glance using Duhamel's principle, this seems to be connected
with whether one can control the commutator $[\oPi_\theta,\Box_{\uA}]$.
Unfortunately, it is not too difficult to see that one runs into
serious difficulties as soon as $\theta \ll 1$. This is not the end
of the story however, because it turns out that modulo a very nice
error term, one can control the commutator with the ``integrated''
form the equations $[\oPi_\theta,\Phi]$. This shows one of the deep
advantages to working with the parametrix as opposed to dealing
directly with the equations themselves\footnote{This also seems
to have far
reaching philosophical consequences for how one should proceed in
lower dimensions. Specifically, it seems to suggest that the correct
``covariant'' $X^{s,\theta}$ space should be defined in terms of the
parametrix and \emph{not} in terms of the symbol of the covariant
equation.}. We proceed as follows.\\

Our first step is to fix a scale $\theta$ and run a cap decomposition
$\mathbb{S}^{n-1} = \cup_\phi \Gamma_\phi$. The next thing we do is
to decompose the parametrix $\Phi(\widehat{f})$ into a sum of three
pieces:
\begin{align}
        \Phi(\widehat{f}) \ &= \ \int_{\mathbb{R}^n} \
        e^{2\pi i \lambda \ou} \ \og_{\bullet\ll\theta}^{-1} \,
        \widehat{f}(\lambda \omega)
        \, \og_{\bullet\ll\theta} \
        \chi_{(\frac{1}{2},2)}(\lambda)\
        \lambda^{n-1} d\lambda d\omega \notag\\
        &\ \ \ \ \ \ \ \ \ + \ \int_{\mathbb{R}^n} \
        e^{2\pi i \lambda \ou} \ \og_{\bullet\ll\theta}^{-1} \,
        \widehat{f}(\lambda \omega)
        \, \og_{\theta\lesssim \bullet}\
        \chi_{(\frac{1}{2},2)}(\lambda)\
        \lambda^{n-1} d\lambda d\omega\notag\\
        &\ \ \ \ \ \ \ \ \ \ \ \ \ \ \ \ \ \ \ \
        + \  \int_{\mathbb{R}^n} \
        e^{2\pi i \lambda \ou} \ \og_{\theta\lesssim \bullet}^{-1} \,
        \widehat{f}(\lambda \omega)
        \, \og \
        \chi_{(\frac{1}{2},2)}(\lambda)\
        \lambda^{n-1} d\lambda d\omega \ , \notag\\
        &= \ I_1 \ +\  I_2 \ +\  I_3 \ . \notag
\end{align}
Here:
\begin{equation}
        \og \ = \ \og_{\bullet\ll\theta} + \og_{\theta\lesssim \bullet} \ = \
        P_{\bullet\ll\theta}(\og) + P_{\theta\lesssim \bullet}(\og)
        \ , \notag
\end{equation}
is a low-high frequency decomposition of the group element $\og$. We
define the decomposition for $\og^{-1}$ similarly. Our goal is now
to prove the following three estimates:
\begin{align}
        \sum_{\substack{\phi\ :\\
     \omega_0 \in \Gamma_{\phi}}}
        \ \lp{{}^{\omega_0}\!\Pi_\theta P_1 (I_1)}{L^2(L^\frac{2(n-1)}{n-3})}^2
        \ &\lesssim \ \lp{\widehat{f}}{L^2}^2 \ ,
        \label{I1_sq_sum_est}\\
        \sum_{\substack{\phi\ :\\
     \omega_0 \in \Gamma_{\phi}}}
        \ \lp{{}^{\omega_0}\!\Pi_\theta P_1 (I_2)}{L^2(L^\frac{2(n-1)}{n-3})}^2
        \ &\lesssim \ \lp{\widehat{f}}{L^2}^2 \ ,
        \label{I2_sq_sum_est}\\
        \sum_{\substack{\phi\ :\\
     \omega_0 \in \Gamma_{\phi}}}
        \ \lp{{}^{\omega_0}\!\Pi_\theta P_1 (I_3)}{L^2(L^\frac{2(n-1)}{n-3})}^2
        \ &\lesssim \ \lp{\widehat{f}}{L^2}^2 \ .
        \label{I3_sq_sum_est}
\end{align}
The proof of the first bound \eqref{I1_sq_sum_est} follows easily
from the plain endpoint Strichartz estimate we have already
established. To see this, first notice that for a fixed angle one
has the identity:
\begin{multline}
        {}^{\omega_0}\!\Pi_\theta P_1\ \int_{\mathbb{R}^n} \
        e^{2\pi i \lambda \ou} \ \og_{\bullet\ll\theta}^{-1} \,
        \widehat{f}(\lambda \omega)
        \, \og_{\bullet\ll\theta} \
        \chi_{(\frac{1}{2},2)}(\lambda)\
        \lambda^{n-1} d\lambda d\omega \\
        = \ \ {}^{\omega_0}\!\Pi_\theta P_1\ \int_{\mathbb{R}^n} \
        e^{2\pi i \lambda \ou} \ \og_{\bullet\ll\theta}^{-1} \,
        \overline{b^{\phi'}}_\theta(\omega) \widehat{f}(\lambda \omega)
        \, \og_{\bullet\ll\theta} \
        \chi_{(\frac{1}{2},2)}(\lambda)\
        \lambda^{n-1} d\lambda d\omega \ , \notag
\end{multline}
where $\Gamma_\phi \subset \frac{1}{2}\Gamma_{\phi'}$ for some fixed 
thickening $\Gamma_{\phi'}$
of the spherical cap that $\omega_0\in\Gamma_\phi$. That this is the
case follows easily from the fact that the Fourier transform of the
function:
\begin{equation}
        e^{2\pi i x\cdot\xi}\ \og_{\bullet\ll\theta}^{-1}(x) \,
        \widehat{f}(\lambda \omega)
        \, \og_{\bullet\ll\theta}(x) \ , \notag
\end{equation}
is a tempered distribution with support contained in an $O(c\theta)$
neighborhood of the point $\xi$ for some small constant $c$, uniform
in the value of $\theta$. Using now the boundedness of the
multiplier ${}^{\omega_0}\!\Pi_\theta P_1$, we only need to
establish that the truncated parametrix $I_1$ obeys the endpoint
Strichartz estimate. We reduce this claim further by writing the
integral in the form:
\begin{equation}
        I_1 \ =\
        \int_{\RR}\int_{\RR}\
        K^{P_{\bullet\ll\theta}}(w) K^{P_{\bullet\ll\theta}}(y)\
        \int_{\mathbb{R}^n} \
        e^{2\pi i \lambda \ou} \ \og^{-1}_w \,
        \widehat{f}(\lambda \omega)
        \, \og_y \
        \chi_{(\frac{1}{2},2)}(\lambda)\
        \lambda^{n-1} d\lambda d\omega\ dw dy \ , \notag
\end{equation}
where $\og^{-1}_w(x)=  \og^{-1}(x-w)$ and  $\og_y(x)= \og(x-y)$
denote the translated group elements. Using the fact that the
convolution kernel $K^{P_{\bullet\ll\theta}}$ has  $O(1)$ $L^1$ norm
uniform in the value of $\theta$, we are left with establishing the
$L^2$ and dispersive estimates of the previous sections for more
general kernels of the form:
\begin{equation}
        \Phi_{g_1,g_2}(\widehat{f}) \ = \
        \int_{\mathbb{R}^n} \
        e^{2\pi i \lambda \ou} \ \og^{-1}_1 \,
        \widehat{f}(\lambda \omega)
        \, \og_2 \
        \chi_{(\frac{1}{2},2)}(\lambda)\
        \lambda^{n-1} d\lambda d\omega \ , \label{general_trans_kernel}
\end{equation}
where $\og_1$ and $\og_2$ are unrelated group elements which are
generated from Hodge systems and connections of the form
\eqref{A_omega_def}--\eqref{C0_system}, which satisfy the general
requirements \eqref{red_conctn_cond} and \eqref{freq_loc_cond} for
$\lambda=1$. This indeed turns out to be the case, and the key
observation is that by using the identity \eqref{Killing_iden}, all
of the $TT^*$ arguments go through just as they did in previous
sections.\\

It remains for us to prove the bounds
\eqref{I2_sq_sum_est}--\eqref{I3_sq_sum_est}. These are essentially
identical to each other so we concentrate on the proof of the first
of these, leaving the other one to the reader. By an application of
Bernstein's inequality and orthogonality, we see that it
suffices for us to show the estimate:
\begin{equation}
    \lp{\theta \ I_2}{L^2(L^2)} \ \lesssim \ \lp{\widehat{f}}{L^2}
     \ . \label{L2_I2_est_to_show}
\end{equation}
At a heuristic level, this estimate is true because one has the
identity $\theta\ \og_{\theta\lesssim\bullet} \approx \nabla_x\og =
g\uoC$. And we see that in this case things would follow easily from
the $D\big(L^2(L^\infty)\big)$ contained in the estimates
\eqref{main_C_decomp_ests1}. To implement this in a rigorous way, we
derive the following elliptic equation for
$\og_{\theta\lesssim\bullet}$ based on the formulas \eqref{C_def}:
\begin{align}
    \og_{\theta\lesssim\bullet} \ &= \
    \nabla^i \Delta^{-1}P_{\theta\lesssim\bullet}(\og \uoC_i)
    \ , \notag\\
    &= \ \sum_{\substack{\lambda \ : \\
    \theta\lesssim \lambda}}\
    \nabla^i \Delta^{-1}P_\lambda(\og \uoC_i) \ . \notag
\end{align}
If we denote the (vector) kernel of operator
$\nabla_x \Delta^{-1}P_{\lambda}$ by
$K_\lambda^{\nabla\Delta^{-1}}$, then we have the uniform $L^1$
bounds:
\begin{equation}
    \lp{K_\lambda^{\nabla\Delta^{-1}}}{L^1} \
    \lesssim \ \lambda^{-1}\ . \notag
\end{equation}
Using this and taking into account the previous reductions
used in the proof of estimate \eqref{I1_sq_sum_est} above
we easily arrive at the bound:
\begin{align}
    \lp{\theta \ I_2}{L^2(L^2)} \ &\lesssim \
    \sum_{\substack{\lambda \ : \\
    \theta\lesssim \lambda}}\ \left(\frac{\theta}{\lambda}\right)
    \ \sup_{w,y}\
    \lp{\td{I}_{w,y}}{L^2(L^2)} \ , \notag\\
    &\lesssim \  \sup_{w,y}\
    \lp{\td{I}_{w,y}}{L^2(L^2)}\ , \notag
\end{align}
where $\td{I}_{w,y}$ is the family of translated kernels:
\begin{equation}
    \td{I}_{w,y} \ = \
        \int_{\mathbb{R}^n} \
        e^{2\pi i \lambda \ou} \ \og^{-1}_w \,
        \widehat{f}(\lambda \omega)
        \, \og_y \uoC_y \
        \chi_{(\frac{1}{2},2)}(\lambda)\
        \lambda^{n-1} d\lambda d\omega\ dw dy \ , \label{tdI_wy_def}
\end{equation}
where we have also now set $\uoC_y(x) = \uoC(x-y)$.
Using the decomposable estimate \eqref{l1_decomp_norm},
we now have that:
\begin{equation}
    \lp{\td{I}_{w,y}}{L^2(L^2)} \ \lesssim \
    \lp{\uoC_y}{D\big(L^2(L^\infty)\big)}
    \cdot\lp{I_{w,y}}{L^\infty(L^2)} \ , \notag
\end{equation}
where the integral $I_{w,y}$ is the same as $\td{I}_{w,y}$ but
with the matrix $\uoC_y$ removed. Using now the nesting:
\begin{equation}
    D\big(L^2(\dot{B}_{2,10n}^{q_\gamma
    ,(2,\frac{n-1}{2})})\big)
    \ \subseteq \ D\big(L^2(L^\infty)\big)
    \ , \label{L2Linfty_decomp_nest}
\end{equation}
the estimate \eqref{main_C_decomp_ests1}, and the
remarks made above about general kernels of the form
\eqref{general_trans_kernel}, we have the pair of estimates:
\begin{align}
    \lp{\uoC_y}{D\big(L^2(L^\infty)\big)} \ &\lesssim \
    \mathcal{E} \ ,
    &\lp{I_{w,y}}{L^\infty(L^2)} \ &\lesssim \
    \lp{\widehat{f}}{L^2} \ , \notag
\end{align}
uniform in the values of $w,y$.
This is enough to prove the estimate \eqref{I2_sq_sum_est}.
This completes our proof of the square sum Strichartz
estimates contained in \eqref{half_wave_est1}.\\

\subsection{Proof of the Differentiated Strichartz Estimates
\eqref{half_wave_est1.5}--\eqref{half_wave_est2} }
To wrap things up for this overall section, we prove the
estimates \eqref{half_wave_est1.5}--\eqref{half_wave_est2}.
This will follows easily from the general list of
decomposable estimates contained in Lemma
\ref{core_decomp_besov_lem}. In what follows, we will only
bother to prove  the time differentiated estimate
\eqref{half_wave_est2}. The proof of the gradient estimate
\eqref{half_wave_est1.5} follows from identical reasoning and
is left to the reader (in fact, one only need apply the plain
Strichartz estimates shown in previous sections followed by a
$D\big(L^\infty(L^\infty)\big)$ estimate for the spatial
potentials $\{\uoC\}$). Time differentiating the parametrix
$\Phi^\pm(\widehat{f})$ we see that:
\begin{align}
    \nabla_t \ \Phi^\pm(\widehat{f})
    \ &= \ \int_{\mathbb{R}^n} \
        (\pm2\pi i\, |\xi|)\, e^{2\pi i \lambda \ou^\pm} \ \og_\pm^{-1} \,
        \widehat{f}(\lambda \omega)
        \, \og_\pm \
        \chi_{(\frac{1}{2},2)}(\lambda)\
        \lambda^{n-1} d\lambda d\omega \notag \\
    &\ \ \ \ \ + \ \ \int_{\mathbb{R}^n} \
        e^{2\pi i \lambda \ou^\pm} \ \big[\og_\pm^{-1} \,
        \widehat{f}(\lambda \omega)
        \, \og_\pm \ , \ \oC_0^\pm \big]\
        \chi_{(\frac{1}{2},2)}(\lambda)\
        \lambda^{n-1} d\lambda d\omega \ , \notag\\
    &=\ \Phi^\pm(\pm2\pi i\, |\xi|\, \widehat{f})
     \ + \ \td{I} \ . \notag
\end{align}
Therefore, our task is to show the pair of estimates:
\begin{align}
    \lp{P_1\td{I}_1}{L^2(SL^\frac{2(n-1)}{n-3})} \ &\lesssim \
    \mathcal{E}\cdot \lp{\widehat{f}}{L^2} \ ,
    \label{dt_Phi_sqsum_est}\\
    \lp{P_1\td{I}_1}{L^\infty(L^2)} \ &\lesssim \
    \mathcal{E}\cdot \lp{\widehat{f}}{L^2} \ .
    \label{dt_Phi_energy_est}
\end{align}
The estimate \eqref{dt_Phi_sqsum_est} follows from
essentially identical reasoning to that employed in the proof
of estimates \eqref{I2_sq_sum_est}--\eqref{I3_sq_sum_est} above.
The main point is to drop to $L^2(L^2)$ via Bernstein, and then
use the $D\big(L^2(L^\infty)\big)$ estimate for the potential
$\oC_0$ contained on line \eqref{main_C_decomp_ests1} above.
The proof of the second  estimate \eqref{dt_Phi_energy_est} above
follows easily from the $D\big(L^\infty(L^\infty)\big)$
estimate for $\oC_0$ contained on line
\eqref{main_C_decomp_ests1} above. Specifically, one has the nesting:
\begin{equation}
    D\big(L^\infty(\dot{B}_{2,10n}^{p_\gamma,(2,\frac{n-2}{2})})\big)
    \ \subseteq \ D\big(L^\infty(L^\infty)\big)
    \ . \notag
\end{equation}
This completes our demonstration of
\eqref{half_wave_est1.5}--\eqref{half_wave_est2} and ends
this section.

\ret

\section{Completion of the proof:\
Controlling the $L^1(L^2)$ Norm of the Differentiated Parametrix}
Our final task here is to prove the estimate \eqref{half_wave_est4}
which guarantees that our parametrix is a good approximation the
covariant wave equation $\Box_{\uA}$. This essentially boils down to
applying the estimates
\eqref{main_A_decomp_ests1}--\eqref{imp_main_C_decomp_ests2} to the
various error terms listed on the right hand side of equation
\eqref{error_terms} above. We will prove the desired estimates
for each of these terms separately.\\

\subsection*{$\bullet$\ \
Decomposing the term\ \  $\uA(\oL^\mp) - \oC^\pm(\oL^\mp)$}
This represents the worst error term which comes out of our
approximation, as well as the main ``renormalization'' which the
parametrix creates. In what follows we will eliminate the $\pm$
notation on favor of the $\uL$ notation introduced on line
\eqref{LLb_wave} above. Using this convention, a short computation
involving the formulas \eqref{uCpm_system}--\eqref{C0_system} and
the structure equation \eqref{red_conctn_cond6} yields the identity:
\begin{align}
    &\uA(\uL) - \oC(\uL) \notag \\
    = \ \ &(I - \ooPi^{(\frac{1}{2} -
    \delta)})\uA(\partial_\omega) \ +
    \ \ooPi^{(\frac{1}{2} - \delta)}
    \Delta_{\omega^\perp}^{-1}\td{\mathcal{P}}
    ([B,H])(\partial_\omega)\notag \\
    &\ \ \ \ \ \ + \ \uL\Delta^{-1}([\uoA,\uoC]) \ - \
    d^*\Delta^{-1}[\uoC, \oC(\uL)] \ , \notag\\
    = \ \ &T_1 \ +\  T_2 \ +\  T_3 \ +\  T_4 \ . \label{T1_T4_defs}
\end{align}
Our goal is prove the following four estimates:
\begin{align}
    \lp{T_1}{D\big(L^2(L^{n-1})\big)} \ &\lesssim \
    \mathcal{E} \ ,
    &\lp{T_2}{D\big(L^1(L^\infty)\big)} \ &\lesssim \
    \mathcal{E}  \ , \label{main_decomp_bL_est1}\\
    \lp{T_3}{D\big(L^1(L^\infty)\big)} \ &\lesssim \
    \mathcal{E} \ ,
    &\lp{T_4}{D\big(L^1(L^{\infty})\big)} \ &\lesssim \
    \mathcal{E}  \ . \label{main_decomp_bL_est2}
\end{align}\ret

To prove the first estimate on line \eqref{main_decomp_bL_est1},
we see from the decomposable version of the Besov nesting
\eqref{Linfty_besov_incl} that is suffices to prove
the following:
\begin{equation}
    \lp{(I - \ooPi^{(\frac{1}{2} -
    \delta)})\uA(\partial_\omega)}{D\big(L^2(\dot{B}
    _{2}^{n-1,(2,\frac{n(n-3)}{2(n-1)})})\big)}
    \ \lesssim \ \mathcal{E} \ . \notag
\end{equation}
By the square sum nature of the Besov and decomposable norms,
and keeping in mind the Besov version of the
endpoint Strichartz estimate contained
in the bootstrapping estimate \eqref{red_conctn_cond5}, we see that
it suffices to prove this estimate at fixed frequency. Thus, we
are trying to prove that:
\begin{equation}
    \lp{(I - \ooPi^{(\frac{1}{2} -
    \delta)})P_\mu\ \uA(\partial_\omega)}{D\big(L^2(
    L^{n-1})\big)} \ \lesssim \
    \lp{P_\mu\uA }{L^2(S
    \dot{B}_2^{\frac{2(n-1)}{n-3},(2,\frac{n-1}{2})}
    )} \ .
    \label{dyadic_trunc_error_est}
\end{equation}
Decomposing the term on the left hand side of this expression into
all dyadic angular regions spread from the direction $\omega$
this is further reduced to showing that:
\begin{equation}
    \lp{\oPi_\theta(I - \ooPi^{(\frac{1}{2} -
    \delta)})P_\mu\ \uA(\partial_\omega)}{D_\theta\big(L^2(
    L^{n-1})\big)} \ \lesssim \
    \theta^\gamma
    \lp{P_\mu\uA }{L^2(S
    \dot{B}_2^{\frac{2(n-1)}{n-3},(2,\frac{n-1}{2})}
    )} \ . \notag
\end{equation}
Notice that we are only trying to show this for
values $\theta \lesssim \ \mu^{\frac{1}{2}-\delta}$.
Further computing the term on the left hand side of this last
expression, and applying the heuristic multiplier bound
(also using the Coulomb savings \eqref{main_coulomb_savings}):
\begin{equation}
    (\theta\nabla_\xi^k)\oPi_\theta(I - \ooPi^{(\frac{1}{2} -
    \delta)})P_\mu\ \uA(\partial_\omega) \ \approx \
    \theta\ \oPi_\theta P_\mu\ \uA \ . \notag
\end{equation}
Plugging this into the definition of the norm $D_\theta\big(L^2(
L^{n-1})\big)$, using the multiplier-sum reductions employed in the
proof of the inequality \eqref{sq_sum_diff_to_prove}, and reverting
back to Besov notation we have the inequality sequence involving
Bernstein's inequality \eqref{Bernstein2} and a simple index
manipulation:
\begin{align}
    \hbox{(L.H.S.)}\eqref{dyadic_trunc_error_est}
    \ &\lesssim \ \theta\ \left(\sum_{\substack{\phi\ :\\
     \omega_0 \in \Gamma_{\phi}} }
    \ \lp{\td{{}^{\omega_0}\!\Pi
        }_\theta P_\mu \uA }{L^2(
        \dot{B}_{2}^{n-1,(2,\frac{n(n-3)}{2(n-1)})} )}^2
    \right)^\frac{1}{2} \ , \notag\\
    &\lesssim \ \theta^\frac{n-3}{2}
    \left(\sum_{\substack{\phi\ :\\
     \omega_0 \in \Gamma_{\phi}} }
    \ \lp{\td{{}^{\omega_0}\!\Pi
        }_\theta P_\mu \uA }{L^2(
        \dot{B}_{2}^{\frac{2(n-1)}{n-3},(2,\frac{n(n-3)}{2(n-1)})} )}^2
    \right)^\frac{1}{2} , \notag\\
    &\lesssim \ \theta^\frac{n-3}{2}\mu^\frac{-n-1}{2(n-1)}\
    \left(\sum_{\substack{\phi\ :\\
     \omega_0 \in \Gamma_{\phi}} }
    \ \lp{\td{{}^{\omega_0}\!\Pi
        }_\theta P_\mu \uA }{L^2(
        \dot{B}_{2}^{\frac{2(n-1)}{n-3},(2,\frac{n-1}{2})} )}^2
    \right)^\frac{1}{2} . \notag
\end{align}
Estimate \eqref{dyadic_trunc_error_est} now follows from the fact
that:
\begin{equation}
    \theta^\frac{n-3}{2}\mu^\frac{-n-1}{2(n-1)} \ \lesssim \
    \theta^\gamma \ , \notag
\end{equation}
which is a consequence of the truncation condition $\theta \lesssim \ \mu^{\frac{1}{2}-\delta}$ and the fact that $6\leqslant n$, and the
fact that we have chosen $\delta,\gamma$ according to
\eqref{dimensional_constant}. This ends our proof of the
first estimate on line \eqref{main_decomp_bL_est1}.\\

Our next step is to prove the second estimate on line
\eqref{main_decomp_bL_est1} above. We will show the somewhat more
regular estimate:
\begin{equation}
        \lp{\ooPi^{(\frac{1}{2} - \delta)}
        \Delta_{\omega^\perp}^{-1}\td{\mathcal{P}}
        ([B,H])(\partial_\omega)}{D\big(L^1(\dot{B}_{1}^{\infty,(n,\frac{n}{2})})\big)}
        \ \lesssim \ \mathcal{E} \ . \label{more_reg_Linfty_decomp}
\end{equation}
Decomposing the term inside the norm on the left hand side of this
last inequality into dyadic angular scales, applying the definition
of the fixed scale decomposable norms
$D_\theta\big(L^1(\dot{B}_{1}^{\infty,(n,\frac{n}{2})})\big)$, using
the (fixed time) fixed frequency heuristic multiplier bound (which
again takes into account the savings \eqref{main_coulomb_savings}):
\begin{equation}
        (\theta\nabla_\xi)^k\ \oPi_\theta\ooPi^{(\frac{1}{2} - \delta)}
        \Delta_{\omega^\perp}^{-1}P_\lambda\ \td{\mathcal{P}}
        ([B,H])(\partial_\omega) \ \approx \ \theta^{-1}\lambda^{-2}\
        \oPi_\theta P_\lambda ([B,H]) \ , \notag
\end{equation}
expanding the resulting expression into a trichotomy, applying the
multiplier square sum reduction used previously in the proof of
estimate \eqref{sq_sum_diff_to_prove} above, and keeping in mind the
bootstrapping structure estimates \eqref{red_conctn_cond7}, we see
that the estimate \eqref{more_reg_Linfty_decomp} reduces to the
demonstration of the following three fixed time bounds:
\begin{align}
    \begin{split}
        &\sum_{\substack{\lambda , \mu_i\ : \\ \mu_1\ll
        \mu_2\sim\lambda}}\ \lambda^{-2}\
        \Big(\sum_{\substack{\phi\ :\\
        \omega_0 \in \Gamma_{\phi}} }
        \ \lp{
        \td{{}^{\omega_0}\!\Pi}_\theta P_\lambda
        ([P_{\mu_1}(B)(t),P_{\mu_2}(H)(t)])
        }{  L^\infty }^2 \Big)^\frac{1}{2} \\
        &\hspace{1.5in}  \lesssim \ \theta^{1+\gamma}\
        \lp{B(t)}{\dot{B}_2^{\frac{2(n-1)}{n-3},(2,\frac{n-1}{2})}}
        \cdot\lp{H(t)}{\dot{B}_2^{\frac{2(n-1)}{n-3},(2,\frac{n-3}{2})}}
        \ ,
    \end{split} \label{bad_LH_Linfty_decomp}\\
    \begin{split}
        &\sum_{\substack{\lambda , \mu_i\ : \\ \mu_2\ll
        \mu_1\sim\lambda}}\ \lambda^{-2}\
        \Big(\sum_{\substack{\phi\ :\\
        \omega_0 \in \Gamma_{\phi}} }
        \ \lp{
        \td{{}^{\omega_0}\!\Pi}_\theta P_\lambda
        ([P_{\mu_1}(B)(t),P_{\mu_2}(H)(t)])
        }{  L^\infty }^2 \Big)^\frac{1}{2} \\
        &\hspace{1.5in}  \lesssim \ \theta^{1+\gamma}\
        \lp{B(t)}{\dot{B}_2^{\frac{2(n-1)}{n-3},(2,\frac{n-1}{2})}}
        \cdot\lp{H(t)}{\dot{B}_2^{\frac{2(n-1)}{n-3},(2,\frac{n-3}{2})}}
        \ ,
    \end{split} \label{HL_Linfty_decomp}\\
    \begin{split}
        &\sum_{\substack{\lambda , \mu_i\ : \\ \lambda\lesssim
        \mu_1\sim \mu_2}}\ \lambda^{-2}\
        \Big(\sum_{\substack{\phi\ :\\
        \omega_0 \in \Gamma_{\phi}} }
        \ \lp{
        \td{{}^{\omega_0}\!\Pi}_\theta P_\lambda
        ([P_{\mu_1}(B)(t),P_{\mu_2}(H)(t)])
        }{  L^\infty }^2 \Big)^\frac{1}{2} \\
        &\hspace{1.5in}  \lesssim \ \theta^{1+\gamma}\
        \lp{B(t)}{\dot{B}_2^{\frac{2(n-1)}{n-3},(2,\frac{n-1}{2})}}
        \cdot\lp{H(t)}{\dot{B}_2^{\frac{2(n-1)}{n-3},(2,\frac{n-3}{2})}}
        \ .
    \end{split} \label{HH_Linfty_decomp}
\end{align}
We begin with the proof of the first estimate
\eqref{bad_LH_Linfty_decomp}. This is the most singular of the
three. Fixing all of the spatial frequencies on the left hand side
of this bound, we see that by an application of Young's inequality,
it suffices to prove the following refinement:
\begin{multline}
        \lambda^{-2}\
        \Big(\sum_{\substack{\phi\ :\\
        \omega_0 \in \Gamma_{\phi}} }
        \ \lp{
        \td{{}^{\omega_0}\!\Pi}_\theta P_\lambda
        ([P_{\mu_1}(B)(t),P_{\mu_2}(H)(t)])
        }{  L^\infty }^2 \Big)^\frac{1}{2} \\
         \lesssim \ \left(\frac{\mu_1}{\mu_2}
         \right)^\gamma \theta^{1+\gamma}\
        \lp{P_{\mu_1}(B)(t)}{\dot{B}_2^{\frac{2(n-1)}{n-3},(2,\frac{n-1}{2})}}
        \cdot\lp{P_{\mu_2}(H)(t)}{\dot{B}_2^{\frac{2(n-1)}{n-3},(2,\frac{n-3}{2})}}
        \ . \label{fixed_freq_bad_LH_Linfty_decomp}
\end{multline}
This bound is scale invariant, so we may assume that
$1=\lambda\sim\mu_2$. To aid in the demonstration, we introduce the
auxiliary index:
\begin{equation}
        \td{r}_\gamma \ = \ \frac{2n(n-1)}{n^2-2n-1 -2(n-1)\gamma}
        \ . \notag
\end{equation}
Notice that this has been chosen precisely so that one has the
identity:
\begin{equation}
        \gamma \ = \ \frac{1}{2} + n(\frac{n-3}{2(n-1)} -
        \frac{1}{\td{r}_\gamma}) \ , \notag
\end{equation}
so that ultimately we can make a reference to the fixed frequency
bound \eqref{gen_besov_fixed_freq}. The problem here is that we have
$2 < \td{r}_\gamma$ (in any dimension), so we are going to run into
orthogonality issues in the square-sum on the left hand side of
\eqref{fixed_freq_bad_LH_Linfty_decomp}. This will end up costing
some extra powers of $\theta^{-1}$, but luckily the Bernstein
inequality will more than make up for this. Applying Bernstein to
each term in the sum on the left hand side of
\eqref{fixed_freq_bad_LH_Linfty_decomp} we arrive at the bound:
\begin{equation}
        \hbox{(L.H.S.)}\eqref{fixed_freq_bad_LH_Linfty_decomp}
        \ \lesssim \ \theta^\frac{n-1}{\td{r}_\gamma}\
        \Big(\sum_{\substack{\phi\ :\\
        \omega_0 \in \Gamma_{\phi}} }
        \ \lp{
        \td{{}^{\omega_0}\!\Pi}_\theta P_1
        ([P_{\mu_1}(B)(t),P_{\mu_2}(H)(t)])
        }{  L^{\td{r}_\gamma} }^2 \Big)^\frac{1}{2} \ .
        \label{bad_LH_after_bernstein}
\end{equation}
To get rid of the square-sum on the right hand side of this last
expression, we introduce the following map from $L^p(\RR)$ to
$\ell^2(L^p(\RR))$:
\begin{equation}
        \mathcal{T}^\theta (A) \ = \
        \left(\td{{}^{\omega_1}\!\Pi}_\theta P_1 (A), \ldots ,
        \td{{}^{\omega_N}\!\Pi}_\theta P_1(A)\right) \ , \notag
\end{equation}
where $(\omega_1,\ldots,\omega_N)$ is some ordering of the
$\Gamma_\phi$ spherical cap ``base-points''. Notice that there are
$N\sim \theta^{1-n}$ of these. By orthogonality, and using the
uniform boundedness of the multipliers $\oPi_\theta P_1$ on $L^\infty$
we have the pair of estimates:
\begin{align}
        \lp{\mathcal{T}^\theta (A)}{\ell^2(L^2)} \ &\lesssim \
        \lp{P_1 (A)}{L^2} \ , \notag\\
        \lp{\mathcal{T}^\theta (A)}{\ell^2(L^\infty)} \ &\lesssim \
        \theta^\frac{1-n}{2}\
        \lp{P_1 (A)}{L^\infty} \ . \notag
\end{align}
By interpolating these to bounds in the pair of spaces
$(\ell^2(L^2),\ell^2(L^\infty))$ and $(L^2,L^\infty)$ (see
\cite{BL_interp}), we have the bound:
\begin{equation}
        \lp{\mathcal{T}^\theta (A)}{\ell^2(L^{\td{r}_\gamma})} \
        \lesssim \ \theta^{(1-n)(\frac{1}{2}-\frac{1}{\td{r}_\gamma})}\
        \lp{P_1 (A)}{L^{\td{r}_\gamma}} \ . \notag
\end{equation}
Plugging this last estimate into the right hand side of
\eqref{bad_LH_after_bernstein} above, and finally applying generic
fixed frequency estimate \eqref{gen_besov_fixed_freq} we have that:
\begin{align}
        &\hbox{(L.H.S.)}\eqref{fixed_freq_bad_LH_Linfty_decomp}
        \ , \notag\\
        \lesssim \ &\theta^{(n-1)(\frac{2}{r_\gamma}-\frac{1}{2})}\
        \lp{  P_1
        ([P_{\mu_1}(B)(t),P_{\mu_2}(H)(t)])
        }{  L^{\td{r}_\gamma} } \ , \notag\\
        \lesssim \ &\left(\frac{\mu_1}{\mu_2}\right)^\gamma
        \theta^{(n-1)(\frac{2}{r_\gamma}-\frac{1}{2})}\
        \lp{P_{\mu_1}(B)(t)}{\dot{B}_2^{\frac{2(n-1)}{n-3},(2,\frac{n-1}{2})}}
        \cdot\lp{P_{\mu_2}(H)(t)}{\dot{B}_2^{\frac{2(n-1)}{n-3},(2,\frac{n-3}{2})}}
        \ . \notag
\end{align}
The estimate \eqref{fixed_freq_bad_LH_Linfty_decomp} now follows
from the bound:
\begin{equation}
        \theta^{(n-1)(\frac{2}{r_\gamma}-\frac{1}{2})} \ \lesssim \
        \theta^{1 + \gamma} \ , \notag
\end{equation}
which holds in dimensions $6\leqslant n$. We leave the verification
of this to the reader. This ends our demonstration of the $Low\times
High$ frequency estimate \eqref{bad_LH_Linfty_decomp}. Notice that
the second estimate \eqref{HL_Linfty_decomp} is simply a less
singular version of this. In fact, repeating the above procedure, we
see that in that case there is an extra factor of
$(\frac{\mu_2}{\mu_1})$ in the analog of the fixed frequency bound
\eqref{fixed_freq_bad_LH_Linfty_decomp}.\\

We have now reduced the proof of the second estimate on line
\eqref{main_decomp_bL_est1} to the $High\times High$ interaction
bound \eqref{HH_Linfty_decomp}. By applying the $L^\infty\to L^2$
version of Bernstein, using orthogonality, and then applying the
general bound \eqref{freq_loc_general_besov_embed2}, we have the
fixed frequency estimate:
\begin{align}
        &\lambda^{-2}\ \Big(\sum_{\substack{\phi\ :\\
        \omega_0 \in \Gamma_{\phi}} }
        \ \lp{
        \td{{}^{\omega_0}\!\Pi}_\theta P_\lambda
        ([P_{\mu_1}(B)(t),P_{\mu_2}(H)(t)])
        }{  L^\infty }^2 \Big)^\frac{1}{2} \ , \notag\\
        \lesssim \
        &\theta^\frac{n-1}{2}\lambda^\frac{n-4}{2}\
        \lp{  P_\lambda
        ([P_{\mu_1}(B)(t),P_{\mu_2}(H)(t)])
        }{  L^2 } \ , \notag\\
        \lesssim \ &\left(\frac{\lambda}{\mu_1}\right)^{\sigma}
        \theta^\frac{n-1}{2}\ \lp{P_{\mu_1}(B)(t)}{\dot{B}_2^{\frac{2(n-1)}{n-3},(2,\frac{n-1}{2})}}
        \cdot\lp{P_{\mu_2}(H)(t)}{\dot{B}_2^{\frac{2(n-1)}{n-3},(2,\frac{n-3}{2})}}
        \ , \notag
\end{align}
where $0 < \sigma = n(\frac{n-3}{n-1}) -2$. By summing this last
line and then applying Cauchy-Schwartz, we easily arrive at the
bound \eqref{HH_Linfty_decomp}.\\

To finish this subsection, we only need to prove the two estimates
on line \eqref{main_decomp_bL_est2} above. To show the first
estimate involving the $T_3$ term, we simply expand the $\uL$
derivative into the product via the Leibniz rule, and then use the
decomposable bounds \eqref{main_A_decomp_ests1} and
\eqref{main_C_decomp_ests1} and
\eqref{imp_main_A_decomp_ests2}--\eqref{imp_main_C_decomp_ests2} in
conjunction with the following instance of the bilinear decomposable
estimate \eqref{general_decomp_besov_embed}:
\begin{equation}
        \Delta^{-1}\ : \
        D\big(L^2_t(\dot{B}_{2,10n}^{p_\gamma,(2,\frac{n-3}{2})})\big)\cdot
        D\big(L^2_t(\dot{B}_{2,10n}^{q_\gamma,(2,\frac{n-1}{2})})\big)
        \ \hookrightarrow \
        D\big(L^1_t(\dot{B}_{1}^{\infty,(2,\frac{n}{2})})\big) \ .
        \notag
\end{equation}
To show the second bound on line \eqref{main_decomp_bL_est2}, we
again use the estimates \eqref{main_A_decomp_ests1} and
\eqref{main_C_decomp_ests1}, this time in conjunction with:
\begin{equation}
        \nabla_x \Delta^{-1}\ : \
        D\big(L^2_t(\dot{B}_{2,10n}^{q_\gamma,(2,\frac{n-1}{2})})\big)\cdot
        D\big(L^2_t(\dot{B}_{2,10n}^{q_\gamma,(2,\frac{n-1}{2})})\big)
        \ \hookrightarrow \
        D\big(L^1_t(\dot{B}_{1}^{\infty,(2,\frac{n}{2})})\big) \ .
        \notag
\end{equation}
This completes our decomposable estimates for the error term
$\uA(\oL^\mp) - \oC^\pm(\oL^\mp)$.\\

\subsection*{$\bullet$\ \
Decomposing the term $D^{\uA}_\alpha \big({\oC^\pm}\big)^\alpha$}
Again dropping the $\pm$ notation and using the equations
\eqref{uCpm_system}--\eqref{C0_system} and the identity
\eqref{LLb_wave} as well as the structure equation
\eqref{red_conctn_cond6}, we can write this as:
\begin{align}
        D^{\uA}_\alpha( {\oC})^\alpha \ &= \
        - \ \ooPi^{(\frac{1}{2} - \delta)}
        \oL\, \Delta_{\omega^\perp}^{-1}\,
        \td{\mathcal{P}}([B,H])(\partial_\omega) \notag\\
        &\ \ \ \ \ \ \ \ \ + \ \ (\pm \uL \mp
        \omega\cdot\nabla_x)\nabla_t\Delta^{-1}
        [\uoA,\uoC]
        \ + \ \nabla_t d^*\Delta^{-1}[\oC_0,\uoC] \notag\\
        &\ \ \ \ \ \ \ \ \ \ \ \ \ \ \ \ \ \
        - \ [\uoA,\uoC] \ + \
        [\uA,\uoC] \ , \notag\\
        &= \ \td{T}_2 \ + \ \td{T}_3 \ + \  \td{T}_4 \ + \ \td{T}_5 \ + \
        \td{T}_6 \ . \notag
\end{align}
We will show that all of these terms obey the estimate:
\begin{align}
        \lp{\td{T}_k}{D\big(L^1(L^\infty)\big)}
        \ &\lesssim \ \mathcal{E} \ ,
        &2 \ \leqslant \ k \ \leqslant \ 6 \ . \label{tdT_k_group_est}
\end{align}
Notice that, for the most part, the terms $\td{T}_k$ represent less
singular versions of the $T_k$ on line \eqref{T1_T4_defs} above. In
fact, they can all be treated using similar embeddings by simply
wasting one derivative. Specifically, the estimate
\eqref{tdT_k_group_est} for the first term $\td{T}_2$ follows
directly from \eqref{more_reg_Linfty_decomp} above once one takes
into account the presence of the truncation \eqref{red_conctn_cond4}
inherent in the projection $\td{\mathcal{P}}$. To prove the estimate
\eqref{tdT_k_group_est} for the portion of term $\td{T}_3$
containing the $\uL$ derivative, we use the same embedding
employed in the proof of the estimate for $T_3$ on line
\eqref{main_decomp_bL_est2} above. This follows because one can
distribute the time derivative and simply waste smoothness in the
estimates \eqref{main_A_decomp_ests2}, \eqref{main_C_decomp_ests2},
and
\eqref{imp_main_A_decomp_ests2}--\eqref{imp_main_C_decomp_ests2}.
Specifically, by taking advantage of the low frequency behavior of
these estimates, we have the bounds:
\begin{align}
        \lp{\nabla_t \uoA}{
        D\big(L^2_t(\dot{B}_{2,9n}^{q_\gamma,(2,\frac{n-1}{2})})\big)
        } \ &\lesssim \ \mathcal{E} \ ,
        &\lp{\nabla_t \uoC}{
        D\big(L^2_t(\dot{B}_{2,9n}^{q_\gamma,(2,\frac{n-1}{2})})\big)
        } \ &\lesssim \ \mathcal{E} \ , \label{der_wasted_bound1}\\
        \lp{ \nabla_t \uL \uoA}{
        D\big(L^2_t(\dot{B}_{2,9n}^{p_\gamma,(2,\frac{n-3}{2})})\big)
        } \ &\lesssim \ \mathcal{E} \ ,
        &\lp{  \nabla_t \uL \uoC}{
        D\big(L^2_t(\dot{B}_{2,9n}^{p_\gamma,(2,\frac{n-3}{2})})\big)
        } \ &\lesssim \ \mathcal{E} \ . \label{der_wasted_bound2}
\end{align}
Using a similar strategy, we can prove the estimate
\eqref{tdT_k_group_est} for the portion of $\td{T}_3$ containing the
$\omega\cdot\nabla_x$ derivative (notice that the functions
$\omega_i$ are trivially decomposable) as well as the term
$\td{T}_4$ in the same way as we showed \eqref{main_decomp_bL_est2}
for the term $T_4$ above. All we need to do is to show the estimate:
\begin{equation}
        \lp{\nabla_t \oC_0}{
        D\big(L^2_t(\dot{B}_{2,9n}^{q_\gamma,(2,\frac{n-1}{2})})\big)
        } \ \lesssim \ \mathcal{E} \ . \notag
\end{equation}
This follows in the same way  we proved the undifferentiated
estimate \eqref{main_C_decomp_ests1} for $\oC_0$ above, but instead
of using the undifferentiated versions of
\eqref{main_A_decomp_ests1}, \eqref{main_C_decomp_ests1}, and
\eqref{imp_main_A_decomp_ests2}--\eqref{imp_main_C_decomp_ests2}, we
simply use \eqref{der_wasted_bound1}--\eqref{der_wasted_bound2}.
Finally, notice that the proof of the estimate
\eqref{tdT_k_group_est} for the terms $\td{T}_5$ and $\td{T}_6$
above follows by simply multiplying (decompose twice!) the
$D\big(L^2(L^\infty)\big)$ estimate which is implied by the bounds
\eqref{main_A_decomp_ests1} and \eqref{main_C_decomp_ests1} above.
This completes our decomposition of the second error term on the
right hand side of \eqref{error_terms} above.\\

\subsection*{$\bullet$\ \
Decomposing the term $ \Big[ \uA^\alpha - (\oC^\pm)^\alpha \ , \ \big[
(\uA)_\alpha - \oC^\pm_\alpha \  , \ \bullet \  \big] \Big]$} Here
we again use the norm $D\big(L^1(L^\infty)\big)$, which we can
achieve as a product of $D\big(L^2(L^\infty)\big)$ estimates, again
making an appeal to \eqref{main_A_decomp_ests1} and
\eqref{main_C_decomp_ests1} above.\\

This completes our proof of the approximation estimate
\eqref{half_wave_est4} and thus, at last, the proof of Proposition
\ref{half_wave_prop} which allows us to close the bootstrapping
begun in Proposition \eqref{critical_bs_prop}. FP.


\ret


\end{document}